# Universal Coefficients and Mayer-Vietoris Sequence for Groupoid Homology

## Proof of a universal coefficient theorem and a Mayer–Vietoris sequence for the homology of ample groupoids

*by*

Luciano Melodia

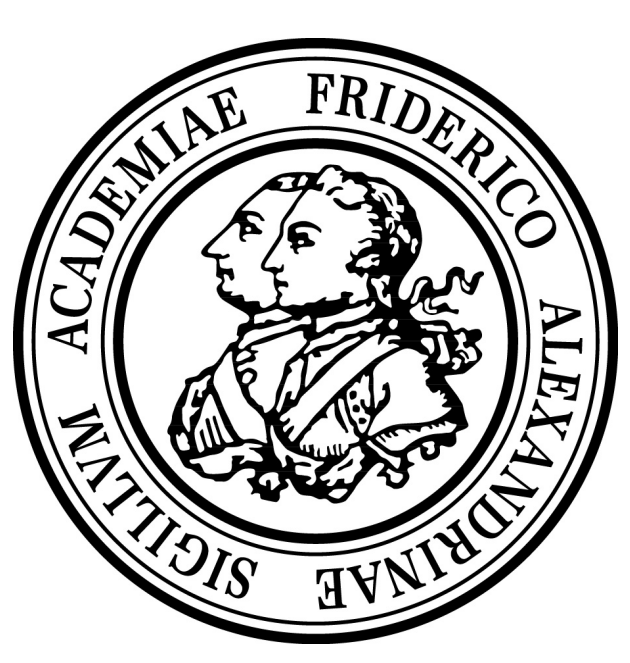

A thesis submitted in
partial fulfillment of the requirements
for the degree

**Master of Science**

in

**Mathematics**

| | |
|---|---|
| Supervisors | Prof. Dr. Kang Li |
| | Prof. Dr. Cathérine Meusburger |
| Institute | Friedrich-Alexander Universität Erlangen-Nürnberg |

# Abstract

We investigate homology of ample groupoids via the compactly supported Moore chain complex of the nerve $\mathcal{G}_\bullet$. Under standing hypotheses that guarantee well behaved compact supports and well defined pushforwards along the face maps, in particular étaleness and local compactness with Hausdorff separation, we define for each $n \geq 0$ the Moore chain group $C_c(\mathcal{G}_n, A)$ of compactly supported continuous $A$-valued functions on $\mathcal{G}_n$, with boundary $\partial_n^A := \sum_{i=0}^{n}(-1)^i (d_i)_*$. The resulting homology groups $H_n(\mathcal{G}; A)$ are functorial for continuous étale homomorphisms and compatible with the standard reduction operations used in computations, including reduction to saturated clopen subsets and, in the ample setting, invariance under Kakutani equivalence. Within this Moore formulation we reprove Matui type long exact sequences and identify the comparison maps already at the chain level, see Theorem 3.1.10.

A central theme is a universal coefficient phenomenon for Moore homology with compact supports. For discrete abelian coefficients $A$ we prove a natural short exact sequence

$$0 \to H_n(\mathcal{G}) \otimes_{\mathbb{Z}} A \xrightarrow{\iota_n^{\mathcal{G}}} H_n(\mathcal{G}; A) \xrightarrow{\kappa_n^{\mathcal{G}}} \mathrm{Tor}_1^{\mathbb{Z}}(H_{n-1}(\mathcal{G}), A) \to 0,$$

natural in both $\mathcal{G}$ and $A$, see Theorem 3.2.3. The key input is the chain level identification $C_c(\mathcal{G}_n, \mathbb{Z}) \otimes_{\mathbb{Z}} A \cong C_c(\mathcal{G}_n, A)$, which reduces the groupoid statement to the classical algebraic universal coefficient theorem applied to the free chain complex $C_c(\mathcal{G}_\bullet, \mathbb{Z})$.

We then isolate the precise obstruction to extending this mechanism beyond discrete coefficients. For a locally compact totally disconnected Hausdorff space $X$ with a basis of compact open sets and a topological abelian group $A$, the image of the canonical comparison map $\Phi_X : C_c(X, \mathbb{Z}) \otimes_{\mathbb{Z}} A \longrightarrow C_c(X, A)$ consists exactly of those compactly supported functions with finite image. Consequently $\Phi_X$ is surjective if and only if every $f \in C_c(X, A)$ has finite image, see Corollary 3.2.4. This shows that the Moore universal coefficient theorem is, in a precise sense, a discrete coefficient phenomenon. Under mild countability hypotheses on $A$ we further construct, for suitable non discrete ample spaces $X$, compactly supported continuous functions $X \to A$ with infinite image, forcing $\Phi_X$ to fail to be surjective.

Finally, we develop a Mayer–Vietoris principle for ample groupoids with discrete coefficients. Given a clopen saturated cover $\mathcal{G}_0 = U_1 \cup U_2$, we construct a short exact sequence of Moore chain complexes and derive a Mayer–Vietoris long exact homology sequence, see Theorem 3.3.10. This sequence is tailored for explicit computations by cutting the unit space into saturated clopen pieces and reconstructing $H_\bullet(\mathcal{G}; A)$ from the corresponding reductions. Combined with the universal coefficient theorem, it cleanly isolates how torsion in integral homology contributes additional homology through $\mathrm{Tor}_1^{\mathbb{Z}}$, a phenomenon illustrated later on examples built from standard ample groupoids such as those arising from shifts of finite type.

## Acknowledgement

The end of my studies closes one chapter and, mercifully, does not open a new one titled "Errata". Writing this thesis took persistence from me and, more importantly, time, patience, and kindness from the people around me. It was a demanding stretch, but also a very happy one. Thank you to everyone who made it possible.

This work sits where algebraic topology meets operator algebras, a place where ideas commute less often than one might hope. I am grateful to my lecturers Cathérine Meusburger, Gandalf Lechner, Peter Fiebig, Thomas Creutzig, and Hermann Schulz-Baldes for excellent courses that provided the background this thesis relies on.

I also thank the department for making it financially possible for me to study mathematics through tutoring positions in linear algebra, analysis, topology, functional analysis, and engineering mathematics. I am especially grateful to Karl-Hermann Neeb, Cathérine Meusburger, Kang Li, and Michael Fried for trusting me with these opportunities.

My special thanks go to Aidan Sims, who visited our university and gave a lecture series on the homology of ample groupoids while I was working on this thesis. Between the discussions during the lectures and our evening conversations over dinner, I learned a great deal and gained real traction on the subject. His stay was as enjoyable as it was enlightening. Thank you.

I am also very grateful to Christian Bönicke. He took the time to discuss my results once they were written down, suggested nontrivial examples, and improved this thesis through careful listening and perfectly timed, practical advice. The evening at Alter Simpl and our conversations were both instructive and great fun. Thank you.

I thank my supervisors Kang Li and, as my mentor in the Master's programme, Cathérine Meusburger, and I also thank Karl-Hermann Neeb once again. I could always come to them with questions and receive quick, thoughtful replies. Their guidance was invaluable, both mathematically and for what comes next. It was a pleasure working with you. Thank you.

I would also like to thank the teaching staff Christiaan van de Ven, Ricardo Correa da Silva, Michael Preeg, Tobias Simon, and Matthias Bauer for a wonderful time.

Finally, my deepest thanks to my partner, Luciana Diaconescu, who had to listen to all of this and still chose to listen. She will soon be writing her own thesis, and I hope I can repay her patience with the same steadiness she showed me and my moods.

# 1 Preliminaries

This chapter fixes notation and collects the background material used throughout the thesis. The central objects are étale groupoids, a common framework that encompasses discrete groups, equivalence relations, and dynamical systems. We begin with the algebraic notion of a groupoid and its orbit decomposition. We then introduce isotropy and principality, which measure stabiliser phenomena. Next we pass to topological groupoids, where all structure maps are continuous. Finally we specialise to étale groupoids, characterised by the source map being a local homeomorphism. In this case the range map is also a local homeomorphism since $r = s \circ i$. This regularity is the input for the analytic and homological constructions in later chapters.

**Setting 1.0.1.** Throughout we write arrows in groupoids as $\gamma, \eta$ and units as $x, y, z$.

## 1.1 Groupoids

A groupoid is a small category in which every morphism is invertible. Concretely, it is a set $\mathcal{G}$ of arrows with partially defined composition $m : \mathcal{G}_2 \to \mathcal{G}$, where

$$\mathcal{G}_2 := \{(\gamma, \eta) \in \mathcal{G} \times \mathcal{G} \mid s(\gamma) = r(\eta)\},$$

together with an object set, the unit space $\mathcal{G}_0$, range and source maps $r, s : \mathcal{G} \to \mathcal{G}_0$, and inversion $\gamma \mapsto \gamma^{-1}$, such that composition is associative on $\mathcal{G}_2$ and units and inverses satisfy the usual identities [21, Remark 2.1.5]. Groupoids are a robust algebraic model for orbit and quotient spaces that can fail to be well behaved topologically. Typical examples are non-$T_1$ quotients such as the shift space modulo the shift, or the circle modulo an irrational rotation [21, p. 1]. As an illustration, let $\mathbb{T} = \mathbb{R}/\mathbb{Z}$ and fix $\alpha \in \mathbb{R} \setminus \mathbb{Q}$. Define $r_\alpha : \mathbb{T} \to \mathbb{T}$ by $r_\alpha(x) = x + \alpha \bmod 1$, and write $x \sim y$ if there exists $n \in \mathbb{Z}$ with $y = r_\alpha^n(x)$. Then $\mathbb{T}/\sim$ is not $T_1$. Indeed, by Dirichlet's Approximation Theorem [9, Theorem 185], for every interval $I \subset \mathbb{T}$ there exists $n \in \mathbb{Z}$ with $n\alpha \bmod 1 \in I$. Hence each orbit $\{x + n\alpha \bmod 1 \mid n \in \mathbb{Z}\}$ is dense and therefore not closed in $\mathbb{T}$. Let $\pi : \mathbb{T} \to \mathbb{T}/\sim$ be the quotient map. Then $\pi^{-1}([x])$ is the non-closed orbit of $x$. If $[x]$ were closed in $\mathbb{T}/\sim$, its preimage $\pi^{-1}([x])$ would be closed in $\mathbb{T}$, a contradiction. We include the proof for completeness.

**Setting 1.1.1.** For $x \in \mathbb{R}$, write its fractional part as $\{x\} := x - \lfloor x \rfloor \in [0, 1)$, where $\lfloor x \rfloor := \max\{n \in \mathbb{Z} \mid n \leq x\}$. The distance to the nearest integer is $\|x\| := \min\{\{x\}, 1 - \{x\}\} \in [0, \frac{1}{2}]$.

**Theorem 1.1.2** (Dirichlet Approximation Theorem [9, Theorem 185]). For any real $\alpha$ and any integer $Q \geq 1$, there exist integers $p, q$ with $1 \leq q \leq Q$ such that

$$|q\alpha - p| \leq \tfrac{1}{Q}, \quad \text{that is} \quad \|q\alpha\| \leq \tfrac{1}{Q}.$$

*Proof.* Consider the $Q+1$ numbers $\{0\}, \{\alpha\}, \ldots, \{Q\alpha\}$ in $[0,1)$. Partition $[0,1)$ into $Q$ half-open intervals $I_j = [\frac{j}{Q}, \frac{j+1}{Q})$ for $j = 0, 1, \ldots, Q-1$. By the pigeonhole principle, two of the $Q+1$ fractional parts, say $\{\ell\alpha\}$ and $\{k\alpha\}$ with $0 \leq k < \ell \leq Q$, lie in the same $I_j$. Hence

$$d := |\{\ell\alpha\} - \{k\alpha\}| < \frac{1}{Q}.$$

Let $q = \ell - k$ with $1 \leq q \leq Q$ and set $p = \lfloor \ell\alpha \rfloor - \lfloor k\alpha \rfloor$. Then

$$q\alpha - p = (\ell\alpha - \lfloor \ell\alpha \rfloor) - (k\alpha - \lfloor k\alpha \rfloor) = \{\ell\alpha\} - \{k\alpha\},$$

so $|q\alpha - p| < 1/Q$. In particular $\|q\alpha\| \leq |q\alpha - p| < 1/Q$. Since $Q$ is arbitrary, for every $\varepsilon > 0$ there exists $q \in \mathbb{N}$ with $\|q\alpha\| < \varepsilon$, hence the set of such $q$ is infinite. If $\alpha$ is rational, then $\|q\alpha\| = 0$ for infinitely many $q$. If $\alpha$ is irrational, continued fractions yield infinitely many coprime $p, q$ with $|q\alpha - p| < 1/q^2$, hence $\|q\alpha\| < 1/q$ for infinitely many $q$ [9, Chapter X]. □

**Example 1.1.3.**

- **Irrational rotation on the circle.** Let $\mathbb{T} = \mathbb{R}/\mathbb{Z}$ with addition modulo 1 and $r_\alpha : \mathbb{T} \to \mathbb{T}$, $r_\alpha(x) = x + \alpha \bmod 1$, for $\alpha \in \mathbb{R} \setminus \mathbb{Q}$. Consider the orbit relation $x \sim y \Leftrightarrow \exists n \in \mathbb{Z} : y = r_\alpha^n(x)$ and the quotient map $\pi : \mathbb{T} \to \mathbb{T}/\sim$. For irrational $\alpha$, each orbit $O(x) = \{r_\alpha^n(x) \mid n \in \mathbb{Z}\}$ is dense in $\mathbb{T}$: let $I \subset [0,1)$ be an interval of length $l > 0$ and write $I = [a, a+l) \bmod 1$. Apply the pigeonhole argument used in Dirichlet's theorem to the $Q+1$ numbers $\{x - a\}$, $\{x + \alpha - a\}$, $\ldots$, $\{x + Q\alpha - a\} \in [0,1)$ and to the partition of $[0,1)$ into $Q$ half-open intervals of length $1/Q$. Then there exist $0 \leq k < \ell \leq Q$ with $\{x + \ell\alpha - a\}$ and $\{x + k\alpha - a\}$ lying in the same subinterval. Without loss of generality assume $\{x + \ell\alpha - a\} \geq \{x + k\alpha - a\}$. Hence

  $$0 \leq \{x + \ell\alpha - a\} - \{x + k\alpha - a\} < \frac{1}{Q}.$$

  Set $q := \ell - k$ and $p := \lfloor x + \ell\alpha - a \rfloor - \lfloor x + k\alpha - a \rfloor$. Then

  $$\{x + q\alpha - a\} = \{(x + \ell\alpha - a) - (x + k\alpha - a)\} = \{x + \ell\alpha - a\} - \{x + k\alpha - a\} < \frac{1}{Q} < l,$$

  so $x + q\alpha \bmod 1 \in I$. As $I$ was arbitrary, $O(x)$ is dense and not closed. Since $\pi$ is continuous and $\pi^{-1}([x]) = O(x)$ is not closed, the singleton $\{[x]\}$ is not closed in $\mathbb{T}/\sim$, thus not $T_1$.

  *Remark* 1.1.4. Dirichlet provides a small step $q$ with $\delta := \|q\alpha\| < l$. Stepping by $\pm q\alpha$, that is along the subsequence $r_\alpha^{n_0 \pm mq}(x)$, forces a hit in any arc of length $> \delta$.

- **Shift space modulo the shift.** An alphabet is a finite, nonempty set $A$ with the discrete topology. The one-sided full shift space over $A$ is $\Sigma := A^{\mathbb{N}}$ with the product topology. A basic open set, also cylinder, determined by a word $w = (w_0, \dots, w_{k-1}) \in A^k$ for $k \in \mathbb{N}$ is

$$[w] := \{x \in \Sigma \mid x_0 = w_0, \dots, x_{k-1} = w_{k-1}\}.$$

The shift map $\sigma : \Sigma \to \Sigma$, $(\sigma(x))_n := x_{n+1}$, is continuous and surjective. Consider the equivalence relation $x \sim y \Leftrightarrow O^+(x) = O^+(y)$, where $O^+(x) := \{\sigma^n(x) \mid n \in \mathbb{N}\}$. Let $\pi_\Sigma : \Sigma \to \Sigma/\sim$ be the quotient map, so $\pi_\Sigma^{-1}([x]) = O^+(x)$ for each $x$. Fix distinct $a, b \in A$ and set $z := (b)_{n\in\mathbb{N}} \in \Sigma$. Define $x \in \Sigma$ by placing $a$ at positions $n_j := 2^j$ and $b$ elsewhere, that is $x_{n_j} = a$ for $j \geq 1$ and $x_n = b$ if $n \notin \{2^j \mid j \geq 1\}$. Then $O^+(x) = \{\sigma^n(x) \mid n \in \mathbb{N}\}$ is not closed and $z \in \overline{O^+(x)} \setminus O^+(x)$. Let $U$ be any cylinder neighbourhood of $z$ determined by its first $k$ coordinates, so $U = [w]$ with $w = b^k$. Choose $j$ so large that $2^{j-1} > k+1$ and set $n := 2^j - (k+1)$. The symbol $a$ at position $2^j$ moves to position $k+1$ in $\sigma^n(x)$, hence the first $k$ coordinates of $\sigma^n(x)$ are all $b$. Moreover, any earlier $a$ at position $2^m$ with $m \leq j-1$ moves to position $2^m - n \leq (k+1) - 2^{j-1} < 0$, so it does not affect these first $k$ coordinates. Thus $\sigma^n(x) \in U$. As $k$ was arbitrary, $z \in \overline{O^+(x)}$. On the other hand, $z \notin O^+(x)$: if $\sigma^n(x) = z$ for some $n$, then $x_{n+i} = b$ for all $i \geq 0$, contradicting the infinitely many occurrences $x_{2^j} = a$. Therefore $O^+(x)$ is not closed in $\Sigma$. Since $\pi_\Sigma^{-1}([\pi_\Sigma(x)]) = O^+(x)$ is not closed and $\pi_\Sigma$ is continuous, the singleton $\{\pi_\Sigma(x)\}$ is not closed in $\Sigma/\sim$. Hence $\Sigma/\sim$ is not $T_1$.

To connect the examples above with the general theory, we begin with a hands-on axiomatisation of a groupoid in the notation $(\mathcal{G}, \mathcal{G}_0, r, s, m, i)$ and with left–to–right composition $\gamma \cdot \eta$, defined when $s(\gamma) = r(\eta)$. This presentation is slightly redundant but keeps all structure maps visible and is useful for intuition and for later constructions, such as nerves and pushforwards. After deriving a few basic consequences, we recast the same object in its concise categorical form, namely a small category in which every arrow is invertible, and verify that the two formulations are equivalent. For étale groupoids we record the equivalent reformulation that $r$ and $s$ are local homeomorphisms, aligning with the simplicial setup in Chapter 2.

**Definition 1.1.5** (Groupoid [21, Remark 2.1.5])**.** A groupoid is a sextuple $(\mathcal{G}, \mathcal{G}_0, r, s, m, i)$ consisting of a set of morphisms $\mathcal{G}$, a distinguished subset $\mathcal{G}_0 \subseteq \mathcal{G}$ of units, range and source maps $r, s : \mathcal{G} \to \mathcal{G}_0$, a partially defined multiplication on composable pairs

$$\cdot : \mathcal{G}_2 := \{(\gamma, \eta) \in \mathcal{G} \times \mathcal{G} \mid s(\gamma) = r(\eta)\} \;\to\; \mathcal{G}, \quad (\gamma, \eta) \mapsto \gamma \cdot \eta,$$

and an inversion map $i : \mathcal{G} \to \mathcal{G}$, $i(\gamma) := \gamma^{-1}$, subject to:
(G1) $r(x) = x = s(x)$ for all $x \in \mathcal{G}_0$.
(G2) $r(\gamma) \cdot \gamma = \gamma = \gamma \cdot s(\gamma)$ for all $\gamma \in \mathcal{G}$.
(G3) $r(\gamma^{-1}) = s(\gamma)$ and $s(\gamma^{-1}) = r(\gamma)$ for all $\gamma \in \mathcal{G}$.
(G4) $\gamma^{-1} \cdot \gamma = s(\gamma)$ and $\gamma \cdot \gamma^{-1} = r(\gamma)$ for all $\gamma \in \mathcal{G}$.
(G5) If $s(\alpha) = r(\beta)$, then $r(\alpha \cdot \beta) = r(\alpha)$ and $s(\alpha \cdot \beta) = s(\beta)$.
(G6) If $s(\alpha) = r(\beta)$ and $s(\beta) = r(\gamma)$, then $(\alpha \cdot \beta) \cdot \gamma = \alpha \cdot (\beta \cdot \gamma)$.

*Remark* 1.1.6. The unit space is

$$\mathcal{G}_0 \coloneqq \{\gamma^{-1}\gamma \mid \gamma \in \mathcal{G}\} = \{\gamma\gamma^{-1} \mid \gamma \in \mathcal{G}\} \subseteq \mathcal{G},$$

and its elements are the unit arrows. Concretely, for $x \in \mathcal{G}_0$ one has $x\gamma = \gamma$ whenever $r(\gamma) = x$ and $\gamma x = \gamma$ whenever $s(\gamma) = x$. The source and range maps $s, r : \mathcal{G} \to \mathcal{G}_0$ decompose $\mathcal{G}$ into the disjoint unions $\mathcal{G} = \bigsqcup_{x\in\mathcal{G}_0} s^{-1}(x) = \bigsqcup_{x\in\mathcal{G}_0} r^{-1}(x)$. The isotropy group at $x \in \mathcal{G}_0$ is

$$\mathcal{G}_x \coloneqq \{\gamma \in \mathcal{G} \mid s(\gamma) = r(\gamma) = x\},$$

a group under multiplication. If $\mathcal{G}$ is a topological groupoid, then $\mathcal{G}_0$, each fibre $s^{-1}(x)$ and $r^{-1}(x)$, and each $\mathcal{G}_x$ inherit the subspace topology. In the special case of a group, a one-object groupoid, $\mathcal{G}_0 = \{e\}$ and $\mathcal{G}_e = \mathcal{G}$.

In general, we will write $\mathcal{G}$ for the arrow space and denote the groupoid by the sextuple or quadruple as in Definition 1.1.5 or Definition 1.1.21. However, sometimes it might be convenient for notational reasons to denote the objects as usual by $\mathcal{G}_0$ and the morphisms by $\mathcal{G}_1$, see Example 1.4.8.

*Remark* 1.1.7. Writing fibre products over $\mathcal{G}_0$ with respect to $s$ and $r$, we have

$$\mathcal{G}_n = \{(\gamma_1, \dots, \gamma_n) \in \mathcal{G}^n \mid s(\gamma_i) = r(\gamma_{i+1}) \text{ for } 1 \le i < n\} \cong \underbrace{\mathcal{G} \,{}_s\!\times_r \mathcal{G} \,{}_s\!\times_r \cdots {}_s\!\times_r \mathcal{G}}_{n \text{ factors}},$$

where $\mathcal{G} \,{}_s\!\times_r \mathcal{G} = \{(\gamma, \eta) \in \mathcal{G} \times \mathcal{G} \mid s(\gamma) = r(\eta)\}$. For $n = 1$ set $\mathcal{G}_1 = \mathcal{G}$. We equip $\mathcal{G}_n$ with the subspace topology inherited from $\mathcal{G}^n$.

**Definition 1.1.8.** Let $(\mathcal{G}, \mathcal{G}_0, r_\mathcal{G}, s_\mathcal{G}, \cdot_\mathcal{G}, -^{-1}_\mathcal{G})$ and $(\mathcal{H}, \mathcal{H}_0, r_\mathcal{H}, s_\mathcal{H}, \cdot_\mathcal{H}, -^{-1}_\mathcal{H})$ be two groupoids. A homomorphism of groupoids is a functor consisting of a map $F : \mathcal{G} \to \mathcal{H}$ and a map $F_0 : \mathcal{G}_0 \to \mathcal{H}_0$ such that:

(F1) $F(\mathcal{G}_0) \subseteq \mathcal{H}_0$ and $F|_{\mathcal{G}_0} = F_0$.
(F2) For all $\gamma \in \mathcal{G}$, $r_\mathcal{H}(F(\gamma)) = F_0(r_\mathcal{G}(\gamma))$ and $s_\mathcal{H}(F(\gamma)) = F_0(s_\mathcal{G}(\gamma))$.
(F3) For all $\alpha, \beta \in \mathcal{G}$ with $s_\mathcal{G}(\alpha) = r_\mathcal{G}(\beta)$, $F(\alpha) \cdot F(\beta)$ is defined and $F(\alpha \cdot \beta) = F(\alpha) \cdot F(\beta)$.

In particular, $F(1_x) = 1_{F_0(x)}$ for all $x \in \mathcal{G}_0$ and $F(\gamma^{-1}) = F(\gamma)^{-1}$ for all $\gamma \in \mathcal{G}$.

*Remark* 1.1.9. We will write $\mathcal{G}$ for the data $(\mathcal{G}, \mathcal{G}_0, r, s, m, i)$ whenever the unit space, the range and source maps, the multiplication, and the inversion are clear from context.

**Definition 1.1.10.** A homomorphism of groupoids $F : \mathcal{G} \to \mathcal{H}$ is an isomorphism of groupoids if there exists a homomorphism $G : \mathcal{H} \to \mathcal{G}$ with $G \circ F = \mathrm{id}_\mathcal{G}$ and $F \circ G = \mathrm{id}_\mathcal{H}$. If $F$ admits an inverse homomorphism, then $F$ is bijective and satisfies (F1)–(F3). Conversely, if $F$ is bijective and satisfies (F1)–(F3), then the set theoretic inverse $F^{-1}$ is again a homomorphism and hence $F$ is an isomorphism of groupoids.

*Remark* 1.1.11. For étale groupoids we will require $F$ and $F_0$ to be continuous. Since $\mathcal{G}_0 \subseteq \mathcal{G}$ carries the subspace topology, continuity of $F$ implies continuity of $F_0 = F|_{\mathcal{G}_0}$.

Principal groupoids capture pure orbit structure without isotropy groups of units, also known as internal stabilisers. Every morphism is determined uniquely by its source and range. For transformation groupoids this corresponds to free actions.

**Definition 1.1.12.** We call $\mathcal{G}$ principal if any of the following equivalent conditions hold:
(P1) **Trivial isotropy:** $\forall x \in \mathcal{G}_0$, the isotropy group $\mathcal{G}_x := \{\gamma \in \mathcal{G} \mid r(\gamma) = s(\gamma) = x\}$ equals $\{x\}$.
(P2) **Injective anchor:** $(r,s) : \mathcal{G} \to \mathcal{G}_0 \times \mathcal{G}_0$ is injective, hence there is at most one arrow $y \to x$.

**Definition 1.1.13.** Let a discrete group $\Gamma$ act on a locally compact Hausdorff space $X$ by homeomorphisms. The action is free if the stabiliser at $x$ given by $\Gamma_x := \{g \in \Gamma \mid g \cdot x = x\}$ is trivial for every $x \in X$.

**Proposition 1.1.14.** The transformation groupoid from Example 1.1.18 $\Gamma \ltimes X$ with objects $X$, morphisms $\Gamma \times X$, $r(g,x) = g \cdot x$, $s(g,x) = x$, is principal if and only if the action is free.

*Remark* 1.1.15. Here, $\Gamma \ltimes X$ is also étale, see Section 1.4. Étaleness does not depend on freeness.

*Proof.* For $x \in X$, the isotropy group of $\Gamma \ltimes X$ at $x$ is $(\Gamma \ltimes X)_x = \{(g,x) \mid g \cdot x = x\}$, which is isomorphic to the stabiliser $\Gamma_x = \{g \in \Gamma \mid g \cdot x = x\}$ via $(g,x) \mapsto g$. Thus $\Gamma \ltimes X$ is principal if and only if every $\Gamma_x$ is trivial, that is the action is free. For étaleness, assume $\Gamma$ is discrete and the action is continuous. Give $\Gamma \times X$ the product topology. For any $(g,x) \in \Gamma \times X$ and any open $U \ni x$, the restrictions $s|_{\{g\}\times U} : \{g\} \times U \to U$ and $r|_{\{g\}\times U} : \{g\} \times U \to g \cdot U$ are homeomorphisms. Hence $s$ and $r$ are local homeomorphisms and $\Gamma \ltimes X$ is étale. □

Minimality singles out groupoids with no proper invariant pieces of the unit space. Every orbit is dense, so the dynamics are indecomposable. For transformation groupoids this corresponds to minimal actions.

**Definition 1.1.16** (Minimal groupoid)**.** Let $\mathcal{G}$ be a topological groupoid. For $x \in \mathcal{G}_0$ write the orbit $O(x) := \{r(\gamma) \mid \gamma \in \mathcal{G},\ s(\gamma) = x\}$. A subset $U \subseteq \mathcal{G}_0$ is invariant, also known as saturated, if it is a union of orbits. This means that $r(s^{-1}(U)) = U$. The same condition can be written as $s(r^{-1}(U)) = U$. We call $\mathcal{G}$ minimal if any of the following hold:

(M1) **Dense orbits:** $\overline{O(x)} = \mathcal{G}_0$ for all $x \in \mathcal{G}_0$.
(M2) There is no nonempty proper open invariant $U \subset \mathcal{G}_0$.

*Remark* 1.1.17. For a left action $\Gamma \triangleright X$, the transformation groupoid $\Gamma \ltimes X$ from Example 1.1.18 is minimal if and only if the action is minimal, hence every $\Gamma$-orbit is dense in $X$.

**Example 1.1.18** (Groupoids [21, Ex. 2.1.7–2.1.11])**.**

- **Groups.** Every group $G$ is a discrete groupoid with one unit: $\mathcal{G} = G$, $\mathcal{G}_0 = \{e\}$, $r(g) = s(g) = e$, multiplication $g \cdot h = gh$, inversion $\gamma \mapsto \gamma^{-1}$.
  A groupoid is a group if and only if $\mathcal{G}_0$ is a singleton.
- **Group bundles.** Let $X$ be a set and $\{G_x\}_{x\in X}$ a family of groups. Put $\mathcal{G} := \bigsqcup_{x\in X}\{x\} \times G_x$, $\mathcal{G}_0 = \{(x, e_{G_x}) \mid x \in X\} \cong X$, $r(x,g) = x = s(x,g)$, $(x,g)\cdot(x,h) = (x,gh)$, $(x,g)^{-1} = (x,g^{-1})$. If each $G_x$ is a topological group and the disjoint union carries the sum topology, this is a topological groupoid, see Definition 1.3.1.

- **Equivalence-relation groupoid.** For a set $X$ and an equivalence relation $R \subseteq X \times X$, let $\mathcal{G} = R$, $\mathcal{G}_0 = X$, $r(x,y) = x$, $s(x,y) = y$, $(x,y)^{-1} = (y,x)$, and $(x,y) \cdot (y,z) = (x,z)$. If $X$ is a topological space and $R \subseteq X \times X$ is equipped with the subspace topology, this is a topological groupoid. It is principal. It is étale whenever the coordinate projections are local homeomorphisms, for instance if $X$ is discrete, or if $R$ is an étale equivalence relation on a Cantor set.
- **Matrix finite pair groupoid.** For $X = \{1, \dots, n\}$ the pair groupoid $\mathcal{G} = X \times X$ with morphisms $(i,j)$, units $X$, and composition $(i,k) \cdot (k,j) = (i,j)$ is discrete, principal, and étale. Its groupoid C*-algebra is $\mathrm{Mat}(n \times n, \mathbb{C})$ with matrix units $E_{ij}$ that correspond to $(i,j)$.
- **Transformation groupoid.** Let a discrete group $\Gamma$ act on an **LCH** space $X$ by homeomorphisms. The transformation groupoid $\mathcal{G} = \Gamma \ltimes X$ has objects $\mathcal{G}_0 = X$, morphisms $\Gamma \times X$, $r(g,x) = g \cdot x$, $s(g,x) = x$, $(g,x)^{-1} = (g^{-1}, g \cdot x)$, and $(h, g \cdot x) \cdot (g,x) = (hg, x)$ whenever $s(h, g \cdot x) = r(g,x)$. Since $\Gamma$ is discrete, $r$ and $s$ are local homeomorphisms. Thus $\Gamma \ltimes X$ is étale. For $r_\alpha$ on $\mathbb{T}$, the groupoid $\mathbb{Z} \ltimes \mathbb{T}$ is Hausdorff, étale, principal, and minimal.
- **Deaconu–Renault groupoid of a local homeomorphism.** Let $f : X \to X$ be a local homeomorphism of an **LCH** space. Define

  $$\mathcal{G}_f := \{(x,k,y) \in X \times \mathbb{Z} \times X \mid \exists m,n \geq 0 : k = n - m \text{ and } f^n(x) = f^m(y)\},$$

  with $r(x,k,y) = x$, $s(x,k,y) = y$, $(x,k,y)^{-1} = (y,-k,x)$, and $(x,k,y) \cdot (y,\ell,z) = (x, k+\ell, z)$. Then $\mathcal{G}_f$ is Hausdorff and étale. It is principal if and only if no point is eventually periodic. It becomes principal on the aperiodic part.

Using the axioms (G1)–(G6) from Definition 1.1.5 we first establish the invertibility of the morphisms within a groupoid.

**Lemma 1.1.19.** $(\gamma^{-1})^{-1} = \gamma$ for all $\gamma \in \mathcal{G}$.

*Proof.* By (G4) applied to $\gamma^{-1}$ we have $(\gamma^{-1})^{-1} \cdot \gamma^{-1} = s(\gamma^{-1})$ and $\gamma^{-1} \cdot (\gamma^{-1})^{-1} = r(\gamma^{-1})$. By (G4) and (G3), $\gamma \cdot \gamma^{-1} = r(\gamma) = s(\gamma^{-1})$ and $\gamma^{-1} \cdot \gamma = s(\gamma) = r(\gamma^{-1})$. Thus both $(\gamma^{-1})^{-1}$ and $\gamma$ are two-sided inverses of $\gamma^{-1}$. In a groupoid, such an inverse is unique: if $\alpha, \beta \in \mathcal{G}$ satisfy $\alpha \cdot \gamma^{-1} = s(\gamma^{-1})$ and $\gamma^{-1} \cdot \beta = r(\gamma^{-1})$, then, using (G2) and associativity (G6), $\alpha = \alpha \cdot r(\gamma^{-1}) = \alpha \cdot (\gamma^{-1} \cdot \beta) = (\alpha \cdot \gamma^{-1}) \cdot \beta = s(\gamma^{-1}) \cdot \beta = \beta$. Hence $(\gamma^{-1})^{-1} = \gamma$. □

Observe that $(\gamma, \gamma^{-1}) \in \mathcal{G}_2$ for all $\gamma \in \mathcal{G}$, since $r(\gamma) = s(\gamma^{-1})$ by (G3).

**Lemma 1.1.20.** Let $\mathcal{G}$ be a groupoid and $\gamma \in \mathcal{G}$. Then

1. for all $(\gamma, \eta) \in \mathcal{G}_2$ one has $\gamma^{-1} \cdot (\gamma \cdot \eta) = \eta$ and $(\gamma \cdot \eta) \cdot \eta^{-1} = \gamma$.
2. $(r(\gamma), \gamma)$ and $(\gamma, s(\gamma))$ lie in $\mathcal{G}_2$.
3. $\gamma^{-1}$ is unique such that $(\gamma, \gamma^{-1}) \in \mathcal{G}_2$, $\gamma \cdot \gamma^{-1} = r(\gamma)$, $(\gamma^{-1}, \gamma) \in \mathcal{G}_2$, and $\gamma^{-1} \cdot \gamma = s(\gamma)$.

*Proof.*

1. Suppose $(\gamma, \eta) \in \mathcal{G}_2$, that is $s(\gamma) = r(\eta)$. Then

   $$\gamma^{-1} \cdot (\gamma \cdot \eta) = (\gamma^{-1} \cdot \gamma) \cdot \eta = s(\gamma) \cdot \eta = r(\eta) \cdot \eta = \eta,$$

where $(\gamma^{-1}, \gamma) \in \mathcal{G}_2$ by (G3) and we used (G6), (G4), and (G2). Similarly,

$$(\gamma \cdot \eta) \cdot \eta^{-1} = \gamma \cdot (\eta \cdot \eta^{-1}) = \gamma \cdot r(\eta) = \gamma,$$

using (G6), (G4), and (G2), and $r(\eta) = s(\gamma)$ to justify the last step.

2. Using (G4) and (G5),

$$s(r(\gamma)) = s(\gamma \cdot \gamma^{-1}) = s(\gamma^{-1}) = r(\gamma),$$

so $(r(\gamma), \gamma) \in \mathcal{G}_2$. Likewise,

$$r(s(\gamma)) = r(\gamma^{-1} \cdot \gamma) = r(\gamma^{-1}) = s(\gamma),$$

so $(\gamma, s(\gamma)) \in \mathcal{G}_2$.

3. We split this proof into three parts:
   - **Existence:** By (G3), $(\gamma, \gamma^{-1})$ and $(\gamma^{-1}, \gamma)$ lie in $\mathcal{G}_2$, and by (G4), $\gamma \cdot \gamma^{-1} = r(\gamma)$ and $\gamma^{-1} \cdot \gamma = s(\gamma)$.
   - **Uniqueness from the right-unit equation:** Assume $(\gamma, \alpha) \in \mathcal{G}_2$ and $\gamma \cdot \alpha = r(\gamma)$. Then

     $$\alpha = s(\gamma) \cdot \alpha = (\gamma^{-1} \cdot \gamma) \cdot \alpha = \gamma^{-1} \cdot (\gamma \cdot \alpha) = \gamma^{-1} \cdot r(\gamma) = \gamma^{-1},$$

     using (G2), (G4), (G6), and $s(\gamma^{-1}) = r(\gamma)$ from (G3).
   - **Uniqueness from the left-unit equation:** Assume $(\alpha, \gamma) \in \mathcal{G}_2$ and $\alpha \cdot \gamma = s(\gamma)$. Then

     $$\alpha = \alpha \cdot r(\gamma) = \alpha \cdot (\gamma \cdot \gamma^{-1}) = (\alpha \cdot \gamma) \cdot \gamma^{-1} = s(\gamma) \cdot \gamma^{-1} = \gamma^{-1},$$

     using (G2), (G4), and (G6), and $r(\gamma^{-1}) = s(\gamma)$ from (G3).

   Thus $\gamma^{-1}$ is unique.

□

For later simplicial constructions, we want the domain of composition to be recorded as part of the data. Presenting a groupoid via its set of composable pairs makes the multiplication explicitly partial, allows us to form the iterated fibre products $\mathcal{G}_n$ cleanly, and aligns with the étale setting in which $\mathcal{G}_2$ is the fibre product over $\mathcal{G}_0$. This viewpoint streamlines nerves and pushforwards with topological arguments.

**Definition 1.1.21** (Groupoid [21, Definition 2.1.1])**.** A groupoid is given by the quadruple $(\mathcal{G}, \mathcal{G}_2, m, i)$ where $\mathcal{G}$ is a set of morphisms, $\mathcal{G}_2 \subseteq \mathcal{G} \times \mathcal{G}$ is the set of composable pairs, $m : \mathcal{G}_2 \to \mathcal{G}$ multiplication, written $m(\alpha, \beta) = \alpha \cdot \beta$, and $i : \mathcal{G} \to \mathcal{G}$ inversion, subject to:

(G1′) $(\gamma^{-1})^{-1} = \gamma$ for all $\gamma \in \mathcal{G}$.

(G2′) If $(\alpha, \beta)$ and $(\beta, \gamma)$ lie in $\mathcal{G}_2$, then $(\alpha \cdot \beta, \gamma)$ and $(\alpha, \beta \cdot \gamma)$ lie in $\mathcal{G}_2$ and $(\alpha \cdot \beta) \cdot \gamma = \alpha \cdot (\beta \cdot \gamma)$.

(G3′) For all $\gamma \in \mathcal{G}$, $(\gamma, \gamma^{-1}) \in \mathcal{G}_2$.

(G4′) For all $(\gamma, \eta) \in \mathcal{G}_2$ one has $\gamma^{-1} \cdot (\gamma \cdot \eta) = \eta$ and $(\gamma \cdot \eta) \cdot \eta^{-1} = \gamma$.

Units are derived rather than primitive in this version. The four axioms encode inversion, associativity with domain control, and the two unit laws via cancellation identities, from which $\mathcal{G}_0$ together with $r$ and $s$ are defined. This presentation is economical for proofs and aligns with categorical practice, morphisms first, objects reconstructed, yet it is equivalent to the expanded axiom list used earlier. The relationship with the previous definition is as follows:

**Proposition 1.1.22** (Equivalence of presentations)**.** Let $\mathcal{G}$ be a groupoid $(\mathcal{G}, \mathcal{G}_0, r, s, m, i)$ and axioms (G1)–(G6) as in Definition 1.1.5. Then one obtains the presentation from Definition 1.1.21 by setting $\mathcal{G}_2 := \{(\alpha, \beta) \in \mathcal{G} \times \mathcal{G} \mid s(\alpha) = r(\beta)\}$ and $m(\alpha, \beta) := \alpha \cdot \beta$, with inversion as before. Conversely, let $(\mathcal{G}, \mathcal{G}_2, m, i)$ satisfy (G1')–(G4'). Define $\mathcal{G}_0 := \{\gamma \cdot \gamma^{-1} \mid \gamma \in \mathcal{G}\} = \{\gamma^{-1} \cdot \gamma \mid \gamma \in \mathcal{G}\}$, $r(\gamma) := \gamma \cdot \gamma^{-1}$, and $s(\gamma) := \gamma^{-1} \cdot \gamma$. Then the axioms (G1)–(G6) hold with left–to–right composition, and the domain of multiplication agrees with source–range matching, that is $(\alpha, \beta) \in \mathcal{G}_2 \Leftrightarrow s(\alpha) = r(\beta)$.

*Proof.*

- **Assume we are given** $(\mathcal{G}, \mathcal{G}_0, r, s, m, i)$ **satisfying (G1)–(G6).** Define $\mathcal{G}_2 := \{(\alpha, \beta) \in \mathcal{G} \times \mathcal{G} \mid s(\alpha) = r(\beta)\}$ and $m(\alpha, \beta) := \alpha \cdot \beta$. Then:
  - (G1') holds by the involutivity lemma $(\gamma^{-1})^{-1} = \gamma$, which follows from (G2)–(G4) and associativity (G6).
  - (G2') holds by (G5) and (G6): if $(\alpha, \beta)$ and $(\beta, \gamma)$ lie in $\mathcal{G}_2$, then $s(\alpha) = r(\beta)$ and $s(\beta) = r(\gamma)$. By (G5), $s(\alpha \cdot \beta) = s(\beta) = r(\gamma)$ and $r(\beta \cdot \gamma) = r(\beta) = s(\alpha)$, hence $(\alpha \cdot \beta, \gamma)$ and $(\alpha, \beta \cdot \gamma)$ lie in $\mathcal{G}_2$, and (G6) yields $(\alpha \cdot \beta) \cdot \gamma = \alpha \cdot (\beta \cdot \gamma)$.
  - (G3') follows from (G3), which gives $(\gamma, \gamma^{-1}) \in \mathcal{G}_2$ for all $\gamma$.
  - (G4') is exactly the cancellation lemma above: for all $(\gamma, \eta) \in \mathcal{G}_2$ one has $\gamma^{-1} \cdot (\gamma \cdot \eta) = \eta$ and $(\gamma \cdot \eta) \cdot \eta^{-1} = \gamma$.
- **Assume we are given** $(\mathcal{G}, \mathcal{G}_2, m, i)$ **satisfying (G1')–(G4').** Define $\mathcal{G}_0 := \{\gamma \cdot \gamma^{-1} \mid \gamma \in \mathcal{G}\} = \{\gamma^{-1} \cdot \gamma \mid \gamma \in \mathcal{G}\}$, $r(\gamma) := \gamma \cdot \gamma^{-1}$, and $s(\gamma) := \gamma^{-1} \cdot \gamma$. We verify (G1)–(G6) and the identification of the domain of composition.
  - By (G3') we have $(\gamma, \gamma^{-1})$ and $(\gamma^{-1}, \gamma)$ in $\mathcal{G}_2$. Taking $\eta = \gamma^{-1}$ in (G4') and using (G1'), we obtain $(\gamma \cdot \gamma^{-1}) \cdot (\gamma^{-1})^{-1} = \gamma$, hence $r(\gamma) \cdot \gamma = \gamma$. Taking $\gamma = \eta^{-1}$ in the second identity of (G4') yields $\eta \cdot (\eta^{-1} \cdot \eta) = \eta$, hence $\gamma \cdot s(\gamma) = \gamma$. Thus (G2) holds.
  - Let $(\alpha, \beta) \in \mathcal{G}_2$ and set $\delta := \beta^{-1} \cdot \alpha^{-1}$. Using (G2') and (G4'),

    $$(\alpha \cdot \beta) \cdot \delta = ((\alpha \cdot \beta) \cdot \beta^{-1}) \cdot \alpha^{-1} = \alpha \cdot \alpha^{-1} = r(\alpha),$$

    and similarly

    $$\delta \cdot (\alpha \cdot \beta) = \beta^{-1} \cdot (\alpha^{-1} \cdot (\alpha \cdot \beta)) = \beta^{-1} \cdot \beta = s(\beta).$$

    Thus $\delta$ is a two-sided inverse of $\alpha \cdot \beta$, hence $(\alpha \cdot \beta)^{-1} = \beta^{-1} \cdot \alpha^{-1}$.
  - For $x \in \mathcal{G}_0$ there is $\gamma$ with $x = \gamma \cdot \gamma^{-1}$. Then $x^{-1} = (\gamma \cdot \gamma^{-1})^{-1} = (\gamma^{-1})^{-1} \cdot \gamma^{-1} = \gamma \cdot \gamma^{-1} = x$, so $x$ is fixed by inversion. Consequently $r(x) = x \cdot x^{-1} = x$ and $s(x) = x^{-1} \cdot x = x$, which is (G1). Also $r(\gamma^{-1}) = \gamma^{-1} \cdot (\gamma^{-1})^{-1} = \gamma^{-1} \cdot \gamma = s(\gamma)$ and $s(\gamma^{-1}) = (\gamma^{-1})^{-1} \cdot \gamma^{-1} = \gamma \cdot \gamma^{-1} = r(\gamma)$, which is (G3). Finally, (G4) holds by definition of $r$ and $s$.

- Let $(\alpha,\beta) \in \mathcal{G}_2$. Then

$$r(\alpha \cdot \beta) = (\alpha \cdot \beta) \cdot (\alpha \cdot \beta)^{-1} = (\alpha \cdot \beta) \cdot (\beta^{-1} \cdot \alpha^{-1}) = \alpha \cdot \alpha^{-1} = r(\alpha),$$

    and

$$s(\alpha \cdot \beta) = (\alpha \cdot \beta)^{-1} \cdot (\alpha \cdot \beta) = (\beta^{-1} \cdot \alpha^{-1}) \cdot (\alpha \cdot \beta) = \beta^{-1} \cdot \beta = s(\beta),$$

    which is (G5).
  - Associativity on composable triples is exactly (G2'), giving (G6).
- $(\boldsymbol{\alpha},\boldsymbol{\beta}) \in \mathcal{G}_2 \Leftrightarrow \mathbf{s}(\boldsymbol{\alpha}) = \mathbf{r}(\boldsymbol{\beta})$. If $(\alpha,\beta) \in \mathcal{G}_2$, then $(\alpha \cdot \beta, \beta^{-1}) \in \mathcal{G}_2$ by (G2') and (G3'), hence by (G5) and (G3),

$$s(\alpha) = s((\alpha \cdot \beta) \cdot \beta^{-1}) = s(\beta^{-1}) = r(\beta).$$

  Conversely, assume $s(\alpha) = r(\beta)$. By (G3') and (G2') applied to $(\alpha,\alpha^{-1})$ and $(\alpha^{-1},\alpha)$, we have $(\alpha, s(\alpha)) \in \mathcal{G}_2$. Similarly, $(r(\beta),\beta) \in \mathcal{G}_2$. Since $s(\alpha) = r(\beta)$, it follows that $(s(\alpha),\beta) \in \mathcal{G}_2$, and another application of (G2') to $(\alpha, s(\alpha))$ and $(s(\alpha),\beta)$ yields $(\alpha,\beta) \in \mathcal{G}_2$. □

## 1.2 Isotropy

In many applications groupoids serve as replacement without singularities for quotients by group actions, and the isotropy encodes precisely where this replacement still carries nontrivial stabiliser symmetry. Understanding the isotropy subgroupoid and the notion of principality is therefore essential for detecting when a groupoid behaves like a genuine equivalence relation, and for relating its analytical and homological invariants to those of classical free actions.

**Setting 1.2.1.** For $x, y \in \mathcal{G}_0$ set $\mathcal{G}(x,y) := \{\gamma \in \mathcal{G} \mid r(\gamma) = x,\ s(\gamma) = y\} = r^{-1}(x) \cap s^{-1}(y)$. For $U, V \subseteq \mathcal{G}$ write $UV := \{\alpha \cdot \beta \mid \alpha \in U,\ \beta \in V,\ s(\alpha) = r(\beta)\}$.

**Definition 1.2.2** (Isotropy subgroupoid [21, Chapter 2.2])**.** $\mathrm{Iso}(\mathcal{G}) := \bigsqcup_{x \in \mathcal{G}_0} \mathcal{G}(x,x) = \{\gamma \in \mathcal{G} \mid r(\gamma) = s(\gamma)\}$. It is a subgroupoid, a group bundle over $\mathcal{G}_0$ with fibres $\mathcal{G}_x := \mathcal{G}(x,x)$. In particular $\mathcal{G}_0 \subseteq \mathrm{Iso}(\mathcal{G})$.

**Lemma 1.2.3** ([21, Lemma 2.2.1])**.** $\mathcal{G}$ is principal if and only if $\mathrm{Iso}(\mathcal{G}) = \mathcal{G}_0$.

*Proof.* If $\mathcal{G}$ is principal, then for every $x \in \mathcal{G}_0$ the isotropy group $\mathcal{G}_x$ equals $\{x\}$. Let $\gamma \in \mathrm{Iso}(\mathcal{G})$. Then $r(\gamma) = s(\gamma) = x$ for some $x \in \mathcal{G}_0$, hence $\gamma \in \mathcal{G}_x = \{x\}$, so $\gamma = x \in \mathcal{G}_0$. Thus $\mathrm{Iso}(\mathcal{G}) \subseteq \mathcal{G}_0$. The reverse inclusion holds because for each $x \in \mathcal{G}_0$ we have $r(x) = x = s(x)$ by (G1). Hence $\mathrm{Iso}(\mathcal{G}) = \mathcal{G}_0$. Conversely, suppose $\mathrm{Iso}(\mathcal{G}) = \mathcal{G}_0$. Fix $x \in \mathcal{G}_0$ and $\gamma \in \mathcal{G}$ with $r(\gamma) = s(\gamma) = x$. Then $\gamma \in \mathrm{Iso}(\mathcal{G}) = \mathcal{G}_0$. By (G1), $r(\gamma) = \gamma$, and since also $r(\gamma) = x$, it follows that $\gamma = x$. Therefore $\mathcal{G}_x = \{x\}$ for all $x \in \mathcal{G}_0$, and $\mathcal{G}$ is principal. □

**Example 1.2.4** (Isotropy in a transformation groupoid [21, Example 2.2.2])**.** Let a discrete group $\Gamma$ act on a set or on an **LCH** space $X$ by homeomorphisms. In $\Gamma \ltimes X$ one has

$$\mathrm{Iso}(\Gamma \ltimes \mathrm{X}) = \bigsqcup_{x \in X} \{(g, x) \mid g \cdot x = x\} \;\cong\; \bigsqcup_{x \in X} \Gamma_x,$$

where $\Gamma_x := \{g \in \Gamma \mid g \cdot x = x\}$ are the stabilisers. Hence $\Gamma \ltimes X$ is principal if and only if the action is free.

**Lemma 1.2.5** ([21, Lemma 2.2.3])**.** For $\gamma \in \mathcal{G}$ the map $\mathrm{Ad}_\gamma : \mathcal{G}_{s(\gamma)} \to \mathcal{G}_{r(\gamma)}$, $\alpha \mapsto \gamma \cdot \alpha \cdot \gamma^{-1}$ is a group isomorphism.

*Proof.*

- **Well definedness:** Let $\alpha \in \mathcal{G}_{s(\gamma)}$. Then $r(\alpha) = s(\gamma)$ and $s(\alpha) = s(\gamma)$. By (G5), $r(\gamma \cdot \alpha) = r(\gamma)$ and $s(\gamma \cdot \alpha) = s(\alpha) = s(\gamma)$. Since $r(\gamma^{-1}) = s(\gamma)$ by (G3), we have $(\gamma \cdot \alpha, \gamma^{-1}) \in \mathcal{G}_2$, and another application of (G5) gives $r((\gamma \cdot \alpha) \cdot \gamma^{-1}) = r(\gamma)$ and $s((\gamma \cdot \alpha) \cdot \gamma^{-1}) = s(\gamma^{-1}) = r(\gamma)$, so $\mathrm{Ad}_\gamma(\alpha) \in \mathcal{G}_{r(\gamma)}$.
- **Homomorphism property:** For $\alpha, \beta \in \mathcal{G}_{s(\gamma)}$ the pairs $(\gamma, \alpha)$, $(\alpha, \beta)$, and $(\beta, \gamma^{-1})$ are composable. By associativity (G6),

  $$\mathrm{Ad}_\gamma(\alpha) \cdot \mathrm{Ad}_\gamma(\beta) = (\gamma \cdot \alpha \cdot \gamma^{-1}) \cdot (\gamma \cdot \beta \cdot \gamma^{-1}) = \gamma \cdot \alpha \cdot (\gamma^{-1} \cdot \gamma) \cdot \beta \cdot \gamma^{-1} = \gamma \cdot (\alpha \cdot \beta) \cdot \gamma^{-1} = \mathrm{Ad}_\gamma(\alpha \cdot \beta),$$

  using (G4) in the third equality.
- **Bijectivity:** Consider $\mathrm{Ad}_{\gamma^{-1}} : \mathcal{G}_{r(\gamma)} \to \mathcal{G}_{s(\gamma)}$. For $\alpha \in \mathcal{G}_{s(\gamma)}$,

  $$\mathrm{Ad}_{\gamma^{-1}}(\mathrm{Ad}_\gamma(\alpha)) = \gamma^{-1} \cdot (\gamma \cdot \alpha \cdot \gamma^{-1}) \cdot \gamma = (\gamma^{-1} \cdot \gamma) \cdot \alpha \cdot (\gamma^{-1} \cdot \gamma) = s(\gamma) \cdot \alpha \cdot s(\gamma) = \alpha,$$

  using (G6), (G4), and (G2). Similarly, for $\beta \in \mathcal{G}_{r(\gamma)}$,

  $$\mathrm{Ad}_\gamma(\mathrm{Ad}_{\gamma^{-1}}(\beta)) = \gamma \cdot (\gamma^{-1} \cdot \beta \cdot \gamma) \cdot \gamma^{-1} = (\gamma \cdot \gamma^{-1}) \cdot \beta \cdot (\gamma \cdot \gamma^{-1}) = r(\gamma) \cdot \beta \cdot r(\gamma) = \beta,$$

  using (G6), (G4), and (G2). Thus $\mathrm{Ad}_{\gamma^{-1}}$ is the inverse of $\mathrm{Ad}_\gamma$.

Therefore $\mathrm{Ad}_\gamma$ is a group isomorphism. □

## 1.3 Topological Groupoids

Topological structure allows groupoids to encode convergence, compactness, and local symmetries, unifying groups, group actions, and equivalence relations within a single analytic–geometric framework. It also captures pathological quotients faithfully – for instance the orbit spaces in Example 1.1.3 – while retaining continuous structure maps that are crucial for later homological constructions.

**Definition 1.3.1** (Topological groupoids [11, §2.1])**.** A topological groupoid is a groupoid $(\mathcal{G}, \mathcal{G}_0, r, s, m, i, u)$ equipped with topologies on $\mathcal{G}$ and $\mathcal{G}_0$ such that:

(T1) $\mathcal{G}_0 \subseteq \mathcal{G}$ is the unit space, endowed with the subspace topology.
(T2) $r, s : \mathcal{G} \to \mathcal{G}_0$ are continuous.
(T3) $\mathcal{G}_2 := \{(\gamma, \eta) \in \mathcal{G} \times \mathcal{G} \mid s(\gamma) = r(\eta)\}$ carries the subspace topology from $\mathcal{G} \times \mathcal{G}$, and multiplication $m : \mathcal{G}_2 \to \mathcal{G}$, $(\gamma, \eta) \mapsto \gamma \cdot \eta$, is continuous.
(T4) The inversion map $i : \mathcal{G} \to \mathcal{G}$, $i(\gamma) := \gamma^{-1}$, is continuous.
(T5) The unit map $u : \mathcal{G}_0 \to \mathcal{G}$, $u(x) := 1_x$, is continuous.

Subgroupoids isolate controlled pieces of a groupoid – reductions to invariant subspaces, isotropy, kernels/images of homomorphisms. They are building blocks for exact constructions.

**Definition 1.3.2** (Subgroupoid [11, §2.1])**.** Let $\mathcal{G}$ be a groupoid. A subgroupoid consists of subsets $\mathcal{H} \subseteq \mathcal{G}$ and $\mathcal{H}_0 \subseteq \mathcal{G}_0$ together with the restricted structure maps $r|_{\mathcal{H}}^{\mathcal{H}_0} : \mathcal{H} \to \mathcal{H}_0, s|_{\mathcal{H}}^{\mathcal{H}_0} : \mathcal{H} \to \mathcal{H}_0, m|_{\mathcal{H}_2}^{\mathcal{H}} : \mathcal{H}_2 \to \mathcal{H}, i|_{\mathcal{H}}^{\mathcal{H}} : \mathcal{H} \to \mathcal{H}$, where $\mathcal{H}_2 := \{(\alpha, \beta) \in \mathcal{H} \times \mathcal{H} \mid s(\alpha) = r(\beta)\}$, such that:
(S1) $\mathcal{H}_0 = \mathcal{H} \cap \mathcal{G}_0$ and $r|_{\mathcal{H}}^{\mathcal{H}_0}(\mathcal{H}) \cup s|_{\mathcal{H}}^{\mathcal{H}_0}(\mathcal{H}) \subseteq \mathcal{H}_0$.
(S2) If $(\alpha, \beta) \in \mathcal{H}_2$, then $m|_{\mathcal{H}_2}^{\mathcal{H}}(\alpha, \beta) \in \mathcal{H}$.
(S3) If $\gamma \in \mathcal{H}$, then $i|_{\mathcal{H}}^{\mathcal{H}}(\gamma) \in \mathcal{H}$.
(S4) For all $x \in \mathcal{H}_0$ one has $r|_{\mathcal{H}}^{\mathcal{H}_0}(x) = x = s|_{\mathcal{H}}^{\mathcal{H}_0}(x)$ and $x \in \mathcal{H}$.

*Remark* 1.3.3. We suppress the notation of the restriction and corestriction, due to readability.

*Remark* 1.3.4. The subgroupoid $\mathcal{H}$ is wide if $\mathcal{H}_0 = \mathcal{G}_0$. It is normal if it is wide and for every $\gamma \in \mathcal{G}$ with $s(\gamma) = x, r(\gamma) = y$ one has $\gamma \mathcal{H}_x \gamma^{-1} \subseteq \mathcal{H}_y$, where $\mathcal{H}_x := \{\eta \in \mathcal{H} \mid r(\eta) = s(\eta) = x\}$. Since the same inclusion holds with $\gamma$ replaced by $\gamma^{-1}$, normality is equivalent to $\gamma \mathcal{H}_x \gamma^{-1} = \mathcal{H}_y$ for all such $\gamma$.

**Proposition 1.3.5** (Quotient groupoid [2, §11.3.1])**.** Let $\mathcal{N} \subseteq \mathcal{G}$ be a wide normal subgroupoid whose morphisms are isotropy elements, that is $\{\eta \in \mathcal{N} \mid r(\eta) = x,\ s(\eta) = y\} = \emptyset$ for $x \neq y$. The quotient groupoid $\mathcal{G}/\mathcal{N}$ is defined by:
(Q1) **Objects:** $(\mathcal{G}/\mathcal{N})_0 = \mathcal{G}_0$.
(Q2) **Morphisms:** for $\gamma \in \mathcal{G}$, write $\mathcal{N}_{r(\gamma)} \gamma \mathcal{N}_{s(\gamma)} := \{n \gamma n' \mid n \in \mathcal{N}_{r(\gamma)},\ n' \in \mathcal{N}_{s(\gamma)}\}$ and denote this double coset by $[\gamma]$. Then $(\mathcal{G}/\mathcal{N})(x, y) := \{[\gamma] \mid \gamma \in \mathcal{G}(x, y)\}$.
(Q3) **Structure maps:** $r([\gamma]) := r(\gamma)$, $s([\gamma]) := s(\gamma)$, $[g]^{-1} := [g^{-1}]$, and if $s(g) = r(h)$ then $[\gamma] \cdot [\eta] := [g \cdot h]$.

*Proof.* Define $g \sim g'$ if $g' \in \mathcal{N}_{r(g)} g \mathcal{N}_{s(g)}$. Then $\sim$ is an equivalence relation on $\mathcal{G}$. Reflexivity follows from $g = 1_{r(g)} g 1_{s(g)}$. If $g' = agb$ with $a \in \mathcal{N}_{r(g)}$ and $b \in \mathcal{N}_{s(g)}$, then $g = a^{-1} g' b^{-1}$, with $a^{-1} \in \mathcal{N}_{r(g)}$ and $b^{-1} \in \mathcal{N}_{s(g)}$; hence the relation $g' \in \mathcal{N}_{r(g)} g \mathcal{N}_{s(g)}$ is symmetric. If $g' = agb$ and $g'' = cg'd$ with $a \in \mathcal{N}_{r(g)}$, $b \in \mathcal{N}_{s(g)}$, $c \in \mathcal{N}_{r(g')}$, $d \in \mathcal{N}_{s(g')}$, then $r(g') = r(g)$ and $s(g') = s(g)$ because $a, b$ are isotropy at $r(g), s(g)$. Hence $c \in \mathcal{N}_{r(g)}$ and $d \in \mathcal{N}_{s(g)}$. Since each $\mathcal{N}_x$ is a subgroup of $\mathcal{G}_x$, we have $ca \in \mathcal{N}_{r(g)}$ and $bd \in \mathcal{N}_{s(g)}$, so $g'' = cg'd = c(agb)d = (ca)g(bd) \in \mathcal{N}_{r(g)} g \mathcal{N}_{s(g)}$. Thus $g \sim g''$ and the transitivity follows.

- **Source, range, and inversion are well defined.** If $g' = agb$ with $a \in \mathcal{N}_{r(g)}$, $b \in \mathcal{N}_{s(g)}$, then $r(g') = r(g)$ and $s(g') = s(g)$. Thus $r([g])$ and $s([g])$ are well defined.
  Moreover $g'^{-1} = b^{-1} g^{-1} a^{-1}$ with $b^{-1} \in \mathcal{N}_{s(g)}$, $a^{-1} \in \mathcal{N}_{r(g)}$, hence $[g'^{-1}] = [g^{-1}]$.

- **Composition is well defined.** Assume $s(g) = r(h)$, and choose representatives $g_1 = agb$, $h_1 = chd$ with $a \in \mathcal{N}_{r(g)}$, $b \in \mathcal{N}_{s(g)} = \mathcal{N}_{r(h)}$, $c \in \mathcal{N}_{r(h)}$, $d \in \mathcal{N}_{s(h)}$. Then $g_1 \cdot h_1$ is defined since $s(g_1) = s(g) = r(h) = r(h_1)$. Set $n := bc \in \mathcal{N}_{s(g)}$. Then $g_1 \cdot h_1 = ag(bc)hd = agnhd$. By normality of $\mathcal{N}$, $gng^{-1} \in \mathcal{N}_{r(g)}$. Set $a' := a(gng^{-1}) \in \mathcal{N}_{r(g)}$. Then $a'(g \cdot h)d = a(gng^{-1})(gh)d = agnhd = g_1 \cdot h_1$, so $[g_1 \cdot h_1] = [g \cdot h]$.
- **Units and associativity.** For $x \in \mathcal{G}_0$, $[1_x]$ is a two-sided unit $[1_{r(g)}] \cdot [g] = [1_{r(g)} \cdot g] = [g]$ and $[g] \cdot [1_{s(g)}] = [g \cdot 1_{s(g)}] = [g]$. Associativity follows from associativity in $\mathcal{G}$ $([g] \cdot [h]) \cdot [k] = [(g \cdot h) \cdot k] = [g \cdot (h \cdot k)] = [g] \cdot ([h] \cdot [k])$.

□

*Remark* 1.3.6. If $\mathcal{G}$ is a topological groupoid and $\mathcal{N}$ is a closed normal wide subgroupoid, endow $\mathcal{G}/\mathcal{N}$ with the quotient topology making the functor $\mathcal{G} \to \mathcal{G}/\mathcal{N}$ continuous.

**Setting 1.3.7.** For units $x, y \in \mathcal{G}_0$ we set $\mathcal{G}(x,y) := \{\gamma \in \mathcal{G} \mid r(\gamma) = x,\ s(\gamma) = y\}$, the set of morphisms from $y$ to $x$. Composition restricts to $\mathcal{G}(x,y) \times \mathcal{G}(y,z) \to \mathcal{G}(x,z)$, $(\alpha, \beta) \mapsto \alpha \cdot \beta$. $\mathcal{G}(x,y)$ and $\mathcal{N}_x$ carry the subspace topology from $\mathcal{G}$.

The following examples form the basic operations for constructing and analysing groupoids. The pair groupoid $X \times X$ encodes the space itself. The transformation groupoid $\Gamma \ltimes X$ models orbit spaces of actions while remaining étale for discrete $\Gamma$. Reductions $\mathcal{G}|_U$ implement localisation along units, and quotients by wide normal isotropy subgroupoids collapse internal symmetries. Such quotients are needed for functoriality, Morita invariance, and homology.

**Example 1.3.8.**

- **Topological groupoid.** Let $X$ be a topological space. Set $\mathcal{G} = X \times X$, $\mathcal{G}_0 = X$, $r(x,y) = x$, $s(x,y) = y$, inversion $(x,y)^{-1} = (y,x)$, and composition $(x,y) \cdot (y,z) = (x,z)$. All structure maps are continuous, so $\mathcal{G}$ is a topological groupoid. It is étale iff $X$ is discrete. Otherwise the projections are not local homeomorphisms. If $X$ is discrete, then for each $(x,y) \in \mathcal{G}$ the singleton $\{(x,y)\}$ is open in $X \times X$, and $r|_{\{(x,y)\}} : \{(x,y)\} \to \{x\}$, $s|_{\{(x,y)\}} : \{(x,y)\} \to \{y\}$, are homeomorphisms onto open subsets of $X$. Hence $r$ and $s$ are local homeomorphisms, so $\mathcal{G}$ is étale. Conversely, suppose $\mathcal{G}$ is étale. Then in particular $r$ is a local homeomorphism. Fix $x \in X$ and consider the arrow $(x,x) \in \mathcal{G}$. By local homeomorphy of $r$ there exist open neighbourhoods $U \subseteq X \times X$ of $(x,x)$ and $V \subseteq X$ of $x$ such that $r|_U : U \to V$ is a homeomorphism. Since $r$ is the projection onto the first coordinate, the restriction $r|_U$ is bijective, so $U$ is the graph of a continuous map $f : V \to X$: $U = \{(v, f(v)) \mid v \in V\}$. Because $U$ is open in $X \times X$ and $(x,x) \in U$, there exist open neighbourhoods $W_1, W_2 \subseteq X$ of $x$ such that $W_1 \times W_2 \subseteq U$. Let $w \in W_2$. Then for every $v \in W_1$ we have $(v,w) \in W_1 \times W_2 \subseteq U$, so $(v,w) = (v, f(v))$ and hence $w = f(v)$. Thus $w$ does not depend on $v$, and in particular for $v = x$ we get $w = f(x)$. Since also $(x,x) \in U$, we have $f(x) = x$, so $w = x$. Therefore every $w \in W_2$ equals $x$, and hence $W_2 = \{x\}$ is open in $X$. As $x \in X$ was arbitrary, every singleton $\{x\}$ is open, so $X$ is discrete.
- **Étale groupoid.** Let a discrete group $\Gamma$ act by homeomorphisms on a locally compact Hausdorff space $X$. Set $\mathcal{G} = \Gamma \ltimes X$, $\mathcal{G}_0 = X$, $r(g,x) = g \cdot x$, $s(g,x) = x$, inversion $(g,x)^{-1} =$ $(g^{-1}, g \cdot x)$, and composition $(h, g \cdot x) \cdot (g, x) = (hg, x)$. Since $\Gamma$ is discrete, $r$ and $s$ are local homeomorphisms. Hence $\Gamma \ltimes X$ is étale.

- **Subgroupoid.** Let $\mathcal{G}$ be a topological groupoid and $U \subseteq \mathcal{G}_0$. Define the reduction

$$\mathcal{G}|_U := \{\gamma \in \mathcal{G} \mid r(\gamma) \in U, s(\gamma) \in U\}$$

  with unit space $U$ and the induced structure maps $r|_{\mathcal{G}|_U}$, $s|_{\mathcal{G}|_U}$, $\cdot|_{(\mathcal{G}|_U)_2}$, and ${}^{-1}|_{\mathcal{G}|_U}$. Then $\mathcal{G}|_U$ is a subgroupoid of $\mathcal{G}$. If $U$ is open or closed, it is an open or closed topological subgroupoid.
- **Quotient groupoid.** View a topological group $G$ as a one-object groupoid with $(\mathcal{G}, \mathcal{G}_0) = (G, \{e\})$. Let $N \trianglelefteq G$ be a closed normal subgroup. Then $N \subseteq G$ is a wide normal subgroupoid whose morphisms are isotropy elements, and the quotient groupoid $\mathcal{G}/N$ is the one-object groupoid $(G/N, \{[e]\})$ with the usual quotient topology and structure maps. Here $[g] \cdot [h] = [gh]$, $[g]^{-1} = [g^{-1}]$, and $r([g]) = s([g]) = [e]$.

*Remark* 1.3.9. For a left group action $\Gamma \triangleright X$ by homeomorphisms, $\Gamma \ltimes X$ is the transformation groupoid. The object is $X$ and the morphisms are $\Gamma \times X$, with structure maps as in Example 1.3.8.

One might expect the unit space to be closed in a topological groupoid. In general this fails. However, it is closed if the arrow space is Hausdorff, but the converse needs not to hold. The converse implication fails for general groupoids in **Top** if one does not assume that $\mathcal{G}_0$ is Hausdorff: the unit groupoid on a non-Hausdorff space $X$ satisfies $\mathcal{G} = \mathcal{G}_0 = X$, so $\mathcal{G}_0$ is closed in $\mathcal{G}$ but $\mathcal{G}$ is not Hausdorff. This example is excluded by the convention in [21], where a topological groupoid is required to have Hausdorff unit space.

**Proposition 1.3.10.** If $\mathcal{G}$ is Hausdorff, then $\mathcal{G}_0$ is closed in $\mathcal{G}$.

*Proof.* Consider the continuous map $H : \mathcal{G} \to \mathcal{G} \times \mathcal{G}$ given by $H(\gamma) := (\gamma, r(\gamma))$. Let $\Delta_{\mathcal{G}} := \{(\eta, \eta) \mid \eta \in \mathcal{G}\}$ be the diagonal. Since $\mathcal{G}$ is Hausdorff, $\Delta_{\mathcal{G}}$ is closed in $\mathcal{G} \times \mathcal{G}$. An arrow $\gamma$ is a unit if and only if $\gamma = r(\gamma)$, hence $\mathcal{G}_0 = H^{-1}(\Delta_{\mathcal{G}})$. Therefore $\mathcal{G}_0$ is closed in $\mathcal{G}$. □

**Proposition 1.3.11.** Let $\mathcal{G}$ be a topological groupoid such that $\mathcal{G}_0 \subseteq \mathcal{G}$ is Hausdorff in the relative topology. Then $\mathcal{G}$ is Hausdorff if and only if $\mathcal{G}_0$ is closed in $\mathcal{G}$.

*Proof.*

- **Assume that $\mathcal{G}$ is Hausdorff.** Let $(x_i)_{i \in I}$ be a net in $\mathcal{G}_0$ such that $x_i \to \gamma \in \mathcal{G}$. Since $r$ is continuous and $r(x_i) = x_i$ for all $i$, we have $x_i = r(x_i) \to r(\gamma) \in \mathcal{G}_0$. By uniqueness of limits in a Hausdorff space, $\gamma = r(\gamma)$, hence $\gamma \in \mathcal{G}_0$. Thus $\mathcal{G}_0$ is closed.
- **Assume that $\mathcal{G}_0$ is closed in $\mathcal{G}$.** To show that $\mathcal{G}$ is Hausdorff, it suffices to show that every convergent net has a unique limit. Let $(\gamma_i)_{i \in I}$ be a net in $\mathcal{G}$ with $\gamma_i \to \alpha$ and $\gamma_i \to \beta$. First, continuity of $r : \mathcal{G} \to \mathcal{G}_0$ gives $r(\gamma_i) \to r(\alpha)$ and $r(\gamma_i) \to r(\beta)$ in the Hausdorff space $\mathcal{G}_0$, so $r(\alpha) = r(\beta)$. Hence $(\alpha^{-1}, \beta) \in \mathcal{G}_2$. Next, since inversion is continuous, $\gamma_i^{-1} \to \alpha^{-1}$. Therefore $(\gamma_i^{-1}, \gamma_i) \to (\alpha^{-1}, \beta)$ in $\mathcal{G} \times \mathcal{G}$, and hence in the subspace $\mathcal{G}_2$. By continuity of multiplication, $\gamma_i^{-1}\gamma_i \to \alpha^{-1}\beta$ in $\mathcal{G}$. But $\gamma_i^{-1}\gamma_i = s(\gamma_i) \in \mathcal{G}_0$ for all $i$. Since $\mathcal{G}_0$ is closed in $\mathcal{G}$, the limit $\alpha^{-1}\beta$ belongs to $\mathcal{G}_0$. Thus $\alpha^{-1}\beta$ is a unit, and hence $\beta = \alpha(\alpha^{-1}\beta) = \alpha$. So limits are unique, and $\mathcal{G}$ is Hausdorff.

□

*Remark* 1.3.12. In [21, Lemma 2.3.2] the statement is proved under the standing hypotheses of [21, Definition 2.3.1], in particular assuming local compactness. However, those uses are inessential for the implication itself. Indeed, the argument only requires that the structure maps $s, r, m, i, u$ are continuous and that $\mathcal{G}_0$ is Hausdorff in the relative topology.

## 1.4 Étale Groupoids

Étale groupoids occupy the middle ground between the purely topological and the smooth settings. They retain enough local regularity to support constructions such as bisections, étale sheaves, nerves and counting measures, while remaining flexible enough to model orbit and equivalence–relation dynamics arising from local homeomorphisms and actions of discrete groups. A groupoid can be summarised by the diagram [4, §1]:

$$\mathcal{G}_2 = \mathcal{G}\ {}_s\!\times_r \mathcal{G} \xrightarrow{m} \mathcal{G} \xrightarrow{i} \mathcal{G} \overset{r}{\underset{s}{\rightrightarrows}} \mathcal{G}_0 \xrightarrow{u} \mathcal{G},$$

where $s, r : \mathcal{G} \to \mathcal{G}_0$ are the source and range maps, $u(x) := 1_x$ inserts the identity morphism at $x \in \mathcal{G}_0$, $i(g) := g^{-1}$ is inversion, and $m(g, h) := g \cdot h$ is defined exactly when $s(g) = r(h)$, that is, $(g, h) \in \mathcal{G}_2$. We write $g : x \to y$ for $s(g) = x$ and $r(g) = y$. A topological groupoid is obtained by equipping $\mathcal{G}_0$ and $\mathcal{G}$ with topologies and requiring all structure maps in the diagram above to be continuous. In the smooth (Lie) case, $\mathcal{G}_0$ and $\mathcal{G}$ are smooth manifolds, the structure maps are smooth, and one assumes that $s$ and $r$ are submersions so that the fibre product $\mathcal{G}\ {}_s\!\times_r \mathcal{G}$ is again a manifold [4, §1]. Étaleness is the topological analogue of this local regularity: we require $s$ and $r$ to be local homeomorphisms. This forces the fibres of $s$ and $r$ to be discrete, aligns étale groupoids with transformation groupoids of actions of discrete groups and equivalence–relation groupoids of local homeomorphisms, and at the same time preserves enough categorical structure for the homological constructions in Chapter 2.

**Definition 1.4.1** (Étale groupoids [4, Definition 1.1])**.** Let $\mathcal{G}$ be a topological groupoid. The groupoid $\mathcal{G}$ is called étale if the source map $s : \mathcal{G} \to \mathcal{G}_0$ is a local homeomorphism. In this case the range map $r = s \circ i$ is also a local homeomorphism.

*Remark* 1.4.2. Since inversion $i$ is a homeomorphism and $r = s \circ i$, this is equivalent to requiring that the range map $r : \mathcal{G} \to \mathcal{G}_0$ is a local homeomorphism. Many authors state the condition as: $r$ and $s$ are local homeomorphisms. In particular, for all $g \in \mathcal{G}$ there is an open neighbourhood $U \subset \mathcal{G}$ with $g \in U$ such that $s|_U : U \xrightarrow{\cong} s(U)$ and $r|_U : U \xrightarrow{\cong} r(U)$ are homeomorphisms.

Such a set $U$ is called a local bisection. In the smooth or Lie setting the same definition applies with local diffeomorphisms in place of local homeomorphisms.

**Lemma 1.4.3** (Discrete fibers of local homeomorphisms)**.** Let $p : X \to Y$ be a local homeomorphism. Then for every $y \in Y$ the fiber $p^{-1}(y)$ is a discrete subspace of $X$. In particular, if $\mathcal{G}$ is an étale groupoid, then for every $x \in \mathcal{G}_0$ the fibers $r^{-1}(x)$ and $s^{-1}(x)$ are discrete subspaces of $\mathcal{G}$.

*Proof.* Fix $y \in Y$. To show that $p^{-1}(y)$ is discrete, it suffices to prove that every singleton $\{x\}$ with $x \in p^{-1}(y)$ is open in the subspace topology on $p^{-1}(y)$.

Fix $x \in p^{-1}(y)$. Since $p$ is a local homeomorphism, there exists an open neighbourhood $U \subseteq X$ of $x$ such that the restriction $p|_U : U \to p(U)$ is a homeomorphism onto an open subset of $Y$. Set $V := U \cap p^{-1}(y)$. Then $V$ is open in the subspace $p^{-1}(y)$. Moreover, $V = \{x\}$. Indeed, if $x' \in V$, then $x' \in U$ and $p(x') = y = p(x)$, hence $x' = x$ because $p|_U$ is injective. Thus $\{x\}$ is open in $p^{-1}(y)$. Since $x$ was arbitrary, the fiber $p^{-1}(y)$ is discrete.

If $\mathcal{G}$ is étale, then by definition $r$ and $s$ are local homeomorphisms. Applying the first part to $p = r$ and $p = s$ gives that $r^{-1}(x)$ and $s^{-1}(x)$ are discrete for every $x \in \mathcal{G}_0$. □

Direct limits assemble compatible local data around a point. Two local maps represent the same element once they agree on some smaller neighbourhood.

**Definition 1.4.4** (Direct limit)**.** Fix $x \in \mathcal{G}_0$ and let $\mathcal{N}(x)$ denote the directed set of open neighbourhoods of $x$, ordered by reverse inclusion. For each $U \in \mathcal{N}(x)$ set $\mathrm{LH}(U, \mathcal{G}_0) := \{\varphi : U \to \mathcal{G}_0 \mid \varphi \text{ is a local homeomorphism}\}$. Whenever $U' \subseteq U$ we have the restriction map $\rho_{U,U'} : \mathrm{LH}(U, \mathcal{G}_0) \to \mathrm{LH}(U', \mathcal{G}_0)$, $\varphi \mapsto \varphi|_{U'}$. This yields a direct system $(\mathrm{LH}(U, \mathcal{G}_0), \rho_{U,U'})_{U \in \mathcal{N}(x)}$ in the category of sets. Its direct limit is the quotient

$$\varinjlim_{U \in \mathcal{N}(x)} \mathrm{LH}(U, \mathcal{G}_0) \cong \Big( \bigsqcup_{U \in \mathcal{N}(x)} \mathrm{LH}(U, \mathcal{G}_0) \Big) \Big/ \sim,$$

where $(U, \varphi) \sim (U', \varphi')$ if there exists an open set $W \subseteq U \cap U'$ with $x \in W$ and $\varphi|_W = \varphi'|_W$.

More generally, a direct system in a category $\mathcal{C}$ consists of a directed set $(I, \leq)$, a family of objects $(A_i)_{i \in I}$, and morphisms $\phi_{ij} : A_i \to A_j$ for all $i \leq j$ such that $\phi_{ii} = \mathrm{id}_{A_i}$ and $\phi_{jk} \circ \phi_{ij} = \phi_{ik}$ whenever $i \leq j \leq k$. In particular, the family $(\mathrm{LH}(U, \mathcal{G}_0))_{U \in \mathcal{N}(x)}$, together with the restriction maps $\rho_{U,U'} : \mathrm{LH}(U, \mathcal{G}_0) \to \mathrm{LH}(U', \mathcal{G}_0)$ for $U' \subseteq U$, is a direct system in **Set** indexed by $\mathcal{N}(x)$. In an étale groupoid, every arrow sufficiently close to a unit arises from a local bisection, so its local effect is recorded by the germ of a local homeomorphism of the unit space, that is, by an element of this direct limit.

**Definition 1.4.5.** Let $\mathcal{G}$ be an étale groupoid and let $g \in \mathcal{G}$ with $s(g) = x$ and $r(g) = y$. Choose an open bisection $B \subseteq \mathcal{G}$ with $g \in B$, and put $U := s(B) \subseteq \mathcal{G}_0$. Let $\sigma : U \to \mathcal{G}$ be the inverse of $s|_B : B \to U$, so that $\sigma(x) = g$ and $s \circ \sigma = \mathrm{id}_U$. Set $V := r(B)$; then $V$ is open in $\mathcal{G}_0$ and $\varphi_{(U,\sigma)} := r \circ \sigma : U \to V$ is a homeomorphism (hence a local homeomorphism) sending $x$ to $y$.

Let $\Sigma_g$ be the set of all such pairs $(U, \sigma)$ arising from open bisections $B \ni g$. Declare $(U, \sigma) \sim (U', \sigma')$ if and only if there exists an open neighbourhood $W \subseteq U \cap U'$ of $x$ such that $\varphi_{(U,\sigma)}|_W = \varphi_{(U',\sigma')}|_W$. The germ of $g$ at $x$ is the equivalence class

$$\tilde{g} := [(U, \sigma)] \in \varinjlim_{W \in \mathcal{N}(x)} \mathrm{LH}(W, \mathcal{G}_0).$$

In the smooth setting replace local homeomorphisms by local diffeomorphisms.

*Remark* 1.4.6. This construction is well defined and depends only on $g$. Since $s$ is a local homeomorphism at $g$, there exist an open neighbourhood $N \subseteq \mathcal{G}$ of $g$ and an open set $W_0 \subseteq \mathcal{G}_0$ such that $s|_N : N \xrightarrow{\cong} W_0$. For any $(U,\sigma),(U',\sigma') \in \Sigma_g$ we may shrink $U$ and $U'$ so that $\sigma(U) \subseteq N$ and $\sigma'(U') \subseteq N$. On $W := U \cap U' \cap W_0$ we then have $\sigma|_W = (s|_N)^{-1}|_W = \sigma'|_W$, hence $(U,\sigma) \sim (U',\sigma')$ and $\varphi_{(U,\sigma)}|_W = \varphi_{(U',\sigma')}|_W$. Thus the germ $\tilde{g}$ is independent of the choice of $(U,\sigma) \in \Sigma_g$.

*Remark* 1.4.7 (Basic properties [4, §1.2]).

(a) Since $\mathcal{G}$ is étale, $s$ is a local homeomorphism. After shrinking $U$, there exists a local section $\sigma : U \to \mathcal{G}$ of $s$ with $\sigma(x) = g$ and $s \circ \sigma = \mathrm{id}_U$. The inversion map $\mathcal{G} \to \mathcal{G}$, $h \mapsto h^{-1}$, is a homeomorphism and $r = s \circ i_{\mathcal{G}}$, so $r$ is a local homeomorphism as well. In particular, $r|_{\sigma(U)} : \sigma(U) \to V$ is a homeomorphism onto the open set $V := r(\sigma(U)) \subseteq \mathcal{G}_0$. Consequently the representative $\varphi_{(U,\sigma)} = r \circ \sigma : U \to V$ of $\tilde{g}$ is a local homeomorphism.

(b) For the unit $1_x \in \mathcal{G}$ we may take $U \subseteq \mathcal{G}_0$ and $\sigma(u) := 1_u$ for $u \in U$. Then $s \circ \sigma = \mathrm{id}_U$ and $r \circ \sigma = \mathrm{id}_U$, so the germ $\widetilde{1_x}$ is represented by the identity map near $x$ and hence $\widetilde{1_x} = 1_x$.

(c) If $g : x \to y$ and $h : y \to z$, choose local sections $\sigma_g : U \to \mathcal{G}$ and $\sigma_h : V \to \mathcal{G}$ with $x \in U$, $y \in V$, $s \circ \sigma_g = \mathrm{id}_U$, $s \circ \sigma_h = \mathrm{id}_V$, $\sigma_g(x) = g$, and $\sigma_h(y) = h$. Shrinking $U$ if necessary, we may assume $\varphi_{(U,\sigma_g)}(U) \subseteq V$. On $U$ the map $\varphi_{(U,\sigma_g)} = r \circ \sigma_g$ represents $\tilde{g}$ and $\varphi_{(V,\sigma_h)} = r \circ \sigma_h$ represents $\tilde{h}$. The product section $U \to \mathcal{G}$, $u \mapsto \sigma_h(\varphi_{(U,\sigma_g)}(u)) \cdot \sigma_g(u)$, represents $hg$, and its range map is the composition $\varphi_{(V,\sigma_h)} \circ \varphi_{(U,\sigma_g)}$ near $x$. Hence the germ of $hg$ at $x$ is the composite germ $\tilde{h} \circ \tilde{g}$.

**Example 1.4.8.**

- **Spaces as étale groupoids [4, §1.3(1)].** Let $X$ be a topological space. Consider the groupoid $(\mathcal{G}, \mathcal{G}_0, r_{\mathcal{G}}, s_{\mathcal{G}}, \cdot_{\mathcal{G}}, i_{\mathcal{G}})$ given by $\mathcal{G} := X$, $\mathcal{G}_0 := X$, $r_{\mathcal{G}} := \mathrm{id}_X$, $s_{\mathcal{G}} := \mathrm{id}_X$, $i_{\mathcal{G}}(g) := g$, and $g \cdot_{\mathcal{G}} h := g$ whenever $g, h \in X$ and $g = h$. Then $(\mathcal{G}, \mathcal{G}_0, r_{\mathcal{G}}, s_{\mathcal{G}}, \cdot_{\mathcal{G}}, i_{\mathcal{G}})$ is a topological groupoid, and since $r_{\mathcal{G}}$ and $s_{\mathcal{G}}$ are homeomorphisms, $\mathcal{G}$ is étale.
- **Translation groupoid of a discrete group action [4, §1.3(2)].** Let a discrete group $\Gamma$ act on a locally compact Hausdorff space $X$ by homeomorphisms via a left action $\Gamma \times X \to X$, $(g,x) \mapsto g \cdot x$. The associated transformation groupoid $\Gamma \ltimes X$ is the étale groupoid with $(\Gamma \ltimes X)_0 := X$, source and range maps $s_{\Gamma \ltimes X}(g,x) := x$, $r_{\Gamma \ltimes X}(g,x) := g \cdot x$, unit map $u_{\Gamma \ltimes X}(x) := (e,x)$, inverse map $(g,x)^{-1} := (g^{-1}, g \cdot x)$, and multiplication $(h,y) \cdot (g,x) := (hg,x)$ whenever $s_{\Gamma \ltimes X}(h,y) = r_{\Gamma \ltimes X}(g,x)$, that is whenever $y = g \cdot x$, so that explicitly $(h, g \cdot x) \cdot (g,x) = (hg,x)$. Since $\Gamma$ is discrete, sets of the form $\{g\} \times U \subseteq \Gamma \times X$ with $U \subseteq X$ open form a basis. On $\{g\} \times U$ the source map is the homeomorphism $(g,x) \mapsto x$, hence $s_{\Gamma \ltimes X}$ is a local homeomorphism. Likewise, $r_{\Gamma \ltimes X}$ restricts to the homeomorphism $(g,x) \mapsto g \cdot x$ from $\{g\} \times U$ onto $g \cdot U$. Therefore $r_{\Gamma \ltimes X}$ is a local homeomorphism as well, and $\Gamma \ltimes X$ is étale. For a right action $(x,\gamma) \mapsto x \cdot \gamma$ one can also use the arrow space $X \times \Gamma$ and define $r_{X \rtimes \Gamma}(x,\gamma) := x$, $s_{X \rtimes \Gamma}(x,\gamma) := x \cdot \gamma$, with $u_{X \rtimes \Gamma}(x) = (x,e)$, $(x,\gamma)^{-1} = (x \cdot \gamma, \gamma^{-1})$, and $(x,\gamma) \cdot (x \cdot \gamma, \eta) = (x, \gamma\eta)$.
- **Action groupoid of a right $\mathcal{G}$-space [4, §1.3(6)].** Let $(\mathcal{G}, \mathcal{G}_0, r, s, \cdot, i)$ be an étale groupoid with unit map $u : \mathcal{G}_0 \to \mathcal{G}$, $x \mapsto 1_x$. A right $\mathcal{G}$-space consists of a topological space $X$, a continuous anchor map $p : X \to \mathcal{G}_0$, and a continuous action map $\cdot : X\ {}_p\!\times_r \mathcal{G} \to X$, $(x,g) \mapsto x \cdot g$, defined on the fibre product $X\ {}_p\!\times_r \mathcal{G} = \{(x,g) \in X \times \mathcal{G} \mid p(x) = r(g)\}$, such that for all $(x,g) \in X\ {}_p\!\times_r \mathcal{G}$ one has $p(x \cdot g) = s(g)$ and $x \cdot 1_{p(x)} = x$, and whenever $(g,h) \in \mathcal{G}_2$ and $p(x) = r(g)$ one has $(x \cdot g) \cdot h = x \cdot (gh)$.

The action groupoid $X \rtimes \mathcal{G}$ is the topological groupoid with objects $(X \rtimes \mathcal{G})_0 = X$ and arrows $(X \rtimes \mathcal{G})_1 = X\ {}_p\!\times_r \mathcal{G}$ with the subspace topology from $X \times \mathcal{G}$, and structure maps $r_{X\rtimes\mathcal{G}}(x,\gamma) = x$, $s_{X\rtimes\mathcal{G}}(x,\gamma) = x \cdot \gamma$, $u_{X\rtimes\mathcal{G}}(x) = (x, 1_{p(x)})$, $(x,\gamma)^{-1} = (x \cdot \gamma, \gamma^{-1})$, and multiplication $(x,\gamma) \cdot (x \cdot \gamma, h) = (x, \gamma\eta)$, defined whenever $(\gamma,\eta) \in \mathcal{G}_2$. The space of composable pairs is $(X \rtimes \mathcal{G})_2 = \{((x,\gamma),(y,\eta)) \in (X \rtimes \mathcal{G})_1^2 \mid y = x \cdot \gamma,\ (\gamma,\eta) \in \mathcal{G}_2\}$.

Fix $(x,g) \in (X \rtimes \mathcal{G})_1$. Choose an open bisection $U \subseteq \mathcal{G}$ with $g \in U$, so $r|_U : U \to r(U)$ and $s|_U : U \to s(U)$ are homeomorphisms onto open subsets of $\mathcal{G}_0$. Choose an open neighbourhood $W \subseteq X$ of $x$ with $p(W) \subseteq r(U)$. Then $W\ {}_p\!\times_r U = \{(x',g') \in W \times U \mid p(x') = r(g')\}$ is an open neighbourhood of $(x,g)$ in $(X \rtimes \mathcal{G})_1$.

On $W\ {}_p\!\times_r U$, the restriction of $r_{X\rtimes\mathcal{G}}$ is a homeomorphism onto $W$, with inverse $W \to W\ {}_p\!\times_r U$, $x' \mapsto (x', (r|_U)^{-1}(p(x')))$. Moreover, the restriction of $s_{X\rtimes\mathcal{G}}$ to $W\ {}_p\!\times_r U$ is a homeomorphism onto its image. Indeed, for $y \in s_{X\rtimes\mathcal{G}}(W\ {}_p\!\times_r U)$ we have $p(y) \in s(U)$, hence $k(y) := (s|_U)^{-1}(p(y)) \in U$ is well defined, and $y \mapsto (y \cdot k(y)^{-1}, k(y))$ is an inverse to $s_{X\rtimes\mathcal{G}}|_{W_p\times_r U}$ since $(y \cdot k(y)^{-1}) \cdot k(y) = y$ and $p(y \cdot k(y)^{-1}) = r(k(y))$. Therefore $r_{X\rtimes\mathcal{G}}$ and $s_{X\rtimes\mathcal{G}}$ are local homeomorphisms, and $X \rtimes \mathcal{G}$ is étale.

**Definition 1.4.9** ([4, §1.4])**.** Let $(\mathcal{K}, \mathcal{K}_0, r_\mathcal{K}, s_\mathcal{K}, \cdot_\mathcal{K}, i_\mathcal{K})$ and $(\mathcal{G}, \mathcal{G}_0, r_\mathcal{G}, s_\mathcal{G}, \cdot_\mathcal{G}, i_\mathcal{G})$ be étale groupoids. A homomorphism of étale groupoids $\varphi : \mathcal{K} \to \mathcal{G}$ is a pair of continuous maps $\varphi_0 : \mathcal{K}_0 \to \mathcal{G}_0$ and $\varphi_1 : \mathcal{K} \to \mathcal{G}$ such that the following identities hold:

$$
\begin{aligned}
&r_\mathcal{G} \circ \varphi_1 = \varphi_0 \circ r_\mathcal{K}, && s_\mathcal{G} \circ \varphi_1 = \varphi_0 \circ s_\mathcal{K},\\
&\varphi_1 \circ i_\mathcal{K} = i_\mathcal{G} \circ \varphi_1, && \varphi_1 \circ u_\mathcal{K} = u_\mathcal{G} \circ \varphi_0,\\
&\varphi_1(g \cdot_\mathcal{K} h) = \varphi_1(g) \cdot_\mathcal{G} \varphi_1(h) && \text{for all } (g,h) \in \mathcal{K} \times_{\mathcal{K}_0} \mathcal{K}.
\end{aligned}
$$

This is precisely the data of a functor of small categories $\mathcal{K} \to \mathcal{G}$.

Isomorphism of étale groupoids is usually too rigid for geometric purposes, since different groupoids can present the same quotient or orbit data. Morita equivalence isolates the underlying geometric object by requiring a homomorphism $\varphi : \mathcal{K} \to \mathcal{G}$ to behave like an equivalence of categories in a way compatible with the topology. Recall that a functor between small categories is an equivalence if and only if it is essentially surjective on objects and fully faithful on morphisms. In the étale setting these two conditions are reformulated in terms of local homeomorphisms of the unit spaces and of the bibundle of arrows, leading to the notion of Morita equivalence for étale groupoids.

**Definition 1.4.10** (Morita equivalence [4, §1.5])**.** Let $\mathcal{K}$ and $\mathcal{G}$ be étale groupoids and let $\varphi : \mathcal{K} \to \mathcal{G}$ be a homomorphism with components $\varphi_0 : \mathcal{K}_0 \to \mathcal{G}_0$ and $\varphi_1 : \mathcal{K}_1 \to \mathcal{G}_1$. We say that $\varphi$ is a Morita equivalence, also known as a weak or essential equivalence, if:

1. **Essential surjectivity (étale surjection).** Consider the fibre product

$$\mathcal{K}_0 \ {}_{\varphi_0}\!\times_{r_\mathcal{G}} \mathcal{G}_1 = \{(y,g) \in \mathcal{K}_0 \times \mathcal{G}_1 \mid \varphi_0(y) = r_\mathcal{G}(g)\},$$

with projections $\pi_1(y,g) = y$ and $\pi_2(y,g) = g$. The map

$$s_\mathcal{G} \circ \pi_2 : \mathcal{K}_0 \ {}_{\varphi_0}\!\times_{r_\mathcal{G}} \mathcal{G}_1 \to \mathcal{G}_0, \qquad (y,g) \mapsto s_\mathcal{G}(g),$$

is a surjective local homeomorphism.

2. **Full faithfulness (pullback square).** The canonical square

$$\begin{array}{ccc} \mathcal{K}_1 & \xrightarrow{\ \varphi_1\ } & \mathcal{G}_1 \\ {\scriptstyle (r_\mathcal{K}, s_\mathcal{K})}\big\downarrow & & \big\downarrow{\scriptstyle (r_\mathcal{G}, s_\mathcal{G})} \\ \mathcal{K}_0 \times \mathcal{K}_0 & \xrightarrow[\varphi_0 \times \varphi_0]{} & \mathcal{G}_0 \times \mathcal{G}_0 \end{array}$$

is a pullback in the sense that the natural map

$$\mathcal{K}_1 \to (\mathcal{K}_0 \times \mathcal{K}_0) \ {}_{\varphi_0 \times \varphi_0}\!\times_{r_\mathcal{G} \times s_\mathcal{G}} \mathcal{G}_1, \qquad k \mapsto ((r_\mathcal{K}(k), s_\mathcal{K}(k)), \varphi_1(k)),$$

is a homeomorphism, where the fibre product is

$$(\mathcal{K}_0 \times \mathcal{K}_0) \ {}_{\varphi_0 \times \varphi_0}\!\times_{r_\mathcal{G} \times s_\mathcal{G}} \mathcal{G}_1 = \{((y_1,y_2),g) \in (\mathcal{K}_0 \times \mathcal{K}_0) \times \mathcal{G}_1 \mid (\varphi_0(y_1), \varphi_0(y_2)) = (r_\mathcal{G}(g), s_\mathcal{G}(g))\}.$$

*Remark* 1.4.11. The essential surjectivity condition says that every object of $\mathcal{G}$ is locally in the image of $s_\mathcal{G} \circ \pi_2$, and hence admits local lifts along $s_\mathcal{G} \circ \pi_2$ after restricting to suitable neighbourhoods in $\mathcal{G}_0$. The pullback condition identifies $\mathcal{K}_1$ with the space of arrows of $\mathcal{G}$ between points in the image of $\mathcal{K}_0$, so that arrows of $\mathcal{G}$ over $\varphi_0$ correspond uniquely and continuously to arrows of $\mathcal{K}$. Together these conditions express that $\varphi$ is fully faithful and essentially surjective in a way compatible with the topologies on $\mathcal{K}$ and $\mathcal{G}$.

When these conditions hold we write $\varphi : \mathcal{K} \xrightarrow{\sim} \mathcal{G}$. Two étale groupoids $\mathcal{H}$ and $\mathcal{G}$ are Morita equivalent if there exist Morita equivalences $\mathcal{H} \xleftarrow{\sim} \mathcal{K} \xrightarrow{\sim} \mathcal{G}$. This generates an equivalence relation, and one often works in the localisation of the category of étale groupoids obtained by formally inverting all Morita equivalences. A morphism $\mathcal{H} \to \mathcal{G}$ there can be represented by a zig–zag $\mathcal{H} \xleftarrow{\sim} \mathcal{K} \xrightarrow{\sim} \mathcal{G}$ [4, §1.5].

**Definition 1.4.12** ([13, §2.1])**.** For an étale groupoid $\mathcal{G}$ and $Y \subseteq \mathcal{G}_0$ the reduction is $\mathcal{G}|_Y := r^{-1}(Y) \cap s^{-1}(Y)$, equipped with the subspace topology. It is an étale subgroupoid with unit space $Y$ whenever $Y$ is open.

**Definition 1.4.13** ([13, §2.1])**.** Let $(\mathcal{G}, \mathcal{G}_0, r_\mathcal{G}, s_\mathcal{G}, \cdot_\mathcal{G}, i_\mathcal{G})$ be an étale groupoid and write $r := r_\mathcal{G}$, $s := s_\mathcal{G}$. A subset $U \subseteq \mathcal{G}$ is a $\mathcal{G}$-bisection if the restrictions $r|_U : U \to r(U)$ and $s|_U : U \to s(U)$ are injective. If in addition $U$ is open, we call it an open $\mathcal{G}$-bisection. For subsets $U, V \subseteq \mathcal{G}$ set

$$U^{-1} := \{g \in \mathcal{G} \mid g^{-1} \in U\}, \qquad UV := \{gh \mid g \in U,\ h \in V,\ s(g) = r(h)\}.$$

If $U$ and $V$ are $\mathcal{G}$-bisections, then $U^{-1}$ and $UV$ are again $\mathcal{G}$-bisections. If $U$ and $V$ are open $\mathcal{G}$-bisections, then $U^{-1}$ and $UV$ are open $\mathcal{G}$-bisections.

*Remark* 1.4.14. In the étale case every arrow $g \in \mathcal{G}$ admits an open neighbourhood $U \subseteq \mathcal{G}$ which is a $\mathcal{G}$-bisection. On such a set, the restrictions $r|_U : U \to r(U)$ and $s|_U : U \to s(U)$ are homeomorphisms onto open subsets of $\mathcal{G}_0$ [21, §2.4].

**Lemma 1.4.15.** Let $\mathcal{G}$ be an étale groupoid. If the unit space $\mathcal{G}_0$ is discrete, then $\mathcal{G}$ is discrete, and for every $n \geq 0$ the space of $n$-composables $\mathcal{G}_n$ is discrete.

*Proof.* Assume $\mathcal{G}_0$ is discrete.

- **$\mathcal{G}$ is discrete.** Since $\mathcal{G}$ is étale, the range map $r : \mathcal{G} \to \mathcal{G}_0$ is a local homeomorphism. Fix $g \in \mathcal{G}$ and set $u := r(g) \in \mathcal{G}_0$. Because $\mathcal{G}_0$ is discrete, the singleton $\{u\}$ is open. As $r$ is a local homeomorphism, there exists an open neighbourhood $U \subset \mathcal{G}$ of $g$ such that $r|_U : U \to r(U)$ is a homeomorphism onto an open subset $r(U) \subset \mathcal{G}_0$. Then $r(U) \cap \{u\}$ is open in $r(U)$, hence $U' := (r|_U)^{-1}(r(U) \cap \{u\})$ is open in $U$, hence open in $\mathcal{G}$. Moreover $g \in U'$ and $r(U') = \{u\}$. Since $r|_U$ is injective, $U'$ contains at most one point, hence $U' = \{g\}$. Thus $\{g\}$ is open in $\mathcal{G}$, so $\mathcal{G}$ is discrete.
- **$\mathcal{G}^n$ is discrete for all $n \geq 0$.** For $n = 0$ we have $\mathcal{G}_0$ discrete by assumption. For $n \geq 1$, the product $\mathcal{G}^n$ is discrete because $\mathcal{G}$ is discrete.

The space $\mathcal{G}_n$ of composable $n$-tuples is a subspace of $\mathcal{G}^n$, hence discrete as well. □

**Definition 1.4.16.** A subset $X \subseteq \mathcal{G}_0$ is called full if $r(s^{-1}(X)) = \mathcal{G}_0$.

**Definition 1.4.17** (Kakutani equivalence [14, Definition 3.8])**.** Let $\mathcal{G}$ and $\mathcal{H}$ be étale groupoids. We say that $\mathcal{G}$ and $\mathcal{H}$ are Kakutani equivalent if there exist full clopen subsets $X \subseteq \mathcal{G}_0$ and $Y \subseteq \mathcal{H}_0$ such that the reductions $\mathcal{G}|_X$ and $\mathcal{H}|_Y$ are isomorphic as étale groupoids.

*Remark* 1.4.18. Let $\mathcal{G}$ be an étale groupoid and let $X \subseteq \mathcal{G}_0$ be clopen and full. The reduction $\mathcal{G}|_X$ has unit space $X$ and arrow space $\mathcal{G}|_X := \{g \in \mathcal{G} \mid r_{\mathcal{G}}(g) \in X,\ s_{\mathcal{G}}(g) \in X\}$. Since $X$ is clopen, the subsets $r_{\mathcal{G}}^{-1}(X)$ and $s_{\mathcal{G}}^{-1}(X)$ are open in $\mathcal{G}$, and hence $\mathcal{G}|_X = r_{\mathcal{G}}^{-1}(X) \cap s_{\mathcal{G}}^{-1}(X)$ is open in $\mathcal{G}$. The structure maps of $\mathcal{G}|_X$ are the restrictions of those of $\mathcal{G}$:

$$r_{\mathcal{G}|_X} = r_{\mathcal{G}}|_{\mathcal{G}|_X}, \quad s_{\mathcal{G}|_X} = s_{\mathcal{G}}|_{\mathcal{G}|_X}, \quad (g,h) \mapsto g \cdot_{\mathcal{G}} h, \quad g \mapsto g^{-1},$$

with composition defined whenever $s_{\mathcal{G}}(g) = r_{\mathcal{G}}(h) \in X$. Since $r_{\mathcal{G}}$ and $s_{\mathcal{G}}$ are local homeomorphisms, their restrictions $r_{\mathcal{G}|_X}$ and $s_{\mathcal{G}|_X}$ are local homeomorphisms as well. Thus $\mathcal{G}|_X$ is again an étale groupoid with unit space $X$. Fullness of $X$, that is $r_{\mathcal{G}}(s_{\mathcal{G}}^{-1}(X)) = \mathcal{G}_0$, means that every unit $x \in \mathcal{G}_0$ is the range of some arrow whose source lies in $X$. This means that every $\mathcal{G}$-orbit meets $X$. Hence the reduction $\mathcal{G}|_X$ still sees the entire orbit structure of $\mathcal{G}$. If $\mathcal{G}|_X \cong \mathcal{H}|_Y$ for full clopen $X \subseteq \mathcal{G}_0$ and $Y \subseteq \mathcal{H}_0$, then $\mathcal{G}$ and $\mathcal{H}$ have the same dynamics up to cutting down to representatives in $X$ and $Y$, which is what Kakutani equivalence is designed to capture.

**Lemma 1.4.19** ([21, Lemma 2.4.9])**.** If $\mathcal{G}$ is second countable, Hausdorff, and étale, then $\mathcal{G}_1$ admits a countable base consisting of open bisections.

*Proof.* Let $\mathcal{B}$ be a countable base for the topology on $\mathcal{G}_1$. Fix $B \in \mathcal{B}$. For every $\gamma \in B$, since $\mathcal{G}$ is étale there exists an open bisection $U_\gamma \subseteq \mathcal{G}_1$ with $\gamma \in U_\gamma \subseteq B$. Thus $\{U_\gamma\}_{\gamma \in B}$ is an open cover of $B$ by open bisections. Since $\mathcal{G}_1$ is second countable, it is Lindelöf, so there exists a countable subcover $\{U_{B,n}\}_{n\in\mathbb{N}}$ of $B$ by open bisections. Define

$$\mathcal{W} := \{U_{B,n} \mid B \in \mathcal{B},\ n \in \mathbb{N}\}.$$

Then $\mathcal{W}$ is countable. To see that $\mathcal{W}$ is a base, let $W \subseteq \mathcal{G}_1$ be open and $\gamma \in W$. Choose $B \in \mathcal{B}$ with $\gamma \in B \subseteq W$. Since $\{U_{B,n}\}_{n\in\mathbb{N}}$ covers $B$, there exists $n$ with $\gamma \in U_{B,n} \subseteq B \subseteq W$. Hence $\mathcal{W}$ is a countable base consisting of open bisections. □

**Corollary 1.4.20** ([21, Corollary 2.4.10])**.** If $\mathcal{G}$ is étale, then for every $x \in \mathcal{G}_0$ the sets $r^{-1}(x)$, $s^{-1}(x)$, and $\mathcal{G}_x := \{\gamma \in \mathcal{G} \mid r(\gamma) = s(\gamma) = x\}$ are discrete in the subspace topology.

*Proof.* Fix $x \in \mathcal{G}_0$ and $\gamma \in r^{-1}(x)$. Since $\mathcal{G}$ is étale, there exists an open bisection $U \subseteq \mathcal{G}$ with $\gamma \in U$. On $U$ the restriction $r|_U : U \to r(U)$ is injective, so $U \cap r^{-1}(x)$ contains at most one point. Since $r(\gamma) = x$, we have $U \cap r^{-1}(x) = \{\gamma\}$, and $\{\gamma\}$ is open in the subspace $r^{-1}(x)$. Thus $r^{-1}(x)$ is discrete. The same argument applied to $s$ shows that $s^{-1}(x)$ is discrete. Now let $\gamma \in \mathcal{G}_x = r^{-1}(x) \cap s^{-1}(x)$. From the first two steps there exist open sets $U_r, U_s \subseteq \mathcal{G}$ such that $U_r \cap r^{-1}(x) = \{\gamma\}$ and $U_s \cap s^{-1}(x) = \{\gamma\}$. Then $(U_r \cap U_s) \cap \mathcal{G}_x = (U_r \cap r^{-1}(x)) \cap (U_s \cap s^{-1}(x)) = \{\gamma\}$, so $\{\gamma\}$ is open in the subspace $\mathcal{G}_x$. Therefore $\mathcal{G}_x$ is discrete. □

**Lemma 1.4.21** ([21, Lemma 2.4.11])**.** Let $\mathcal{G}$ be a topological groupoid. If the range map $r : \mathcal{G} \to \mathcal{G}_0$ is open, then the multiplication $m : \mathcal{G}_2 \to \mathcal{G}$ is an open map. In particular, $m$ is open for any étale groupoid.

*Proof.* Let $U, V \subseteq \mathcal{G}$ be open and let $(\alpha, \beta) \in \mathcal{G}_2 \cap (U \times V)$, so $\alpha \in U$, $\beta \in V$ and the product $\alpha\beta$ is defined. We show that $\alpha\beta$ is an interior point of

$$UV := \{\mu\nu \mid \mu \in U,\ \nu \in V,\ s(\mu) = r(\nu)\} = m((U \times V) \cap \mathcal{G}_2).$$

Fix a decreasing neighbourhood base $(U_j)_{j\in J}$ of $\alpha$ in $\mathcal{G}$, with each $U_j \subseteq U$. Since $r$ is open, each $r(U_j)$ is an open neighbourhood of $r(\alpha) = r(\alpha\beta)$ in $\mathcal{G}_0$. Let $(\gamma_i)_{i\in I}$ be a net in $\mathcal{G}$ with $\gamma_i \to \alpha\beta$. Then $r(\gamma_i) \to r(\alpha\beta) = r(\alpha)$. For each $j \in J$, the set $r(U_j)$ is an open neighbourhood of $r(\alpha)$, so there exists $i_0(j) \in I$ such that $r(\gamma_i) \in r(U_j)$ for all $i \geq i_0(j)$. For each pair $(i, j)$ with $i \geq i_0(j)$, choose $\alpha_{(i,j)} \in U_j$ such that $r(\alpha_{(i,j)}) = r(\gamma_i)$. Let $D := \{(i, j) \in I \times J \mid i \geq i_0(j)\}$, directed by $(i, j) \leq (i', j')$ if $i \leq i'$ and $j \leq j'$. Then $\alpha_{(i,j)} \to \alpha$ as $(i, j)$ tends to infinity in $D$. Moreover, inversion is continuous, hence $\alpha^{-1}_{(i,j)} \to \alpha^{-1}$. For each $(i, j) \in D$ the pair $(\alpha^{-1}_{(i,j)}, \gamma_i)$ lies in $\mathcal{G}_2$ because $s(\alpha^{-1}_{(i,j)}) = r(\alpha_{(i,j)}) = r(\gamma_i)$. By continuity of multiplication,

$$\alpha^{-1}_{(i,j)}\gamma_i \to \alpha^{-1}(\alpha\beta) = \beta$$

in $\mathcal{G}$ along the net indexed by $D$. Since $V$ is an open neighbourhood of $\beta$, there exists $(i_1, j_1) \in D$ such that $\alpha^{-1}_{(i,j)}\gamma_i \in V$ for all $(i, j) \geq (i_1, j_1)$. For such $(i, j)$ we have $\alpha_{(i,j)} \in U$, $\alpha^{-1}_{(i,j)}\gamma_i \in V$, and the pair $(\alpha_{(i,j)}, \alpha_{(i,j)}^{-1}\gamma_i)$ is composable. Therefore $\gamma_i = \alpha_{(i,j)}(\alpha_{(i,j)}^{-1}\gamma_i) \in UV$ for all sufficiently large $(i,j)$. Since $(\gamma_i)$ was an arbitrary net converging to $\alpha\beta$, this shows that $\alpha\beta$ lies in the interior of $UV$. Hence $UV$ is open and $m$ is an open map. If $\mathcal{G}$ is étale, then $r$ is a local homeomorphism, hence an open map, so $m$ is open. □

Principal $\mathcal{G}$-bundles enter naturally when one wants to express Morita equivalence in geometric rather than purely functorial terms. In our setting, Morita equivalent étale groupoids are related by principal bibundles, and these bibundles are built from right and left principal $\mathcal{G}$-bundles as in the definition below. Introducing principal bundles therefore prepares the ground for describing Morita equivalence via correspondences, and for proving that invariants such as Moore homology are preserved under these geometric equivalences of groupoids.

**Definition 1.4.22** (Principal $\mathcal{G}$-bundle [4, §1.6])**.** Let $X$ be a topological space and $\mathcal{G}$ an étale groupoid. A right principal $\mathcal{G}$-bundle over $X$ consists of a topological space $P$, a surjective open map $\pi : P \to X$, an anchor $\epsilon : P \to \mathcal{G}_0$, and a continuous right action $P\ {}_s{\times}_r\ \mathcal{G}_1 \to P$, $(p,g) \mapsto p \cdot g$, defined when $\epsilon(p) = r(g)$, such that $\epsilon(p \cdot g) = s(g), p \cdot 1_{\epsilon(p)} = p, (p \cdot g) \cdot h = p \cdot (gh)$, and the canonical map of fibre products $P\ {}_s{\times}_r\ \mathcal{G}_1 \to P\ {}_\pi{\times}_\pi\ P$, $(p,g) \mapsto (p, p \cdot g)$, is a homeomorphism.

In the context of Morita equivalence we do not only need right principal $\mathcal{G}$-bundles, but bundles that carry in addition a compatible action of a second groupoid $\mathcal{K}$. A $\mathcal{K}$-equivariant principal $\mathcal{G}$-bundle is precisely a principal $\mathcal{G}$-bundle together with a left $\mathcal{K}$-action which commutes with the right $\mathcal{G}$-action and respects the projection to the base. Such objects organise into $(\mathcal{K}, \mathcal{G})$-bibundles and give the geometric realisation of Morita equivalences between $\mathcal{K}$ and $\mathcal{G}$. In particular, every Morita equivalence can be encoded by a $\mathcal{K}$-equivariant principal $\mathcal{G}$-bundle, and conversely these bibundles provide the correct framework to transport structures such as homology and cohomology functorially along equivalences of étale groupoids.

**Definition 1.4.23** ([4, §1.6])**.** Let $\mathcal{K}$ and $\mathcal{G}$ be étale groupoids.

A $\mathcal{K}$-equivariant principal $\mathcal{G}$-bundle consists of:

- a topological space $P$ together with continuous maps $\pi : P \to \mathcal{K}_0, \varepsilon : P \to \mathcal{G}_0$,
- a continuous right $\mathcal{G}$-action with anchor $\varepsilon$, $P\ {}_\varepsilon{\times}_{r_\mathcal{G}}\ \mathcal{G}_1 \to P$, $(p,g) \mapsto p \cdot g$, defined whenever $\varepsilon(p) = r_\mathcal{G}(g)$, such that $\varepsilon(p \cdot g) = s_\mathcal{G}(g), p \cdot 1_{\varepsilon(p)} = p, (p \cdot g) \cdot h = p \cdot (gh)$ for all composable $g, h \in \mathcal{G}_1$, and the canonical map of fibre products $\Theta :\ P\ {}_s{\times}_r \mathcal{G}_1 \to P\ {}_\pi{\times}_\pi P$, $(p,g) \mapsto (p, p{\cdot}g)$, is a homeomorphism. In particular, $\pi : P \to \mathcal{K}_0$ is a surjective open map and $(P, \pi, \varepsilon)$ is a principal right $\mathcal{G}$-bundle over $\mathcal{K}_0$ in the sense defined above,
- a continuous left $\mathcal{K}$-action with anchor $\pi$, $\mathcal{K}_1\ {}_\pi{\times}_\pi\ P \to P$, $(k,p) \mapsto k \cdot p$, defined whenever $s_\mathcal{K}(k) = \pi(p)$, such that $\pi(k \cdot p) = r_\mathcal{K}(k), 1_{\pi(p)} \cdot p = p, (k\ell) \cdot p = k \cdot (\ell \cdot p)$ for all composable $k, \ell \in \mathcal{K}_1$, and $\varepsilon(k \cdot p) = \varepsilon(p)$ for all $(k,p) \in \mathcal{K}_1\ {}_\pi{\times}_\pi\ P$,
- the left $\mathcal{K}$-action and the right $\mathcal{G}$-action commute, $(k \cdot p) \cdot g = k \cdot (p \cdot g)$ whenever defined.

This says that the left $\mathcal{K}$-action is by bundle automorphisms of the principal right $\mathcal{G}$-bundle $(P, \pi, \varepsilon)$: the map $\Theta$ is $\mathcal{K}$-equivariant for the diagonal $\mathcal{K}$-action on $P\ {}_\pi{\times}_\pi\ P$ and the induced $\mathcal{K}$-action on $P\ {}_s{\times}_r\ \mathcal{G}_1$. A morphism of $\mathcal{K}$-equivariant principal $\mathcal{G}$-bundles is a homeomorphism $f : P \to P'$ such that $\pi' = \pi \circ f$, $\varepsilon' = \varepsilon \circ f$, $f(k \cdot p) = k \cdot f(p)$ and $f(p \cdot g) = f(p) \cdot g$. Isomorphism classes of such objects are the generalised morphisms $\mathcal{K} \dashrightarrow \mathcal{G}$ [4, §1.6]. To apply the homology theories of [5, 13, 14] we adopt the setting of [4, §1.8]. Throughout, all groupoids are étale and their unit and morphism spaces are Hausdorff, locally compact and second countable. Under these standing assumptions the nerves $\mathcal{G}_n$ inherit compatible topologies, the bisection bases of Lemma 1.4.19 exist, the fibres in Corollary 1.4.20 are discrete and the multiplication is open by Lemma 1.4.21. These properties enter both the construction of convolution algebras [21] and the homology long exact sequences [14].

# 2 Moore Homology and Cohomology

Let $A$ be a topological abelian group, written additively with neutral element $0_A$. A subset of a topological space is called clopen if it is both closed and open. A topological space is totally disconnected if each of its connected components is a singleton. A Cantor set is a compact metrizable totally disconnected perfect space. A compact metrizable totally disconnected space is a Cantor set if and only if it has no isolated points. Any two Cantor sets are homeomorphic. Unless stated otherwise, all groupoids $\mathcal{G}$ considered below are second countable locally compact and Hausdorff. When we say that $\mathcal{G}$ is étale we mean that the source map $s : \mathcal{G} \to \mathcal{G}_0$ is a local homeomorphism. Then $r = s \circ i$ is a local homeomorphism as well, where $i(\gamma) = \gamma^{-1}$. This convention is used in [21, §2] and [13, §2.1]. For a locally compact Hausdorff space $X$ we write $C_c(X, A)$ for the abelian group of continuous functions $f : X \to A$ with compact support, where

$$\operatorname{supp}(f) := \overline{\{x \in X \mid f(x) \neq 0_A\}}.$$

If $X$ is compact, we abbreviate $C_c(X, A)$ by $C(X, A)$, since every continuous $f : X \to A$ then has compact support. With pointwise addition, both $C_c(X, A)$ and $C(X, A)$ are abelian groups. This section relies on the recent theory of Matui [13, 14, 15] and Sims et. al. [6, 20, 21].

## 2.1 Covariant Pushforward

Pushing forward along a local homeomorphism lets one transport compactly supported $A$-valued functions from $X$ to $Y$ by summing over fibres. This construction will be used to define maps induced by the structure maps of étale groupoids on compactly supported chains.

**Definition 2.1.1** (Pushforward)**.** Let $\phi : X \to Y$ be a local homeomorphism between locally compact Hausdorff spaces. For $f \in C_c(X, A)$ define $\phi_* f : Y \to A$ by

$$(\phi_* f)(y) := \sum_{x \in \phi^{-1}(y)} f(x) \quad \text{for all } y \in Y.$$

This is well defined: for each $y \in Y$ the fibre $\phi^{-1}(y)$ is discrete, hence $\operatorname{supp}(f) \cap \phi^{-1}(y)$ is a compact discrete space and therefore finite, so the sum has only finitely many nonzero terms. The empty sum is $0_A$. We next show that taking the pushforward is compatible with composition of local homeomorphisms.

**Proposition 2.1.2.** Let $A$ be a topological abelian group. Let $\phi\colon X \to Y$ and $\phi'\colon Y \to Z$ be local homeomorphisms between locally compact Hausdorff spaces. For every $f \in C_c(X, A)$,

$$(\phi' \circ \phi)_* f = \phi'_*(\phi_* f) \in C_c(Z, A).$$

We need a technical lemma as preparation.

**Lemma 2.1.3.** Let $\phi\colon X \to Y$ be a local homeomorphism between locally compact Hausdorff spaces, and let $K \subseteq X$ be compact. For every $y_0 \in Y$ there exist an open neighborhood $V$ of $y_0$, an integer $m \in \mathbb{N}_0$, and pairwise disjoint open sets $U_1, \dots, U_m \subseteq X$ such that

$$\phi^{-1}(V) \cap K \subseteq \bigcup_{i=1}^{m} U_i \quad \text{and} \quad \phi|_{U_i}\colon U_i \xrightarrow{\cong} V$$

is a homeomorphism for each $i \in \{1, \dots, m\}$.

*Proof.* If $\phi^{-1}(y_0) \cap K = \emptyset$, then $y_0 \notin \phi(K)$. Since $\phi(K)$ is compact and $Y$ is Hausdorff, $\phi(K)$ is closed. Hence there is an open neighborhood $V$ of $y_0$ with $V \cap \phi(K) = \emptyset$. Then $\phi^{-1}(V) \cap K = \emptyset$, and the conclusion holds with $m = 0$.

Otherwise, for each $x \in \phi^{-1}(y_0) \cap K$ choose open sets $U_x \ni x$ and $V_x \ni y_0$ such that $\phi|_{U_x}\colon U_x \to V_x$ is a homeomorphism. The fibre $\phi^{-1}(y_0)$ is discrete, hence $\phi^{-1}(y_0) \cap K$ is a discrete subspace of the compact Hausdorff space $K$, and therefore finite; enumerate it as $\{x_1, \dots, x_m\}$. Since $X$ is Hausdorff and $\{x_1, \dots, x_m\}$ is finite, after shrinking we may assume that the $U_{x_i}$ are pairwise disjoint. Set $V := \bigcap_{i=1}^{m} V_{x_i}$ and replace each $U_{x_i}$ by $U_{x_i} := (\phi|_{U_{x_i}})^{-1}(V)$. Then $\phi|_{U_{x_i}}\colon U_{x_i} \xrightarrow{\cong} V$ is a homeomorphism for each $i$, and the $U_{x_i}$ remain pairwise disjoint. We claim that, after possibly shrinking $V$, one has $\phi^{-1}(V) \cap K \subseteq \bigcup_{i=1}^{m} U_{x_i}$. If not, there exist a net $(y_\lambda)_\lambda$ in $Y$ with $y_\lambda \to y_0$ and points $z_\lambda \in (\phi^{-1}(y_\lambda) \cap K) \setminus \bigcup_{i=1}^{m} U_{x_i}$. By compactness of $K$, pass to a subnet with $z_\lambda \to z \in K$. Continuity of $\phi$ gives $\phi(z) = y_0$. Thus $z \in \phi^{-1}(y_0) \cap K = \{x_1, \dots, x_m\} \subseteq \bigcup_{i=1}^{m} U_{x_i}$. Since $\bigcup_{i=1}^{m} U_{x_i}$ is open, this implies $z_\lambda \in \bigcup_{i=1}^{m} U_{x_i}$ eventually, a contradiction. □

*Proof of Proposition 2.1.2.* Let $f \in C_c(X, A)$ and set $K := \operatorname{supp}(f)$.

1. **Pointwise equality:** For $z \in Z$,

$$\begin{aligned}(\phi'_*(\phi_* f))(z) &= \sum_{y \in \phi'^{-1}(z)} (\phi_* f)(y) = \sum_{y \in \phi'^{-1}(z)} \sum_{x \in \phi^{-1}(y)} f(x) \\ &= \sum_{x \in (\phi' \circ \phi)^{-1}(z)} f(x) = ((\phi' \circ \phi)_* f)(z).\end{aligned}$$

   These sums are well defined and finite: for each $y \in Y$, the fibre $\phi^{-1}(y)$ is discrete, hence $\phi^{-1}(y) \cap K$ is a discrete subspace of the compact Hausdorff space $K$, and therefore finite, so only finitely many $x$ contribute. For fixed $z \in Z$, the fibre $\phi'^{-1}(z)$ is discrete and $\phi(K)$ is compact, so $\phi'^{-1}(z) \cap \phi(K)$ is finite, hence only finitely many $y$ contribute.

2. $\boldsymbol{\phi_* f \in C_c(Y, A)}$: Fix $y_0 \in Y$. By Lemma 2.1.3 applied to $K$, there exist an open neighborhood $V \ni y_0$ and finitely many pairwise disjoint open sets $U_1, \dots, U_m \subseteq X$ such that $\phi|_{U_i}\colon U_i \to V$ is a homeomorphism and $\phi^{-1}(V) \cap K \subseteq \bigcup_{i=1}^m U_i$. Let $s_i\colon V \to U_i$ be the inverse $(\phi|_{U_i})^{-1}$. For $y \in V$, $(\phi_* f)(y) = \sum_{i=1}^m f(s_i(y))$, where the sum is finite (and is $0_A$ if $m = 0$). This is a finite sum of continuous maps, hence continuous. Since $y_0$ was arbitrary, $\phi_* f$ is continuous on $Y$. If $y \notin \phi(K)$, then $\phi^{-1}(y) \cap K = \emptyset$, hence $(\phi_* f)(y) = 0_A$. Therefore $\operatorname{supp}(\phi_* f) \subseteq \phi(K)$. Since $\phi(K)$ is compact, we have $\phi_* f \in C_c(Y, A)$.

3. $\boldsymbol{\phi'_*(\phi_* f) \in C_c(Z, A)}$: Apply 2. to $\phi'$ and the compact set $K' := \operatorname{supp}(\phi_* f) \subseteq Y$ to obtain $\phi'_*(\phi_* f) \in C_c(Z, A)$.

Combining 1.–3. yields $(\phi' \circ \phi)_* f = \phi'_*(\phi_* f)$ in $C_c(Z, A)$. □

**Corollary 2.1.4.** Let $A$ be a topological abelian group. Then $C_c(-, A)$ is a covariant functor from the category **LCH** of locally compact Hausdorff spaces with local homeomorphisms as morphisms to the category **Ab** of abelian groups.

*Proof.*

- **Objects and arrows.** For a locally compact Hausdorff space $X$, the set $C_c(X, A)$ of continuous, compactly supported maps $f : X \to A$ is an abelian group under pointwise addition. If $\phi\colon X \to Y$ is a local homeomorphism, define

$$(\phi_* f)(y) := \sum_{x \in \phi^{-1}(y)} f(x) \quad \text{for } y \in Y.$$

- **Well-definedness of $\phi_*$.** By Lemma 2.1.3, each fibre $\phi^{-1}(y)$ is discrete and meets $\operatorname{supp}(f)$ in only finitely many points, so the sum is finite in $A$. The continuity of $\phi_* f$ on $Y$ and the inclusion $\operatorname{supp}(\phi_* f) \subseteq \phi(\operatorname{supp}(f))$, hence compactness, were established in Proposition 2.1.2. $\phi_*\colon C_c(X, A) \to C_c(Y, A)$ is well-defined.
- **Homomorphism property.** For $f, g \in C_c(X, A)$ and $y \in Y$,

$$\phi_*(f + g)(y) = \sum_{x \in \phi^{-1}(y)} (f(x) + g(x)) = \sum_{x \in \phi^{-1}(y)} f(x) + \sum_{x \in \phi^{-1}(y)} g(x) = \phi_* f(y) + \phi_* g(y),$$

  so $\phi_*$ is a group homomorphism.
- **Functoriality.** For the identity, if $\mathrm{id}_X\colon X \to X$ then $(\mathrm{id}_X)^{-1}(y) = \{y\}$, hence $(\mathrm{id}_X)_* f(y) = f(y)$ and $(\mathrm{id}_X)_* = \mathrm{id}_{C_c(X,A)}$. If $\psi\colon Y \to Z$ is another local homeomorphism, then for $z \in Z$,

$$(\psi_* \phi_* f)(z) = \sum_{y \in \psi^{-1}(z)} \sum_{x \in \phi^{-1}(y)} f(x) = \sum_{x \in (\psi \circ \phi)^{-1}(z)} f(x) = ((\psi \circ \phi)_* f)(z),$$

  where interchanging the sums is valid because only finitely many terms are nonzero as above. In particular, Proposition 2.1.2 yields $(\psi \circ \phi)_* = \psi_* \circ \phi_*$.

Therefore $C_c(-, A)\colon \mathbf{LCH} \to \mathbf{Ab}$ is a covariant functor. □

Let $A$ be a topological abelian group. The constructions above show that $C_c(-, A)$ assigns to every local homeomorphism a pushforward group homomorphism, and that these maps are compatible with composition, see Proposition 2.1.2. Consequently, any diagram of locally compact Hausdorff spaces with local homeomorphisms as arrows gives rise, after applying $C_c(-, A)$ and pushforward, to a diagram of abelian groups. This functoriality is exactly what is needed to turn simplicial structure maps into boundary operators.

A simplicial space $X_\bullet$ consists of spaces $X_n$ for $n \geq 0$, face maps $d_i\colon X_n \to X_{n-1}$ and degeneracy maps $\sigma_i\colon X_n \to X_{n+1}$ satisfying the simplicial identities. Assume that all $d_i$ and $\sigma_i$ are local homeomorphisms between locally compact Hausdorff spaces. Then each $(d_i)_*\colon C_c(X_n, A) \to C_c(X_{n-1}, A)$ is well-defined, and we obtain a chain complex of abelian groups by setting $\partial_n := \sum_{i=0}^{n} (-1)^i (d_i)_*\colon C_c(X_n, A) \to C_c(X_{n-1}, A)$. The simplicial relations $d_i d_j = d_{j-1} d_i$ for $i < j$, together with functoriality of pushforward, imply $\partial_{n-1} \circ \partial_n = 0$. We write $H_n(X; A) := H_n(C_c(X_\bullet, A), \partial_\bullet)$ for the resulting homology groups. Our primary application will be to étale groupoids. If $\mathcal{G} \rightrightarrows \mathcal{G}_0$ is étale, its nerve $\mathcal{G}_\bullet$ has $\mathcal{G}_n$ the space of composable $n$-tuples, with face and degeneracy maps induced by the structure maps $r, s, u, i$ and multiplication. In the étale setting these induced maps are local homeomorphisms.[1] Therefore the preceding construction applies to $\mathcal{G}_\bullet$, and we define $H_n(\mathcal{G}; A) := H_n(C_c(\mathcal{G}_\bullet, A), \partial_\bullet)$, the Moore homology of $\mathcal{G}$ with coefficients in $A$, which we develop next.

## 2.2 Simplicial Spaces

For an étale groupoid $\mathcal{G}$ we need a canonical way to package finite composable strings of arrows so that the groupoid operations (units, inversion, and composition) are reflected by structure maps in a functorial topological object. The standard construction is the nerve $\mathcal{G}_\bullet$, a simplicial space whose $n$-simplices are composable $n$-tuples in $\mathcal{G}$ and whose face and degeneracy maps are induced by the groupoid structure.

- Set $C_n := C_c(\mathcal{G}_n, A)$ and use the face maps $d_i\colon \mathcal{G}_n \to \mathcal{G}_{n-1}$ to define the boundary $\partial_n := \sum_{i=0}^{n} (-1)^i (d_i)_*\colon C_n \to C_{n-1}$. If $\mathcal{G}$ is étale, then each $d_i$ is a local homeomorphism, hence $(d_i)_*$ is well defined on $C_c(-, A)$ by Definition 2.1.1. The simplicial identities $d_i d_j = d_{j-1} d_i$ for $i < j$, together with compatibility of pushforward with composition from Proposition 2.1.2, imply $\partial_{n-1} \circ \partial_n = 0$. Thus $(C_c(\mathcal{G}_\bullet, A), \partial_\bullet)$ is a chain complex, and its homology defines the Moore homology of $\mathcal{G}$ with coefficients in $A$.
- The simplicial space $\mathcal{G}_\bullet$ also has a geometric realization $B\mathcal{G}$, which provides a standard model for a classifying space of $\mathcal{G}$. This gives a conceptual link to classical constructions, but the Moore complex $C_c(\mathcal{G}_\bullet, A)$ is generally not the singular chain complex of $B\mathcal{G}$, which can be seen in Example 2.3.16. The nerve viewpoint is used to organize functoriality and exactness for compactly supported chains.
- An étale homomorphism of étale groupoids $\varphi\colon \mathcal{K} \to \mathcal{G}$ induces a simplicial map $N\varphi\colon \mathcal{K}_\bullet \to \mathcal{G}_\bullet$ whose components are local homeomorphisms. Applying pushforward degreewise yields a chain map $(\varphi_n)_*\colon C_c(\mathcal{K}_n, A) \to C_c(\mathcal{G}_n, A)$, hence induced homomorphisms $H_n(\mathcal{K}; A) \to H_n(\mathcal{G}; A)$. In the ample setting, reductions to full clopen subsets and, more generally, Kakutani equivalences will therefore yield canonical isomorphisms on Moore homology (see Theorem 2.5.5).

[1] The unit map $u\colon \mathcal{G}_0 \to \mathcal{G}$ is a continuous section of the local homeomorphism $s$, hence a homeomorphism onto an open subset of $\mathcal{G}$. The inversion map is a homeomorphism. The multiplication map $m\colon \mathcal{G}_2 \to \mathcal{G}$ is a local homeomorphism: for each $(g, h) \in \mathcal{G}_2$ one can choose open bisections $U, V \subseteq \mathcal{G}$ with $g \in U$, $h \in V$ and $s(U) = r(V)$, and then $m\colon U\,{}_s{\times}_r\, V \to UV$ is a homeomorphism onto an open subset of $\mathcal{G}$. The remaining face and degeneracy maps are obtained from these maps by finite products and pullbacks, which preserve local homeomorphisms.

- The simplicial structure also supplies standard homological algebra: normalized complexes, spectral sequences from filtrations of the nerve, long exact sequences associated to short exact sequences of chain complexes, and Mayer–Vietoris type constructions arising from clopen saturated covers and reductions.

To formalize nerves and, more generally, simplicial objects, we use the simplex category $\Delta$, which encodes the combinatorics of faces and degeneracies of simplices. Contravariant functors from $\Delta$ to a category $\mathcal{C}$ are simplicial objects in $\mathcal{C}$, and the relations in $\Delta$ yield the simplicial identities, including those needed to ensure $\partial_{n-1} \circ \partial_n = 0$ in the associated Moore complex.

**Definition 2.2.1** (Simplex Category)**.** Define the set of objects and morphisms between objects as

$$\begin{gathered} \mathrm{Ob}(\Delta) := \{[n] \mid n \in \mathbb{N},\ [n] := \{0, 1, \cdots, n\}\}, \\ \mathrm{Hom}_\Delta([m], [n]) := \{\theta\colon [m] \to [n] \mid \theta \text{ is order-preserving}\}. \end{gathered}$$

The morphism set is $\mathrm{Mor}(\Delta) := \bigsqcup_{m,n \in \mathbb{N}} \mathrm{Hom}_\Delta([m], [n])$. Composition is function composition

$$\circ\colon \mathrm{Hom}_\Delta([n], [p]) \times \mathrm{Hom}_\Delta([m], [n]) \to \mathrm{Hom}_\Delta([m], [p]), \quad (\eta, \theta) \mapsto \eta \circ \theta,$$

and the identity on $[n]$ is $\mathrm{id}_{[n]} \in \mathrm{Hom}_\Delta([n], [n])$. This makes $\Delta$ a small category.

1. **Generators and relations.** For $n \geq 1$ and $0 \leq i \leq n$, the coface map $\delta^i\colon [n-1] \to [n]$ is the injective order-preserving map and for $n \geq 0$ and $0 \leq j \leq n$, the codegeneracy map $\sigma^j\colon [n+1] \to [n]$ is the surjective order-preserving map for $0 \leq k \leq n$

$$\delta^i(k) = \begin{cases} k, & \text{for } k < i, \\ k+1, & \text{for } k \geq i. \end{cases} \quad \sigma^j(k) = \begin{cases} k, & \text{for } k \leq j, \\ k-1, & \text{for } k \geq j+1. \end{cases}$$

   These maps generate all morphisms of $\Delta$ and satisfy

$$\begin{aligned} \delta^j \delta^i &= \delta^i \delta^{j-1} \text{ for } i < j, \\ \sigma^j \sigma^i &= \sigma^i \sigma^{j+1} \text{ for } i \leq j, \\ \sigma^j \delta^i &= \begin{cases} \delta^i \sigma^{j-1}, & \text{for } i < j, \\ \mathrm{id}_{[n]}, & \text{for } i = j \text{ or } i = j+1, \\ \delta^{i-1} \sigma^j, & \text{for } i > j+1. \end{cases} \end{aligned}$$

2. **Opposite category $\Delta^{\mathrm{op}}$.** $\mathrm{Ob}(\Delta^{\mathrm{op}}) = \mathrm{Ob}(\Delta)$, $\mathrm{Hom}_{\Delta^{\mathrm{op}}}([n], [m]) = \mathrm{Hom}_\Delta([m], [n])$, with composition reversed. Writing $d_i := (\delta^i)^{\mathrm{op}}\colon [n] \to [n-1], s_j := (\sigma^j)^{\mathrm{op}}\colon [n] \to [n+1]$, the above relations become the simplicial identities for the $d_i$ and $s_j$.

To turn the combinatorics of simplices into topological input for chain complexes, we package the spaces of "$n$-simplices" together with their face and degeneracy operations into a single object. This is the natural setting for nerves of groupoids and for the Moore chain complex built from pushforward along face maps.

**Definition 2.2.2** (Simplicial space)**.** A family $(X_n, (d_i)_{i=0}^n, (s_j)_{j=0}^n)_{n\geq 0}$ of topological spaces with continuous face maps $d_i\colon X_n \to X_{n-1}$ and degeneracy maps $s_j\colon X_n \to X_{n+1}$ is called a simplicial space if the simplicial identities hold:

$$d_i d_j = d_{j-1} d_i \quad \text{for } i < j, \qquad s_i s_j = s_{j+1} s_i \quad \text{for } i \leq j, \qquad d_i s_j = \begin{cases} s_{j-1} d_i, & \text{for } i < j, \\ \mathrm{id}_{X_n}, & \text{for } i = j \text{ or } i = j+1, \\ s_j d_{i-1}, & \text{for } i > j+1. \end{cases}$$

*Remark* 2.2.3. In other words, a simplicial space is a functor $X_\bullet\colon \Delta^{\mathrm{op}} \to \mathbf{Top}$. In our applications we assume $X_n \in \mathbf{LCH}$ for all $n$, so it suffices to regard $X_\bullet$ as a functor $X_\bullet\colon \Delta^{\mathrm{op}} \to \mathbf{LCH}$. For the definition of the opposite simplex category $\Delta^{\mathrm{op}}$ we refer to Definition 2.2.1.

Simplicial spaces are best viewed as functors $X_\bullet : \Delta^{\mathrm{op}} \to \mathbf{Top}$. Accordingly, a map $f_\bullet : X_\bullet \to Y_\bullet$ should be a morphism of such functors, that is, a natural transformation. Concretely, it is a family $(f_n)_{n\geq 0}$ commuting with all face and degeneracy maps.

**Definition 2.2.4** (Maps of simplicial spaces)**.** A map $f_\bullet\colon X_\bullet \to Y_\bullet$ is a natural transformation, i.e. a family of continuous maps $f_n\colon X_n \to Y_n$ for $n \geq 0$ such that the following diagrams commute. For every $n \geq 1$ and $0 \leq i \leq n$, and for every $n \geq 0$ and $0 \leq j \leq n$,

$$\begin{array}{ccc} X_n & \xrightarrow{f_n} & Y_n \\ {\scriptstyle d_i^X}\downarrow & & \downarrow{\scriptstyle d_i^Y} \\ X_{n-1} & \xrightarrow[f_{n-1}]{} & Y_{n-1}, \end{array} \qquad \begin{array}{ccc} X_n & \xrightarrow{f_n} & Y_n \\ {\scriptstyle s_j^X}\downarrow & & \downarrow{\scriptstyle s_j^Y} \\ X_{n+1} & \xrightarrow[f_{n+1}]{} & Y_{n+1}. \end{array}$$

To define homology from an étale groupoid $\mathcal{G}$, we need a systematic way to package all composable strings of arrows together with the structural operations (source, range, units, inversion, multiplication) so that the resulting spaces vary functorially under étale homomorphisms of groupoids. The standard tool for this is the nerve $\mathcal{G}_\bullet$: its $n$-simplices are composable $n$-tuples, and the simplicial operators are obtained by composing adjacent arrows or inserting units. This provides the simplicial space on which the Moore chain complex $C_c(\mathcal{G}_n, A)$ is built. Let $(\mathcal{G}, \mathcal{G}_0, r, s, m, {}^{-1})$ be an étale groupoid with unit space $\mathcal{G}_0$ and range/source maps $r, s\colon \mathcal{G} \to \mathcal{G}_0$. Since $\mathcal{G}$ is étale, $s$ is a local homeomorphism, hence $r = s \circ (-^{-1})$ is a local homeomorphism as well. Write $\mathcal{G}_2 \coloneqq \{(g,h) \in \mathcal{G} \times \mathcal{G} \mid s(g) = r(h)\} = \mathcal{G}\,{}_s\!\times_r \mathcal{G}$ for the space of composable pairs, with projections $p_1, p_2$. Composition is the continuous map $m\colon \mathcal{G}_2 \to \mathcal{G}$, written left–to–right: $g \cdot h := m(g,h)$ is defined when $s(g) = r(h)$, with $s(g \cdot h) = s(h)$ and $r(g \cdot h) = r(g)$. Diagrammatically,

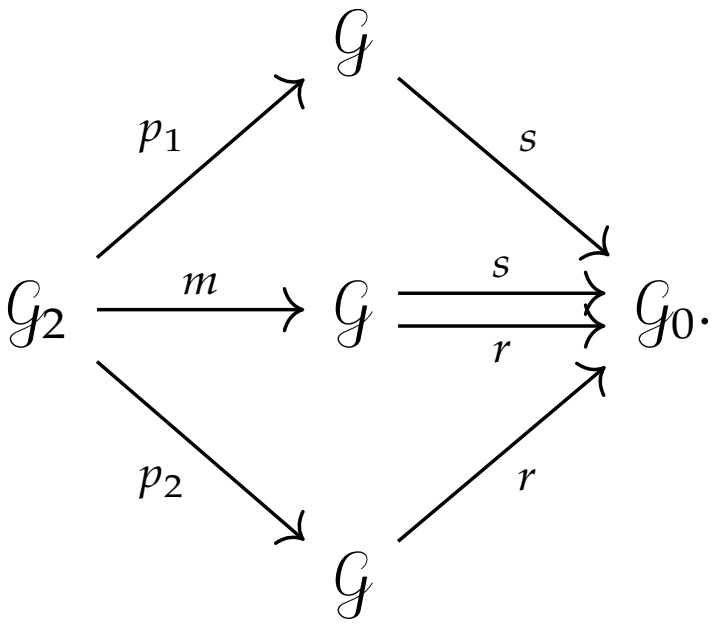

For $n \geq 0$ define the levels of the nerve, see Definition 2.3.1, by

$$\mathcal{G}_0 := \mathcal{G}_0, \quad \mathcal{G}_n := \{(g_1, \dots, g_n) \in \mathcal{G}^n \mid s(g_i) = r(g_{i+1}) \text{ for } i \in \{1, \dots, n-1\}\} \quad \text{for } n \geq 1,$$

with the subspace topology from $\mathcal{G}^n$. In particular, $\mathcal{G}_1 = \mathcal{G}$.

When convenient, we abbreviate an $n$-tuple $(g_1, \dots, g_n) \in \mathcal{G}_n$ by $\mathbf{g}$, and we write $g_1 \cdots g_n$ for its composite. This is unambiguous: since $\mathbf{g} \in \mathcal{G}_n$, all intermediate products are defined, and associativity of the groupoid multiplication implies that every parenthesization yields the same arrow in $\mathcal{G}$. In particular, $g_1 \cdots g_n \in \mathcal{G}$ satisfies $r(g_1 \cdots g_n) = r(g_1)$ and $s(g_1 \cdots g_n) = s(g_n)$.

We fix the indexing of the groupoid once and for all. While the unit space is often written $\mathcal{G}^{(0)}$, we write $\mathcal{G}_0$ to emphasize that we organize $\mathcal{G}$ via its nerve $\mathcal{G}_\bullet$: the 0-simplices are exactly the units, $\mathcal{G}_1 = \mathcal{G}$, and $\mathcal{G}_n$ is the space of composable $n$-tuples. With this convention the face maps $d_i : \mathcal{G}_n \to \mathcal{G}_{n-1}$ and degeneracies $s_j : \mathcal{G}_n \to \mathcal{G}_{n+1}$ always change the simplicial degree by $\pm 1$, and we avoid additional parenthesized notations for composability spaces. This alignment is tailored to the Moore–complex viewpoint. In degree $n$ the chain group is $C_c(\mathcal{G}_n, A)$, and the boundary is the alternating sum $\partial_n = \sum_{i=0}^{n} (-1)^i (d_i)_*$. Thus the single index $n$ simultaneously records the simplicial level, the underlying space $\mathcal{G}_n$, and the corresponding chain group, which keeps later degreewise constructions readable (pushforwards, functoriality, and exact sequences), all of them living on the simplicial object $\mathcal{G}_\bullet$, see Definition 2.3.1.

*Remark* 2.2.5. For $n \geq 2$ one may identify $\mathcal{G}_n$ with the iterated fibre product

$$\underbrace{\mathcal{G} \,{}_s\!\times_r \mathcal{G} \,{}_s\!\times_r \cdots {}_s\!\times_r \mathcal{G}}_{n \text{ factors}};$$

the composability condition $s(g_i) = r(g_{i+1})$ is then built into the pullback. Under this description the endpoint faces $d_0, d_n$ forget the first, respectively last, arrow, while the inner faces $d_i$ compose the $(i, i+1)$-entries. Finally, in the étale setting the structure maps of the nerve are local homeomorphisms, hence have discrete fibres; combined with compact support, this ensures that the fibrewise sums defining pushforwards $(d_i)_*$ are finite, and therefore that the Moore boundary formula is well defined.

**Definition 2.2.6** (Face and degeneracy maps [14, p. 4])**.** Let $\mathcal{G}$ be an étale groupoid with unit map $u\colon \mathcal{G}_0 \to \mathcal{G}_1 = \mathcal{G}$. For $n \geq 1$ and $i \in \{0, 1, \dots, n\}$ define $d_i\colon \mathcal{G}_n \to \mathcal{G}_{n-1}$, and for $n \geq 0$ and $j \in \{0, 1, \dots, n\}$ define $s_j\colon \mathcal{G}_n \to \mathcal{G}_{n+1}$, by

$$d_i(\mathbf{g}) = \begin{cases} s(g_1), & \text{if } n = 1,\ i = 0,\\ r(g_1), & \text{if } n = 1,\ i = 1,\\ (g_2, \dots, g_n), & \text{if } n \geq 2,\ i = 0,\\ (g_1, \dots, g_i \cdot g_{i+1}, \dots, g_n), & \text{if } n \geq 2,\ 1 \leq i \leq n-1,\\ (g_1, \dots, g_{n-1}), & \text{if } n \geq 2,\ i = n, \end{cases}$$

$$s_j(\mathbf{g}) = \begin{cases} u(x), & \text{if } n = 0,\ \mathbf{g} = x \in \mathcal{G}_0,\\ (u(r(g_1)), g_1, \dots, g_n), & \text{if } n \geq 1,\ j = 0,\\ (g_1, \dots, g_j,\ u(r(g_{j+1})),\ g_{j+1}, \dots, g_n), & \text{if } n \geq 2,\ 1 \leq j \leq n-1,\\ (g_1, \dots, g_n,\ u(s(g_n))), & \text{if } n \geq 1,\ j = n, \end{cases}$$

where $\mathbf{g} = (g_1, \dots, g_n)$ for $n \geq 1$.

**Proposition 2.2.7.** $\mathcal{G}_\bullet := (\mathcal{G}_n, (d_i)_{i=0}^n, (s_j)_{j=0}^n)_{n\geq 0}$ is a simplicial space.

*Proof.* We verify simplicial identities by computations of the formulas from Definition 2.2.6.

- **Continuity:** For $n \geq 2$, the maps $d_0$ and $d_n$ are coordinate projections $d_0(g_1, \dots, g_n) = (g_2, \dots, g_n)$ and $d_n(g_1, \dots, g_n) = (g_1, \dots, g_{n-1})$, hence continuous. For $1 \leq i \leq n-1$, $d_i(g_1, \dots, g_n) = (g_1, \dots, g_i \cdot g_{i+1}, \dots, g_n)$, which is obtained by composing a projection $\mathcal{G}_n \to \mathcal{G}_2$, $(g_1, \dots, g_n) \mapsto (g_i, g_{i+1})$, with $m : \mathcal{G}_2 \to \mathcal{G}$, and then inserting the product back into the tuple, hence $d_i$ is continuous. For $n = 1$, $d_0 = s$ and $d_1 = r$, hence continuous. Each $s_j$ is defined by inserting a unit $u(\cdot)$ into a tuple and leaving all other entries unchanged, so $s_j$ is a product of coordinate projections with $u$, hence continuous.
- **Face–face identities:** Fix $n \geq 2$, $0 \leq i < j \leq n$, and $\mathbf{g} = (g_1, \dots, g_n) \in \mathcal{G}_n$.
  We show $d_i d_j(\mathbf{g}) = d_{j-1} d_i(\mathbf{g})$ by cases.
  - **Case** $n = 2$**:** Then $\mathbf{g} = (g_1, g_2) \in \mathcal{G}_2$.
    Recall that after applying one face map we land in $\mathcal{G}_1 = \mathcal{G}$, and then $d_0 = s$, $d_1 = r$.
    - If $(i, j) = (0, 1)$, then

      $$\begin{aligned} d_1(g_1, g_2) &= g_1 \cdot g_2, & d_0 d_1(g_1, g_2) &= s(g_1 \cdot g_2) = s(g_2),\\ d_0(g_1, g_2) &= g_2, & d_0 d_0(g_1, g_2) &= s(g_2). \end{aligned}$$

    - If $(i, j) = (0, 2)$, then

      $$\begin{aligned} d_2(g_1, g_2) &= g_1, & d_0 d_2(g_1, g_2) &= s(g_1),\\ d_0(g_1, g_2) &= g_2, & d_1 d_0(g_1, g_2) &= r(g_2) = s(g_1), \end{aligned}$$

      using $s(g_1) = r(g_2)$.

- If $(i,j) = (1,2)$, then

$$d_2(g_1,g_2) = g_1, \qquad d_1 d_2(g_1,g_2) = r(g_1),$$
$$d_1(g_1,g_2) = g_1 \cdot g_2, \quad d_1 d_1(g_1,g_2) = r(g_1 \cdot g_2) = r(g_1).$$

- **Case $n = 3$:** Then $\mathbf{g} = (g_1, g_2, g_3) \in \mathcal{G}_3$. First record the faces:

$$d_0(g_1,g_2,g_3) = (g_2,g_3),$$
$$d_1(g_1,g_2,g_3) = (g_1 \cdot g_2, g_3),$$
$$d_2(g_1,g_2,g_3) = (g_1, g_2 \cdot g_3),$$
$$d_3(g_1,g_2,g_3) = (g_1,g_2).$$

Now check all $(i,j)$ with $0 \leq i < j \leq 3$:

- If $(i,j) = (0,1)$, then

$$d_0 d_1(g_1,g_2,g_3) = d_0(g_1 \cdot g_2, g_3) = g_3,$$
$$d_0 d_0(g_1,g_2,g_3) = d_0(g_2,g_3) = g_3.$$

- If $(i,j) = (0,2)$, then

$$d_0 d_2(g_1,g_2,g_3) = d_0(g_1, g_2 \cdot g_3) = g_2 \cdot g_3,$$
$$d_1 d_0(g_1,g_2,g_3) = d_1(g_2,g_3) = g_2 \cdot g_3.$$

- If $(i,j) = (0,3)$, then

$$d_0 d_3(g_1,g_2,g_3) = d_0(g_1,g_2) = g_2,$$
$$d_2 d_0(g_1,g_2,g_3) = d_2(g_2,g_3) = g_2.$$

- If $(i,j) = (1,2)$, then

$$d_1 d_2(g_1,g_2,g_3) = d_1(g_1, g_2 \cdot g_3) = g_1 \cdot (g_2 \cdot g_3),$$
$$d_1 d_1(g_1,g_2,g_3) = d_1(g_1 \cdot g_2, g_3) = (g_1 \cdot g_2) \cdot g_3,$$

  which are equal by associativity.

- If $(i,j) = (1,3)$, then

$$d_1 d_3(g_1,g_2,g_3) = d_1(g_1,g_2) = g_1 \cdot g_2,$$
$$d_2 d_1(g_1,g_2,g_3) = d_2(g_1 \cdot g_2, g_3) = g_1 \cdot g_2.$$

- If $(i,j) = (2,3)$, then

$$d_2 d_3(g_1,g_2,g_3) = d_2(g_1,g_2) = g_1,$$
$$d_2 d_2(g_1,g_2,g_3) = d_2(g_1, g_2 \cdot g_3) = g_1.$$

- **Case** $n \geq 4$**:** Fix $0 \leq i < j \leq n$.
  - **Subcase** $i = 0$**:**
    - If $j = 1$, then

$$\begin{aligned} d_1(\mathbf{g}) &= (g_1 \cdot g_2, g_3, \cdots, g_n), \\ d_0 d_1(\mathbf{g}) &= (g_3, \cdots, g_n), \\ d_0(\mathbf{g}) &= (g_2, \cdots, g_n), \\ d_0 d_0(\mathbf{g}) &= (g_3, \cdots, g_n). \end{aligned}$$

    - If $2 \leq j \leq n-1$, then

$$\begin{aligned} d_j(\mathbf{g}) &= (g_1, \cdots, g_{j-1}, g_j \cdot g_{j+1}, g_{j+2}, \cdots, g_n), \\ d_0 d_j(\mathbf{g}) &= (g_2, \cdots, g_{j-1}, g_j \cdot g_{j+1}, g_{j+2}, \cdots, g_n), \\ d_0(\mathbf{g}) &= (g_2, \cdots, g_n), \\ d_{j-1} d_0(\mathbf{g}) &= (g_2, \cdots, g_{j-1}, g_j \cdot g_{j+1}, g_{j+2}, \cdots, g_n). \end{aligned}$$

    - If $j = n$, then

$$\begin{aligned} d_n(\mathbf{g}) &= (g_1, \cdots, g_{n-1}), \\ d_0 d_n(\mathbf{g}) &= (g_2, \cdots, g_{n-1}), \\ d_0(\mathbf{g}) &= (g_2, \cdots, g_n), \\ d_{n-1} d_0(\mathbf{g}) &= (g_2, \cdots, g_{n-1}). \end{aligned}$$

  - **Subcase** $1 \leq i \leq n-1$ **and** $j = n$**:**
    - If $1 \leq i \leq n-2$, then

$$\begin{aligned} d_n(\mathbf{g}) &= (g_1, \cdots, g_{n-1}), \\ d_i d_n(\mathbf{g}) &= (g_1, \cdots, g_i \cdot g_{i+1}, \cdots, g_{n-1}), \\ d_i(\mathbf{g}) &= (g_1, \cdots, g_i \cdot g_{i+1}, \cdots, g_n), \\ d_{n-1} d_i(\mathbf{g}) &= (g_1, \cdots, g_i \cdot g_{i+1}, \cdots, g_{n-1}). \end{aligned}$$

    - If $i = n-1$, then

$$\begin{aligned} d_n(\mathbf{g}) &= (g_1, \cdots, g_{n-1}), \\ d_{n-1} d_n(\mathbf{g}) &= (g_1, \cdots, g_{n-2}), \\ d_{n-1}(\mathbf{g}) &= (g_1, \cdots, g_{n-2}, g_{n-1} \cdot g_n), \\ d_{n-1} d_{n-1}(\mathbf{g}) &= (g_1, \cdots, g_{n-2}). \end{aligned}$$

- **Subcase $1 \le i < j \le n-1$:**
  - If $j \ge i+2$, then $d_j$ composes $(g_j, g_{j+1})$ and $d_i$ composes $(g_i, g_{i+1})$, which are disjoint pairs. Hence

$$\begin{aligned} d_i d_j(\mathbf{g}) &= (g_1, \dots, g_{i-1}, g_i \cdot g_{i+1}, g_{i+2}, \dots, g_{j-1}, g_j \cdot g_{j+1}, g_{j+2}, \dots, g_n) \\ &= d_{j-1} d_i(\mathbf{g}). \end{aligned}$$

  - If $j = i+1$, then

$$\begin{aligned} d_i d_{i+1}(\mathbf{g}) &= (g_1, \dots, g_{i-1}, g_i \cdot (g_{i+1} \cdot g_{i+2}), g_{i+3}, \dots, g_n), \\ d_i d_i(\mathbf{g}) &= (g_1, \dots, g_{i-1}, (g_i \cdot g_{i+1}) \cdot g_{i+2}, g_{i+3}, \dots, g_n), \end{aligned}$$

    which are equal by associativity.

Thus $d_i d_j = d_{j-1} d_i$ for all $n \ge 2$ and $0 \le i < j \le n$.

- **Degeneracy–degeneracy identities:** Fix $n \ge 0$, $0 \le i \le j \le n$, and $\mathbf{g} \in \mathcal{G}_n$.
  We show $s_i s_j(\mathbf{g}) = s_{j+1} s_i(\mathbf{g})$ by cases.
  - **Case $n = 0$:**
    Then $i = j = 0$ and $\mathbf{g} = x \in \mathcal{G}_0$. Using $r(u(x)) = x = s(u(x))$,

$$\begin{aligned} s_0 s_0(x) &= s_0(u(x)) = (u(r(u(x))), u(x)) = (u(x), u(x)), \\ s_1 s_0(x) &= s_1(u(x)) = (u(x), u(s(u(x)))) = (u(x), u(x)). \end{aligned}$$

  - **Case $n \ge 1$:** Write $\mathbf{g} = (g_1, \dots, g_n)$.
    - If $j = 0$, then $i = 0$ and

$$\begin{aligned} s_0(\mathbf{g}) &= (u(r(g_1)), g_1, \dots, g_n), \\ s_0 s_0(\mathbf{g}) &= (u(r(u(r(g_1)))), u(r(g_1)), g_1, \dots, g_n) = (u(r(g_1)), u(r(g_1)), g_1, \dots, g_n), \\ s_1 s_0(\mathbf{g}) &= (u(r(g_1)), u(r(g_1)), g_1, \dots, g_n). \end{aligned}$$

    - If $1 \le j \le n-1$, then $s_j(\mathbf{g}) = (g_1, \dots, g_j, u(r(g_{j+1})), g_{j+1}, \dots, g_n)$.
      - If $i = 0$, then

$$s_0 s_j(\mathbf{g}) = (u(r(g_1)), g_1, \dots, g_j, u(r(g_{j+1})), g_{j+1}, \dots, g_n) = s_{j+1} s_0(\mathbf{g}).$$

      - If $1 \le i \le j$, then

$$\begin{aligned} s_i s_j(\mathbf{g}) &= (g_1, \dots, g_i, u(r(g_{i+1})), g_{i+1}, \dots, g_j, u(r(g_{j+1})), g_{j+1}, \dots, g_n), \\ s_{j+1} s_i(\mathbf{g}) &= (g_1, \dots, g_i, u(r(g_{i+1})), g_{i+1}, \dots, g_j, u(r(g_{j+1})), g_{j+1}, \dots, g_n). \end{aligned}$$

    - If $j = n$, then $s_n(\mathbf{g}) = (g_1, \dots, g_n, u(s(g_n)))$.
      - If $0 \le i \le n-1$, then

$$s_i s_n(\mathbf{g}) = (g_1, \dots, g_i, u(r(g_{i+1})), g_{i+1}, \dots, g_n, u(s(g_n))) = s_{n+1} s_i(\mathbf{g}).$$

- If $i = n$, then using $s(u(x)) = x$,

$$s_n s_n(\mathbf{g}) = (g_1, \dots, g_n, u(s(g_n)), u(s(g_n))) = s_{n+1} s_n(\mathbf{g}).$$

- **Face–degeneracy identities:** Fix $n \geq 0$, $0 \leq j \leq n$, and $\mathbf{g} \in \mathcal{G}_n$. We verify that for all $0 \leq i \leq n+1$,

$$d_i s_j = \begin{cases} s_{j-1} d_i, & \text{for } i < j, \\ \mathrm{id}_{\mathcal{G}_n}, & \text{for } i = j \text{ or } i = j+1, \\ s_j d_{i-1}, & \text{for } i > j+1. \end{cases}$$

  - **Case $n = 0$:** Then $j = 0$ and $\mathbf{g} = x \in \mathcal{G}_0$, so $s_0(x) = u(x)$. Thus

$$d_0 s_0(x) = s(u(x)) = x,$$
$$d_1 s_0(x) = r(u(x)) = x.$$

  - **Case $n \geq 1$:** Write $\mathbf{g} = (g_1, \dots, g_n)$.
    - **Subcase $j = 0$:** Then $s_0(\mathbf{g}) = (u(r(g_1)), g_1, \dots, g_n)$.
      - If $i = 0$, then $d_0 s_0(\mathbf{g}) = (g_1, \dots, g_n) = \mathbf{g}$.
      - If $i = 1$, then

$$d_1 s_0(\mathbf{g}) = (u(r(g_1)) \cdot g_1, g_2, \dots, g_n) = (g_1, \dots, g_n) = \mathbf{g},$$

        by the left unit law $u(r(g_1)) \cdot g_1 = g_1$.
      - If $2 \leq i \leq n$, then

$$d_i s_0(\mathbf{g}) = (u(r(g_1)), g_1, \dots, g_{i-2}, g_{i-1} \cdot g_i, g_{i+1}, \dots, g_n) = s_0 d_{i-1}(\mathbf{g}).$$

      - If $i = n+1$, then

$$d_{n+1} s_0(\mathbf{g}) = (u(r(g_1)), g_1, \dots, g_{n-1}) = s_0 d_n(\mathbf{g}).$$

    - **Subcase $1 \leq j \leq n-1$:** Then $s_j(\mathbf{g}) = (g_1, \dots, g_j, u(r(g_{j+1})), g_{j+1}, \dots, g_n)$.
      - If $i < j$, then

$$d_i s_j(\mathbf{g}) = (g_1, \dots, g_{i-1}, g_i \cdot g_{i+1}, g_{i+2}, \dots, g_j, u(r(g_{j+1})), g_{j+1}, \dots, g_n) = s_{j-1} d_i(\mathbf{g}).$$

      - If $i = j$, then

$$d_j s_j(\mathbf{g}) = (g_1, \dots, g_{j-1}, \; g_j \cdot u(r(g_{j+1})), \; g_{j+1}, \dots, g_n) = \mathbf{g},$$

        using $s(g_j) = r(g_{j+1})$ and $g_j \cdot u(s(g_j)) = g_j$.

· If $i = j+1$, then

$$d_{j+1}s_j(\mathbf{g}) = (g_1, \cdots, g_j, u(r(g_{j+1})) \cdot g_{j+1}, g_{j+2}, \cdots, g_n) = \mathbf{g}.$$

· If $j+2 \leq i \leq n$, then

$$d_i s_j(\mathbf{g}) = (g_1, \dots, g_j, u(r(g_{j+1})), g_{j+1}, \cdots, g_{i-2}, g_{i-1} \cdot g_i, g_{i+1}, \cdots, g_n) = s_j d_{i-1}(\mathbf{g}).$$

· If $i = n+1$, then

$$d_{n+1}s_j(\mathbf{g}) = (g_1, \dots, g_j, u(r(g_{j+1})), g_{j+1}, \cdots, g_{n-1}) = s_j d_n(\mathbf{g}).$$

- **Subcase $j = n$:** Then $s_n(\mathbf{g}) = (g_1, \dots, g_n, u(s(g_n)))$.
  · If $0 \leq i < n$, then

$$d_i s_n(\mathbf{g}) = (g_1, \dots, g_{i-1}, g_i \cdot g_{i+1}, g_{i+2}, \cdots, g_n, u(s(g_n))) = s_{n-1}d_i(\mathbf{g}).$$

  · If $i = n$, then

$$d_n s_n(\mathbf{g}) = (g_1, \dots, g_{n-1}, g_n \cdot u(s(g_n))) = \mathbf{g}.$$

  · If $i = n+1$, then $d_{n+1}s_n(\mathbf{g}) = (g_1, \dots, g_n) = \mathbf{g}$.

Therefore all simplicial identities hold, and $(\mathcal{G}_n, (d_i)_{i=0}^n, (s_j)_{j=0}^n)_{n\geq 0}$ is a simplicial space. □

Starting from the simplex category $\Delta$ from Definition 2.2.1, a simplicial space is, by definition, a functor $X_\bullet\colon \Delta^{\mathrm{op}} \to \mathbf{Top}$: the objects $[n]$ give the levels $X_n$, and the generating maps in $\Delta$, cofaces $\delta^i$ and codegeneracies $\sigma^j$, correspond in $\Delta^{\mathrm{op}}$ to face maps $d_i$ and degeneracy maps $s_j$ satisfying the simplicial identities. For an étale groupoid $(\mathcal{G}, \mathcal{G}_0, r, s, m, {}^{-1})$, the nerve $\mathcal{G}_\bullet$ is obtained by taking $\mathcal{G}_n$ to be the space of composable $n$-tuples and by defining $d_i$ and $s_j$ via the structure maps $r, s, m, u$ as in Definition 2.2.6. Proposition 2.2.7 verifies by explicit computation that these maps satisfy the simplicial identities, hence $\mathcal{G}_\bullet$ is a simplicial space. Since $\mathcal{G}$ is étale, each face map $d_i\colon \mathcal{G}_n \to \mathcal{G}_{n-1}$ is a local homeomorphism. Therefore the push-forward construction from Definition 2.1.1 applies and yields well-defined homomorphisms $(d_i)_*\colon C_c(\mathcal{G}_n, A) \to C_c(\mathcal{G}_{n-1}, A)$, compatible with composition by Proposition 2.1.2. Moreover, for $f \in C_c(\mathcal{G}_n, A)$ one has $\mathrm{supp}((d_i)_* f) \subseteq d_i(\mathrm{supp}(f))$, and continuity and compact support are ensured by Lemma 2.1.3.

## 2.3 Homology Groups

Groupoid homology in this work is built from compactly supported chains on the spaces of composable strings of arrows. The nerve construction packages these spaces, together with the algebraic operations of a groupoid, into a single simplicial space in a functorial way, so that later chain-level constructions are automatically compatible with functors of groupoids.

**Definition 2.3.1** (Nerve functor for groupoids)**.** Let $\mathcal{G}$**Top** be the category of topological groupoids with continuous functors, and let **sTop** $:= \mathrm{Fun}(\Delta^{\mathrm{op}}, \mathbf{Top})$ be the category of simplicial spaces with simplicial maps. Define an assignment $N\colon \mathcal{G}\mathbf{Top} \to \mathbf{sTop}$ as follows.

- **On objects.** Let $\mathcal{G}$ be a groupoid with source and range maps $s, r\colon \mathcal{G}_1 \to \mathcal{G}_0$, unit map $u\colon \mathcal{G}_0 \to \mathcal{G}_1$, and multiplication written left to right: $\gamma \cdot \eta$ is defined when $s(\gamma) = r(\eta)$. Set $N(\mathcal{G}) = \mathcal{G}_\bullet$, where $\mathcal{G}_0$ is the unit space and, for $n \geq 1$,

$$\mathcal{G}_n := \{(\gamma_1, \dots, \gamma_n) \in \mathcal{G}_1^n \mid s(g_i) = r(g_{i+1}) \text{ for } 1 \leq i < n\},$$

  endowed with the subspace topology from $\mathcal{G}_1^n$.
  For $n \geq 1$ and $0 \leq i \leq n$, define the face maps $d_i\colon \mathcal{G}_n \to \mathcal{G}_{n-1}$ by

$$d_i(\gamma_1, \dots, \gamma_n) := \begin{cases} s(\gamma_1), & n = 1,\ i = 0, \\ r(\gamma_1), & n = 1,\ i = 1, \\ (\gamma_2, \dots, \gamma_n), & n \geq 2,\ i = 0, \\ (\gamma_1, \dots, g_i \cdot g_{i+1}, \dots, \gamma_n), & n \geq 2,\ 1 \leq i \leq n-1, \\ (\gamma_1, \dots, g_{n-1}), & n \geq 2,\ i = n. \end{cases}$$

  For $n \geq 0$ and $0 \leq j \leq n$, define the degeneracy maps $s_j\colon \mathcal{G}_n \to \mathcal{G}_{n+1}$ by

$$s_j(\gamma_1, \dots, \gamma_n) := \begin{cases} u(x), & n = 0,\ j = 0,\ x \in \mathcal{G}_0, \\ (u(r(\gamma_1)), \gamma_1, \dots, \gamma_n), & n \geq 1,\ j = 0, \\ (\gamma_1, \dots, g_j,\ u(r(g_{j+1})),\ g_{j+1}, \dots, \gamma_n), & n \geq 2,\ 1 \leq j \leq n-1, \\ (\gamma_1, \dots, \gamma_n,\ u(s(\gamma_n))), & n \geq 1,\ j = n. \end{cases}$$

  These maps agree with Definition 2.2.6, written in left-to-right convention $s(g_i) = r(g_{i+1})$.
- **On morphisms.** For a continuous functor $\varphi := (\varphi_0, \varphi_1)\colon \mathcal{H} \to \mathcal{G}$, define $N\varphi$ by setting $\varphi_0\colon \mathcal{H}_0 \to \mathcal{G}_0$ and $\varphi_n\colon \mathcal{H}_n \to \mathcal{G}_n$, $(h_1, \dots, h_n) \mapsto (\varphi_1(h_1), \dots, \varphi_1(h_n))$ for $n \geq 1$.

**Proposition 2.3.2.** The assignment $N$ from Definition 2.3.1 defines a functor $N\colon \mathcal{G}\mathbf{Top} \to \mathbf{sTop}$.

*Proof.* Let $\varphi = (\varphi_0, \varphi_1)\colon \mathcal{H} \to \mathcal{G}$ be a continuous functor of topological groupoids.

- **Well-definedness and continuity on each level.** If $(h_1, \dots, h_n) \in \mathcal{H}_n$, then $s(h_i) = r(h_{i+1})$ for $1 \leq i < n$. Using $s \circ \varphi_1 = \varphi_0 \circ s$ and $r \circ \varphi_1 = \varphi_0 \circ r$, we obtain

$$s(\varphi_1(h_i)) = \varphi_0(s(h_i)) = \varphi_0(r(h_{i+1})) = r(\varphi_1(h_{i+1})),$$

  so $\varphi_n(h_1, \dots, h_n) \in \mathcal{G}_n$. Continuity follows because $\varphi_n$ is the restriction of the continuous map $\varphi_1^{\times n}\colon \mathcal{H}_1^n \to \mathcal{G}_1^n$ to the subspaces $\mathcal{H}_n \subseteq \mathcal{H}_1^n$ and $\mathcal{G}_n \subseteq \mathcal{G}_1^n$. The case $n = 0$ is the continuity of $\varphi_0$.

- **Commutation with face maps.** For $n = 1$ and $h \in \mathcal{H}_1$,

$$d_0(\varphi_1(h)) = s(\varphi_1(h)) = \varphi_0(s(h)) = \varphi_0(d_0(h)),$$
$$d_1(\varphi_1(h)) = r(\varphi_1(h)) = \varphi_0(r(h)) = \varphi_0(d_1(h)).$$

For $n \geq 2$ and $(h_1, \dots, h_n) \in \mathcal{H}_n$, the cases $i = 0$ and $i = n$ are immediate from the coordinate descriptions. For $1 \leq i \leq n-1$,

$$\begin{aligned} d_i(\varphi_n(h_1, \dots, h_n)) &= (\varphi_1(h_1), \dots, \varphi_1(h_i) \cdot \varphi_1(h_{i+1}), \dots, \varphi_1(h_n)) \\ &= (\varphi_1(h_1), \dots, \varphi_1(h_i \cdot h_{i+1}), \dots, \varphi_1(h_n)) \\ &= \varphi_{n-1}(d_i(h_1, \dots, h_n)), \end{aligned}$$

using $\varphi_1(h_i \cdot h_{i+1}) = \varphi_1(h_i) \cdot \varphi_1(h_{i+1})$.
- **Commutation with degeneracy maps.** For $n = 0$ and $x \in \mathcal{H}_0$,

$$\varphi_1(s_0(x)) = \varphi_1(u(x)) = u(\varphi_0(x)) = s_0(\varphi_0(x)).$$

For $n \geq 1$ and $(h_1, \dots, h_n) \in \mathcal{H}_n$, the cases $j = 0$ and $j = n$ follow from $\varphi_1 \circ u = u \circ \varphi_0$ and the identities $r \circ \varphi_1 = \varphi_0 \circ r$, $s \circ \varphi_1 = \varphi_0 \circ s$. For $1 \leq j \leq n-1$ and $n \geq 2$,

$$\begin{aligned} s_j(\varphi_n(h_1, \dots, h_n)) &= (\varphi_1(h_1), \dots, \varphi_1(h_j),\, u(r(\varphi_1(h_{j+1}))),\, \varphi_1(h_{j+1}), \dots, \varphi_1(h_n)) \\ &= (\varphi_1(h_1), \dots, \varphi_1(h_j),\, u(\varphi_0(r(h_{j+1}))),\, \varphi_1(h_{j+1}), \dots, \varphi_1(h_n)) \\ &= (\varphi_1(h_1), \dots, \varphi_1(h_j),\, \varphi_1(u(r(h_{j+1}))),\, \varphi_1(h_{j+1}), \dots, \varphi_1(h_n)) \\ &= \varphi_{n+1}(s_j(h_1, \dots, h_n)). \end{aligned}$$

- **Functoriality.** For the identity functor $\mathrm{id}_{\mathcal{G}}$ one has $N\mathrm{id}_{\mathcal{G}} = \mathrm{id}_{\mathcal{G}_\bullet}$ levelwise. If $\psi\colon \mathcal{G} \to \mathcal{K}$ and $\varphi\colon \mathcal{H} \to \mathcal{G}$, then for $n \geq 1$,

$$\begin{aligned} N(\psi \circ \varphi)_n(h_1, \dots, h_n) &= (\psi_1(\varphi_1(h_1)), \dots, \psi_1(\varphi_1(h_n))) \\ &= (N\psi_n \circ N\varphi_n)(h_1, \dots, h_n), \end{aligned}$$

and similarly for $n = 0$. Hence $N(\psi \circ \varphi) = N\psi \circ N\varphi$.

$\square$

*Remark* 2.3.3. For a topological groupoid, these structure maps agree with Definition 2.2.6, and Proposition 2.2.7 verifies the simplicial identities in this setting.

### 2.3.1 Classifying Spaces

The nerve $\mathcal{G}_\bullet$ records the algebraic structure of an étale groupoid $\mathcal{G}$ in a simplicial object, hence in a form amenable to homological constructions via compactly supported chains. At the same time, the nerve has an intrinsic topological meaning: it presents a canonical space whose homotopy type captures the global geometry encoded by $\mathcal{G}$. This space is the classifying space of the groupoid. There are two reasons why the classifying space is central in our context.

1. First, it provides the conceptual bridge between groupoid homology defined from the simplicial space $\mathcal{G}_\bullet$ and familiar invariants from algebraic topology. In particular, it clarifies why the simplicial machinery is the correct setting: the simplicial structure maps of the nerve are precisely what is needed to build a space by gluing simplices along faces, and the resulting homotopy type is functorial in $\mathcal{G}$.
2. Second, the classifying space is the natural recipient of Morita invariance: for étale groupoids, Morita equivalences induce weak homotopy equivalences on nerves, hence do not change the homotopy type of the associated classifying spaces. In principle Morita equivalent groupoids represent the same underlying geometric object.

In this section we recall the construction of the geometric realization $B\mathcal{G}$ and define the classifying space $B\mathcal{G}$ of a topological groupoid $\mathcal{G}$ as the realization of its nerve. We also record the functoriality properties needed later: continuous functors of groupoids induce continuous maps on classifying spaces, and the passage $\mathcal{G} \mapsto B\mathcal{G}$ is compatible with the nerve functor $N$ from Definition 2.3.1.

**Definition 2.3.4** (Classifying space [4, §1.7])**.** Let $\mathcal{G}_\bullet = (\mathcal{G}_n, (d_i)_{i=0}^n, (s_j)_{j=0}^n)_{n\in\mathbb{N}_0}$ be the nerve of a topological groupoid $\mathcal{G}$. The geometric realization of $\mathcal{G}_\bullet$ is the quotient

$$B\mathcal{G} = \Big(\bigsqcup_{n\geq 0} \mathcal{G}_n \times \Delta^n\Big)\Big/\sim,$$

where $\Delta^n := \{(t_0,\dots,t_n) \in [0,1]^{n+1} \mid \sum_{j=0}^n t_j = 1\}$ is the standard $n$-simplex, and $\sim$ is the smallest equivalence relation such that

$$\begin{aligned}
(d_i x, t) &\sim (x, \delta^i t) \quad \text{for all } n \geq 1,\ 0 \leq i \leq n,\ x \in \mathcal{G}_n,\ t \in \Delta^{n-1},\\
(s_i x, t) &\sim (x, \sigma^i t) \quad \text{for all } n \geq 0,\ 0 \leq i \leq n,\ x \in \mathcal{G}_n,\ t \in \Delta^{n+1},
\end{aligned}$$

where the coface and codegeneracy maps on simplices are

$$\begin{aligned}
\delta^i &\colon \Delta^{n-1} \to \Delta^n, \quad (t_0,\dots,t_{n-1}) \mapsto (t_0,\dots,t_{i-1},0,t_i,\dots,t_{n-1}),\\
\sigma^i &\colon \Delta^{n+1} \to \Delta^n, \quad (t_0,\dots,t_{n+1}) \mapsto (t_0,\dots,t_{i-1},\ t_i+t_{i+1},\ t_{i+2},\dots,t_{n+1}).
\end{aligned}$$

The space $B\mathcal{G}$ is called the classifying space of $\mathcal{G}$.

The classifying space $B\mathcal{G}$ classifies principal $\mathcal{G}$-bundles over paracompact spaces [4, §1.7]. In particular, if $\varphi : \mathcal{H} \xrightarrow{\sim} \mathcal{G}$ is a Morita equivalence, then $\mathcal{H}$ and $\mathcal{G}$ present the same quotient geometry, so one expects $B\mathcal{H}$ and $B\mathcal{G}$ to have the same weak homotopy type. Theorem 2.3.5 makes this precise by giving an explicit proof that the induced map $B\varphi : B\mathcal{H} \to B\mathcal{G}$ is a weak homotopy equivalence. The argument is organised around Quillen's Theorem A. For each $y \in \mathcal{G}_0$ we consider the comma groupoid $\varphi \downarrow y$, whose objects are arrows $\alpha : \varphi_0(x) \mapsto y$ in $\mathcal{G}$ and whose morphisms are the arrows in $\mathcal{H}$ making the evident triangles commute. We show that $\varphi \downarrow y$ has a terminal object, obtained from essential surjectivity and full faithfulness of $\varphi$ in Definition 1.4.10. From this we construct a simplicial strong deformation retraction of the nerve $N\varphi \downarrow y$ onto a vertex: letting $p : \varphi \downarrow y \to \mathbf{1}$ be the unique functor to the terminal category and choosing a section $s : \mathbf{1} \to \varphi \downarrow y$ picking the terminal object, we build a natural transformation $\eta : \mathrm{id} \Rightarrow s \circ p$, which induces a contraction after geometric realization. Hence $B\varphi \downarrow y$ is contractible for every $y$, and Quillen's Theorem A implies that $B\varphi$ is a weak homotopy equivalence. This has an immediate homological consequence: for every discrete abelian group $A$, the induced map $B\varphi$ yields isomorphisms

$$H_n^{\mathrm{sing}}(B\mathcal{H}; A) \xrightarrow{\cong} H_n^{\mathrm{sing}}(B\mathcal{G}; A) \quad \text{for all } n \geq 0.$$

Accordingly, any invariant extracted functorially from $B\mathcal{G}$ and invariant under weak homotopy equivalence is automatically Morita invariant. In the present thesis this provides the homotopical backdrop for Morita invariance: although our Moore homology is defined on the nerve via compactly supported chains rather than via singular chains on $B\mathcal{G}$, the weak equivalence $B\mathcal{H} \simeq B\mathcal{G}$ explains why Morita equivalent groupoids should encode the same homological information and guides the chain-level constructions establishing invariance in our setting.

**Theorem 2.3.5** ([6, Theorem 3.12, Lemma 4.3])**.** Let $\varphi : \mathcal{H} \xrightarrow{\sim} \mathcal{G}$ be a Morita equivalence in the sense of Definition 1.4.10. Then the induced map on classifying spaces $B\varphi : B\mathcal{H} \to B\mathcal{G}$ is a weak homotopy equivalence. This means that for every $k \geq 0$ and every basepoint $b \in B\mathcal{H}$, the induced map $\pi_k(B\mathcal{H}, b) \xrightarrow{\cong} \pi_k(B\mathcal{G}, B\varphi(b))$ is an isomorphism.

*Proof.* The strategy is to apply Quillen's Theorem A to the functor $\varphi : \mathcal{H} \to \mathcal{G}$. Concretely, for each $y \in \mathcal{G}_0$ we study the comma groupoid $\varphi \downarrow y$ and show that its classifying space is contractible; the conclusion then follows formally from Quillen's Theorem A.

- **Well-definedness and continuity of $B\varphi$.** Let $q_{\mathcal{H}} \colon \bigsqcup_{n\geq 0} \mathcal{H}_n \times \Delta^n \to B\mathcal{H}$ and $q_{\mathcal{G}} \colon \bigsqcup_{n\geq 0} \mathcal{G}_n \times \Delta^n \to B\mathcal{G}$ be the quotient maps. Define

  $$\Phi := \bigsqcup_{n\geq 0} (\varphi_n \times \mathrm{id}_{\Delta^n}) \colon \bigsqcup_{n\geq 0} \mathcal{H}_n \times \Delta^n \to \bigsqcup_{n\geq 0} \mathcal{G}_n \times \Delta^n,$$

  which is continuous by the coproduct topology. We show that $\Phi$ respects the generating relations of the geometric realizations. We use that $N\varphi$ is simplicial, this is verified next, and is built into Definition 2.3.1. If $(d_i x, t) \sim (x, \delta^i t)$ with $x \in \mathcal{H}_n$ and $t \in \Delta^{n-1}$, then

  $$\Phi(d_i x, t) = (\varphi_{n-1}(d_i x), t) = (d_i\, \varphi_n(x), t) \sim (\varphi_n(x), \delta^i t) = \Phi(x, \delta^i t).$$

  Similarly, if $(s_i x, t) \sim (x, \sigma^i t)$ with $x \in \mathcal{H}_n$ and $t \in \Delta^{n+1}$, then

  $$\Phi(s_i x, t) = (\varphi_{n+1}(s_i x), t) = (s_i\, \varphi_n(x), t) \sim (\varphi_n(x), \sigma^i t) = \Phi(x, \sigma^i t).$$

  Since $\sim$ is the smallest equivalence relation containing these generating relations, it follows that $q_{\mathcal{G}} \circ \Phi$ is constant on $q_{\mathcal{H}}$-equivalence classes. Hence $\Phi$ descends to a unique map $B\varphi \colon B\mathcal{H} \to B\mathcal{G}$ such that $B\varphi \circ q_{\mathcal{H}} = q_{\mathcal{G}} \circ \Phi$.
  By the universal property of quotient maps $B\varphi$ is continuous.

- **Simplicial maps $N\varphi$.** Let $\varphi = (\varphi_0, \varphi_1) : \mathcal{H} \to \mathcal{G}$ be a continuous functor of topological groupoids. For $n \geq 1$ define $\varphi_n : \mathcal{H}_n \to \mathcal{G}_n$, $(h_1, \dots, h_n) \mapsto (\varphi_1(h_1), \dots, \varphi_1(h_n))$.
  - **Well-definedness:** If $(h_1, \dots, h_n) \in \mathcal{H}_n$, then $s_{\mathcal{H}}(h_i) = r_{\mathcal{H}}(h_{i+1})$ for $1 \leq i < n$. Using $s_{\mathcal{G}} \circ \varphi_1 = \varphi_0 \circ s_{\mathcal{H}}$ and $r_{\mathcal{G}} \circ \varphi_1 = \varphi_0 \circ r_{\mathcal{H}}$, we obtain $s_{\mathcal{G}}(\varphi_1(h_i)) = \varphi_0(s_{\mathcal{H}}(h_i)) = \varphi_0(r_{\mathcal{H}}(h_{i+1})) = r_{\mathcal{G}}(\varphi_1(h_{i+1}))$, so $(\varphi_1(h_1), \dots, \varphi_1(h_n)) \in \mathcal{G}_n$. Hence $\varphi_n$ is well defined.
  - **Continuity:** The map $\varphi_1^{\times n} : \mathcal{H}_1^n \to \mathcal{G}_1^n$ is continuous. Since $\mathcal{H}_n \subseteq \mathcal{H}_1^n$ and $\mathcal{G}_n \subseteq \mathcal{G}_1^n$ carry the subspace topologies and $\varphi_1^{\times n}(\mathcal{H}_n) \subseteq \mathcal{G}_n$, the restriction $\varphi_n = \varphi_1^{\times n}|_{\mathcal{H}_n} : \mathcal{H}_n \to \mathcal{G}_n$ is continuous. The map $\varphi_0$ is continuous by assumption.
  - **Simpliciality.** By the formulas in Definition 2.2.6, and since $\varphi$ preserves multiplication and units, $\varphi_1(h \cdot_{\mathcal{H}} h') = \varphi_1(h) \cdot_{\mathcal{G}} \varphi_1(h')$, $\varphi_1(u_{\mathcal{H}}(x)) = u_{\mathcal{G}}(\varphi_0(x))$, one checks for all $n \geq 1$ and $0 \leq i \leq n$ that
$$d_i^{\mathcal{G}} \circ \varphi_n = \varphi_{n-1} \circ d_i^{\mathcal{H}},$$
$$s_i^{\mathcal{G}} \circ \varphi_n = \varphi_{n+1} \circ s_i^{\mathcal{H}}.$$
    - Case $i = 0$:
$$\begin{aligned} d_0^{\mathcal{G}}(\varphi_n(h_1, \dots, h_n)) &= d_0^{\mathcal{G}}(\varphi(h_1), \dots, \varphi(h_n)) \\ &= (\varphi(h_2), \dots, \varphi(h_n)) \\ &= \varphi_{n-1}(h_2, \dots, h_n) \\ &= \varphi_{n-1}(d_0^{\mathcal{H}}(h_1, \dots, h_n)). \end{aligned}$$
    - Case $1 \leq i \leq n-1$:
$$\begin{aligned} d_i^{\mathcal{G}}(\varphi_n(h_1, \dots, h_n)) &= d_i^{\mathcal{G}}(\varphi(h_1), \dots, \varphi(h_n)) \\ &= (\varphi(h_1), \dots, \varphi(h_{i+1}) \circ \varphi(h_i), \dots, \varphi(h_n)) \\ &= (\varphi(h_1), \dots, \varphi(h_{i+1} \circ h_i), \dots, \varphi(h_n)) \\ &= \varphi_{n-1}(h_1, \dots, h_{i+1} \circ h_i, \dots, h_n) \\ &= \varphi_{n-1}(d_i^{\mathcal{H}}(h_1, \dots, h_n)), \end{aligned}$$

      where the third equality uses functoriality of $\varphi$ for composition, namely $\varphi(h_{i+1} \circ h_i) = \varphi(h_{i+1}) \circ \varphi(h_i)$.
    - Case $i = n$:
$$\begin{aligned} d_n^{\mathcal{G}}(\varphi_n(h_1, \dots, h_n)) &= d_n^{\mathcal{G}}(\varphi(h_1), \dots, \varphi(h_n)) \\ &= (\varphi(h_1), \dots, \varphi(h_{n-1})) \\ &= \varphi_{n-1}(h_1, \dots, h_{n-1}) \\ &= \varphi_{n-1}(d_n^{\mathcal{H}}(h_1, \dots, h_n)). \end{aligned}$$

    Hence $d_i^{\mathcal{G}} \circ \varphi_n = \varphi_{n-1} \circ d_i^{\mathcal{H}}$ for all $0 \leq i \leq n$.

Then

$$
\begin{aligned}
s_i^{\mathcal{G}}(\varphi_n(h_1,\dots,h_n)) &= s_i^{\mathcal{G}}(\varphi(h_1),\dots,\varphi(h_n)) \\
&= (\varphi(h_1),\dots,\varphi(h_i),1_{\varphi(x_i)},\varphi(h_{i+1}),\dots,\varphi(h_n)) \\
&= (\varphi(h_1),\dots,\varphi(h_i),\varphi(1_{x_i}),\varphi(h_{i+1}),\dots,\varphi(h_n)) \\
&= \varphi_{n+1}(h_1,\dots,h_i,1_{x_i},h_{i+1},\dots,h_n) \\
&= \varphi_{n+1}(s_i^{\mathcal{H}}(h_1,\dots,h_n)),
\end{aligned}
$$

where the third equality uses functoriality of $\varphi$ for identities, namely $\varphi(1_{x_i}) = 1_{\varphi(x_i)}$. Thus $N\varphi = (\varphi_n)_{n\geq 0}$ is a simplicial map $N\mathcal{H} \to N\mathcal{G}$.

- **Induced map on geometric realizations:** Passing to geometric realizations yields a continuous map $B\varphi : B\mathcal{H} \to B\mathcal{G}$.
  1. **Definition on the disjoint unions.** Each $\varphi_n$ is continuous since it is the restriction of the continuous product map $\varphi_1^n : \mathcal{H}_1^n \to \mathcal{G}_1^n$. Define $F : \coprod_{n\geq 0} \mathcal{H}_n \times \Delta^n \to \coprod_{n\geq 0} \mathcal{G}_n \times \Delta^n$ by the restrictions $F_n : \mathcal{H}_n \times \Delta^n \to \mathcal{G}_n \times \Delta^n, F_n(x,t) := (\varphi_n(x), t)$. Each $F_n$ is continuous, hence $F$ is continuous.
  2. **Quotient model and target statement.** Let

$$
q_{\mathcal{H}} : \coprod_{n\geq 0} \mathcal{H}_n \times \Delta^n \to B\mathcal{H}, \quad q_{\mathcal{G}} : \coprod_{n\geq 0} \mathcal{G}_n \times \Delta^n \to B\mathcal{G}
$$

     be the quotient maps, so

$$
B\mathcal{H} = \Big(\coprod_{n\geq 0} \mathcal{H}_n \times \Delta^n\Big)/ \sim_{\mathcal{H}}, \quad B\mathcal{G} = \Big(\coprod_{n\geq 0} \mathcal{G}_n \times \Delta^n\Big)/ \sim_{\mathcal{G}},
$$

     where $\sim_{\mathcal{H}}$ is the smallest equivalence relation such that for all $n \geq 1$, all $0 \leq i \leq n$, all $x \in \mathcal{H}_n$, all $t \in \Delta^{n-1}$, $(d_i^{\mathcal{H}}x, t) \sim_{\mathcal{H}} (x, \delta^i t)$, and for all $n \geq 0$, all $0 \leq i \leq n$, all $x \in \mathcal{H}_n$, all $t \in \Delta^{n+1}$, $(s_i^{\mathcal{H}}x, t) \sim_{\mathcal{H}} (x, \sigma^i t)$, with $\delta^i$ and $\sigma^i$ as in Definition 2.3.1. The relation $\sim_{\mathcal{G}}$ is defined analogously with $d_i^{\mathcal{G}}$ and $s_i^{\mathcal{G}}$.
     It suffices to show that $F$ preserves $\sim_{\mathcal{H}}$ in the sense that $a \sim_{\mathcal{H}} b \Rightarrow F(a) \sim_{\mathcal{G}} F(b)$. Then there exists a unique map $B\varphi : B\mathcal{H} \to B\mathcal{G}$ such that $B\varphi \circ q_{\mathcal{H}} = q_{\mathcal{G}} \circ F$, hence the diagram commutes

$$
\begin{array}{ccc}
\coprod_{n\geq 0} \mathcal{H}_n \times \Delta^n & \xrightarrow{F} & \coprod_{n\geq 0} \mathcal{G}_n \times \Delta^n \\
{\scriptstyle q_{\mathcal{H}}}\downarrow & & \downarrow{\scriptstyle q_{\mathcal{G}}} \\
B\mathcal{H} & \xrightarrow[B\varphi]{} & B\mathcal{G}.
\end{array}
$$

  3. **Verification on generators of the relation.** Fix $n \geq 1$ and $0 \leq i \leq n$.
     - **For face generators:** Let $x \in \mathcal{H}_n$ and $t \in \Delta^{n-1}$. Then

$$
F(d_i^{\mathcal{H}}x, t) = (\varphi_{n-1}(d_i^{\mathcal{H}}x), t).
$$

Using the simplicial identities, $d_i^{\mathcal{G}} \circ \varphi_n = \varphi_{n-1} \circ d_i^{\mathcal{H}}$, this equals

$$(d_i^{\mathcal{G}}(\varphi_n(x)), t) \sim_{\mathcal{G}} (\varphi_n(x), \delta^i t) = F(x, \delta^i t),$$

hence $F(d_i^{\mathcal{H}} x, t) \sim_{\mathcal{G}} F(x, \delta^i t)$.

- **For degeneracy generators:** Fix $n \geq 0$ and $0 \leq i \leq n$, let $x \in \mathcal{H}_n$ and $t \in \Delta^{n+1}$. Then

$$F(s_i^{\mathcal{H}} x, t) = (\varphi_{n+1}(s_i^{\mathcal{H}} x), t).$$

Using the simplicial identities, $s_i^{\mathcal{G}} \circ \varphi_n = \varphi_{n+1} \circ s_i^{\mathcal{H}}$, hence

$$(s_i^{\mathcal{G}}(\varphi_n(x)), t) \sim_{\mathcal{G}} (\varphi_n(x), \sigma^i t) = F(x, \sigma^i t).$$

Therefore $F(s_i^{\mathcal{H}} x, t) \sim_{\mathcal{G}} F(x, \sigma^i t)$.

Thus $F$ preserves the generating relations of $\sim_{\mathcal{H}}$, hence $F$ preserves $\sim_{\mathcal{H}}$.

4. **Existence and continuity of the induced map.** By 3. the composite $q_{\mathcal{G}} \circ F$ is constant on $\sim_{\mathcal{H}}$ equivalence classes. Therefore there exists a unique map $B\varphi : B\mathcal{H} \to B\mathcal{G}$ such that $B\varphi \circ q_{\mathcal{H}} = q_{\mathcal{G}} \circ F$. This identity implies commutativity of the diagram in 2. The map $B\varphi$ is continuous since $q_{\mathcal{H}}$ is a quotient map and $q_{\mathcal{G}} \circ F$ is continuous.

- **The comma groupoid $\varphi \downarrow y$:** To apply Quillen's Theorem A to the map on classifying spaces induced by $\varphi$, we need a concrete model for the homotopy fibre of $B\varphi$ over each unit $y \in \mathcal{G}_0$. This model is provided by the comma groupoid $\varphi \downarrow y$, whose objects are arrows of $\mathcal{G}$ landing in $y$ with source in the image of $\varphi_0$, and whose morphisms record how these arrows vary along arrows of $\mathcal{H}$. Let $\varphi = (\varphi_0, \varphi_1) : \mathcal{H} \to \mathcal{G}$ be a functor of small categories and fix $y \in \mathcal{G}_0$. Define

$$\mathrm{Ob}\, \varphi \downarrow y := \{(x, \alpha) \in \mathcal{H}_0 \times \mathcal{G}_1 \mid s(\alpha) = \varphi_0(x),\ r(\alpha) = y\}.$$

Thus an object of the comma groupoid is an object $x \in \mathcal{H}_0$ together with an arrow $\alpha : \varphi_0(x) \mapsto y$ in $\mathcal{G}$. For $(x, \alpha), (x', \alpha') \in \mathrm{Ob}\, \varphi \downarrow y$ set

$$\mathrm{Hom}_{\varphi \downarrow y}((x, \alpha), (x', \alpha')) := \{h \in \mathcal{H}_1 \mid s(h) = x,\ r(h) = x',\ \alpha' \cdot \varphi_1(h) = \alpha\}.$$

Thus morphisms are those arrows $h : x \mapsto x'$ in $\mathcal{H}$ for which the triangle commutes in $\mathcal{G}$. The identity at $(x, \alpha)$ is $1_{(x,\alpha)} := 1_x \in \mathcal{H}_1$. Indeed $s(1_x) = x = r(1_x)$ and, since $\varphi$ preserves units, $\alpha \cdot \varphi_1(1_x) = \alpha \cdot 1_{\varphi_0(x)} = \alpha$, so $1_x \in \mathrm{Hom}_{\varphi \downarrow y}((x, \alpha), (x, \alpha))$. Given $h \in \mathrm{Hom}_{\varphi \downarrow y}((x, \alpha), (x', \alpha'))$ and $h' \in \mathrm{Hom}_{\varphi \downarrow y}((x', \alpha'), (x'', \alpha''))$, define composition by transport from $\mathcal{H}$, $h' \cdot_{\varphi\downarrow y} h := h' \cdot_{\mathcal{H}} h \in \mathcal{H}_1$. This is well defined because $s(h' \cdot h) = s(h) = x$ and $r(h' \cdot h) = r(h') = x''$, and functoriality of $\varphi_1$ gives

$$\begin{aligned}\alpha'' \cdot \varphi_1(h' \cdot h) &= \alpha'' \cdot (\varphi_1(h') \cdot \varphi_1(h)) \\ &= (\alpha'' \cdot \varphi_1(h')) \cdot \varphi_1(h) \\ &= \alpha' \cdot \varphi_1(h) \\ &= \alpha.\end{aligned}$$

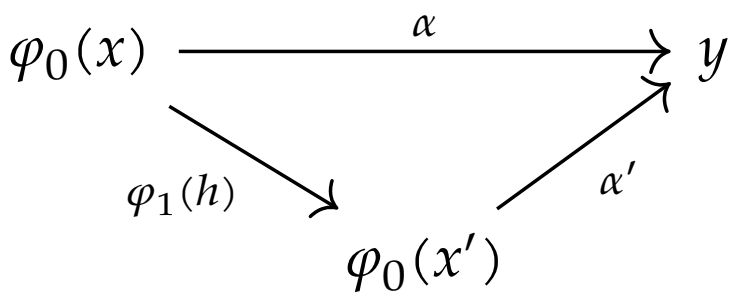

Associativity and unit laws hold because they hold in $\mathcal{H}$. If $\mathcal{H}$ and $\mathcal{G}$ are groupoids, then $\varphi \downarrow y$ is again a groupoid: if $\alpha' \cdot \varphi_1(h) = \alpha$, then multiplying on the right by $\varphi_1(h)^{-1}$ yields

$$\begin{aligned}\alpha \cdot \varphi_1(h^{-1}) &= (\alpha' \cdot \varphi_1(h)) \cdot \varphi_1(h)^{-1} \\ &= \alpha' \cdot 1_{\varphi_0(x')} = \alpha',\end{aligned}$$

so $h^{-1} \in \mathrm{Hom}_{\varphi\downarrow y}((x', \alpha'), (x, \alpha))$.
Finally, $\varphi \downarrow y$ is small since $\mathrm{Ob}\,\varphi \downarrow y \subseteq \mathcal{H}_0 \times \mathcal{G}_1$ and $\mathrm{Hom}_{\varphi\downarrow y}((x, \alpha), (x', \alpha')) \subseteq \mathcal{H}_1$.

- **Terminal object in $\varphi \downarrow y$.** Fix $y \in \mathcal{G}_0$. To show that the relevant comma groupoid has contractible classifying space, we produce a terminal object and then contract its nerve. Since $\varphi : \mathcal{H} \xrightarrow{\sim} \mathcal{G}$ is a Morita equivalence as in Definition 1.4.10, the essential surjectivity condition yields an element $x_0 \in \mathcal{H}_0$ and an arrow $g \in \mathcal{G}_1$ such that $r(g) = \varphi_0(x_0)$, $s(g) = y$. Because $\mathcal{G}$ is a groupoid, $g$ is invertible. Then $s(g^{-1}) = r(g) = \varphi_0(x_0)$ and $r(g^{-1}) = s(g) = y$, so $(x_0, g^{-1}) \in \mathrm{Ob}\,\varphi \downarrow y$. We claim that $(x_0, g^{-1})$ is terminal in $\varphi \downarrow y$. Let $(x, \alpha) \in \mathrm{Ob}\,\varphi \downarrow y$, so $s(\alpha) = \varphi_0(x)$ and $r(\alpha) = y$. Observe that $s((g^{-1})^{-1}) = r(g^{-1}) = y = r(\alpha)$. Moreover, using the source/range rules for left-to-right multiplication,

$$s((g^{-1})^{-1} \cdot \alpha) = s(\alpha) = \varphi_0(x), \quad r((g^{-1})^{-1} \cdot \alpha) = r((g^{-1})^{-1}) = s(g^{-1}) = \varphi_0(x_0),$$

so $(g^{-1})^{-1} \cdot \alpha \in \mathcal{G}_1(\varphi_0(x), \varphi_0(x_0))$. By the fully faithful condition in Definition 1.4.10, the map $\varphi_1$ induces a bijection $\mathcal{H}_1(x, x_0) \xrightarrow{\cong} \mathcal{G}_1(\varphi_0(x), \varphi_0(x_0))$. Hence there exists a unique $h \in \mathcal{H}_1(x, x_0)$ with $\varphi_1(h) = (g^{-1})^{-1} \cdot \alpha$. Then

$$\begin{aligned}g^{-1} \cdot \varphi_1(h) &= g^{-1} \cdot ((g^{-1})^{-1} \cdot \alpha) \\ &= (g^{-1} \cdot (g^{-1})^{-1}) \cdot \alpha \\ &= 1_y \cdot \alpha \\ &= \alpha,\end{aligned}$$

so $h \in \mathrm{Hom}_{\varphi\downarrow y}((x, \alpha), (x_0, g^{-1}))$. If $k : (x, \alpha) \mapsto (x_0, g^{-1})$ is any other morphism in $\varphi \downarrow y$, then $g^{-1} \cdot \varphi_1(k) = \alpha = g^{-1} \cdot \varphi_1(h)$. Left-multiplying by $(g^{-1})^{-1}$ gives $\varphi_1(k) = \varphi_1(h)$, and by injectivity of $\varphi_1 : \mathcal{H}_1(x, x_0) \to \mathcal{G}_1(\varphi_0(x), \varphi_0(x_0))$ we conclude $k = h$. Thus every object admits a unique morphism to $(x_0, g^{-1})$, and $(x_0, g^{-1})$ is terminal.

- **A simplicial deformation retraction of $N\varphi{\downarrow}y$ onto a point.** To apply Quillen's Theorem A to $B\varphi$, we must understand the comma groupoids $\varphi{\downarrow}y$. Since we have exhibited a terminal object $T$ in $\varphi \downarrow y$, it is enough to show directly that $B\varphi \downarrow y$ is contractible. We do this constructively by writing down an explicit simplicial homotopy which contracts the nerve $N\varphi{\downarrow}y$ to the 0-simplex corresponding to $T$. Fix $y \in \mathcal{G}_0$ and let $T := (x_0, g^{-1}) \in \operatorname{Ob}\varphi{\downarrow}y$ be the terminal object constructed above. Let $\mathbf{1}$ be the terminal category, $p : \varphi{\downarrow}y \to \mathbf{1}$ the unique functor, and $s : \mathbf{1} \to \varphi{\downarrow}y$ the section with $s(\star) = T$, where $\star$ is the unique object of $\mathbf{1}$. There is a natural transformation $\eta : \mathrm{id}_{\varphi\downarrow y} \Rightarrow s \circ p$ whose component $\eta_{(x,\alpha)} : (x, \alpha) \to T$ is the unique morphism; in particular $\eta_T = \mathrm{id}_T$. We identify $N\mathbf{1} = \Delta^0$ and set

$$r := Np : N\varphi{\downarrow}y \to \Delta^0, \quad i := Ns : \Delta^0 \to N\varphi{\downarrow}y.$$

Then $r \circ i = \mathrm{id}_{\Delta^0}$ and $i \circ r = N(s \circ p)$. For each $n \geq 0$ and $0 \leq k \leq n$ define a map $H_k : N\varphi{\downarrow}y_n \to N\varphi{\downarrow}y_{n+1}$ as follows: given an $n$-simplex $\mathbf{a}$ we define $H_k(\mathbf{a})$ by

$$\begin{aligned} \mathbf{a} &= (a_0 \xrightarrow{h_1} a_1 \xrightarrow{h_2} \cdots \xrightarrow{h_n} a_n) \quad \in N\varphi{\downarrow}y_n, \\ H_k(\mathbf{a}) &= (a_0 \xrightarrow{h_1} \cdots \xrightarrow{h_k} a_k \xrightarrow{\eta_{a_k}} T \xrightarrow{\mathrm{id}_T} \cdots \xrightarrow{\mathrm{id}_T} T). \end{aligned}$$

Concretely, we insert $\eta_{a_k} : a_k \to T$ after the $k$-th arrow and then fill the remaining arrows with identities of $T$, which is consistent because $(s \circ p)$ sends every morphism to $\mathrm{id}_T$. We now verify that $(H_k)_{0 \leq k \leq n}$ is a simplicial homotopy from $i \circ r$ to $N\mathrm{id}_{\varphi\downarrow y}$. The required identities are: $d_0 H_0 = i \circ r$ and $d_{n+1} H_n = N\mathrm{id}_{\varphi\downarrow y}$, and for all $n \geq 0$ and $0 \leq k \leq n$,

$$d_\ell H_k = \begin{cases} H_{k-1} d_\ell, & \text{for } 0 < \ell \leq k, \\ H_k d_{\ell-1}, & \text{for } k < \ell \leq n+1, \end{cases} \qquad s_\ell H_k = \begin{cases} H_{k+1} s_\ell, & \text{for } 0 \leq \ell \leq k, \\ H_k s_{\ell-1}, & \text{for } k < \ell \leq n. \end{cases}$$

  - **Faces.** For $\mathbf{b} = (b_0 \xrightarrow{g_1} \cdots \xrightarrow{g_{n+1}} b_{n+1}) \in N\varphi{\downarrow}y_{n+1}$ and for $1 \leq \ell \leq n$,

    $$\begin{aligned} d_0(\mathbf{b}) &= (b_1 \xrightarrow{g_2} \cdots \xrightarrow{g_{n+1}} b_{n+1}), \\ d_{n+1}(\mathbf{b}) &= (b_0 \xrightarrow{g_1} \cdots \xrightarrow{g_n} b_n), \\ d_\ell(\mathbf{b}) &= (b_0 \xrightarrow{g_1} \cdots \xrightarrow{g_\ell g_{\ell+1}} \cdots \xrightarrow{g_{n+1}} b_{n+1}). \end{aligned}$$

    Let $\mathbf{a} \in N\varphi{\downarrow}y_n$. First,

    $$\begin{aligned} d_0 H_0(\mathbf{a}) &= d_0(a_0 \xrightarrow{\eta_{a_0}} T \xrightarrow{\mathrm{id}_T} \cdots \xrightarrow{\mathrm{id}_T} T) = (T \xrightarrow{\mathrm{id}_T} \cdots \xrightarrow{\mathrm{id}_T} T) = (i \circ r)(\mathbf{a}), \\ d_{n+1} H_n(\mathbf{a}) &= d_{n+1}(a_0 \xrightarrow{h_1} \cdots \xrightarrow{h_n} a_n \xrightarrow{\eta_{a_n}} T) = (a_0 \xrightarrow{h_1} \cdots \xrightarrow{h_n} a_n) = \mathbf{a}. \end{aligned}$$

    Now fix $0 \leq k \leq n$ and $1 \leq \ell \leq n+1$.
    If $1 \leq \ell < k$, $d_\ell$ composes two arrows among $h_1, \ldots, h_k$ and does not involve $\eta_{a_k}$, hence

    $$d_\ell H_k(\mathbf{a}) = H_{k-1}\, d_\ell(\mathbf{a}).$$

If $\ell = k$, then $d_k$ composes $h_k$ with the inserted arrow $\eta_{a_k}\colon a_{k-1} \xrightarrow{h_k} a_k \xrightarrow{\eta_{a_k}} T$. We have $\eta_{a_k} \circ h_k = \eta_{a_{k-1}}$, hence

$$d_k H_k(\mathbf{a}) = (a_0 \xrightarrow{h_1} \cdots \xrightarrow{h_{k-1}} a_{k-1} \xrightarrow{\eta_{a_{k-1}}} T \xrightarrow{\mathrm{id}_T} \cdots \xrightarrow{\mathrm{id}_T} T) = H_{k-1}\, d_k(\mathbf{a}).$$

If $\ell = k+1$, then $d_{k+1}$ composes $\eta_{a_k}$ with the next arrow $\mathrm{id}_T$, so $\mathrm{id}_T \circ \eta_{a_k} = \eta_{a_k}$ and

$$d_{k+1} H_k(\mathbf{a}) = H_k\, d_k(\mathbf{a}).$$

If $k+1 < \ell \leq n+1$, then $d_\ell$ only acts in the tail of $H_k(\mathbf{a})$, where all arrows are $\mathrm{id}_T$, so $d_\ell H_k(\mathbf{a}) = H_k(\mathbf{a})$. Also $d_{\ell-1}$ acts on $\mathbf{a}$ only at indices strictly larger than $k$, hence it does not change $a_k$, so

$$d_\ell H_k(\mathbf{a}) = H_k\, d_{\ell-1}(\mathbf{a}).$$

- **Degeneracies.** For $\mathbf{b} = (b_0 \xrightarrow{g_1} \cdots \xrightarrow{g_{n+1}} b_{n+1}) \in N\varphi\!\downarrow\! y_{n+1}$, the degeneracy $s_\ell$ inserts an identity at $b_\ell$:

$$s_\ell(\mathbf{b}) = (b_0 \xrightarrow{g_1} \cdots \xrightarrow{g_\ell} b_\ell \xrightarrow{\mathrm{id}_{b_\ell}} b_\ell \xrightarrow{g_{\ell+1}} \cdots \xrightarrow{g_{n+1}} b_{n+1}), \quad \text{for } 0 \leq \ell \leq n+1.$$

Let $\mathbf{a} \in N\varphi\!\downarrow\! y_n$. If $0 \leq \ell \leq k$, then $s_\ell$ inserts an identity arrow somewhere in the initial segment $a_0 \xrightarrow{h_1} \cdots \xrightarrow{h_k} a_k$ before the arrow $\eta_{a_k}$. Thus the vertex at which we insert $\eta$ shifts from position $k$ to position $k+1$, and we get

$$s_\ell H_k(\mathbf{a}) = H_{k+1}\, s_\ell(\mathbf{a}).$$

If $k < \ell \leq n$, then $s_\ell$ inserts an identity arrow in the part of $\mathbf{a}$ strictly after $a_k$. In $H_k(\mathbf{a})$ this entire part has been replaced by the constant tail of identities of $T$, so inserting such an identity commutes with inserting $\eta_{a_k}$ and only changes indices in the evident way:

$$s_\ell H_k(\mathbf{a}) = H_k\, s_{\ell-1}(\mathbf{a}).$$

- **Contractibility of $B\varphi\!\downarrow\! y$.** For each $n \geq 0$, the simplicial set $\Delta^0$ has a unique $n$-simplex. Its image under $i = Ns : \Delta^0 \to N\varphi\!\downarrow\! y$ is the degenerate $n$-simplex

$$T \xrightarrow{\mathrm{id}_T} \cdots \xrightarrow{\mathrm{id}_T} T \quad \in N\varphi\!\downarrow\! y_n.$$

Since $\eta_T = \mathrm{id}_T$, inserting $\eta$ after the $k$-th arrow does not change this simplex, except for adding one more identity arrow. Concretely, for all $n \geq 0$ and $0 \leq k \leq n$, $H_k \circ i_n = s_k \circ i_n$. In particular, the simplicial homotopy $H$ is stationary on the image of $i$. Therefore $(r, i, H)$ is a simplicial strong deformation retraction of $N\varphi\!\downarrow\! y$ onto the vertex $i(\star) = Ns(\star)$, that is, the 0-simplex corresponding to $T$. Thus, $B\varphi\!\downarrow\! y$ is contractible.

- **Strong deformation retraction** $(Br, Bi, BH)$. To pass from the simplicial homotopy $(H_k)_{0\le k\le n}$ to an ordinary homotopy on geometric realizations, we use the standard prism decomposition of $\Delta^n \times [0,1]$, see [17, Chapter 14]. For $n \ge 0$ and $0 \le k \le n$ define

$$\lambda_k : \Delta^n \times [0,1] \to \Delta^{n+1},$$
$$(u,\tau) \mapsto ((1-\tau)u_0, \dots, (1-\tau)u_k,\ \tau,\ (1-\tau)u_{k+1}, \dots, (1-\tau)u_n),$$

and let

$$R_k := \Big\{(u,\tau) \in \Delta^n \times [0,1] \mid \sum_{j<k} u_j \le 1-\tau \le \sum_{j\le k} u_j\Big\}.$$

Then $(R_k)_{k=0}^n$ covers $\Delta^n \times [0,1]$, because for each $(u,\tau)$ the number $1-\tau \in [0,1]$ lies in the interval $\big[\sum_{j<k} u_j,\ \sum_{j\le k} u_j\big]$ for some $k$, since the partial sums $\sum_{j\le k} u_j$ form an increasing chain from 0 to 1. Also $\lambda_k(u,\tau) \in \Delta^{n+1}$ for every $(u,\tau)$, since all coordinates are nonnegative and

$$\sum_{i=0}^{n+1} \lambda_k(u,\tau)_i = (1-\tau)\sum_{i=0}^{n} u_i + \tau = (1-\tau)\cdot 1 + \tau = 1.$$

Let

$$q : \bigsqcup_{n\ge 0} N\varphi{\downarrow}y_n \times \Delta^n \to B\varphi{\downarrow}y$$

be the realization quotient map. For each $n \ge 0$ define a map

$$\widehat{H}_n : N\varphi{\downarrow}y_n \times \Delta^n \times [0,1] \to N\varphi{\downarrow}y_{n+1} \times \Delta^{n+1}$$

by the piecewise rule

$$\widehat{H}_n(x,u,\tau) := (H_k(x),\ \lambda_k(u,\tau)) \quad \text{whenever} \quad (u,\tau) \in R_k.$$

Set $\widehat{H} := \bigsqcup_{n\ge 0} \widehat{H}_n$.

  - **Compatibility on overlaps in the simplex factor by explicit expansion.** The only overlaps are $R_k \cap R_{k+1}$. On such an overlap we compare the two affine formulas $\lambda_k$ and $\lambda_{k+1}$ using the standard degeneracy map

$$\sigma^{k+1} : \Delta^{n+1} \to \Delta^n, \quad \sigma^{k+1}(t_0, \dots, t_{n+1}) := (t_0, \dots, t_k,\ t_{k+1} + t_{k+2},\ t_{k+3}, \dots, t_{n+1}).$$

    For every $(u,\tau) \in \Delta^n \times [0,1]$ one computes

$$\begin{aligned}
\lambda_k(u,\tau) &= ((1-\tau)u_0, \dots, (1-\tau)u_k,\ \tau,\ (1-\tau)u_{k+1},\ (1-\tau)u_{k+2}, \dots, (1-\tau)u_n),\\
\sigma^{k+1}\lambda_k(u,\tau) &= ((1-\tau)u_0, \dots, (1-\tau)u_k,\ \tau + (1-\tau)u_{k+1},\ (1-\tau)u_{k+2}, \dots, (1-\tau)u_n),\\
\lambda_{k+1}(u,\tau) &= ((1-\tau)u_0, \dots, (1-\tau)u_{k+1},\ \tau,\ (1-\tau)u_{k+2}, \dots, (1-\tau)u_n),\\
\sigma^{k+1}\lambda_{k+1}(u,\tau) &= ((1-\tau)u_0, \dots, (1-\tau)u_k,\ (1-\tau)u_{k+1} + \tau,\ (1-\tau)u_{k+2}, \dots, (1-\tau)u_n).
\end{aligned}$$

The two $n$-tuples are identical term by term, so

$$\sigma^{k+1}\lambda_k(u,\tau) = \sigma^{k+1}\lambda_{k+1}(u,\tau) \quad \text{for all} \quad (u,\tau) \in \Delta^n \times [0,1]. \tag{2.3.1}$$

- **Compatibility on overlaps in the nerve factor by explicit expansion.** Let $x \in N\varphi{\downarrow}y_n$ be written as $x = (a_0 \xrightarrow{h_1} a_1 \xrightarrow{h_2} \cdots \xrightarrow{h_n} a_n)$. Then

$$H_k(x) = (a_0 \xrightarrow{h_1} \cdots \xrightarrow{h_k} a_k \xrightarrow{\eta_{a_k}} T \xrightarrow{\mathrm{id}_T} \cdots \xrightarrow{\mathrm{id}_T} T),$$
$$H_{k+1}(x) = (a_0 \xrightarrow{h_1} \cdots \xrightarrow{h_k} a_k \xrightarrow{h_{k+1}} a_{k+1} \xrightarrow{\eta_{a_{k+1}}} T \xrightarrow{\mathrm{id}_T} \cdots \xrightarrow{\mathrm{id}_T} T).$$

Apply the face operator $d_{k+1}$ on $N\varphi{\downarrow}y_{n+1}$, which composes the $(k+1)$-st and $(k+2)$-nd arrows. For $H_k(x)$ the relevant consecutive pair is $a_k \xrightarrow{\eta_{a_k}} T \xrightarrow{\mathrm{id}_T} T$, so

$$\begin{aligned} d_{k+1}H_k(x) &= (a_0 \xrightarrow{h_1} \cdots \xrightarrow{h_k} a_k \xrightarrow{\mathrm{id}_T \circ \eta_{a_k}} T \xrightarrow{\mathrm{id}_T} \cdots \xrightarrow{\mathrm{id}_T} T) \\ &= (a_0 \xrightarrow{h_1} \cdots \xrightarrow{h_k} a_k \xrightarrow{\eta_{a_k}} T \xrightarrow{\mathrm{id}_T} \cdots \xrightarrow{\mathrm{id}_T} T). \end{aligned}$$

For $H_{k+1}(x)$ the relevant consecutive pair is $a_{k+1} \xrightarrow{\eta_{a_{k+1}}} T \xrightarrow{\mathrm{id}_T} T$, so

$$\begin{aligned} d_{k+1}H_{k+1}(x) &= (a_0 \xrightarrow{h_1} \cdots \xrightarrow{h_k} a_k \xrightarrow{h_{k+1}} a_{k+1} \xrightarrow{\mathrm{id}_T \circ \eta_{a_{k+1}}} T \xrightarrow{\mathrm{id}_T} \cdots \xrightarrow{\mathrm{id}_T} T) \\ &= (a_0 \xrightarrow{h_1} \cdots \xrightarrow{h_k} a_k \xrightarrow{h_{k+1}} a_{k+1} \xrightarrow{\eta_{a_{k+1}}} T \xrightarrow{\mathrm{id}_T} \cdots \xrightarrow{\mathrm{id}_T} T). \end{aligned}$$

Now apply the simplicial homotopy identity already verified for faces in the special case $\ell = k+1$, $d_{k+1}H_{k+1} = H_k d_{k+1}$, $d_{k+1}H_k = H_k d_k$, and expand $d_{k+1}(x)$ and $d_k(x)$ as explicit strings. The only nontrivial composition that appears is the composition of $h_{k+1} : a_k \mapsto a_{k+1}$ with $\eta_{a_{k+1}} : a_{k+1} \mapsto T$. By the naturality identity $\eta_{a_k} = \eta_{a_{k+1}} \circ h_{k+1}$, this composite equals $\eta_{a_k}$. Therefore the two $(n+1)$-simplices $d_{k+1}H_k(x)$ and $d_{k+1}H_{k+1}(x)$ become literally identical strings after expanding the unique composite, hence

$$d_{k+1}H_k(x) = d_{k+1}H_{k+1}(x). \tag{2.3.2}$$

- **Descent of $\widehat{H}$ to the realization.** The geometric realization identifies $(d_i z, u)$ with $(z, \delta^i u)$ and $(s_i z, u)$ with $(z, \sigma^i u)$. On $R_k \cap R_{k+1}$, the two representatives

$$(H_k(x), \lambda_k(u,\tau)), \quad (H_{k+1}(x), \lambda_{k+1}(u,\tau))$$

map under $\sigma^{k+1}$ in the simplex factor to the same $n$-tuple by (2.3.1). They also map under $d_{k+1}$ in the nerve factor to the same $(n+1)$-simplex by (2.3.2). Thus they represent the same point in the quotient defining $B\varphi{\downarrow}y$. Consequently $\widehat{H}$ respects the generating relations of the realization and descends to a unique continuous map

$$BH : B\varphi{\downarrow}y \times [0,1] \to B\varphi{\downarrow}y \quad \text{with} \quad BH \circ (q \times \mathrm{id}) = q \circ \widehat{H}.$$

- **Endpoints and stationarity.** For $\tau = 0$ one has

  $$\lambda_k(u, 0) = (u_0, \dots, u_k, 0, u_{k+1}, \dots, u_n),$$

  so $\lambda_k(u, 0)$ is obtained from $u \in \Delta^n$ by inserting a zero coordinate. In the realization this is identified with $u$ via the standard face inclusion that inserts 0 in the $(k+1)$-st coordinate. Therefore $BH(-, 0) = B(i \circ r) = Bi \circ Br$.
  For $\tau = 1$ one has $\lambda_k(u, 1) = e_{k+1} \in \Delta^{n+1}$, hence $\widehat{H}_n(x, u, 1)$ represents the $(k+1)$-st vertex of the simplex supporting the $(n+1)$-simplex $H_k(x)$. By the defining simplicial homotopy identities, the realization $BH$ is a homotopy with endpoints $BH(-, 0) = Bi \circ Br, BH(-, 1) = \mathrm{id}_{B\varphi\downarrow y}$. Also $Br \circ Bi = \mathrm{id}_{B\Delta^0}$ since $p \circ s = \mathrm{id}_{\mathbf{1}}$. Stationarity of $H$ at $i(\star)$ means $H_k(i(\star)) = s_k(i(\star))$ for all $k$, hence $\widehat{H}$ fixes $i(\star)$ for all $\tau$, and therefore $BH(i(\star), \tau) = i(\star)$. Thus $BH$ is a strong deformation retraction of $B\varphi\downarrow y$ onto $i(\star)$.
- **Continuous homotopy.** The preceding construction is the standard realization of a simplicial homotopy, see [17, Chapter 14]. Alternatively, one may appeal to the realization lemma for simplicial homotopies, see [8, Chapter I, §4] or [16, §11]. We record the endpoint computations.
  Let $BH$ be the unique map with $BH \circ (q \times \mathrm{id}) = q \circ \widehat{H}$. Fix $n \geq 0$, $x \in N\varphi\downarrow y_n$, $u \in \Delta^n$. At $\tau = 0$ we have $(u, 0) \in R_n$, hence

  $$BH(q(x, u), 0) = q(H_n(x), \lambda_n(u, 0)) = q(H_n(x), (u_0, \dots, u_n, 0)) = q(H_n(x), \delta^{n+1}u).$$

  By the defining realization relation $(d_{n+1}z, t) \sim (z, \delta^{n+1}t)$ we obtain

  $$q(H_n(x), \delta^{n+1}u) = q(d_{n+1}H_n(x), u).$$

  Since $d_{n+1}H_n = \mathrm{id}_{N\varphi\downarrow y_n}$, this gives $BH(q(x, u), 0) = q(x, u)$, hence $BH(-, 0) = \mathrm{id}_{B\varphi\downarrow y}$.
  At $\tau = 1$ we have $\lambda_0(u, 1) = e_1$, and $(u, 1) \in R_0$, hence

  $$BH(q(x, u), 1) = q(H_0(x), \lambda_0(u, 1)) = q(H_0(x), e_1).$$

  The vertex $e_1 \in \Delta^{n+1}$ is $\delta^0(\star)$ where $\star \in \Delta^0$. Thus, again by the defining realization relation $(d_0z, t) \sim (z, \delta^0 t)$, $q(H_0(x), e_1) = q(H_0(x), \delta^0\star) = q(d_0H_0(x), \star)$. Since $d_0H_0 = i \circ r$,

  $$BH(q(x, u), 1) = q((i \circ r)(x), \star) = (Bi \circ Br)(q(x, u)), \quad \text{hence} \quad BH(-, 1) = Bi \circ Br.$$

  Since $p \circ s = \mathrm{id}_{\mathbf{1}}$, we have $r \circ i = \mathrm{id}_{\Delta^0}$ and hence $Br \circ Bi = \mathrm{id}_{B\Delta^0}$. Moreover, $\eta_T = \mathrm{id}_T$ implies $H_k \circ i_n = s_k \circ i_n$ for all $n \geq 0$ and $0 \leq k \leq n$, so $BH(i(\star), \tau) = i(\star)$ for all $\tau \in [0, 1]$. Consequently, $(Br, Bi, BH^{\mathrm{rev}})$ is a strong deformation retraction of $B\varphi\downarrow y$ onto $Bi(B\Delta^0) = \{i(\star)\}$, where $BH^{\mathrm{rev}}(z, \tau) := BH(z, 1 - \tau)$.

- **Quillen's Theorem A.** Let $\varphi : \mathcal{H} \xrightarrow{\sim} \mathcal{G}$ be a Morita equivalence in the sense of Definition 1.4.10. By the previous step, $\varphi$ induces a simplicial map $N\varphi : N\mathcal{H} \to N\mathcal{G}$, hence a continuous map on classifying spaces

$$B\varphi : B\mathcal{H} \to B\mathcal{G}.$$

Fix $y \in \mathcal{G}_0$. Consider the comma groupoid $\varphi \downarrow y$ constructed above: its objects are pairs $(x, \alpha)$ with $x \in \mathcal{H}_0$ and $\alpha : \varphi_0(x) \to y$ in $\mathcal{G}$, and its morphisms are those $h : x \mapsto x'$ in $\mathcal{H}$ satisfying $\alpha' \cdot \varphi_1(h) = \alpha$. We have shown that $\varphi \downarrow y$ admits a terminal object $T = (x_0, \alpha_0)$, obtained from essential surjectivity and full faithfulness of $\varphi$. We then constructed an explicit simplicial strong deformation retraction of $N(\varphi \downarrow y)$ onto the vertex corresponding to $T$, hence $B(\varphi \downarrow y)$ is contractible. Since $y \in \mathcal{G}_0$ was arbitrary, $B(\varphi \downarrow y)$ is contractible for every $y \in \mathcal{G}_0$. Quillen's Theorem A [19, Theorem A] applies to the functor $\varphi : \mathcal{H} \to \mathcal{G}$ and yields that the induced map on classifying spaces $B\varphi : B\mathcal{H} \to B\mathcal{G}$ is a weak homotopy equivalence. This means that for every $k \geq 0$ and every basepoint $b \in B\mathcal{H}$ the induced map on homotopy groups is an isomorphism

$$\pi_k(B\mathcal{H}, b) \xrightarrow{\cong} \pi_k(B\mathcal{G}, B\varphi(b)).$$

□

Theorem 2.3.5 identifies Morita equivalence as the correct homotopical notion of equivalence at the level of classifying spaces: if $\varphi : \mathcal{H} \xrightarrow{\sim} \mathcal{G}$ is Morita, then $B\mathcal{H} \simeq B\mathcal{G}$. In particular, any invariant that is functorially extracted from $B\mathcal{G}$ and invariant under weak homotopy equivalence is automatically Morita invariant. This is a useful benchmark for our later constructions: it clarifies that Morita equivalent groupoids are indistinguishable from the viewpoint of the homotopy type of the quotient geometry encoded by the nerve. At the same time, the Moore chain complex $C_c(\mathcal{G}_\bullet, A)$ is built from compactly supported, for discrete $A$ locally constant $A$-valued functions on the levels $\mathcal{G}_n$, and it is generally not the singular chain complex of $B\mathcal{G}$, see Example 2.3.16. Consequently, Morita invariance of $B\mathcal{G}$ does not by itself imply Morita invariance of Moore homology: there is no a priori mechanism that identifies $C_c(\mathcal{H}_\bullet, A)$ with $C_c(\mathcal{G}_\bullet, A)$ from the weak equivalence $B\mathcal{H} \simeq B\mathcal{G}$. What Theorem 2.3.5 does provide is the correct conceptual interpretation: any discrepancy between Moore homology and singular homology of $B\mathcal{G}$ is not an artefact of the Morita representative, but an intrinsic feature of the compactly supported Moore model. In particular, if one seeks a comparison map from Moore homology to $H_\bullet^{\mathrm{sing}}(B\mathcal{G}; A)$, then the theorem forces such a comparison, when it exists, to be Morita invariant on the $B\mathcal{G}$-side, while the Moore side requires independent invariance arguments at chain level.

Theorem 2.3.5 relates to Moore homology in two ways:

1. It justifies Morita equivalence as the natural notion of sameness for ample groupoids: it preserves the quotient geometry, so any homology theory intended to reflect that geometry should be Morita invariant.

2. It clarifies what is specific to Moore homology: since our chains live in $C_c(\mathcal{G}_\bullet, A)$ rather than in the singular complex of $B\mathcal{G}$, Morita invariance cannot be deduced from $B\mathcal{H} \simeq B\mathcal{G}$ alone, but must be proved by explicit chain-level constructions. This is precisely the point of the pushforward formalism and the exact-sequence machinery developed later.

The Moore complex is tailored to ample groupoids. If $\mathcal{G}$ is ample and étale, then each $\mathcal{G}_n$ is locally compact and totally disconnected, hence locally constant compactly supported functions are generated by characteristic functions of compact open subsets. Consequently, $C_c(\mathcal{G}_n, \mathbb{Z})$ admits a concrete description in terms of compact open pieces of $\mathcal{G}_n$, and the boundary maps are given by pushforward along the face maps, which are local homeomorphisms in the étale setting. This yields a homology theory that is computationally tractable and stable under the operations that are natural for ample groupoids: restriction to clopen saturated reductions, compatibility with long exact sequences, and Mayer–Vietoris sequences for clopen covers. In this sense, the Moore complex encodes the locally compact, totally disconnected geometry that is central in the ample setting but largely invisible to singular chains on $B\mathcal{G}$.

This leads to concrete open questions:

*Question* 1. What are criteria under which Moore homology agrees with a classical homology of $B\mathcal{G}$, for instance, under hypotheses on $\mathcal{G}$ or on the coefficient group $A$?

*Question* 2. Even when agreement fails, when is there a natural transformation from Moore homology to $H_\bullet^{\mathrm{sing}}(B\mathcal{G}; A)$ or to other invariants of $\mathcal{G}$?

*Question* 3. Since our theory is designed for topological coefficient groups $A$, to what extent does Morita invariance persist beyond the discrete case? In particular, which additional hypotheses on $A$ (or on $\mathcal{G}$) ensure that Morita equivalent groupoids have isomorphic Moore homology with coefficients in $A$, given that this is not a formal consequence of Theorem 2.3.5?

### 2.3.2 The Moore Complex

Let $\mathcal{G}$ be a topological groupoid. A functor $[n] \to \mathcal{G}$ is the same as a composable $n$-tuple $(g_1, \dots, g_n) \in \mathcal{G}_n$. Assume $\mathcal{G}$ is étale, so the structure maps $s, r\colon \mathcal{G} \to \mathcal{G}_0$ are local homeomorphisms. Then the unit map $u\colon \mathcal{G}_0 \to \mathcal{G}$ is a local homeomorphism. Indeed, for $x \in \mathcal{G}_0$ choose an open bisection $U \ni u(x)$. Then $s|_U\colon U \to s(U)$ is a homeomorphism onto an open neighborhood of $x$, and $u = (s|_U)^{-1}$ on $s(U)$. Inversion $(-)^{-1}\colon \mathcal{G} \to \mathcal{G}$ is a homeomorphism. Multiplication $m\colon \mathcal{G}_2 = \mathcal{G} \,{}_s{\times}_r\, \mathcal{G} \to \mathcal{G}$ is a local homeomorphism. For $(g, h) \in \mathcal{G}_2$ choose open bisections $U \ni g$, $V \ni h$ with $s(U) = r(V)$. Then $U \,{}_s{\times}_r\, V$ is open in $\mathcal{G}_2$, $UV$ is open in $\mathcal{G}$, and $m|_{U \,{}_s{\times}_r V}\colon U \,{}_s{\times}_r\, V \to UV$ is a homeomorphism.

Hence all simplicial structure maps of the nerve $\mathcal{G}_\bullet$ are local homeomorphisms. For $n \geq 1$, the endpoint faces $d_0, d_n\colon \mathcal{G}_n \to \mathcal{G}_{n-1}$ are the projections forgetting $g_1$ or $g_n$. Since $\mathcal{G}_n$ is an iterated fibre product along $s$ and $r$, these maps are base changes of $r$ or $s$, hence local homeomorphisms. For $1 \leq i \leq n-1$, the face map $d_i$ composes $(g_i, g_{i+1})$ and fixes the other coordinates. It is the restriction of $\mathrm{id}_{\mathcal{G}_1} \,{}_s{\times}_r \cdots {}_s{\times}_r\, m \,{}_s{\times}_r \cdots {}_s{\times}_r\, \mathrm{id}_{\mathcal{G}_1}$ to the relevant fibre product, hence a local homeomorphism because $m$ is. Each degeneracy $s_j\colon \mathcal{G}_n \to \mathcal{G}_{n+1}$ inserts a unit via $u$ in the $j$-th slot and is a local homeomorphism by the same product and base-change reasoning. Thus the maps $d_i$ and $s_j$ assemble to the nerve $\mathcal{G}_\bullet = (\mathcal{G}_n, (d_i)_{i=0}^n, (s_j)_{j=0}^n)_{n\geq 0}$.

We will use only the face maps to define boundary operators and therefore work with the corresponding semi-simplicial structure.

**Definition 2.3.6.** A semi-simplicial abelian group consists of abelian groups $(A_n)_{n\geq 0}$ together with group homomorphisms $d_i \colon A_n \to A_{n-1}$ for $n \geq 1$ and $0 \leq i \leq n$, such that

$$d_i d_j = d_{j-1} d_i \quad \text{for all } n \geq 2 \text{ and } 0 \leq i < j \leq n.$$

*Remark* 2.3.7. In other words, a semi-simplicial abelian group is a functor $A_\bullet \colon \Delta_{\mathrm{inj}}^{\mathrm{op}} \to \mathbf{Ab}$, where $\Delta_{\mathrm{inj}}$ denotes the wide subcategory of the simplex category $\Delta$ with the same objects $[n]$ and whose morphisms are the injective order-preserving maps. In particular, $\Delta_{\mathrm{inj}}$ is an ordinary small category, and we do not view it as an internal category in **Top**.

**Definition 2.3.8.** A simplicial abelian group consists of abelian groups $(A_n)_{n\geq 0}$, face maps $d_i \colon A_n \to A_{n-1}$ for $n \geq 1$ and $0 \leq i \leq n$, and degeneracy maps $s_j \colon A_n \to A_{n+1}$ for $n \geq 0$ and $0 \leq j \leq n$, all group homomorphisms, such that for all valid indices,

$$\begin{aligned}
d_i d_j &= d_{j-1} d_i \quad \text{for } i < j,\\
s_i s_j &= s_{j+1} s_i \quad \text{for } i \leq j,\\
d_i s_j &= \begin{cases} s_{j-1} d_i, & \text{for } i < j,\\ \mathrm{id}_{A_n}, & \text{for } i = j \text{ or } i = j+1,\\ s_j d_{i-1}, & \text{for } i > j+1. \end{cases}
\end{aligned}$$

*Remark* 2.3.9. In other words, a simplicial abelian group is a functor $A_\bullet \colon \Delta^{\mathrm{op}} \to \mathbf{Ab}$, where $\Delta$ is the simplex category.

**Corollary 2.3.10** ([8, Example 1.4])**.** Let $\mathcal{G}_\bullet = (\mathcal{G}_n, (d_i)_{i=0}^n, (s_j)_{j=0}^n)_{n\geq 0}$ be the nerve of an étale groupoid $\mathcal{G}$, as in Definition 2.3.1. Then

$$C_c(\mathcal{G}_\bullet, A) := (C_c(\mathcal{G}_n, A))_{n\geq 0}$$

is a simplicial abelian group with face maps and degeneracy maps

$$\begin{aligned}
(d_i)_* \colon C_c(\mathcal{G}_n, A) &\to C_c(\mathcal{G}_{n-1}, A) \quad \text{for } n \geq 1,\ 0 \leq i \leq n,\\
(s_j)_* \colon C_c(\mathcal{G}_n, A) &\to C_c(\mathcal{G}_{n+1}, A) \quad \text{for } n \geq 0,\ 0 \leq j \leq n,
\end{aligned}$$

defined by pushforward along the local homeomorphisms $d_i$ and $s_j$.

*Proof.* Since $\mathcal{G}$ is étale, each face map $d_i$ and each degeneracy map $s_j$ of $\mathcal{G}_\bullet$ is a local homeomorphism. Therefore pushforward is defined in every degree. More precisely, for a local homeomorphism $\phi \colon X \to Y$, Definition 2.1.1 gives a homomorphism $\phi_*$, and Proposition 2.1.2 shows that $\phi_*$ maps $C_c(X, A)$ to $C_c(Y, A)$ and preserves compact support. Moreover, push-forward is functorial for composition, so for composable local homeomorphisms $\psi, \phi$ one has $(\psi \circ \phi)_* = \psi_* \circ \phi_*$. Applying this to the simplicial identities among the maps $d_i$ and $s_j$ in $\mathcal{G}_\bullet$, we obtain the corresponding identities for $(d_i)_*$ and $(s_j)_*$. Hence $C_c(\mathcal{G}_\bullet, A)$ is a simplicial abelian group in the sense of Definition 2.3.8. □

**Proposition 2.3.11** (Moore complex)**.** For $n \geq 1$ define

$$\partial_n := \sum_{i=0}^{n} (-1)^i \, (d_i)_* \colon C_c(\mathcal{G}_n, A) \to C_c(\mathcal{G}_{n-1}, A),$$

and set $\partial_0 := 0$. Then $(C_c(\mathcal{G}_n, A), \partial_n)_{n \geq 0}$ is a chain complex. We call it the Moore complex of $\mathcal{G}$ with coefficients in $A$.

*Proof.* For each $n \geq 1$ and each $0 \leq i \leq n$, the map $d_i \colon \mathcal{G}_n \to \mathcal{G}_{n-1}$ is a local homeomorphism. Thus Lemma 2.1.3 ensures that the fibrewise sum defining $(d_i)_*$ is finite on compact support, and Proposition 2.1.2 shows that $(d_i)_*$ maps $C_c(\mathcal{G}_n, A)$ into $C_c(\mathcal{G}_{n-1}, A)$. Hence $\partial_n$ is a well-defined homomorphism.

It remains to show $\partial_{n-1}\partial_n = 0$ for $n \geq 1$. By functoriality of pushforward,

$$(d_j)_*(d_i)_* = (d_j d_i)_* \quad \text{for all composable face maps } d_i, d_j.$$

Using the simplicial identity $d_j d_i = d_{i-1} d_j$ for $j < i$, we compute

$$\begin{aligned}
\partial_{n-1}\partial_n &= \sum_{i=0}^{n}\sum_{j=0}^{n-1} (-1)^{i+j} (d_j d_i)_* \\
&= \sum_{0 \leq j < i \leq n} (-1)^{i+j} (d_j d_i)_* + \sum_{0 \leq i \leq j \leq n-1} (-1)^{i+j} (d_j d_i)_* \\
&= \sum_{0 \leq j < i \leq n} (-1)^{i+j} (d_j d_i)_* + \sum_{0 \leq i < j \leq n} (-1)^{i+j-1} (d_{j-1} d_i)_* \\
&= \sum_{0 \leq j < i \leq n} (-1)^{i+j} (d_j d_i)_* - \sum_{0 \leq i < j \leq n} (-1)^{i+j} (d_{j-1} d_i)_* \\
&= \sum_{0 \leq j < i \leq n} (-1)^{i+j} (d_j d_i)_* - \sum_{0 \leq j < i \leq n} (-1)^{i+j} (d_{i-1} d_j)_* \\
&= \sum_{0 \leq j < i \leq n} (-1)^{i+j} (d_j d_i)_* - \sum_{0 \leq j < i \leq n} (-1)^{i+j} (d_j d_i)_* \\
&= 0.
\end{aligned}$$

Thus $\partial_{n-1}\partial_n = 0$ for all $n \geq 1$, and the claim follows. □

**Corollary 2.3.12.** For all $n \geq 1$ one has $\partial_{n-1} \circ \partial_n = 0$ as a consequence of the simplicial identities. In particular, $(C_c(\mathcal{G}_\bullet, A), \partial_\bullet)$ is a chain complex.

**Definition 2.3.13** (Homology with constant coefficients [5, Section 3.1], [13, Definition 3.1])**.** Let $(C_c(\mathcal{G}_\bullet, A), \partial_\bullet)$ be the Moore complex of the nerve $\mathcal{G}_\bullet$ from Proposition 2.3.11, that is $\partial_n\colon C_c(\mathcal{G}_n, A) \to C_c(\mathcal{G}_{n-1}, A)$. For $n \geq 0$ define

$$H_n(\mathcal{G}; A) := H_n(C_c(\mathcal{G}_\bullet, A), \partial_\bullet) = \frac{\ker(\partial_n)}{\operatorname{im}(\partial_{n+1})}.$$

*Remark* 2.3.14. We call $H_n(\mathcal{G}; A)$ the $n$-th homology group of $\mathcal{G}$ with constant coefficients in $A$. When $A = \mathbb{Z}$, we simply write $H_n(\mathcal{G}) := H_n(\mathcal{G}; \mathbb{Z})$. This functoriality is a basic input for the exact sequences used later, in particular for Mayer–Vietoris arguments and for computations by decomposing the unit space.

**Theorem 2.3.15.** Let $\varphi : \mathcal{H} \to \mathcal{G}$ be an étale functor between étale groupoids, and let $A$ be a topological abelian group.

1. For each $n \geq 0$, the induced map on nerves $\varphi_n : \mathcal{H}_n \to \mathcal{G}_n$ is a local homeomorphism. Here $\varphi_0$ is the object map, and for $n \geq 1$, $\varphi_n(h_1, \dots, h_n) := (\varphi_1(h_1), \dots, \varphi_1(h_n))$.
2. For each $n \geq 0$, pushforward along $\varphi_n$ defines a group homomorphism

$$(\varphi_n)_* : C_c(\mathcal{H}_n, A) \to C_c(\mathcal{G}_n, A), \quad (\varphi_n)_*(f)(\mathbf{g}) := \sum_{\mathbf{h} \in \varphi_n^{-1}(\mathbf{g})} f(\mathbf{h}),$$

   and the maps $(\varphi_n)_*$ form a chain map between Moore complexes

$$(\varphi_\bullet)_* : C_c(\mathcal{H}_\bullet, A) \to C_c(\mathcal{G}_\bullet, A),\ \partial_n^{\mathcal{G}} \circ (\varphi_n)_* = (\varphi_{n-1})_* \circ \partial_n^{\mathcal{H}} \quad \text{for all } n \geq 1.$$

   Hence $\varphi$ induces homomorphisms on homology $H_n(\varphi) : H_n(\mathcal{H}; A) \to H_n(\mathcal{G}; A)$ for all $n \geq 0$, and $\varphi \mapsto H_\bullet(\varphi)$ is functorial. If $\psi : \mathcal{G} \to \mathcal{K}$ is another étale functor, then $H_\bullet(\psi \circ \varphi) = H_\bullet(\psi) \circ H_\bullet(\varphi)$ and $H_\bullet(\mathrm{id}_{\mathcal{G}}) = \mathrm{id}_{H_\bullet(\mathcal{G}; A)}$.
3. Let $(\mathcal{G}_k)_{k \in \mathbb{N}}$ be an increasing sequence of open subgroupoids of $\mathcal{G}$ with $\mathcal{G} = \bigcup_{k \in \mathbb{N}} \mathcal{G}_k$. Then for each $n \geq 0$, the canonical map is an isomorphism:

$$\varinjlim_k H_n(\mathcal{G}_k; A) \to H_n(\mathcal{G}; A).$$

*Proof.*

1. Since $\varphi$ is an étale functor, both $\varphi_0 : \mathcal{H}_0 \to \mathcal{G}_0$ and $\varphi_1 : \mathcal{H}_1 \to \mathcal{G}_1$ are local homeomorphisms. For $n \geq 2$, write the $n$-simplices as iterated fibre products

$$\mathcal{H}_n = \mathcal{H}_1 \; {}_{s_{\mathcal{H}}}\!\times_{r_{\mathcal{H}}} \cdots \; {}_{s_{\mathcal{H}}}\!\times_{r_{\mathcal{H}}} \mathcal{H}_1, \quad \mathcal{G}_n = \mathcal{G}_1 \; {}_{s_{\mathcal{G}}}\!\times_{r_{\mathcal{G}}} \cdots \; {}_{s_{\mathcal{G}}}\!\times_{r_{\mathcal{G}}} \mathcal{G}_1.$$

   The identities $s_{\mathcal{G}} \circ \varphi_1 = \varphi_0 \circ s_{\mathcal{H}}$ and $r_{\mathcal{G}} \circ \varphi_1 = \varphi_0 \circ r_{\mathcal{H}}$ make the defining pullback squares commute, so $\varphi_n$ is the induced map between these iterated fibre products. Local homeomorphisms are stable under products and base change, hence $\varphi_n$ is a local homeomorphism for all $n \geq 0$.

2. Since $\varphi_n$ is a local homeomorphism, $(\varphi_n)_*$ is well defined by Definition 2.1.1. For $f \in C_c(\mathcal{H}_n, A)$ and $\mathbf{g} \in \mathcal{G}_n$, the sum $\sum_{\mathbf{h}\in\varphi_n^{-1}(\mathbf{g})} f(\mathbf{h})$ is finite because $\varphi_n^{-1}(\mathbf{g})$ is discrete and $\varphi_n^{-1}(\mathbf{g}) \cap \operatorname{supp}(f)$ is finite by Lemma 2.1.3. We show that $(\varphi_\bullet)_*$ is a chain map. For each $n \geq 1$ and $0 \leq i \leq n$, simpliciality of $\varphi_\bullet$ gives

$$d_i^{\mathcal{G}} \circ \varphi_n = \varphi_{n-1} \circ d_i^{\mathcal{H}}.$$

Pushforward and compatibility with composition from Proposition 2.1.2 gives

$$(d_i^{\mathcal{G}})_* \circ (\varphi_n)_* = (\varphi_{n-1})_* \circ (d_i^{\mathcal{H}})_*.$$

Summing with signs yields

$$\partial_n^{\mathcal{G}} \circ (\varphi_n)_* = \sum_{i=0}^{n} (-1)^i (d_i^{\mathcal{G}})_* \circ (\varphi_n)_* = \sum_{i=0}^{n} (-1)^i (\varphi_{n-1})_* \circ (d_i^{\mathcal{H}})_* = (\varphi_{n-1})_* \circ \partial_n^{\mathcal{H}}.$$

Thus $(\varphi_\bullet)_*$ induces $H_n(\varphi)$ on homology. If $\psi : \mathcal{G} \to \mathcal{K}$ is another étale functor, then $(\psi \circ \varphi)_n = \psi_n \circ \varphi_n$, hence $((\psi\circ\varphi)_n)_* = (\psi_n)_* \circ (\varphi_n)_*$ by Proposition 2.1.2, and functoriality on homology follows. The identity functor induces the identity chain map

$$(\mathrm{id}_{\mathcal{G}_n})_*(f)(\mathbf{g}) = \sum_{(\mathrm{id}_{\mathcal{G}_n})(\mathbf{h})=\mathbf{g}} f(\mathbf{h}) = f(\mathbf{g}),$$

so $(\mathrm{id}_{\mathcal{G}_n})_* = \mathrm{id}_{C_c(\mathcal{G}_n,A)}$. Hence $(\mathrm{id}_{\mathcal{G}_\bullet})_* = \mathrm{id}_{C_c(\mathcal{G}_\bullet,A)}$ and $H_n(\mathrm{id}_{\mathcal{G}}) = \mathrm{id}_{H_n(\mathcal{G};A)}$ for all $n \geq 0$.

3. For each $n \geq 0$, the subsets $(\mathcal{G}_k)_n \subseteq \mathcal{G}_n$ form an increasing sequence of open subsets with $\mathcal{G}_n = \bigcup_k (\mathcal{G}_k)_n$. Indeed, if $\mathbf{g} = (g_1, \dots, g_n) \in \mathcal{G}_n$, then each $g_i$ lies in some $(\mathcal{G}_{k_i})_1$, and for $k \geq \max_i k_i$ one has $g_i \in (\mathcal{G}_k)_1$ for all $i$. Since $\mathcal{G}_k$ is a subgroupoid, composability is inherited, hence $\mathbf{g} \in (\mathcal{G}_k)_n$. Consider the directed system $\{C_c((\mathcal{G}_k)_n, A)\}_k$ with transition maps given by extension by $0_A$ along $(\mathcal{G}_k)_n \subseteq (\mathcal{G}_\ell)_n$ for $k \leq \ell$. The canonical map

$$\varinjlim_k C_c((\mathcal{G}_k)_n, A) \to C_c(\mathcal{G}_n, A) \tag{2.3.3}$$

is an isomorphism for each $n$.

- **Surjectivity.** Let $f \in C_c(\mathcal{G}_n, A)$ and set $K = \operatorname{supp}(f)$, which is compact. The open cover $\{(\mathcal{G}_k)_n\}_k$ of $\mathcal{G}_n$ admits a finite subcover of $K$, hence $K \subseteq (\mathcal{G}_k)_n$ for some $k$ because $(\mathcal{G}_k)_n$ is increasing. Then $f \in C_c((\mathcal{G}_k)_n, A)$, and its extension by $0_A$ is $f$.
- **Injectivity.** If $f \in C_c((\mathcal{G}_k)_n, A)$ and $g \in C_c((\mathcal{G}_\ell)_n, A)$ have the same extension to $C_c(\mathcal{G}_n, A)$, let $h$ denote this common extension and set $K = \operatorname{supp}(h)$. As above there exists $m \geq k, \ell$ with $K \subseteq (\mathcal{G}_m)_n$. Restricting $h$ to $(\mathcal{G}_m)_n$ shows that the images of $f$ and $g$ agree in $C_c((\mathcal{G}_m)_n, A)$, hence they define the same element in the colimit.

This proves (2.3.3). Moreover, since each face map $d_i : \mathcal{G}_n \to \mathcal{G}_{n-1}$ restricts to $d_i : (\mathcal{G}_k)_n \to (\mathcal{G}_k)_{n-1}$, the extension-by-zero maps commute with the Moore boundaries $\partial_n$. Therefore $C_c(\mathcal{G}_\bullet, A)$ is the filtered colimit of the subcomplexes $C_c((\mathcal{G}_k)_\bullet, A)$. Filtered colimits are exact in **Ab** [12, p. 67f], so homology commutes with this colimit:

$$H_n(\mathcal{G}; A) \cong H_n\Big(\varinjlim_k C_c((\mathcal{G}_k)_\bullet, A)\Big) \cong \varinjlim_k H_n(\mathcal{G}_k; A) \quad \text{for all } n \geq 0.$$

□

The classifying space $B\mathcal{G}$ captures the homotopy-theoretic content of a topological groupoid. In many standard settings it classifies principal $\mathcal{G}$-bundles over paracompact spaces and records Morita invariance at the level of weak homotopy type. By contrast, the Moore complex $C_c(\mathcal{G}_\bullet, A)$ is built from compactly supported $A$-valued functions on the nerve levels $\mathcal{G}_n$, with boundary given by alternating sums of pushforwards along face maps. In the ample setting with discrete coefficients, these functions are typically locally constant with finite image, so the Moore complex is designed for cut-and-paste arguments and explicit computations using compact open decompositions, reductions, and exact sequences. It is therefore not intended to recover singular homology of $B\mathcal{G}$ in general, and its invariance and exactness properties must be established directly at the chain level.

The following example isolates this distinction in the most elementary case. Even for a principal groupoid, the Moore complex can be chain-contractible in positive degrees while $B\mathcal{G}$ retains nontrivial singular homology.

**Example 2.3.16** (Moore chains vs. singular chains)**.** The Moore complex $C_c(\mathcal{G}_\bullet, A)$ is in general not the singular chain complex of $B\mathcal{G}$. The simplest instance already occurs for the unit groupoid on a space.

Let $X$ be a locally compact Hausdorff space and let $\mathcal{G}$ be the unit groupoid on $X$, so

$$\mathcal{G}_0 := X, \quad \mathcal{G}_1 := X, \quad r = s = \mathrm{id}_X,$$

and every arrow is a unit. For $n \geq 1$, the space of composable $n$-tuples is the diagonal

$$\mathcal{G}_n = \{(x, \dots, x) \in X^n \mid x \in X\} \subseteq X^n.$$

Let

$$\iota_n : X \to \mathcal{G}_n, \quad \iota_n(x) := (x, \dots, x), \quad \pi_n : \mathcal{G}_n \to X, \quad \pi_n(x, \dots, x) := x.$$

Then $\iota_n$ and $\pi_n$ are inverse homeomorphisms, so $\mathcal{G}_n \cong X$ via $\pi_n$. For each $n \geq 1$ and $0 \leq i \leq n$, the face map $d_i : \mathcal{G}_n \to \mathcal{G}_{n-1}$ deletes one coordinate, hence satisfies

$$d_i = \iota_{n-1} \circ \pi_n.$$

In particular, under the identifications $\pi_n$ and $\pi_{n-1}$, each $d_i$ corresponds to $\mathrm{id}_X$. Since pushforward along a homeomorphism is an isomorphism, it follows that

$$(d_i)_* = \mathrm{id}_{C_c(\mathcal{G}_n, A)} \quad \text{for all } n \geq 1 \text{ and } 0 \leq i \leq n,$$

after identifying $C_c(\mathcal{G}_n, A) \cong C_c(X, A)$.

Therefore, for $n \geq 1$ the Moore boundary is

$$\partial_n = \sum_{i=0}^{n} (-1)^i (d_i)_* = \sum_{i=0}^{n} (-1)^i \mathrm{id}_{C_c(\mathcal{G}_n, A)} = \begin{cases} 0, & n \text{ odd}, \\ \mathrm{id}_{C_c(\mathcal{G}_n, A)}, & n \text{ even}, \end{cases}$$

and we set $\partial_0 := 0$. Hence $C_c(\mathcal{G}_\bullet, A)$ is isomorphic to the chain complex

$$\cdots \xrightarrow{\mathrm{id}} C_c(X, A) \xrightarrow{0} C_c(X, A) \xrightarrow{\mathrm{id}} C_c(X, A) \xrightarrow{0} C_c(X, A) \xrightarrow{0} 0,$$

with homology groups

$$H_0(\mathcal{G}; A) \cong C_c(X, A), \quad H_n(\mathcal{G}; A) = 0 \text{ for all } n \geq 1.$$

On the other hand, the simplicial space $\mathcal{G}_\bullet$ is levelwise homeomorphic to $X$ via $\pi_n$, and its geometric realization satisfies $B\mathcal{G} \cong X$. Consequently,

$$H_n^{\mathrm{sing}}(B\mathcal{G}; A) \cong H_n^{\mathrm{sing}}(X; A),$$

which is typically nonzero in higher degrees. For example, if $X = S^1$ and $A = \mathbb{Z}$, then

$$H_1^{\mathrm{sing}}(B\mathcal{G}; \mathbb{Z}) \cong H_1^{\mathrm{sing}}(S^1; \mathbb{Z}) \cong \mathbb{Z}, \quad \text{while} \quad H_1(\mathcal{G}; \mathbb{Z}) = 0.$$

Thus $C_c(\mathcal{G}_\bullet, A)$ cannot agree with the singular chain complex of $B\mathcal{G}$ in general.

## 2.4 Cohomology Groups

The Moore homology groups $H_n(\mathcal{G}; A)$ are computable in the ample setting. They are defined from compactly supported chains on the nerve and fit into long exact sequences for reductions and clopen saturated decompositions, as in [14] and Theorem 3.3.10. Cohomology is dual: many natural groupoid invariants are cocycle-based, detect obstructions, and classify extensions.

In this chapter we define cohomology on the same nerve $\mathcal{G}_\bullet$ by passing to cochains, and we write $H^\bullet(\mathcal{G}; A)$ for its cohomology. Evaluation induces a canonical pairing $H_n(\mathcal{G}; \mathbb{Z}) \times H^n(\mathcal{G}; A) \to A$ compatible with boundaries, hence with the long exact sequences, see [15]. Under the usual hypotheses, in particular that $C_c(\mathcal{G}_n, \mathbb{Z})$ is free for all $n$, the classical universal coefficient theorem identifies $H^n(\mathcal{G}; A)$ in terms of $H_n(\mathcal{G}; \mathbb{Z})$ via $\mathrm{Hom}_{\mathbb{Z}}$ and $\mathrm{Ext}^1_{\mathbb{Z}}$, see Section 3.2. We emphasize that this coefficient description is specific to the compact-support model and, as in homology, the behaviour changes substantially once one allows genuinely topological coefficients.

**Definition 2.4.1** (Cochains and coboundary)**.** Let $\mathcal{G}$ be an étale groupoid and let $A$ be an abelian group. For $n \geq 0$ set $C^n(\mathcal{G}; A) := C_c(\mathcal{G}_n, A)$. For $\xi \in C^n(\mathcal{G}; A)$ define

$$\delta^n \colon C^n(\mathcal{G}; A) \to C^{n+1}(\mathcal{G}; A), \quad \delta^n(\xi) := \sum_{i=0}^{n+1} (-1)^i \, \xi \circ d_i,$$

where $d_i : \mathcal{G}_{n+1} \to \mathcal{G}_n$ are the face maps of the nerve $\mathcal{G}_\bullet$. We set $\delta^{-1} : 0 \to C^0(\mathcal{G}; A)$.

**Proposition 2.4.2.** Assume that for each $n \geq 0$ and each $0 \leq i \leq n+1$, the face map $d_i : \mathcal{G}_{n+1} \to \mathcal{G}_n$ is proper. Then $\delta^n$ is well defined on $C_c(\mathcal{G}_n, A)$, and $\delta^{n+1} \circ \delta^n = 0$. for all $n \geq 0$. In particular, $(C^\bullet(\mathcal{G}; A), \delta^\bullet)$ is a cochain complex.

*Proof.* Let $n \geq 0$ and $\xi \in C^n(\mathcal{G}; A)$. Each $d_i$ is continuous, hence $\xi \circ d_i$ is continuous. Moreover,

$$\operatorname{supp}(\xi \circ d_i) \subseteq d_i^{-1}(\operatorname{supp}(\xi)).$$

Since $\operatorname{supp}(\xi)$ is compact and $d_i$ is proper, the set $d_i^{-1}(\operatorname{supp}(\xi))$ is compact. Thus $\xi \circ d_i \in C_c(\mathcal{G}_{n+1}, A)$, and $\delta^n(\xi) \in C^{n+1}(\mathcal{G}; A)$ is well defined.

To prove $\delta^{n+1}\delta^n = 0$, fix $g \in \mathcal{G}_{n+2}$. We compute

$$\delta^{n+1}\delta^n(\xi)(g) = \sum_{j=0}^{n+2} (-1)^j \, \delta^n(\xi)(d_j(g)) = \sum_{j=0}^{n+2} \sum_{i=0}^{n+1} (-1)^{i+j} \, \xi(d_i d_j(g)).$$

Using the simplicial identity $d_i d_j = d_{j-1} d_i$ for $i < j$, the summands cancel in pairs:

For $0 \leq i < j \leq n+2$, $(-1)^{i+j} \, \xi(d_i d_j(g)) + (-1)^{i+j-1} \, \xi(d_{j-1} d_i(g)) = 0$. Every term occurs in exactly one such pair, hence $\delta^{n+1}\delta^n(\xi)(g) = 0$. Since $g$ was arbitrary, $\delta^{n+1} \circ \delta^n = 0$. □

*Remark* 2.4.3*.* The properness hypothesis in Proposition 2.4.2 is a genuine restriction in the compact-support model. For étale groupoids the face maps are local homeomorphisms, but need not be proper, so pullback does not preserve compact support in general. There are standard alternatives, for instance defining cochains without compact support or working with sheaf or module coefficients, as in the module-theoretic and sheaf-theoretic approaches developed in recent work [4, 5]. We keep the compact-support convention here because it interacts well with the cut-and-paste technology and the exact sequences used throughout this thesis, and we isolate precisely where additional hypotheses are needed.

**Definition 2.4.4** (Cohomology with constant coefficients)**.** Assume the hypotheses of Proposition 2.4.2. For $n \geq 0$ define the $n$-th cohomology group of $\mathcal{G}$ with coefficients in $A$ by

$$H^n(\mathcal{G}; A) := H^n(C^\bullet(\mathcal{G}; A), \delta^\bullet) = \frac{\ker(\delta^n)}{\operatorname{im}(\delta^{n-1})}.$$

*Remark* 2.4.5*.* Elements of $\ker(\delta^n)$ are called $n$-cocycles, and elements of $\operatorname{im}(\delta^{n-1})$ are called $n$-coboundaries. For $\xi \in \ker(\delta^n)$ its class in $H^n(\mathcal{G}; A)$ is denoted by $[\xi]$. When $A = \mathbb{Z}$ we write $H^n(\mathcal{G}) := H^n(\mathcal{G}; \mathbb{Z})$.

Compact support on cochains is not needed for the basic evaluation pairing against compactly supported chains. If $X$ is a space, $f \in C_c(X, \mathbb{Z})$, and $\xi \in C(X, A)$, define

$$(f \cdot \xi)(x) := f(x)\, \xi(x),$$

where $f(x)$ acts on $A$ by repeated addition. Since $f$ is continuous with values in the discrete group $\mathbb{Z}$, it is locally constant, hence $f \cdot \xi$ is continuous. Moreover $\operatorname{supp}(f \cdot \xi) \subseteq \operatorname{supp}(f)$, so $f \cdot \xi \in C_c(X, A)$ even when $\xi$ has no compact support.

The following lemma allows us to combine compactly supported integer-valued chains $f \in C_c(X, \mathbb{Z})$ with $A$-valued cochains $\xi \in C(X, A)$ to obtain compactly supported $A$-valued functions $f \cdot \xi \in C_c(X, A)$. Since pushforwards along local homeomorphisms are defined on $C_c(-, A)$ by finite fibrewise sums, see Proposition 2.1.2, this makes it possible to transport $A$-valued data along the same maps that appear in the Moore complex. This chain-level compatibility is used later to verify that the pairings and the universal coefficient constructions respect boundaries and compact supports.

**Lemma 2.4.6.** Let $X$ be an object of **LCH** and let $A$ be a topological abelian group. For $f \in C_c(X, \mathbb{Z})$ and $\xi \in C(X, A)$ define a map $(f \cdot \xi)(x) := f(x) \cdot \xi(x)$ for all $x \in X$, where $n \cdot a$ denotes the $n$-fold sum of $a \in A$ and its inverse for $n < 0$. Then $f \cdot \xi \in C_c(X, A)$.

*Proof.*

- $f \cdot \xi$ **is continuous.** Define $F\colon X \to \mathbb{Z} \times A$ by $F(x) := (f(x), \xi(x))$. The map $F$ is continuous because $f$ and $\xi$ are continuous and $\mathbb{Z}$ carries the discrete topology. Define $m\colon \mathbb{Z} \times A \to A$ by $m(n, a) := n \cdot a$. For each fixed $n \in \mathbb{Z}$ the restriction $m|_{\{n\} \times A}\colon A \to A$ is the map $a \mapsto n \cdot a$, which is continuous as a finite composition of addition and inversion in $A$. Since $\mathbb{Z}$ is discrete, the subsets $\{n\} \times A$ are open, hence $m$ is continuous on $\mathbb{Z} \times A$. We have $f \cdot \xi = m \circ F$, so $f \cdot \xi$ is continuous.
- $f \cdot \xi$ **has compact support.** If $x \notin \operatorname{supp}(f)$, there exists an open neighbourhood $U$ of $x$ such that $f|_U = 0$. For every $y \in U$ we then have $(f \cdot \xi)(y) = f(y) \cdot \xi(y) = 0 \cdot \xi(y) = 0_A$, so $x \notin \operatorname{supp}(f \cdot \xi)$. Thus $\operatorname{supp}(f \cdot \xi) \subseteq \operatorname{supp}(f)$.

Since $\operatorname{supp}(f)$ is compact, so is $\operatorname{supp}(f \cdot \xi)$. Therefore $f \cdot \xi \in C_c(X, A)$. □

Developing the UCT, we frequently pushforward compactly supported $\mathbb{Z}$-valued functions along local homeomorphisms and then multiply the result with $A$-valued cochains. Lemma 2.4.7 shows this agrees with first multiplying on the source and then pushing forward.

**Lemma 2.4.7.** Let $\pi\colon X \to Y$ be a local homeomorphism in **LCH**, and let $A$ be a topological abelian group. For $f \in C_c(X, \mathbb{Z})$ and $\xi \in C(Y, A)$ one has

$$\pi_*(f) \cdot \xi = \pi_*(f \cdot (\xi \circ \pi))$$

as elements of $C_c(Y, A)$, where $\cdot$ denotes the pointwise product from Lemma 2.4.6 and $\pi_*$ is the pushforward from Proposition 2.1.2.

*Proof.*

- **Well-definedness.** By Lemma 2.4.6, the function $f \cdot (\xi \circ \pi) \colon X \to A$ belongs to $C_c(X, A)$, so both sides are well defined elements of $C_c(Y, A)$.
- **Product and pushforward.** Fix $y \in Y$ and write $E_y := \{x \in X \mid \pi(x) = y\}$. By definition of $\pi_*$ we have

$$\pi_*(f)(y) = \sum_{x \in E_y} f(x), \quad \pi_*(f \cdot (\xi \circ \pi))(y) = \sum_{x \in E_y} f(x) \cdot \xi(\pi(x)).$$

  Only finitely many terms are nonzero. The fibre $E_y$ is discrete because $\pi$ is a local homeomorphism, and $E_y \cap \operatorname{supp}(f)$ is finite since $\operatorname{supp}(f)$ is compact. Using $\pi(x) = y$ for all $x \in E_y$ we obtain

$$\pi_*(f \cdot (\xi \circ \pi))(y) = \sum_{x \in E_y} f(x) \cdot \xi(y).$$

  For fixed $a \in A$, the map $\mu_a \colon \mathbb{Z} \to A$, $n \mapsto n \cdot a$, is a group homomorphism. Applying this with $a = \xi(y)$ gives

$$\Big( \sum_{x \in E_y} f(x) \Big) \cdot \xi(y) = \sum_{x \in E_y} f(x) \cdot \xi(y).$$

  Hence

$$(\pi_*(f) \cdot \xi)(y) = \pi_*(f)(y) \cdot \xi(y) = \sum_{x \in E_y} f(x) \cdot \xi(y) = \pi_*(f \cdot (\xi \circ \pi))(y).$$

  Since $y \in Y$ was arbitrary, the two functions agree on $Y$.

$\square$

*Remark* 2.4.8. We write $\operatorname{Bis}(\mathcal{G})$ for the set of compact open bisections of $\mathcal{G}$. For $U \subseteq \mathcal{G}$ we denote by $\chi_U \in C_c(\mathcal{G}, \mathbb{Z})$ its characteristic function.

**Lemma 2.4.9.** Let $\mathcal{G}$ be an ample étale groupoid. Then every $f \in C_c(\mathcal{G}, \mathbb{Z})$ is locally constant and has compact open support. Moreover, $C_c(\mathcal{G}, \mathbb{Z})$ is generated, as an abelian group, by the characteristic functions $\chi_U$ with $U \in \operatorname{Bis}(\mathcal{G})$.

*Proof.* Let $f \in C_c(\mathcal{G}, \mathbb{Z})$.

- **Local constancy.** Fix $g \in \mathcal{G}$ and set $n := f(g)$. Since $\{n\} \subseteq \mathbb{Z}$ is open in the discrete topology, the set $f^{-1}(\{n\})$ is an open neighbourhood of $g$ on which $f$ is constant.
- **Compact open support.** For each $n \in \mathbb{Z}$, the set $f^{-1}(\{n\})$ is clopen. Put $U := \bigcup_{n \neq 0} f^{-1}(\{n\}) = \{g \in \mathcal{G} \mid f(g) \neq 0\}$. Then $U$ is open, and $\operatorname{supp}(f) = \overline{U}$. Since $U$ is closed, $\operatorname{supp}(f) = U$. In particular, $\operatorname{supp}(f)$ is compact and open.
- **Generation by characteristic functions of bisections.** Since $K := \operatorname{supp}(f)$ is compact and $f$ is locally constant, the image $f(K) \subseteq \mathbb{Z}$ is finite. Let $n_1, \dots, n_m$ be the distinct nonzero values of $f$ on $K$, and set $S_j := f^{-1}(\{n_j\})$ for $1 \le j \le m$. Then each $S_j$ is clopen, contained in $K$, the sets $S_1, \dots, S_m$ are pairwise disjoint, and

$$K = \bigcup_{j=1}^{m} S_j, \quad f = \sum_{j=1}^{m} n_j \chi_{S_j}.$$

Fix $j$. The set $S_j$ is compact open in $\mathcal{G}$. Since $\mathcal{G}$ is ample, $\mathrm{Bis}(\mathcal{G})$ is a basis of compact open bisections, hence for every $g \in S_j$ there exists $U_g \in \mathrm{Bis}(\mathcal{G})$ with $g \in U_g \subseteq S_j$. By compactness of $S_j$ we may choose $U_{j,1}, \dots, U_{j,\ell_j} \in \mathrm{Bis}(\mathcal{G})$ with $S_j \subseteq \bigcup_{k=1}^{\ell_j} U_{j,k}$ and $U_{j,k} \subseteq S_j$ for all $k$. We refine this finite cover to a finite disjoint family of compact open bisections. For $1 \le k \le \ell_j$ define

$$W_{j,k} := U_{j,k} \setminus \bigcup_{r<k} U_{j,r}.$$

Each $W_{j,k}$ is compact open, being the difference of compact open sets. Moreover $W_{j,k} \subseteq U_{j,k}$, hence $W_{j,k}$ is a bisection. The sets $W_{j,1}, \dots, W_{j,\ell_j}$ are pairwise disjoint and satisfy $S_j = \bigsqcup_{k=1}^{\ell_j} W_{j,k}$, so

$$\chi_{S_j} = \sum_{k=1}^{\ell_j} \chi_{W_{j,k}}.$$

Substituting into the previous decomposition yields

$$f = \sum_{j=1}^{m} n_j \chi_{S_j} = \sum_{j=1}^{m} \sum_{k=1}^{\ell_j} n_j \chi_{W_{j,k}},$$

a finite $\mathbb{Z}$-linear combination of characteristic functions of compact open bisections.

□

This lemma is the basic ample input behind the chain-level freeness used later in the universal coefficient theorem for discrete coefficients, and it is also the point where the restriction to compact open data becomes visible.

**Definition 2.4.10.** Let $\mathcal{G}$ be an étale groupoid and, for $u \in \mathcal{G}_0$, set $\mathcal{G}_u := \{g \in \mathcal{G} \mid s(g) = u\}$. For $f_1, f_2 \in C_c(\mathcal{G}, \mathbb{Z})$ define the convolution product $f_1 * f_2 \in C_c(\mathcal{G}, \mathbb{Z})$ by

$$(f_1 * f_2)(g) := \sum_{h \in \mathcal{G}_{r(g)}} f_1(h^{-1}) f_2(hg) \quad \text{for all } g \in \mathcal{G}.$$

The sum is finite because $\mathcal{G}$ is étale and $f_1, f_2$ have compact support.

**Lemma 2.4.11.** Let $\mathcal{G}$ be an ample étale groupoid. With pointwise addition and convolution $*$ as in Definition 2.4.10, the triple $(C_c(\mathcal{G}, \mathbb{Z}), +, *)$ is an associative ring with local units: for every finite subset $\{f_1, \dots, f_n\} \subset C_c(\mathcal{G}, \mathbb{Z})$ there exists an idempotent $e \in C_c(\mathcal{G}, \mathbb{Z})$ such that $e * f_i = f_i = f_i * e$ for all $i$.

*Proof.* Let $f_1, f_2 \in C_c(\mathcal{G}, \mathbb{Z})$ and define $(f_1 * f_2)(g) := \sum_{h \in \mathcal{G}_{r(g)}} f_1(h^{-1}) f_2(hg)$.

- **Well-definedness and compact support.** Fix $g \in \mathcal{G}$. A term is nonzero only if $h^{-1} \in \operatorname{supp}(f_1)$ and $hg \in \operatorname{supp}(f_2)$. Thus $h$ lies in
$$E_g := \{h \in \mathcal{G}_{r(g)} \mid h^{-1} \in \operatorname{supp}(f_1),\ hg \in \operatorname{supp}(f_2)\}.$$
Since $\mathcal{G}$ is étale, the fibre $\mathcal{G}_{r(g)} = s^{-1}(r(g))$ is discrete. The set
$$H_1 := \{h \in \mathcal{G}_{r(g)} \mid h^{-1} \in \operatorname{supp}(f_1)\} = \mathcal{G}_{r(g)} \cap (\operatorname{supp}(f_1))^{-1}$$
is compact as a closed subset of the compact set $(\operatorname{supp}(f_1))^{-1}$, hence finite in the discrete space $\mathcal{G}_{r(g)}$. Moreover, right multiplication by $g$ restricts to a homeomorphism
$$R_g \colon \mathcal{G}_{r(g)} \to r^{-1}(s(g)), \quad h \mapsto hg,$$
hence
$$H_2 := \{h \in \mathcal{G}_{r(g)} \mid hg \in \operatorname{supp}(f_2)\} = R_g^{-1}(\operatorname{supp}(f_2) \cap r^{-1}(s(g)))$$
is compact in $\mathcal{G}_{r(g)}$, hence finite. Therefore $E_g = H_1 \cap H_2$ is finite, and $(f_1 * f_2)(g)$ is a finite sum. If $(f_1 * f_2)(g) \neq 0$, then there exists $h \in \mathcal{G}_{r(g)}$ with $h^{-1} \in \operatorname{supp}(f_1)$ and $hg \in \operatorname{supp}(f_2)$. Put $k_1 := h^{-1}$ and $k_2 := hg$. Then $(k_1, k_2) \in \mathcal{G}_2$ and $k_1 k_2 = g$, so
$$\operatorname{supp}(f_1 * f_2) \subseteq \operatorname{supp}(f_1)\operatorname{supp}(f_2),$$
$$\operatorname{supp}(f_1)\operatorname{supp}(f_2) := \{k_1 k_2 \mid (k_1, k_2) \in \mathcal{G}_2,\ k_1 \in \operatorname{supp}(f_1),\ k_2 \in \operatorname{supp}(f_2)\}.$$
The set $\mathcal{G}_2 \cap (\operatorname{supp}(f_1) \times \operatorname{supp}(f_2))$ is compact, and $m \colon \mathcal{G}_2 \to \mathcal{G}$ is continuous, hence $\operatorname{supp}(f_1)\operatorname{supp}(f_2)$ is compact. Thus $\operatorname{supp}(f_1 * f_2)$ is compact and $f_1 * f_2 \in C_c(\mathcal{G}, \mathbb{Z})$.
- **Local constancy.** The group $\mathbb{Z}$ is discrete, hence every element of $C_c(\mathcal{G}, \mathbb{Z})$ is locally constant. Therefore $f_1 * f_2$, being an element of $C_c(\mathcal{G}, \mathbb{Z})$, is locally constant as well.
- **Bilinearity.** Follows from distributivity of $\mathbb{Z}$-addition and linearity of finite sums.
- **Associativity.** It is convenient to use the equivalent form
$$(f_1 * f_2)(g) = \sum_{\substack{(h_1, h_2) \in \mathcal{G}_2 \\ h_1 h_2 = g}} f_1(h_1) f_2(h_2), \quad g \in \mathcal{G},$$
obtained via the bijection $\mathcal{G}_{r(g)} \to \{(h_1, h_2) \in \mathcal{G}_2 \mid h_1 h_2 = g\}$, $h \mapsto (h^{-1}, hg)$.

Fix $f_1, f_2, f_3 \in C_c(\mathcal{G}, \mathbb{Z})$ and $g \in \mathcal{G}$. Then

$$
\begin{aligned}
((f_1 * f_2) * f_3)(g) &= \sum_{\substack{(k,h_3)\in\mathcal{G}_2\\ kh_3=g}} (f_1 * f_2)(k) f_3(h_3) \\
&= \sum_{\substack{(k,h_3)\in\mathcal{G}_2\\ kh_3=g}} \left( \sum_{\substack{(h_1,h_2)\in\mathcal{G}_2\\ h_1h_2=k}} f_1(h_1) f_2(h_2) \right) f_3(h_3) \\
&= \sum_{\substack{(k,h_3)\in\mathcal{G}_2\\ kh_3=g}} \sum_{\substack{(h_1,h_2)\in\mathcal{G}_2\\ h_1h_2=k}} f_1(h_1) f_2(h_2) f_3(h_3).
\end{aligned}
$$

The map

$$
\begin{aligned}
\left\{(h_1,h_2,h_3) \in \mathcal{G}_3 \mid h_1h_2h_3 = g\right\} \to \left\{(k,h_3) \in \mathcal{G}_2 \mid kh_3 = g\right\} &\times \left\{(h_1,h_2) \in \mathcal{G}_2 \mid h_1h_2 = k\right\}, \\
(h_1,h_2,h_3) \mapsto (h_1h_2, h_3, h_1, h_2)&,
\end{aligned}
$$

is a bijection, hence

$$
((f_1 * f_2) * f_3)(g) = \sum_{\substack{(h_1,h_2,h_3)\in\mathcal{G}_3\\ h_1h_2h_3=g}} f_1(h_1) f_2(h_2) f_3(h_3).
$$

Similarly,

$$
\begin{aligned}
(f_1 * (f_2 * f_3))(g) &= \sum_{\substack{(h_1,k)\in\mathcal{G}_2\\ h_1k=g}} f_1(h_1)\, (f_2 * f_3)(k) \\
&= \sum_{\substack{(h_1,k)\in\mathcal{G}_2\\ h_1k=g}} f_1(h_1) \left( \sum_{\substack{(h_2,h_3)\in\mathcal{G}_2\\ h_2h_3=k}} f_2(h_2) f_3(h_3) \right) \\
&= \sum_{\substack{(h_1,k)\in\mathcal{G}_2\\ h_1k=g}} \sum_{\substack{(h_2,h_3)\in\mathcal{G}_2\\ h_2h_3=k}} f_1(h_1) f_2(h_2) f_3(h_3).
\end{aligned}
$$

The map

$$
\begin{aligned}
\left\{(h_1,h_2,h_3) \in \mathcal{G}_3 \mid h_1h_2h_3 = g\right\} \to \left\{(h_1,k) \in \mathcal{G}_2 \mid h_1k = g\right\} &\times \left\{(h_2,h_3) \in \mathcal{G}_2 \mid h_2h_3 = k\right\}, \\
(h_1,h_2,h_3) \mapsto (h_1, h_2h_3, h_2, h_3)&,
\end{aligned}
$$

is a bijection, hence

$$
(f_1 * (f_2 * f_3))(g) = \sum_{\substack{(h_1,h_2,h_3)\in\mathcal{G}_3\\ h_1h_2h_3=g}} f_1(h_1) f_2(h_2) f_3(h_3).
$$

Therefore $(f_1 * f_2) * f_3 = f_1 * (f_2 * f_3)$, so $*$ is associative.

- **Local units.** Let $\{f_1, \dots, f_n\} \subset C_c(\mathcal{G}, \mathbb{Z})$ be finite and set

$$K := \bigcup_{i=1}^{n} \operatorname{supp}(f_i) \subset \mathcal{G}, \quad K_0 := r(K) \cup s(K) \subset \mathcal{G}_0.$$

  Since $r$ and $s$ are continuous and $K$ is compact, the set $K_0$ is compact. As $\mathcal{G}$ is ample, $\mathcal{G}_0$ has a basis of compact open subsets, hence there exists a compact open set $U \subset \mathcal{G}_0$ with $K_0 \subseteq U$. Let $e := \chi_U \in C_c(\mathcal{G}_0, \mathbb{Z}) \subseteq C_c(\mathcal{G}, \mathbb{Z})$, viewed as a function supported on units. Then $e * e = e$, and we claim $e * f_i = f_i = f_i * e$ for all $i$. Fix $f \in \{f_1, \dots, f_n\}$ and $g \in \mathcal{G}$. If $g \notin K$, then $f(g) = 0$, and both identities are trivial. Assume $g \in K$, so $r(g), s(g) \in K_0 \subseteq U$. Then

$$(e * f)(g) = \sum_{h \in \mathcal{G}_{r(g)}} e(h^{-1}) f(hg).$$

  The term $e(h^{-1})$ is nonzero only if $h^{-1} \in \mathcal{G}_0$, that is, $h \in \mathcal{G}_0$, and the unique unit in $\mathcal{G}_{r(g)}$ is $r(g)$. Hence $(e * f)(g) = e(r(g)) f(r(g)g) = 1 \cdot f(g) = f(g)$. Similarly,

$$(f * e)(g) = \sum_{h \in \mathcal{G}_{r(g)}} f(h^{-1}) e(hg).$$

  The term $e(hg)$ is nonzero only if $hg \in \mathcal{G}_0$, that is, $h = g^{-1}$, the unique such element of $\mathcal{G}_{r(g)}$. Thus $(f * e)(g) = f(g) e(s(g)) = f(g)$. Therefore $e * f = f = f * e$ for all $f$ in the finite set, so $C_c(\mathcal{G}, \mathbb{Z})$ has local units.

□

## 2.5 Invariance under Kakutani Equivalence

Up to this point, we have developed Moore homology for ample groupoids as a computable invariant built from compactly supported chains on the nerve. In many situations, in particular for ample groupoids arising from Cantor dynamics [7, 10], Bratteli–Vershik models [10], and related symbolic constructions, the groupoid presentation is far from unique. The same orbit structure can be modelled by different groupoids obtained by passing to full clopen reductions, refining cross sections, or changing the chosen model. These operations preserve the orbit picture but can drastically simplify the combinatorics of the nerve, and are therefore indispensable for explicit computations and for comparison with other invariants.

Kakutani equivalence formalizes precisely this flexibility. It is weaker than isomorphism, but strong enough to preserve the orbit structure relevant in the ample setting, and it is the standard equivalence relation in the study of Cantor groupoids and orbit equivalence. Establishing invariance under Kakutani equivalence therefore serves two purposes:

1. It shows that Moore homology depends only on the intrinsic orbit geometry of an ample groupoid and not on incidental choices of presentation.

2. It provides a practical tool: one may replace $\mathcal{G}$ by any Kakutani equivalent model, typically a full clopen reduction, before applying the long exact sequences, Moore–Mayer–Vietoris arguments, and UCT results developed in Sections 3.1, 3.2, and 3.3.

This point of view also connects directly to the broader Morita-invariance philosophy for groupoids and their invariants: Kakutani equivalence is implemented by full clopen reductions and can be seen as a concrete ample analogue of Morita equivalence at the level of orbit geometry. Accordingly, the invariance proof below is carried out at the chain level, using functoriality and reduction exact sequences, rather than appealing to classifying-space homotopy type.

Let $\mathcal{G}$ and $\mathcal{H}$ be Kakutani equivalent ample étale groupoids, and let $A$ be an abelian group, viewed with the discrete topology when forming $C_c(-, A)$. We will prove that for every $n \geq 0$ there is a natural isomorphism $H_n(\mathcal{G}; A) \cong H_n(\mathcal{H}; A)$. We now recall the definition of Kakutani equivalence in the ample setting. Let $\mathcal{G}$ and $\mathcal{H}$ be étale groupoids whose unit spaces are compact and totally disconnected.

**Definition 2.5.1** (Kakutani equivalence [3, §3])**.** Étale groupoids $\mathcal{G}$ and $\mathcal{H}$ are Kakutani equivalent if there exist clopen, full subsets $F \subseteq \mathcal{G}_0$ and $E \subseteq \mathcal{H}_0$ and an isomorphism of étale groupoids $\phi\colon \mathcal{G}|_F \xrightarrow{\cong} \mathcal{H}|_E$, where $\mathcal{G}|_F$ denotes the reduction with arrow space $r^{-1}(F) \cap s^{-1}(F)$ and units $F$.

*Remark* 2.5.2. Here full means that the saturation equals the whole unit space, that is, $r(s^{-1}(F)) = \mathcal{G}_0$ and $r(s^{-1}(E)) = \mathcal{H}_0$.

**Definition 2.5.3** (Homological similarity [13, Definition 3.4])**.** Let $\mathcal{G}$ and $\mathcal{H}$ be étale groupoids.

- Functors $\rho, \sigma\colon \mathcal{G} \to \mathcal{H}$ are similar if there exists a continuous map $\theta\colon \mathcal{G}_0 \to \mathcal{H}$ such that, for all $x \in \mathcal{G}_0$,
$$s_{\mathcal{H}}(\theta(x)) = \rho_0(x), \quad r_{\mathcal{H}}(\theta(x)) = \sigma_0(x),$$
and for all $g \in \mathcal{G}$,
$$\theta(r_{\mathcal{G}}(g))\, \rho(g) = \sigma(g)\, \theta(s_{\mathcal{G}}(g)).$$
Thus $\theta\colon \rho \Rightarrow \sigma$ is a natural transformation.
- The groupoids $\mathcal{G}$ and $\mathcal{H}$ are homologically similar if there exist functors $\rho\colon \mathcal{G} \to \mathcal{H}$ and $\sigma\colon \mathcal{H} \to \mathcal{G}$ together with natural transformations
$$\theta_{\mathcal{G}}\colon \mathrm{id}_{\mathcal{G}} \Rightarrow \sigma \circ \rho, \quad \theta_{\mathcal{H}}\colon \mathrm{id}_{\mathcal{H}} \Rightarrow \rho \circ \sigma.$$

*Remark* 2.5.4. For $x \in \mathcal{G}_0$, the component $\theta(x)$ is an arrow in $\mathcal{H}$ from $\rho_0(x)$ to $\sigma_0(x)$. Identity arrows are written $u(y) \in \mathcal{H}$ for $y \in \mathcal{H}_0$.

We now state Kakutani invariance of Moore homology and then prove it by reducing to chain-level functoriality, reduction exact sequences, and homological similarity.

**Theorem 2.5.5** (Kakutani invariance of Moore homology [14, Theorem 4.8])**.** If $\mathcal{G}$ and $\mathcal{H}$ are Kakutani equivalent, then for any topological abelian group $A$ there are natural isomorphisms
$$H_n(\mathcal{G}; A) \cong H_n(\mathcal{H}; A) \quad \text{for all } n \geq 0.$$

For $A = \mathbb{Z}$ there is an isomorphism $\pi\colon H_0(\mathcal{G}) \to H_0(\mathcal{H})$ with $\pi(H_0(\mathcal{G})^+) = H_0(\mathcal{H})^+$.

*Remark* 2.5.6 (Positive cone in $H_0(\mathcal{G})$). Let $q\colon C_c(\mathcal{G}_0, \mathbb{Z}) \to H_0(\mathcal{G})$ be the quotient map $q(f) = [f]$. Define $H_0(\mathcal{G})^+ := q(C_c(\mathcal{G}_0, \mathbb{Z}_{\geq 0}))$. In other words, $H_0(\mathcal{G})^+$ is the submonoid of $H_0(\mathcal{G})$ generated by the classes $[\chi_U]$ of characteristic functions of clopen subsets $U \subseteq \mathcal{G}_0$. When $\mathcal{G}_0$ is totally disconnected, every nonnegative integer-valued compactly supported continuous function is a finite sum of characteristic functions of clopen sets, so the two descriptions agree. The pair $(H_0(\mathcal{G}), H_0(\mathcal{G})^+)$ is an ordered abelian group, and the maps induced on $H_0$ by functors of étale groupoids are positive with respect to these cones.

**Lemma 2.5.7.** Let $\rho : \mathcal{G} \to \mathcal{H}$ be an étale functor of étale groupoids and let $A$ be a topological abelian group. For $n \geq 0$ define

$$\rho_n : \mathcal{G}_n \to \mathcal{H}_n, \quad (g_1, \dots, g_n) \mapsto (\rho_1(g_1), \dots, \rho_1(g_n)),$$

with $\rho_0$ the map on units. Then each $\rho_n$ is a local homeomorphism. Moreover, pushforward along $\rho_n$ defines group homomorphisms

$$(\rho_n)_* : C_c(\mathcal{G}_n, A) \to C_c(\mathcal{H}_n, A), \quad (\rho_n)_*(f)(\mathbf{h}) := \sum_{\mathbf{g} \in \rho_n^{-1}(\mathbf{h})} f(\mathbf{g}),$$

and the maps $(\rho_n)_*$ form a chain map between Moore complexes,

$$\partial_n^{\mathcal{H}} \circ (\rho_n)_* = (\rho_{n-1})_* \circ \partial_n^{\mathcal{G}} \quad \text{for all } n \geq 1.$$

*Proof.* We divide the argument into three steps.

- **The maps $\rho_n$ are local homeomorphisms.** For $n \geq 2$ write the $n$-simplices as iterated fibre products

$$\mathcal{G}_n = \mathcal{G}_1 \,{}_{s_{\mathcal{G}}}\!\times_{r_{\mathcal{G}}} \cdots \,{}_{s_{\mathcal{G}}}\!\times_{r_{\mathcal{G}}} \mathcal{G}_1, \quad \mathcal{H}_n = \mathcal{H}_1 \,{}_{s_{\mathcal{H}}}\!\times_{r_{\mathcal{H}}} \cdots \,{}_{s_{\mathcal{H}}}\!\times_{r_{\mathcal{H}}} \mathcal{H}_1.$$

  Since $\rho$ is an étale functor, the maps $\rho_0 : \mathcal{G}_0 \to \mathcal{H}_0$ and $\rho_1 : \mathcal{G}_1 \to \mathcal{H}_1$ are local homeomorphisms and satisfy $s_{\mathcal{H}} \circ \rho_1 = \rho_0 \circ s_{\mathcal{G}}$ and $r_{\mathcal{H}} \circ \rho_1 = \rho_0 \circ r_{\mathcal{G}}$. Thus $\rho_n$ is the induced map between these iterated pullbacks. Local homeomorphisms are stable under finite products and base change, hence $\rho_n$ is a local homeomorphism for all $n \geq 0$.
- **The maps $\rho_\bullet$ form a simplicial map.** Let $n \geq 1$ and $\mathbf{g} = (g_1, \dots, g_n) \in \mathcal{G}_n$.
  For the endpoint faces,

$$d_0^{\mathcal{H}}(\rho_n(\mathbf{g})) = (\rho_1(g_2), \dots, \rho_1(g_n)) = \rho_{n-1}(d_0^{\mathcal{G}}(\mathbf{g})),$$
$$d_n^{\mathcal{H}}(\rho_n(\mathbf{g})) = (\rho_1(g_1), \dots, \rho_1(g_{n-1})) = \rho_{n-1}(d_n^{\mathcal{G}}(\mathbf{g})).$$

  For $1 \leq i \leq n-1$, using functoriality of $\rho$,

$$\begin{aligned} d_i^{\mathcal{H}}(\rho_n(\mathbf{g})) &= (\rho_1(g_1), \dots, \rho_1(g_i)\rho_1(g_{i+1}), \dots, \rho_1(g_n)) \\ &= (\rho_1(g_1), \dots, \rho_1(g_i g_{i+1}), \dots, \rho_1(g_n)) = \rho_{n-1}(d_i^{\mathcal{G}}(\mathbf{g})). \end{aligned}$$

Hence $d_i^{\mathcal{H}} \circ \rho_n = \rho_{n-1} \circ d_i^{\mathcal{G}}$ for all $n \geq 1$ and $0 \leq i \leq n$.

- **Pushforward yields a chain map.** Since $\rho_n$ is a local homeomorphism, pushforward $(\rho_n)_*$ is well defined on $C_c(\mathcal{G}_n, A)$. For each $n \geq 1$ and $0 \leq i \leq n$, the previous step gives $d_i^{\mathcal{H}} \circ \rho_n = \rho_{n-1} \circ d_i^{\mathcal{G}}$. By functoriality of pushforward from Proposition 2.1.2, $(d_i^{\mathcal{H}})_* \circ (\rho_n)_* = (\rho_{n-1})_* \circ (d_i^{\mathcal{G}})_*$.

Therefore, for $c \in C_c(\mathcal{G}_n, A)$,

$$\begin{aligned} \partial_n^{\mathcal{H}}((\rho_n)_* c) &= \sum_{i=0}^{n} (-1)^i (d_i^{\mathcal{H}})_* (\rho_n)_* c = \sum_{i=0}^{n} (-1)^i (\rho_{n-1})_* (d_i^{\mathcal{G}})_* c \\ &= (\rho_{n-1})_* \left( \sum_{i=0}^{n} (-1)^i (d_i^{\mathcal{G}})_* c \right) = (\rho_{n-1})_* (\partial_n^{\mathcal{G}} c), \end{aligned}$$

which is the desired chain map identity. For $n = 0$ we have $\partial_0^{\mathcal{G}} = 0, \partial_0^{\mathcal{H}} = 0$. □

**Proposition 2.5.8** (Chain homotopy from a similarity [13, Proposition 3.5])**.** Let $\rho, \sigma : \mathcal{G} \to \mathcal{H}$ be étale functors between étale groupoids, and suppose they are similar via $\theta : \mathcal{G}_0 \to \mathcal{H}_1$, meaning $s_{\mathcal{H}}(\theta(x)) = \rho_0(x)$, $r_{\mathcal{H}}(\theta(x)) = \sigma_0(x)$ for all $x \in \mathcal{G}_0$, and $\theta(r(g)) \cdot_{\mathcal{H}} \rho_1(g) = \sigma_1(g) \cdot_{\mathcal{H}} \theta(s(g))$ for all $g \in \mathcal{G}_1$. Assume moreover that $\theta$ is étale, that is, $\theta(\mathcal{G}_0)$ is a bisection of $\mathcal{H}_1$, so $s_{\mathcal{H}} \circ \theta$ and $r_{\mathcal{H}} \circ \theta$ are local homeomorphisms. Define $k_0 : \mathcal{G}_0 \to \mathcal{H}_1$ by $k_0 := \theta$.

For $n \geq 1$ and $0 \leq j \leq n$ define

$$k_j : \mathcal{G}_n \to \mathcal{H}_{n+1},$$

$$k_j(\mathbf{g}) := \begin{cases} (\theta(r(g_1)), \rho_1(g_1), \dots, \rho_1(g_n)), & \text{for } j = 0, \\ (\sigma_1(g_1), \dots, \sigma_1(g_j), \theta(s(g_j)), \rho_1(g_{j+1}), \dots, \rho_1(g_n)), & \text{for } 1 \leq j \leq n-1, \\ (\sigma_1(g_1), \dots, \sigma_1(g_n), \theta(s(g_n))), & \text{for } j = n, \end{cases}$$

which are local homeomorphisms. Set $h_0 := (k_0)_* : C_c(\mathcal{G}_0, A) \to C_c(\mathcal{H}_1, A)$, and for $n \geq 1$ set

$$h_n := \sum_{j=0}^{n} (-1)^j (k_j)_* : C_c(\mathcal{G}_n, A) \to C_c(\mathcal{H}_{n+1}, A).$$

Then $h_\bullet$ is a chain homotopy $h_\bullet : (\rho_n)_* \Rightarrow (\sigma_n)_*$, meaning

$$\partial_1^{\mathcal{H}} \circ h_0 = (\rho_0)_* - (\sigma_0)_*, \quad \partial_{n+1}^{\mathcal{H}} \circ h_n + h_{n-1} \circ \partial_n^{\mathcal{G}} = (\rho_n)_* - (\sigma_n)_* \quad \text{for } n \geq 1.$$

*Proof.* For $n \geq 1$ and $0 \leq j \leq n$, the maps $k_j : \mathcal{G}_n \to \mathcal{H}_{n+1}$ are local homeomorphisms, hence $(k_j)_*$ is defined and preserves compact supports. They satisfy the face identities

$$d_0 k_0 = \rho_n, \tag{2.5.1}$$

$$d_{n+1} k_n = \sigma_n, \tag{2.5.2}$$

$$d_i k_j = k_{j-1} d_i \quad \text{for } 1 \leq i \leq j-1, \tag{2.5.3}$$

$$d_i k_j = k_j d_{i-1} \quad \text{for } j+2 \le i \le n+1, \tag{2.5.4}$$

$$d_{j+1} k_j = d_j k_{j+1} \quad \text{for } 0 \le j \le n-1. \tag{2.5.5}$$

The identities (2.5.1) and (2.5.2) are immediate by deleting the inserted $\theta$ at the front or back. The identities (2.5.3) and (2.5.4) use functoriality of $\sigma$ and $\rho$ to merge adjacent factors away from the $\theta$-slot. The identity (2.5.5) is the similarity relation $\theta(r(g)) \cdot_{\mathcal{H}} \rho_1(g) = \sigma_1(g) \cdot_{\mathcal{H}} \theta(s(g))$ applied at the junction $s(g_j) = r(g_{j+1})$.

Let $n \geq 1$. Using functoriality of pushforward from Proposition 2.1.2,

$$\partial^{\mathcal{H}}_{n+1} h_n = \sum_{i=0}^{n+1} \sum_{j=0}^{n} (-1)^{i+j} (d_i k_j)_*.$$

We use (2.5.1) and (2.5.2) to isolate the edge terms $i = 0$ and $i = n+1$, then (2.5.3) and (2.5.4) to move faces past $k_j$, and obtain

$$\partial^{\mathcal{H}}_{n+1} h_n = (\rho_n)_* - (\sigma_n)_* - \sum_{j=0}^{n-1} \sum_{q=0}^{n} (-1)^{q+j} (k_j d_q)_* + \sum_{j=0}^{n-1} \big( (k_j d_{j+1})_* - (d_j k_{j+1})_* \big).$$

By functoriality of pushforward, the double sum equals $-h_{n-1}\, \partial^{\mathcal{G}}_n$, and the bracket cancels termwise by (2.5.5). Hence

$$\partial^{\mathcal{H}}_{n+1} h_n + h_{n-1}\, \partial^{\mathcal{G}}_n = (\rho_n)_* - (\sigma_n)_* \quad \text{for } n \geq 1.$$

For $n = 0$, we have $h_0 = \theta_*$, and $d_0 \theta = \rho_0$, $d_1 \theta = \sigma_0$, hence

$$\partial^{\mathcal{H}}_1 h_0 = (d_0)_* \theta_* - (d_1)_* \theta_* = (\rho_0)_* - (\sigma_0)_*.$$

□

**Corollary 2.5.9.** If $\rho, \sigma : \mathcal{G} \to \mathcal{H}$ are similar, then the induced maps on homology agree, $H_n(\rho) = H_n(\sigma)$ for all $n \geq 0$.

*Proof.* By Proposition 2.5.8 the induced chain maps $(\rho_\bullet)_*$ and $(\sigma_\bullet)_*$ on Moore complexes are chain homotopic. Chain-homotopic maps induce the same morphisms on homology, hence $H_n(\rho) = H_n(\sigma)$ for all $n \geq 0$. □

**Proposition 2.5.10** ([13, Theorem 4.8])**.** If $\mathcal{G}$ and $\mathcal{H}$ are homologically similar, see Definition 2.5.3, then for any topological abelian group $A$ the induced maps $H_n(\rho) : H_n(\mathcal{G}; A) \to H_n(\mathcal{H}; A)$ are isomorphisms with inverse $H_n(\sigma)$ for all $n \geq 0$. If $A = \mathbb{Z}$, then $H_0(\rho)$ carries the positive cone $H_0(\mathcal{G})^+$ onto $H_0(\mathcal{H})^+$.

*Proof.* By Definition 2.5.3 there exist étale functors $\rho : \mathcal{G} \to \mathcal{H}$ and $\sigma : \mathcal{H} \to \mathcal{G}$ such that $\sigma \circ \rho$ is similar to $\mathrm{id}_{\mathcal{G}}$ and $\rho \circ \sigma$ is similar to $\mathrm{id}_{\mathcal{H}}$. By Corollary 2.5.9,

$$H_n(\sigma \circ \rho) = H_n(\mathrm{id}_{\mathcal{G}}) = \mathrm{id}_{H_n(\mathcal{G};A)}, \quad H_n(\rho \circ \sigma) = H_n(\mathrm{id}_{\mathcal{H}}) = \mathrm{id}_{H_n(\mathcal{H};A)}.$$

Hence $H_n(\rho)$ is an isomorphism with inverse $H_n(\sigma)$ for all $n \geq 0$.

Let $f \in C_c(\mathcal{G}_0, \mathbb{Z})$ be pointwise nonnegative. Then for $y \in \mathcal{H}_0$,

$$(\rho_0)_* f(y) = \sum_{x \in \rho_0^{-1}(y)} f(x) \in \mathbb{Z}_{\geq 0},$$

so $H_0(\rho)$ maps $H_0(\mathcal{G})^+$ into $H_0(\mathcal{H})^+$. The same argument shows that $H_0(\sigma)$ maps $H_0(\mathcal{H})^+$ into $H_0(\mathcal{G})^+$. Since $H_0(\sigma)$ is the inverse of $H_0(\rho)$, we conclude $H_0(\rho)(H_0(\mathcal{G})^+) = H_0(\mathcal{H})^+$. □

**Lemma 2.5.11** ([13, Theorem 3.6(1)]). Let $\mathcal{G}$ be an étale groupoid and let $F \subset \mathcal{G}_0$ be an open $\mathcal{G}$-full subset. Suppose there exists a continuous map $\theta : \mathcal{G}_0 \to \mathcal{G}_1$ such that $r_\mathcal{G}(\theta(x)) = x, s_\mathcal{G}(\theta(x)) \in F$ for all $x \in \mathcal{G}_0$, and whose image $\theta(\mathcal{G}_0)$ is a bisection, that is, $r_\mathcal{G} \circ \theta = \mathrm{id}_{\mathcal{G}_0}$ and $s_\mathcal{G} \circ \theta : \mathcal{G}_0 \to F$ are local homeomorphisms. Define functors

$$\rho : \mathcal{G} \to \mathcal{G}|_F, \quad \rho_0(x) := s_\mathcal{G}(\theta(x)), \quad \rho_1(g) := \theta(r_\mathcal{G}(g))^{-1} \cdot_\mathcal{G} g \cdot_\mathcal{G} \theta(s_\mathcal{G}(g)),$$

and the inclusion functor $\sigma : \mathcal{G}|_F \hookrightarrow \mathcal{G}$ given by

$$(\mathcal{G}|_F)_0 = F, \quad (\mathcal{G}|_F)_1 = \{h \in \mathcal{G}_1 \mid r_\mathcal{G}(h) \in F,\ s_\mathcal{G}(h) \in F\},$$

$$\sigma_0 : F \to \mathcal{G}_0, \quad \sigma_0(x) = x, \quad \sigma_1 : (\mathcal{G}|_F)_1 \to \mathcal{G}_1, \quad \sigma_1(h) = h.$$

Then $\rho$ and $\sigma$ are étale functors and $\mathcal{G}$ and $\mathcal{G}|_F$ are homologically similar in the sense of Definition 2.5.3. More precisely:

- $\sigma \circ \rho$ is similar to $\mathrm{id}_\mathcal{G}$ via $\theta$, meaning $\theta(r_\mathcal{G}(g)) \cdot_\mathcal{G} (\sigma \circ \rho)_1(g) = g \cdot_\mathcal{G} \theta(s_\mathcal{G}(g))$ for all $g \in \mathcal{G}_1$.
- $\rho \circ \sigma$ is similar to $\mathrm{id}_{\mathcal{G}|_F}$ via $\theta|_F$, meaning $\theta(r_\mathcal{G}(h)) \cdot_\mathcal{G} (\rho \circ \sigma)_1(h) = h \cdot_\mathcal{G} \theta(s_\mathcal{G}(h))$ for all $h \in (\mathcal{G}|_F)_1$.

*Proof.*

- **Well-definedness and functoriality of $\rho$.** First, $\rho_0(x) = s_\mathcal{G}(\theta(x)) \in F$, so $\rho_0 : \mathcal{G}_0 \to F = (\mathcal{G}|_F)_0$ is well defined. For arrows $g \in \mathcal{G}_1$,

$$r_\mathcal{G}(\rho_1(g)) = r_\mathcal{G}(\theta(r_\mathcal{G}(g))^{-1} \cdot_\mathcal{G} g \cdot_\mathcal{G} \theta(s_\mathcal{G}(g))) = s_\mathcal{G}(\theta(r_\mathcal{G}(g))) = \rho_0(r_\mathcal{G}(g)) \in F,$$
$$s_\mathcal{G}(\rho_1(g)) = s_\mathcal{G}(\theta(s_\mathcal{G}(g))) = \rho_0(s_\mathcal{G}(g)) \in F,$$

hence $\rho_1(g) \in (\mathcal{G}|_F)_1$ and $\rho$ respects range and source. If $g, h$ are composable in $\mathcal{G}$ with $s_\mathcal{G}(g) = r_\mathcal{G}(h)$, then

$$\begin{aligned}\rho_1(g) \cdot_\mathcal{G} \rho_1(h) &= \theta(r_\mathcal{G}(g))^{-1} \cdot_\mathcal{G} g \cdot_\mathcal{G} \theta(s_\mathcal{G}(g)) \cdot_\mathcal{G} \theta(r_\mathcal{G}(h))^{-1} \cdot_\mathcal{G} h \cdot_\mathcal{G} \theta(s_\mathcal{G}(h)) \\ &= \theta(r_\mathcal{G}(g))^{-1} \cdot_\mathcal{G} g \cdot_\mathcal{G} (\theta(s_\mathcal{G}(g)) \cdot_\mathcal{G} \theta(r_\mathcal{G}(h))^{-1}) \cdot_\mathcal{G} h \cdot_\mathcal{G} \theta(s_\mathcal{G}(h)) \\ &= \theta(r_\mathcal{G}(g))^{-1} \cdot_\mathcal{G} g \cdot_\mathcal{G} h \cdot_\mathcal{G} \theta(s_\mathcal{G}(h)) = \rho_1(g \cdot_\mathcal{G} h),\end{aligned}$$

since $s_{\mathcal{G}}(g) = r_{\mathcal{G}}(h)$ implies $\theta(s_{\mathcal{G}}(g)) = \theta(r_{\mathcal{G}}(h))$, hence $\theta(s_{\mathcal{G}}(g)) \cdot_{\mathcal{G}} \theta(r_{\mathcal{G}}(h))^{-1} = u_{\mathcal{G}}(r_{\mathcal{G}}(\theta(s_{\mathcal{G}}(g))))$. Units are preserved because

$$\rho_1(u_{\mathcal{G}}(x)) = \theta(x)^{-1} \cdot_{\mathcal{G}} u_{\mathcal{G}}(x) \cdot_{\mathcal{G}} \theta(x) = u_{\mathcal{G}}(s_{\mathcal{G}}(\theta(x))) = u_{\mathcal{G}}(\rho_0(x)).$$

- **Continuity.** Follows from continuity of $r_{\mathcal{G}}, s_{\mathcal{G}}$, inversion, multiplication, and $\theta$.
- **Étaleness.** The map $\rho_0 = s_{\mathcal{G}} \circ \theta$ is a local homeomorphism because $\theta(\mathcal{G}_0)$ is a bisection and $s_{\mathcal{G}}$ is a local homeomorphism. The map $\rho_1$ is a local homeomorphism because it is built from the local homeomorphisms $r_{\mathcal{G}}, s_{\mathcal{G}}$, inversion, multiplication restricted to products of bisections, and $\theta$, using stability of local homeomorphisms under finite products and base change. The inclusion $\sigma : \mathcal{G}|_F \hookrightarrow \mathcal{G}$ is étale because $(\mathcal{G}|_F)_1 = s_{\mathcal{G}}^{-1}(F) \cap r_{\mathcal{G}}^{-1}(F)$ is open in $\mathcal{G}_1$ and $\sigma_0, \sigma_1$ are open embeddings.
- **Similarity $\sigma \circ \rho \sim \mathrm{id}_{\mathcal{G}}$ via $\theta$.** For objects $x \in \mathcal{G}_0$,

$$r_{\mathcal{G}}(\theta(x)) = x = \mathrm{id}_{\mathcal{G}_0}(x), \quad s_{\mathcal{G}}(\theta(x)) = \rho_0(x) = (\sigma \circ \rho)_0(x).$$

  For $g \in \mathcal{G}_1$,

$$\theta(r_{\mathcal{G}}(g)) \cdot_{\mathcal{G}} (\sigma \circ \rho)_1(g) = \theta(r_{\mathcal{G}}(g)) \cdot_{\mathcal{G}} (\theta(r_{\mathcal{G}}(g))^{-1} \cdot_{\mathcal{G}} g \cdot_{\mathcal{G}} \theta(s_{\mathcal{G}}(g))) = g \cdot_{\mathcal{G}} \theta(s_{\mathcal{G}}(g)),$$

  so $\sigma \circ \rho$ is similar to $\mathrm{id}_{\mathcal{G}}$ via $\theta$.
- **Similarity $\rho \circ \sigma \sim \mathrm{id}_{\mathcal{G}|_F}$ via $\theta|_F$.** Let $h \in (\mathcal{G}|_F)_1$. Then $r_{\mathcal{G}}(h), s_{\mathcal{G}}(h) \in F$ and

$$\theta(r_{\mathcal{G}}(h)) \cdot_{\mathcal{G}} (\rho \circ \sigma)_1(h) = \theta(r_{\mathcal{G}}(h)) \cdot_{\mathcal{G}} (\theta(r_{\mathcal{G}}(h))^{-1} \cdot_{\mathcal{G}} h \cdot_{\mathcal{G}} \theta(s_{\mathcal{G}}(h))) = h \cdot_{\mathcal{G}} \theta(s_{\mathcal{G}}(h)).$$

  Moreover, for $x \in F$ we have $r_{\mathcal{G}}(\theta(x)) = x$ and $s_{\mathcal{G}}(\theta(x)) \in F = (\rho \circ \sigma)_0(x)$.

□

In the Kakutani invariance argument we repeatedly pass to full open reductions $\mathcal{G}|_F$ and compare $\mathcal{G}$ with $\mathcal{G}|_F$ by explicit étale functors. At chain level this comparison is implemented by pushforwards along the induced maps on nerves, so we need a concrete way to transport compactly supported chains on $\mathcal{G}_\bullet$ into compactly supported chains on $(\mathcal{G}|_F)_\bullet$ and back.

The key technical input is a continuous choice of arrows that move points into $F$: we seek a map $\theta : \mathcal{G}_0 \to \mathcal{G}_1$ with $r(\theta(x)) = x$ and $s(\theta(x)) \in F$. Such a map is a global section of the range map $r : \mathcal{G}_1 \to \mathcal{G}_0$ with values landing in the open subset $s^{-1}(F) \subset \mathcal{G}_1$. Fullness of $F$ guarantees existence of arrows pointwise, but a continuous choice requires patching local bisections coherently.

This is exactly where $\sigma$-compactness enters. Since $\mathcal{G}$ is étale and $\mathcal{G}_0$ is totally disconnected, around each $x \in \mathcal{G}_0$ one can find a compact open bisection $U_x \subset s^{-1}(F)$ with $x \in r(U_x)$ and $s(U_x) \subset F$. A continuous global section can then be obtained by selecting a countable family of such bisections covering $\mathcal{G}_0$ and refining it to a disjoint cover of $\mathcal{G}_0$ by compact open subsets in the range. The refinement step uses only total disconnectedness and Hausdorffness, but the reduction to a countable cover uses $\sigma$-compactness in an essential way, by applying compactness on an exhausting sequence of compact subsets. Once a disjoint compact open cover $\mathcal{G}_0 = \bigsqcup_n r(V_n)$ with bisections $V_n \subset s^{-1}(F)$ is available, the section is forced on each piece by $\theta|_{r(V_n)} = (r|_{V_n})^{-1}$, and continuity follows because the pieces are open.

Lemma 2.5.12 formalizes this construction. It is the mechanism that turns the geometric fullness hypothesis into an explicit, continuous arrow selection, which in turn produces the functors needed for homological similarity and hence for Kakutani invariance.

**Lemma 2.5.12** ([13, Lemma 4.3])**.** Let $\mathcal{G}$ be an étale groupoid with $\mathcal{G}_0$ $\sigma$-compact and totally disconnected, and let $F \subset \mathcal{G}_0$ be an open $\mathcal{G}$-full subset. Then there exists a continuous map $\theta : \mathcal{G}_0 \to \mathcal{G}_1$ such that $r(\theta(x)) = x$ and $s(\theta(x)) \in F$ for all $x \in \mathcal{G}_0$.

*Proof.* For each $x \in \mathcal{G}_0$ choose $g_x \in \mathcal{G}_1$ with $r(g_x) = x$ and $s(g_x) \in F$, which is possible since $F$ is $\mathcal{G}$-full. As $\mathcal{G}$ is étale, there exists an open bisection $B_x \ni g_x$. Since $F$ is open and $s(g_x) \in F$, replacing $B_x$ by $B_x \cap s^{-1}(F)$ we may assume $s(B_x) \subseteq F$.

Because $\mathcal{G}_0$ is totally disconnected and locally compact, we can choose a compact open neighbourhood $W_x \subset r(B_x)$ of $x$ and set $U_x := (r|_{B_x})^{-1}(W_x)$. Then $U_x$ is a compact open bisection, $r(U_x) = W_x$ is compact open, and $s(U_x) \subseteq s(B_x) \subseteq F$.

Since $\mathcal{G}_0$ is $\sigma$-compact, choose an increasing sequence of compact sets $(K_m)_{m\geq 1}$ with $\mathcal{G}_0 = \bigcup_{m\geq 1} K_m$. Fix $m \geq 1$. The family $\{r(U_x) \mid x \in K_m\}$ covers $K_m$, hence by compactness there exist points $x_{m,1}, \dots, x_{m,N_m} \in K_m$ such that $K_m \subset \bigcup_{i=1}^{N_m} r(U_{x_{m,i}})$. Write $U_{m,i} := U_{x_{m,i}}$. The index set $\bigsqcup_{m\geq 1}\{1, \dots, N_m\}$ is countable, hence we can enumerate the corresponding bisections as a sequence $(U_n)_{n\geq 1}$. Then

$$\bigcup_{n\geq 1} r(U_n) = \mathcal{G}_0, \quad s(U_n) \subseteq F \text{ for all } n \geq 1.$$

Since $\mathcal{G}_0$ is Hausdorff, every compact subset is closed, hence every compact open subset of $\mathcal{G}_0$ is clopen. In particular, finite unions and complements of compact open subsets are again compact open. Define compact open bisections $(V_n)_{n\geq 1}$ inductively by

$$V_1 := U_1, \quad V_n := U_n \setminus r^{-1}\Big(r\Big(\bigcup_{j=1}^{n-1} V_j\Big)\Big) = (r|_{U_n})^{-1}\left(r(U_n) \setminus \bigcup_{j=1}^{n-1} r(V_j)\right) \quad \text{for } n \geq 2.$$

Here $r|_{U_n} : U_n \to r(U_n)$ is a homeomorphism and $r(U_n) \setminus \bigcup_{j=1}^{n-1} r(V_j)$ is compact open in $\mathcal{G}_0$, hence $V_n$ is compact open in $U_n$, in particular open in $\mathcal{G}_1$. Being an open subset of the bisection $U_n$, each $V_n$ is again a bisection and $s(V_n) \subseteq s(U_n) \subseteq F$. Moreover, the ranges $r(V_n)$ are pairwise disjoint and

$$\bigcup_{n\geq 1} r(V_n) = \bigcup_{n\geq 1}\left(r(U_n) \setminus \bigcup_{j=1}^{n-1} r(V_j)\right) = \bigcup_{n\geq 1} r(U_n) = \mathcal{G}_0.$$

Define $\theta : \mathcal{G}_0 \to \mathcal{G}_1$ by $\theta(x) := (r|_{V_n})^{-1}(x)$ for the unique $n$ with $x \in r(V_n)$. This is well defined since $\{r(V_n)\}_{n\geq 1}$ is a pairwise disjoint cover of $\mathcal{G}_0$, and each $r|_{V_n}$ is a homeomorphism.

Since each $r(V_n)$ is open and $\theta|_{r(V_n)} = (r|_{V_n})^{-1}$ is continuous, $\theta$ is continuous. Finally, $r(\theta(x)) = x$ and $s(\theta(x)) \in s(V_n) \subseteq F$ for all $x \in \mathcal{G}_0$. □

This construction is the basic compact open partition argument that underlies several chain-level tools later on, such as reduction to full clopen subsets and the resulting exact sequences in the ample setting.

In the proof of Kakutani invariance in Theorem 2.5.5 one repeatedly replaces a groupoid by a full reduction $\mathcal{G}|_F$ to a suitable open, typically clopen, subset $F \subseteq \mathcal{G}_0$. Geometrically, such a reduction does not change the orbit picture: every $\mathcal{G}$-orbit meets $F$, so $\mathcal{G}|_F$ still sees all orbits. To make this usable on the level of Moore chains, one needs a concrete way to move units into $F$ by arrows depending continuously on the unit. This is the role of the next step: Lemma 2.5.12 constructs a continuous section of the range map with image in $s^{-1}(F)$, and Lemma 2.5.11 turns such a section into étale functors $\rho : \mathcal{G} \to \mathcal{G}|_F$ and $\sigma : \mathcal{G}|_F \hookrightarrow \mathcal{G}$ whose composites are similar to the identities. Since similarity yields chain homotopies by Proposition 2.5.8, this produces the chain-level comparison needed for invariance.

**Proposition 2.5.13.** Let $\mathcal{G}$ be étale with $\mathcal{G}_0$ $\sigma$-compact and totally disconnected, and let $F \subseteq \mathcal{G}_0$ be open and $\mathcal{G}$-full. Then $\mathcal{G}$ and $\mathcal{G}|_F$ are homologically similar.

*Proof.* By Lemma 2.5.12 there exists a continuous map $\theta : \mathcal{G}_0 \to \mathcal{G}_1$ such that $r(\theta(x)) = x$ and $s(\theta(x)) \in F$ for all $x \in \mathcal{G}_0$. Moreover, the construction in Lemma 2.5.12 produces $\theta(\mathcal{G}_0)$ as a union of pairwise disjoint compact open bisections, hence $\theta(\mathcal{G}_0)$ is itself a bisection. Therefore $r \circ \theta = \mathrm{id}_{\mathcal{G}_0}$ and $s \circ \theta : \mathcal{G}_0 \to F$ are local homeomorphisms.

Applying Lemma 2.5.11 to this $\theta$ yields étale functors $\rho : \mathcal{G} \to \mathcal{G}|_F$ and $\sigma : \mathcal{G}|_F \hookrightarrow \mathcal{G}$ such that $\sigma \circ \rho$ is similar to $\mathrm{id}_\mathcal{G}$ via $\theta$ and $\rho \circ \sigma$ is similar to $\mathrm{id}_{\mathcal{G}_F}$ via $\theta|_F$. Hence $\mathcal{G}$ and $\mathcal{G}|_F$ are homologically similar in the sense of Definition 2.5.3. □

*Proof of Theorem 2.5.5.*

- **Reduction to full clopen subsets.** By Kakutani equivalence in Definition 2.5.1, there exist clopen, $\mathcal{G}$-full $F_\mathcal{G} \subseteq \mathcal{G}_0$ and clopen, $\mathcal{H}$-full $F_\mathcal{H} \subseteq \mathcal{H}_0$, and an isomorphism of étale groupoids $\phi : \mathcal{G}|_{F_\mathcal{G}} \xrightarrow{\cong} \mathcal{H}|_{F_\mathcal{H}}$. Since $\mathcal{G}_0$ and $\mathcal{H}_0$ are compact, they are $\sigma$-compact. Proposition 2.5.13 yields homological similarities $\mathcal{G} \sim_h \mathcal{G}|_{F_\mathcal{G}}, \mathcal{H} \sim_h \mathcal{H}|_{F_\mathcal{H}}$.
- **Identify the reductions.** Set $\rho := \phi$ and $\sigma := \phi^{-1}$. Then $\rho$ and $\sigma$ are étale functors and $\sigma \circ \rho = \mathrm{id}_{\mathcal{G}_{F_\mathcal{G}}}, \rho \circ \sigma = \mathrm{id}_{\mathcal{H}_{F_\mathcal{H}}}$. In particular, $\mathcal{G}|_{F_\mathcal{G}} \sim_h \mathcal{H}|_{F_\mathcal{H}}$.
- **Transitivity and conclusion.** Combining the previous steps gives $\mathcal{G} \sim_h \mathcal{G}|_{F_\mathcal{G}} \sim_h \mathcal{H}|_{F_\mathcal{H}} \sim_h \mathcal{H}$, hence $\mathcal{G} \sim_h \mathcal{H}$. By Proposition 2.5.10, for every topological abelian group $A$ the induced maps are isomorphisms $H_n(\mathcal{G}; A) \xrightarrow{\cong} H_n(\mathcal{H}; A)$ for all $n \geq 0$. For $A = \mathbb{Z}$, Proposition 2.5.10 also yields that the induced isomorphism on $H_0$ carries $H_0(\mathcal{G})^+$ onto $H_0(\mathcal{H})^+$.

□

We define a homology theory for étale groupoids by using the simplicial geometry of the nerve $\mathcal{G}_\bullet = (\mathcal{G}_n, (d_i)_{i=0}^n, (s_j)_{j=0}^n)_{n \geq 0}$. If $\mathcal{G}$ is étale, then every face map $d_i : \mathcal{G}_n \to \mathcal{G}_{n-1}$ is a local homeomorphism. Hence, for every topological abelian group $A$, pushforward along $d_i$ is defined on compactly supported continuous functions. Applying $C_c(-, A)$ levelwise yields the simplicial abelian group $C_c(\mathcal{G}_\bullet, A) = (C_c(\mathcal{G}_n, A))_{n\geq 0}$, and the Moore boundary

$$\partial_0 := 0 : \; C_c(\mathcal{G}_0, A) \to 0,$$
$$\partial_n := \sum_{i=0}^{n} (-1)^i (d_i)_* : \; C_c(\mathcal{G}_n, A) \to C_c(\mathcal{G}_{n-1}, A).$$

Functoriality of pushforward and the simplicial identities imply $\partial_{n-1}\partial_n = 0$. Thus $C_c(\mathcal{G}_\bullet, A) := (C_c(\mathcal{G}_n, A), \partial_n)_{n\geq 0}$ is a chain complex, and we set $H_n(\mathcal{G}; A) := \ker(\partial_n) / \operatorname{im}(\partial_{n+1})$.

An étale functor $\varphi : \mathcal{H} \to \mathcal{G}$ induces a simplicial map $\varphi_\bullet : \mathcal{H}_\bullet \to \mathcal{G}_\bullet$, hence a chain map $(\varphi_\bullet)_* : C_c(\mathcal{H}_\bullet, A) \to C_c(\mathcal{G}_\bullet, A)$. Similarity of functors $\rho, \sigma : \mathcal{G} \to \mathcal{H}$ produces an explicit chain homotopy $h_\bullet$ with $(\rho_n)_* - (\sigma_n)_* = \partial^{\mathcal{H}}_{n+1} h_n + h_{n-1} \partial^{\mathcal{G}}_n$. In particular, similar functors induce the same morphisms on homology. If $\mathcal{G}$ and $\mathcal{H}$ are homologically similar, then $H_n(\mathcal{G}; A) \cong H_n(\mathcal{H}; A)$ for all $n$, and for $A = \mathbb{Z}$ the induced isomorphism on $H_0$ preserves the positive cone.

To compare a groupoid with a full reduction $\mathcal{G}|_F$ one needs a continuous choice of arrows moving each unit into $F$. If $\mathcal{G}_0$ is $\sigma$-compact and totally disconnected and $F \subset \mathcal{G}_0$ is open and $\mathcal{G}$-full, a countable refinement by compact open bisections yields a continuous map $\theta : \mathcal{G}_0 \to \mathcal{G}_1$ with $r(\theta(x)) = x$ and $s(\theta(x)) \in F$. From $\theta$ one constructs étale functors $\rho : \mathcal{G} \to \mathcal{G}|_F$ and $\sigma : \mathcal{G}|_F \hookrightarrow \mathcal{G}$ such that $\sigma \circ \rho$ is similar to $\mathrm{id}_{\mathcal{G}}$ and $\rho \circ \sigma$ is similar to $\mathrm{id}_{\mathcal{G}_F}$. Hence $\mathcal{G}$ and $\mathcal{G}|_F$ are homologically similar.

Finally, if $\mathcal{G}$ and $\mathcal{H}$ are Kakutani equivalent, there exist full clopen subsets $F_{\mathcal{G}} \subset \mathcal{G}_0$ and $F_{\mathcal{H}} \subset \mathcal{H}_0$ and an isomorphism $\mathcal{G}|_{F_{\mathcal{G}}} \cong \mathcal{H}|_{F_{\mathcal{H}}}$. Combining invariance under full reductions with invariance under isomorphism gives natural isomorphisms $H_n(\mathcal{G}; A) \cong H_n(\mathcal{H}; A)$ for all $n \geq 0$, and for $A = \mathbb{Z}$ an identification of ordered groups $(H_0(\mathcal{G}), H_0(\mathcal{G})^+) \cong (H_0(\mathcal{H}), H_0(\mathcal{H})^+)$. This allows us to replace $\mathcal{G}$ by convenient full clopen reductions without changing Moore homology, and to compare different groupoid models of the same orbit structure.

# 3 Computing Moore Homology

In the preceding chapters we defined Moore homology $H_\bullet(\mathcal{G}; A)$ for étale groupoids from the simplicial geometry of the nerve $\mathcal{G}_\bullet$. For each $n \geq 1$, every face map $d_i : \mathcal{G}_n \to \mathcal{G}_{n-1}$ is a local homeomorphism, hence admits pushforward on compactly supported functions by a finite fibre sum over each point. The alternating sum of these pushforwards defines the Moore boundary on $C_c(\mathcal{G}_n, A)$. We also established the functorial and invariance mechanisms that make $H_\bullet(\mathcal{G}; A)$ workable in the ample setting, in particular invariance under full reductions and under Kakutani equivalence.

The purpose of this chapter is to turn this structural framework into a computational toolkit. In the ample world, the natural cut and paste operations are reductions along open subsets of the unit space, typically clopen and saturated, and gluings along such pieces. We therefore do not pursue excision in full generality. Instead we develop exact sequence technology tailored to compact support and to the orbit geometric operations used throughout the thesis. The basic mechanism is simple: for an open inclusion, extension by zero and restriction produce degreewise short exact sequences of compactly supported chain groups, and the étale pushforward formulas ensure compatibility with the Moore differentials. This is exactly what is needed for reductions, clopen saturated decompositions, and Mayer–Vietoris arguments in the totally disconnected setting.

The chapter develops three complementary tools:

1. **Long exact Moore homology sequence.** Section 3.1 constructs the long exact sequence associated to inclusions of subgroupoids and to reductions. In Section 3.1.1 we prove the chain level short exact sequence using extension by zero and restriction, verify compatibility with the Moore boundary by explicit pushforward computations, and describe the connecting morphism at chain level in the sense of Matui's theory [13]. In Section 3.1.2 we recast the same mechanism in a quotient language by chains supported on complements.
2. **Universal coefficient sequences.** Section 3.2 addresses coefficient changes. When the Moore complex $C_c(\mathcal{G}_\bullet, \mathbb{Z})$ is degreewise free abelian, the classical universal coefficient theorem for chain complexes yields short exact sequences for homology and cohomology. In Section 3.2.1 we obtain the UCT for $H_n(\mathcal{G}; A)$ in terms of $H_n(\mathcal{G}) = H_n(\mathcal{G}; \mathbb{Z})$, $\otimes_{\mathbb{Z}}$, and $\mathrm{Tor}_1^{\mathbb{Z}}$. In Section 3.2.2 we obtain the dual statement for the cohomology of $\mathrm{Hom}_{\mathbb{Z}}(C_c(\mathcal{G}_\bullet, \mathbb{Z}), A)$, involving $\mathrm{Hom}_{\mathbb{Z}}$ and $\mathrm{Ext}^1_{\mathbb{Z}}$. The substantive point here is structural: the Moore complex is built from compactly supported functions on an ample nerve, so freeness must be verified from compact open partitions, and the restriction to discrete coefficients is essential in this compact support model since the tensor comparison map can fail to be surjective for non-discrete $A$.

3. **Moore–Mayer–Vietoris.** Section 3.3 develops a Mayer–Vietoris principle adapted to compactly supported Moore chains. In Definition 3.3.1 we isolate admissible covers by clopen saturated pieces for which compact support and total disconnectedness give a clean decomposition of chains. In Section 3.3.2 we build the Moore–Mayer–Vietoris sequence at chain level and verify exactness by explicit control of supports and compatibility with pushforward along the face maps. Passing to homology yields the Moore–Mayer–Vietoris long exact sequence in Theorem 3.3.10, which is the main gluing tool for computations from clopen saturated decompositions.

Taken together, these results provide a practical calculus for $H_\bullet(\mathcal{G};A)$. One first replaces $\mathcal{G}$ by a convenient full clopen model without changing homology, then decomposes the unit space into admissible pieces, applies the long exact and Moore–Mayer–Vietoris sequences to relate the pieces, and finally uses the universal coefficient sequences to change coefficients whenever the freeness hypothesis holds.

**Setting 3.0.1.** Throughout this chapter, unless explicitly stated otherwise, $\mathcal{G}$ denotes a second countable locally compact Hausdorff ample étale groupoid. Its unit space $\mathcal{G}_0$ is totally disconnected. Coefficient groups $A$ are discrete abelian.

## 3.1 Long Exact Moore Homology Sequence

We study Moore homology for pairs of étale groupoids with totally disconnected unit spaces. The goal is to compare their homology groups via exact sequences coming from open inclusions and from proper quotient maps.

**Definition 3.1.1.** Let $\mathcal{G}$ be an étale groupoid and let $\mathcal{G}' \subseteq \mathcal{G}$ be an open subgroupoid. Set $\triangle := \mathcal{G} \setminus \mathcal{G}'$ and view $r(\triangle) \subseteq \mathcal{G}_0$ with the subspace topology.

The inclusion $\mathcal{G}' \subseteq \mathcal{G}$ is called regular if the restriction $r|_\triangle : \triangle \to r(\triangle)$ is an open map. In particular, $r|_\triangle$ is a quotient map onto its image.

The pair $(\mathcal{G}', \mathcal{G})$ is called a regular pair if $\mathcal{G}' \subseteq \mathcal{G}$ is a regular open subgroupoid.

Regularity means that applying the range map does not create artificial boundary phenomena when one passes from $\triangle$ to the subset $r(\triangle)$ of units. In practice, it ensures that openness of subsets of $\triangle$ is detected after applying $r$, so supports of compactly supported functions project to open subsets of units in a controlled way. This is exactly what one needs for extension by zero and restriction to interact cleanly with the pushforwards along the face maps.

We work in two complementary settings.

- Subgroupoid case: $\mathcal{G}' \subseteq \mathcal{G}$ is an open subgroupoid and $(\mathcal{G}', \mathcal{G})$ is a regular pair.

- Quotient case: $\mathcal{G}'$ is a quotient of $\mathcal{G}$ with a proper[1] and regular[2] factor map $\pi : \mathcal{G} \to \mathcal{G}'$.

### 3.1.1 Subgroupoid Case

In this section we recall the subgroupoid framework of Matui [14, Setting 3.1], which is the starting point for the long exact sequence in the subgroupoid case. Throughout, $\mathcal{G}$ denotes an ample groupoid with arrow space $\mathcal{G}_1$ and unit space $\mathcal{G}_0$.

**Setting 3.1.2.** Let $\mathcal{G}$ be an étale groupoid and let $\mathcal{G}' \subseteq \mathcal{G}$ be an open subgroupoid with the same unit space, $(\mathcal{G}')_0 = \mathcal{G}_0$. Set $\triangle := \mathcal{G} \setminus \mathcal{G}'$ for the complement of $\mathcal{G}'$ inside $\mathcal{G}$ on the level of arrows. Thus $\triangle$ consists exactly of those arrows of $\mathcal{G}$ that do not belong to $\mathcal{G}'$. In particular, $\triangle \cap \mathcal{G}_0 = \emptyset$, since all units lie in $\mathcal{G}'$. Since $\mathcal{G}'$ is a subgroupoid, it is closed under inversion, hence $\triangle$ is also closed under inversion: if $\gamma \in \triangle$ and $\gamma^{-1} \in \mathcal{G}'$, then $\gamma = (\gamma^{-1})^{-1} \in \mathcal{G}'$, a contradiction. Consequently, for every $\gamma \in \triangle$ both endpoints lie in $r(\triangle)$: $r(\gamma) \in r(\triangle), s(\gamma) = r(\gamma^{-1}) \in r(\triangle)$. In particular, the subset $r(\triangle) \subseteq \mathcal{G}_0$ is invariant for $\triangle$ in the sense that every arrow in $\triangle$ both starts and ends in $r(\triangle)$. We now restrict $\mathcal{G}$ and $\mathcal{G}'$ to the unit space $r(\triangle)$. Define $\mathcal{H} := \mathcal{G}|_{r(\triangle)}, \mathcal{H}' := \mathcal{G}'|_{r(\triangle)}$. Thus $\mathcal{H}$ has unit space $\mathcal{H}_0 = r(\triangle)$ and arrow space $\mathcal{H}_1 = \{\gamma \in \mathcal{G}_1 \mid r(\gamma), s(\gamma) \in r(\triangle)\}$, while $\mathcal{H}'$ has the same unit space $\mathcal{H}'_0 = r(\triangle)$ and arrow space $\mathcal{H}'_1 = \{\gamma \in \mathcal{G}'_1 \mid r(\gamma), s(\gamma) \in r(\triangle)\} = \mathcal{H}_1 \cap \mathcal{G}'_1$.

In particular, $\mathcal{H}'$ is a subgroupoid of $\mathcal{H}$ with the same unit space. Since every $\gamma \in \triangle$ satisfies $r(\gamma), s(\gamma) \in r(\triangle)$, we have $\triangle \subseteq \mathcal{H}_1$. Moreover, for $h \in \mathcal{H}_1$ either $h \in \mathcal{G}'_1$ or $h \in \triangle$. If $h \in \mathcal{G}'_1$, then $h \in \mathcal{H}'_1$ by definition. Hence, on the level of arrows,

$$\mathcal{H}_1 = \mathcal{H}'_1 \sqcup \triangle, \qquad \mathcal{H}_0 = \mathcal{H}'_0 = r(\triangle).$$

This decomposition is the input for a short exact sequence of chain groups built from extension by zero on $\mathcal{H}'_1 \subseteq \mathcal{H}_1$ and restriction to the complement $\triangle$. To make this compatible with the étale pushforwards in the Moore boundary, one needs a mild regularity condition ensuring that openness in $\triangle$ is reflected on the unit space through the range map.

**Lemma 3.1.3** ([14, Remark 3.3])**.** Assume that $\triangle$ is first countable, for example if $\mathcal{G}_1$ is second countable. Then the following conditions are equivalent:

1. $r|_\triangle : \triangle \to r(\triangle)$ is open when $r(\triangle)$ is equipped with the final topology induced by $r|_\triangle$.
2. For every open subset $U \subseteq \triangle$, the set $r^{-1}(r(U)) \cap \triangle$ is open in $\triangle$.
3. For every open subset $U \subseteq \triangle$ and every sequence $(\gamma_k)_k$ in $\triangle$ converging to some $\gamma \in \triangle$ with $r(\gamma) \in r(U)$, one has $r(\gamma_k) \in r(U)$ for all sufficiently large $k$.

*Proof.* The point is that regularity can be tested on the arrow space $\triangle$ by asking whether the range of an open set in $\triangle$ is again open in the quotient topology on $r(\triangle)$.

[1] A continuous map is proper if the preimage of every compact set is compact.

[2] Here regular means that the topology on $\mathcal{G}'$ is the quotient topology for $\pi$, and the restriction of $\pi$ to the relevant complement behaves as an open map onto its image, in the same sense as in Definition 3.1.1 after replacing $r|_\triangle$ by the corresponding restriction of $\pi$.

- 1. $\Rightarrow$ 2. Let $U \subseteq \triangle$ be open. If $r|_\triangle$ is open, then $r(U)$ is open in $r(\triangle)$ for the final topology. By definition of the final topology, $(r|_\triangle)^{-1}(r(U)) = r^{-1}(r(U)) \cap \triangle$ is open in $\triangle$.
- 2. $\Rightarrow$ 1. Let $U \subseteq \triangle$ be open. Condition 2. says that $(r|_\triangle)^{-1}(r(U))$ is open in $\triangle$. By definition of the final topology, this implies that $r(U)$ is open in $r(\triangle)$. Hence $r|_\triangle$ is open, which is 1.
- 2. $\Rightarrow$ 3. Let $U \subseteq \triangle$ be open and $(\gamma_k)_k \to g$ in $\triangle$ with $r(\gamma) \in r(U)$. Then $\gamma \in r^{-1}(r(U)) \cap \triangle$, which is open by 2. Hence $\gamma_k \in r^{-1}(r(U)) \cap \triangle$ for large $k$, so $r(\gamma_k) \in r(U)$ eventually.
- 3. $\Rightarrow$ 2. Fix an open $U \subseteq \triangle$ and set $A := r^{-1}(r(U)) \cap \triangle$. Let $g \in A$ and let $(\gamma_k)_k$ be a sequence in $\triangle$ with $\gamma_k \to g$. Since $r(\gamma) \in r(U)$, condition 3. implies $r(\gamma_k) \in r(U)$ for all sufficiently large $k$, hence $\gamma_k \in A$ eventually. Thus $A$ is sequentially open in $\triangle$. Since $\triangle$ is first countable, sequential openness agrees with openness, so $A$ is open in $\triangle$.

□

**Lemma 3.1.4** ([14, Lemma 3.5])**.** Assume that $\mathcal{G}' \subset \mathcal{G}$ is regular. Then the multiplication map

$$m\colon (\triangle \times \triangle) \cap m^{-1}(\mathcal{G}') \to \mathcal{H}_1', \qquad (g,h) \mapsto gh,$$

is open when $\mathcal{H}_1'$ is equipped with the final topology induced by $m$.

*Proof.* Let $U_1, U_2 \subseteq \mathcal{G}_1$ be open and set $V := (U_1 \times U_2) \cap (\triangle \times \triangle) \cap m^{-1}(\mathcal{G}')$. We show that $m(V)$ is open in $\mathcal{H}_1'$ for the final topology. Since $\triangle \times \triangle$ is first countable, it suffices to prove that $m(V)$ is sequentially open in $\mathcal{H}_1'$ in the sense of Lemma 3.1.3.

Let $(\gamma_k, h_k)$ be a sequence in $(\triangle \times \triangle) \cap m^{-1}(\mathcal{G}')$ converging to $(g,h)$, and assume that $gh \in m(V)$. Choose $(a,b) \in V$ with $ab = gh$ and set $c := bh^{-1} = a^{-1}g$. Then $c \in \mathcal{G}_1$ and the products $ch = b$ and $ac = g$ are defined. Moreover $b \in \triangle$, hence $c = (bh^{-1}) \in \triangle$ because $\triangle$ is closed under inversion and multiplication by units.

Since $\mathcal{G}$ is étale, the map $r \times s\colon \mathcal{G}_1 \to \mathcal{G}_0 \times \mathcal{G}_0$ is a local homeomorphism on a bisection neighbourhood of $c$. Shrinking around $c$ if necessary, we can choose a sequence $(c_k)_k$ in $\mathcal{G}_1$ such that $c_k \to c, r(c_k) = r(h_k), s(c_k) = r(h)$ for all $k$ large. The indices split into two sets

$$A := \{k \mid c_k \in \mathcal{G}_1'\}, \qquad B := \{k \mid c_k \in \triangle\}.$$

Passing to a subsequence, we may assume either $c_k \in \mathcal{G}_1'$ for all $k$ or $c_k \in \triangle$ for all $k$.

- **Case $c_k \in \mathcal{G}_1'$ for all $k$.** Since $\triangle$ is open in $\mathcal{G}_1$ and $ch = b \in \triangle$, we have $c_k h_k \in \triangle$ for all $k$ sufficiently large. Similarly, $\gamma_k c_k^{-1} \to gc^{-1} = a \in \triangle$ implies $\gamma_k c_k^{-1} \in \triangle$ for large $k$. Moreover,

$$(\gamma_k c_k^{-1})(c_k h_k) = \gamma_k h_k, \qquad c_k h_k \in U_2, \qquad \gamma_k c_k^{-1} \in U_1$$

  for all $k$ sufficiently large because $U_1, U_2$ are open and $a \in U_1$, $b \in U_2$. Hence $(\gamma_k c_k^{-1}, c_k h_k) \in V$ for large $k$, and therefore $\gamma_k h_k \in m(V)$ eventually.
- **Case $c_k \in \triangle$ for all $k$.** By [18, Lemma 6.4], the set $(\triangle \times \triangle) \cap m^{-1}(\triangle)$ is open in $\triangle \times \triangle$. Since $(c_k, h_k) \to (c,h)$ and $ch = b \in \triangle$, we obtain $c_k h_k \in \triangle$ for all $k$ sufficiently large. As above, $\gamma_k c_k^{-1} \to a \in \triangle$ gives $\gamma_k c_k^{-1} \in \triangle$ eventually, and openness of $U_1, U_2$ yields $\gamma_k c_k^{-1} \in U_1$, $c_k h_k \in U_2$ for large $k$. Thus $(\gamma_k c_k^{-1}, c_k h_k) \in V$ eventually, hence $\gamma_k h_k \in m(V)$ eventually.

In either case, whenever $(\gamma_k, h_k) \to (g, h)$ with $gh \in m(V)$, we have $\gamma_k h_k \in m(V)$ for all sufficiently large $k$. Thus $m(V)$ is sequentially open in $\mathcal{H}_1'$ for the final topology induced by $m$, hence open. $\square$

Now we introduce topologies on $\mathcal{H}$ and $\mathcal{H}'$.

**Definition 3.1.5** ([18, Definition 6.5], [14, Definition 3.6]). Consider the surjection

$$m\colon (\triangle \times \triangle) \cap m^{-1}(\mathcal{G}') \to \mathcal{H}_1', \qquad (g, h) \mapsto gh,$$

given by groupoid multiplication. We equip $\mathcal{H}'$ with the final topology induced by $m$. We equip $\mathcal{H} = \mathcal{H}' \cup \triangle$ with the disjoint union topology, meaning $\mathcal{H}'$ carries the final topology, $\triangle$ carries the subspace topology from $\mathcal{G}$, and both $\mathcal{H}'$ and $\triangle$ are clopen in $\mathcal{H}$.

Regularity is designed so that these quotient-type topologies still interact well with the étale structure maps.

**Theorem 3.1.6** ([14, Theorem 3.7], [18, Theorem 6.8]). Let $\mathcal{G}' \subset \mathcal{G}$ and $\mathcal{H}' \subset \mathcal{H}$ be as above, and suppose that $\mathcal{G}' \subset \mathcal{G}$ is regular. Then:

1. $\mathcal{H}$ and $\mathcal{H}'$ are étale groupoids,
2. $\mathcal{H}$ and $\mathcal{H}'$ are totally disconnected,
3. $\mathcal{H}'$ is a clopen subgroupoid of $\mathcal{H}$.

*Proof.*

1. Set

$$W := ((\triangle \times \triangle) \cap \mathcal{G}_2) \cap m^{-1}(\mathcal{G}') \subseteq \mathcal{G}_2 \subseteq \mathcal{G}_1 \times \mathcal{G}_1.$$

Then $W$ is locally compact Hausdorff as an open subset of the closed subspace $\mathcal{G}_2$. Equip $\mathcal{H}_1'$ with the final topology for the surjection $m : W \to \mathcal{H}_1'$. We first show that $\mathcal{H}_1'$ is Hausdorff. Define

$$\Phi\colon W \times W \to \mathcal{G}_1 \times \mathcal{G}_1, \qquad ((g, h), (g', h')) \mapsto (gh, g'h').$$

Let $\triangle_{\mathcal{G}_1} := \{(x, x) \mid x \in \mathcal{G}_1\}$, which is closed since $\mathcal{G}_1$ is Hausdorff. The equivalence relation of the quotient map $m$ is

$$R := \{(w, w') \in W \times W \mid m(w) = m(w')\} = \Phi^{-1}(\triangle_{\mathcal{G}_1}),$$

hence $R$ is closed in $W \times W$. Therefore the quotient $\mathcal{H}_1' = W/R$ is Hausdorff.

We next show that $\mathcal{H}_1'$ is locally compact. Fix $x \in \mathcal{H}_1'$ and choose $w \in W$ with $m(w) = x$. Since $W$ is locally compact, there exists an open neighbourhood $U \subseteq W$ of $w$ with compact closure $\overline{U} \subseteq W$. By Lemma 3.1.4, the map $m$ is open, hence $m(U)$ is open in $\mathcal{H}_1'$. The set $m(\overline{U})$ is compact, hence closed in the Hausdorff space $\mathcal{H}_1'$, and it contains $m(U)$, so $\overline{m(U)} \subseteq m(\overline{U})$. Thus $x$ has an open neighbourhood with compact closure, and $\mathcal{H}_1'$ is locally compact.

It remains to see that $\mathcal{H}$ and $\mathcal{H}'$ are étale. Since the topologies on $\triangle \subseteq \mathcal{H}$ and on $\mathcal{G}_1' \cap \mathcal{H}_1' \subseteq \mathcal{H}'$ restrict to the given subspace topologies from $\mathcal{G}_1$, it suffices to show that the range map is open on each of the two clopen pieces $\triangle$ and $\mathcal{H}_1'$. Let $U \subseteq \mathcal{G}_1$ be an open bisection. Then $U \cap \triangle$ is open in $\triangle$ and

$$r(U \cap \triangle) = (U \cap \triangle)(U \cap \triangle)^{-1} = m\Big(((U \cap \triangle) \times (U \cap \triangle)) \cap W\Big),$$

so $r(U \cap \triangle)$ is open in $\mathcal{H}_1'$ because $m$ is open. Since $\mathcal{H}_0 = r(\triangle) \subseteq \mathcal{H}_1'$, this shows $r|_\triangle \colon \triangle \to \mathcal{H}_0$ is open. For $\mathcal{H}'$, consider the commutative identity $r \circ m = r \circ p$ on $W$, where $p \colon W \to \triangle$ is the restriction of the projection $(g, h) \mapsto g$. The map $p$ is open, and we have just shown that $r|_\triangle$ is open, hence $r \circ p$ is open. Since $m$ is a quotient map by definition of the final topology, it follows that $r|_{\mathcal{H}_1'}$ is open. The same argument applies to $s$, using $s(gh) = s(h)$ and the projection $(g, h) \mapsto h$. Therefore $r$ and $s$ are local homeomorphisms on $\mathcal{H}$ and on $\mathcal{H}'$, so both are étale groupoids.

2. The topologies on $\mathcal{H}$ and $\mathcal{H}'$ are finer than the subspace topology from $\mathcal{G}_1$ on each clopen piece, so any subset that is connected in $\mathcal{H}$ or $\mathcal{H}'$ is connected in $\mathcal{G}_1$. Since $\mathcal{G}_1$ is totally disconnected, every connected subset is a singleton. Hence $\mathcal{H}, \mathcal{H}'$ are totally disconnected.
3. By Definition 3.1.5, $\mathcal{H} = \mathcal{H}' \sqcup \triangle$ is a disjoint union of clopen subsets. The set $\mathcal{H}'$ is closed under multiplication, inversion, and units by construction, hence is a subgroupoid of $\mathcal{H}$, and it is clopen by definition of the disjoint union topology.

□

Next, we show that $\mathcal{G}_n \setminus \mathcal{G}_n'$ and $\mathcal{H}_n \setminus \mathcal{H}_n'$ are canonically homeomorphic for all $n \geq 0$. This identifies the relative chain complexes

$$C_c(\mathcal{G}_n \setminus \mathcal{G}_n', A) \quad \text{and} \quad C_c(\mathcal{H}_n \setminus \mathcal{H}_n', A),$$

and hence shows that the inclusions $\mathcal{G}' \subseteq \mathcal{G}$ and $\mathcal{H}' \subseteq \mathcal{H}$ have the same relative homology. For this we need the following technical lemma.

**Lemma 3.1.7** ([14, Lemma 3.8])**.** Let $\mathcal{G}' \subseteq \mathcal{G}$ and $\mathcal{H}' \subseteq \mathcal{H}$ be as in Setting 3.1.2, and suppose that $\mathcal{G}' \subseteq \mathcal{G}$ is regular. Let $(\gamma_k)_k$ be a sequence in $\mathcal{H}$ and let $\gamma \in \mathcal{H}$. Then the following conditions are equivalent:

1. $\gamma_k \to g$ in $\mathcal{H}$,
2. $r(\gamma_k) \to r(\gamma)$ in $\mathcal{H}_0$, $s(\gamma_k) \to s(\gamma)$ in $\mathcal{H}_0$, and $\gamma_k \to g$ in $\mathcal{G}$,
3. $r(\gamma_k) \to r(\gamma)$ in $\mathcal{H}_0$ and $\gamma_k \to g$ in $\mathcal{G}$,
4. $s(\gamma_k) \to s(\gamma)$ in $\mathcal{H}_0$ and $\gamma_k \to g$ in $\mathcal{G}$.

*Proof.* The implications 1.⇒2., 2.⇒3., and 2.⇒4. are immediate. It remains to show 3.⇒1.; the implication 4.⇒1. is analogous by symmetry of $r$ and $s$.

Assume $r(\gamma_k) \to r(\gamma)$ in $\mathcal{H}_0$ and $\gamma_k \to g$ in $\mathcal{G}$.

- **Case** $\gamma \in \triangle$**.** Let $U \subseteq \mathcal{G}$ be an open bisection with $\gamma \in U$. Then $\gamma_k \in U$ for all large $k$. By Theorem 3.1.6, the map $r|_{\triangle} \colon \triangle \to \mathcal{H}_0$ is open, hence $r(U \cap \triangle)$ is an open neighbourhood of $r(\gamma)$ in $\mathcal{H}_0$. Thus $r(\gamma_k) \in r(U \cap \triangle)$ for all large $k$. For such $k$, there exists a unique $h_k \in U \cap \triangle$ with $r(h_k) = r(\gamma_k)$, because $r|_U \colon U \to r(U)$ is a homeomorphism and $r(U \cap \triangle) = r(U) \cap r(\triangle)$. Since also $\gamma_k \in U$ and $r(\gamma_k) = r(h_k)$, injectivity of $r|_U$ gives $\gamma_k = h_k \in U \cap \triangle$ for all large $k$. The $\mathcal{H}$-topology on $\triangle$ is the subspace topology from $\mathcal{G}$, and sets of the form $U \cap \triangle$ form a neighbourhood basis at $g$. Hence $\gamma_k \to g$ in $\mathcal{H}$.
- **Case** $\gamma \in \mathcal{H}$**.** By surjectivity of $m : W \to \mathcal{H}'$ from Definition 3.1.5, choose $a, b \in \triangle$ with $g = ab$. Let $U \subseteq \mathcal{G}$ be an open bisection with $a \in U$. As in the first case, $r(U \cap \triangle)$ is an open neighbourhood of $r(a) = r(\gamma)$ in $\mathcal{H}_0$, so $r(\gamma_k) \in r(U \cap \triangle)$ for all large $k$. For each such $k$, let $a_k \in U \cap \triangle$ be the unique element with $r(a_k) = r(\gamma_k)$. Then $a_k \to a$ in $\mathcal{H}$, because $r|_{U\cap\triangle} \colon U \cap \triangle \to r(U \cap \triangle)$ is a homeomorphism and $r(a_k) = r(\gamma_k) \to r(\gamma) = r(a)$ in $\mathcal{H}_0$. On $\triangle$ the $\mathcal{H}$-topology is the subspace topology from $\mathcal{G}$, hence $a_k \to a$ in $\mathcal{G}$.
  Define $b_k := a_k^{-1}\gamma_k$. By continuity of inversion and multiplication in $\mathcal{G}$, we have $b_k \to a^{-1}g = b$ in $\mathcal{G}$. Moreover,

  $$r(b_k) = r(a_k^{-1}\gamma_k) = r(a_k^{-1}) = s(a_k) \to s(a) = r(b) \quad \text{in } \mathcal{H}_0,$$

  where the convergence $s(a_k) \to s(a)$ holds because $s$ is continuous on the étale groupoid $\mathcal{H}$ and $a_k \to a$ in $\mathcal{H}$. Thus $(b_k)_k$ and $b \in \triangle$ satisfy condition 3. of the lemma. Applying the first case to $b_k \to b$, we obtain $b_k \to b$ in $\mathcal{H}$. Finally, multiplication in $\mathcal{H}$ is continuous by Theorem 3.1.6, hence $\gamma_k = a_k b_k \to ab = g$ in $\mathcal{H}$.

This proves 3.⇒1., hence all four conditions are equivalent. ☐

Before comparing the long exact sequences for the pairs $(\mathcal{G}, \mathcal{G}')$ and $(\mathcal{H}, \mathcal{H}')$, we need to understand how the Moore chain complexes decompose along a regular subgroupoid $\mathcal{G}' \subseteq \mathcal{G}$ and how this decomposition behaves under the equivalence from Setting 3.1.2. The next proposition identifies, for discrete coefficients, the quotient complex $C_c(\mathcal{G}_\bullet, A)/C_c(\mathcal{G}'_\bullet, A)$ with the chain complex on the complement $\mathcal{G} \setminus \mathcal{G}'$, and shows that this identification is compatible with the corresponding pair $(\mathcal{H}, \mathcal{H}')$. It is the analogue of [14, Proposition 3.10] in our notation.

**Proposition 3.1.8** ([14, Proposition 3.9])**.** *Assume Setting 3.1.2 and let $\mathcal{G}' \subseteq \mathcal{G}$ be regular. Let $A$ be a discrete abelian group. For each $n \geq 0$ write $\triangle_n := \mathcal{G}_n \setminus \mathcal{G}'_n$. Then $\mathcal{G}'_n \subseteq \mathcal{G}_n$ is clopen and there is a canonical short exact sequence*

$$0 \longrightarrow C_c(\mathcal{G}'_n, A) \xrightarrow{(\iota_n)_*} C_c(\mathcal{G}_n, A) \xrightarrow{(\rho_n)_*} C_c(\triangle_n, A) \longrightarrow 0,$$

*where $\iota_n : \mathcal{G}'_n \hookrightarrow \mathcal{G}_n$ is the inclusion, $(\iota_n)_*$ is extension by zero, and $(\rho_n)_*$ is restriction to $\triangle_n$.*

*Moreover, under Definition 3.1.5 we have equalities of sets $\mathcal{H}_n = \mathcal{H}'_n \sqcup \triangle_n$, $\mathcal{H}'_n = \mathcal{G}'_n$, and $\triangle_n$ is clopen in $\mathcal{H}_n$. Hence there is a second canonical short exact sequence*

$$0 \longrightarrow C_c(\mathcal{H}'_n, A) \xrightarrow{(\jmath_n)_*} C_c(\mathcal{H}_n, A) \xrightarrow{(\tilde{\rho}_n)_*} C_c(\triangle_n, A) \longrightarrow 0,$$

where $\jmath_n : \mathcal{H}'_n \hookrightarrow \mathcal{H}_n$ is the inclusion, $(\jmath_n)_*$ is extension by zero, and $(\tilde{\rho}_n)_*$ is restriction to $\triangle_n$.
For each $n \geq 0$, the universal property of cokernels yields a unique isomorphism

$$(\Theta_n)_* : \; C_c(\mathcal{G}_n, A)/(\iota_n)_*(C_c(\mathcal{G}'_n, A)) \xrightarrow{\cong} C_c(\mathcal{H}_n, A)/(\jmath_n)_*(C_c(\mathcal{H}'_n, A))$$

characterized by the commutative diagram with exact rows

$$\begin{array}{ccccccccc} 0 & \longrightarrow & C_c(\mathcal{G}'_n, A) & \xrightarrow{(\iota_n)_*} & C_c(\mathcal{G}_n, A) & \longrightarrow & C_c(\mathcal{G}_n, A)/(\iota_n)_*(C_c(\mathcal{G}'_n, A)) & \longrightarrow & 0 \\ & & & & & & \downarrow {\scriptstyle (\Theta_n)_*} \; {\scriptstyle \cong} & & \\ 0 & \longrightarrow & C_c(\mathcal{H}'_n, A) & \xrightarrow{(\jmath_n)_*} & C_c(\mathcal{H}_n, A) & \longrightarrow & C_c(\mathcal{H}_n, A)/(\jmath_n)_*(C_c(\mathcal{H}'_n, A)) & \longrightarrow & 0 \end{array}$$

together with the requirement that, under the canonical identifications

$$C_c(\mathcal{G}_n, A)/(\iota_n)_*(C_c(\mathcal{G}'_n, A)) \xrightarrow[\cong]{(\overline{\rho}_n)_*} C_c(\triangle_n, A),$$
$$C_c(\mathcal{H}_n, A)/(\jmath_n)_*(C_c(\mathcal{H}'_n, A)) \xrightarrow[\cong]{(\overline{\tilde{\rho}}_n)_*} C_c(\triangle_n, A),$$

the map $(\Theta_n)_*$ corresponds to $\mathrm{id}_{C_c(\triangle_n, A)}$, meaning $(\overline{\tilde{\rho}}_n)_* \circ (\Theta_n)_* = (\overline{\rho}_n)_*$. Finally, the maps $(\Theta_n)_*$ are compatible with the Moore differentials, hence $(\Theta_\bullet)_*$ is a natural chain isomorphism $(\Theta_\bullet)_* : \; C_c(\mathcal{G}_\bullet, A)/C_c(\mathcal{G}'_\bullet, A) \xrightarrow{\cong} C_c(\mathcal{H}_\bullet, A)/C_c(\mathcal{H}'_\bullet, A)$.

*Proof.* Fix $n \geq 0$.

- **Clopen decomposition for $\mathcal{G}_n$.** Since $\mathcal{G}' \subseteq \mathcal{G}$ is regular, the inclusion $\mathcal{G}'_1 \subseteq \mathcal{G}_1$ is clopen. For $n \geq 1$ we have

$$\mathcal{G}_n = \{(g_1, \dots, g_n) \in \mathcal{G}_1^n \mid s(g_i) = r(g_{i+1}) \text{ for all } i\}$$

as a subspace of $\mathcal{G}_1^n$, and similarly $\mathcal{G}'_n$ is the intersection of $\mathcal{G}_n$ with $(\mathcal{G}'_1)^n$. Because $(\mathcal{G}'_1)^n \subseteq \mathcal{G}_1^n$ is clopen, its intersection with $\mathcal{G}_n$ is clopen in $\mathcal{G}_n$. Thus $\mathcal{G}'_n \subseteq \mathcal{G}_n$ is clopen, so $\triangle_n = \mathcal{G}_n \setminus \mathcal{G}'_n$ is clopen as well. For $n = 0$ we have $\mathcal{G}'_0 = \mathcal{G}_0$, hence $\triangle_0 = \emptyset$.
- **Exact sequence for $\mathcal{G}_n$ and explicit maps.** Let $\iota_n : \mathcal{G}'_n \hookrightarrow \mathcal{G}_n$ be the inclusion. Since $\mathcal{G}'_n$ is open, extension by zero defines an injective homomorphism

$$(\iota_n)_* : C_c(\mathcal{G}'_n, A) \to C_c(\mathcal{G}_n, A), \qquad (\iota_n)_*(f)(g) := \begin{cases} f(g), & \text{for } g \in \mathcal{G}'_n, \\ 0, & \text{for } \gamma \in \triangle_n. \end{cases}$$

Define restriction to the clopen complement

$$(\rho_n)_* : C_c(\mathcal{G}_n, A) \to C_c(\triangle_n, A), \qquad (\rho_n)_*(f) := f|_{\triangle_n}.$$

Then $(\rho_n)_*$ is surjective since $\triangle_n$ is clopen in $\mathcal{G}_n$, so any $h \in C_c(\triangle_n, A)$ extends by zero to an element of $C_c(\mathcal{G}_n, A)$. Moreover,

$$\begin{aligned}\ker((\rho_n)_*) &= \{f \in C_c(\mathcal{G}_n, A) \mid f|_{\triangle_n} = 0\}\\ &= \{f \in C_c(\mathcal{G}_n, A) \mid \operatorname{supp}(f) \subseteq \mathcal{G}'_n\}\\ &= (\iota_n)_*(C_c(\mathcal{G}'_n, A)),\end{aligned}$$

where the last equality holds because $\mathcal{G}'_n$ is open and every $f$ supported in $\mathcal{G}'_n$ restricts to $C_c(\mathcal{G}'_n, A)$ and is recovered by extension by zero. This proves exactness of

$$0 \to C_c(\mathcal{G}'_n, A) \xrightarrow{(\iota_n)_*} C_c(\mathcal{G}_n, A) \xrightarrow{(\rho_n)_*} C_c(\triangle_n, A) \to 0.$$

Let $(\overline{\rho}_n)_*$ denote the induced isomorphism

$$(\overline{\rho}_n)_* : \ C_c(\mathcal{G}_n, A)/(\iota_n)_*(C_c(\mathcal{G}'_n, A)) \xrightarrow{\cong} C_c(\triangle_n, A), \qquad [f] \mapsto f|_{\triangle_n}.$$

- **Exact sequence for $\mathcal{H}_n$ and the same complement.** By Definition 3.1.5 we have $\mathcal{H}_n = \mathcal{H}'_n \sqcup \triangle_n$ as sets and $\triangle_n$ is clopen in $\mathcal{H}_n$. The same argument as in the previous step yields exactness of

$$0 \to C_c(\mathcal{H}'_n, A) \xrightarrow{(\jmath_n)_*} C_c(\mathcal{H}_n, A) \xrightarrow{(\tilde{\rho}_n)_*} C_c(\triangle_n, A) \to 0,$$

where $(\tilde{\rho}_n)_*$ is restriction to $\triangle_n$. Let $(\overline{\tilde{\rho}}_n)_*$ denote the induced isomorphism

$$(\overline{\tilde{\rho}}_n)_* : \ C_c(\mathcal{H}_n, A)/(\jmath_n)_*(C_c(\mathcal{H}'_n, A)) \xrightarrow{\cong} C_c(\triangle_n, A), \qquad [f] \mapsto f|_{\triangle_n}.$$

- **Definition and uniqueness of $(\Theta_n)_*$.** Define $(\Theta_n)_* := (\overline{\tilde{\rho}}_n)_*^{-1} \circ (\overline{\rho}_n)_*$. This is an isomorphism and it is the unique homomorphism such that $(\overline{\tilde{\rho}}_n)_* \circ (\Theta_n)_* = (\overline{\rho}_n)_*$. $(\Theta_n)_*$ sends the class of $f \in C_c(\mathcal{G}_n, A)$ to the unique class in the bottom quotient whose restriction to $\triangle_n$ equals $f|_{\triangle_n}$.
- **Compatibility with Moore differentials.** Let $\partial_n^{\mathcal{G}} : C_c(\mathcal{G}_n, A) \to C_c(\mathcal{G}_{n-1}, A)$ and $\partial_n^{\mathcal{H}} : C_c(\mathcal{H}_n, A) \to C_c(\mathcal{H}_{n-1}, A)$ be the Moore differentials, defined as alternating sums of pushforwards along the face maps. Since each face map restricts to the corresponding face map on the subgroupoids, the extension-by-zero maps are chain maps, meaning

$$\partial_n^{\mathcal{G}} \circ (\iota_n)_* = (\iota_{n-1})_* \circ \partial_n^{\mathcal{G}'}, \qquad \partial_n^{\mathcal{H}} \circ (\jmath_n)_* = (\jmath_{n-1})_* \circ \partial_n^{\mathcal{H}'}.$$

Hence $\partial_n^{\mathcal{G}}$ and $\partial_n^{\mathcal{H}}$ descend to differentials on the quotients, denoted

$$\begin{aligned}\overline{\partial}_n^{\mathcal{G}} &: \ C_c(\mathcal{G}_n, A)/(\iota_n)_*(C_c(\mathcal{G}'_n, A)) \to C_c(\mathcal{G}_{n-1}, A)/(\iota_{n-1})_*(C_c(\mathcal{G}'_{n-1}, A)),\\ \overline{\partial}_n^{\mathcal{H}} &: \ C_c(\mathcal{H}_n, A)/(\jmath_n)_*(C_c(\mathcal{H}'_n, A)) \to C_c(\mathcal{H}_{n-1}, A)/(\jmath_{n-1})_*(C_c(\mathcal{H}'_{n-1}, A)).\end{aligned}$$

To show that $(\Theta_\bullet)_*$ is a chain map, it suffices to show $(\Theta_{n-1})_* \circ \overline{\partial}^{\mathcal{G}}_n = \overline{\partial}^{\mathcal{H}}_n \circ (\Theta_n)_*$. Apply $(\bar{\tilde{\rho}}_{n-1})_*$ to both sides. Using $(\bar{\tilde{\rho}}_{n-1})_* \circ (\Theta_{n-1})_* = (\overline{\rho}_{n-1})_*$ and $(\bar{\tilde{\rho}}_n)_* \circ (\Theta_n)_* = (\overline{\rho}_n)_*$, this reduces to $(\overline{\rho}_{n-1})_* \circ \overline{\partial}^{\mathcal{G}}_n = (\bar{\tilde{\rho}}_{n-1})_* \circ \overline{\partial}^{\mathcal{H}}_n \circ (\bar{\tilde{\rho}}_n)^{-1}_* \circ (\overline{\rho}_n)_*$. But $(\overline{\rho}_{n-1})_* \circ \overline{\partial}^{\mathcal{G}}_n$ is, by definition of the induced quotient differential, the map $[f] \longmapsto (\partial^{\mathcal{G}}_n f)|_{\triangle_{n-1}}$, and similarly the right-hand side is $[f] \longmapsto (\partial^{\mathcal{H}}_n \tilde{f})|_{\triangle_{n-1}}$, where $\tilde{f} \in C_c(\mathcal{H}_n, A)$ is any function with $\tilde{f}|_{\triangle_n} = f|_{\triangle_n}$. Choosing $\tilde{f}$ to be the extension by zero from $\triangle_n \subset \mathcal{H}_n$, we have $\tilde{f}|_{\triangle_n} = f|_{\triangle_n}$ and $\tilde{f}|_{\mathcal{H}'_n} = 0$. Since $\mathcal{H}'_n = \mathcal{G}'_n$ and $\triangle_n$ agree as subsets of the underlying set, and since the face maps and their pushforwards agree on $\triangle_n$ under the identifications fixed in Definition 3.1.5, we get

$$(\partial^{\mathcal{G}}_n f)|_{\triangle_{n-1}} = (\partial^{\mathcal{H}}_n \tilde{f})|_{\triangle_{n-1}}.$$

□

The preceding proposition identifies, in the regular subgroupoid setting, the correct quotient complex: after passing from $\mathcal{G}$ and $\mathcal{G}'$ to the auxiliary pair $\mathcal{H} \supseteq \mathcal{H}'$, the complement $\triangle$ becomes a clopen summand in $\mathcal{H}$, so the quotient by chains supported on $\mathcal{H}'$ is literally the chain complex of compactly supported functions on the complement.

Invariance statements are proved by the same two ingredients. First, any étale functor $\varphi : \mathcal{H} \to \mathcal{G}$ induces maps on nerves $\varphi_n : \mathcal{H}_n \to \mathcal{G}_n$, hence pushforwards $(\varphi_n)_*$ on compactly supported chains, and these assemble to a chain map by functoriality of pushforward. Second, when $\varphi$ implements an equivalence of groupoids, the equivalence data provide a quasi-inverse up to similarity, and similarity yields a chain homotopy. Thus an equivalence gives a chain homotopy equivalence of Moore complexes and therefore an isomorphism on Moore homology. The next proposition records this in the discrete-coefficient, finite-fibre situation.

**Proposition 3.1.9.** Let $\mathcal{G}$ and $\mathcal{H}$ be ample groupoids and let $A$ be a discrete abelian group. Assume there is an étale functor $\varphi : \mathcal{H} \to \mathcal{G}$ which is a groupoid equivalence and such that, for every $n \geq 0$, the induced map on nerves $\varphi_n : \mathcal{H}_n \to \mathcal{G}_n$ is a local homeomorphism with finite fibres. For each $n \geq 0$ define $(\varphi_n)_* : C_c(\mathcal{H}_n, A) \to C_c(\mathcal{G}_n, A)$ by

$$((\varphi_n)_*\xi)(g) := \sum_{\substack{h \in \mathcal{H}_n \\ \varphi_n(h) = g}} \xi(h), \qquad g \in \mathcal{G}_n,\ \xi \in C_c(\mathcal{H}_n, A),$$

and write $(\varphi_\bullet)_* := ((\varphi_n)_*)_{n \geq 0}$. Then:

1. $(\varphi_\bullet)_*$ is a chain map, that is $\partial^{\mathcal{G}}_n \circ (\varphi_n)_* = (\varphi_{n-1})_* \circ \partial^{\mathcal{H}}_n$ for all $n \geq 1$.
2. The induced maps in homology $H_n((\varphi_\bullet)_*) : H_n(\mathcal{H}; A) \to H_n(\mathcal{G}; A)$ are isomorphisms for all $n \geq 0$, natural in $A$.

In particular, $\mathcal{G}$ and $\mathcal{H}$ have canonically isomorphic Moore homology with discrete coefficients.

*Proof.* Since $A$ is discrete and $\mathcal{G}, \mathcal{H}$ are ample, each $\mathcal{G}_n$ and $\mathcal{H}_n$ is locally compact, Hausdorff, and totally disconnected. For fixed $n$ and $\gamma \in \mathcal{G}_n$ the fibre $\varphi_n^{-1}(\gamma)$ is finite by assumption, hence the sum defining $(\varphi_n)_*\xi$ is finite. Let us verify that $(\varphi_n)_*$ maps $C_c(\mathcal{H}_n, A)$ to $C_c(\mathcal{G}_n, A)$.

Since $A$ is discrete, every $\xi \in C_c(\mathcal{H}_n, A)$ is locally constant. Fix $\gamma \in \mathcal{G}_n$ and choose an open neighbourhood $W \subseteq \mathcal{G}_n$ such that $\varphi_n^{-1}(W)$ is a disjoint union of open sets on which $\varphi_n$ restricts to a homeomorphism onto $W$. On each such sheet, $\xi$ is constant on a neighbourhood of every point, hence the finite sum defining $(\varphi_n)_*\xi$ is locally constant on $W$. Thus $(\varphi_n)_*\xi$ is continuous. If $K \subseteq \mathcal{H}_n$ is compact with $\mathrm{supp}(\xi) \subseteq K$, then $\mathrm{supp}((\varphi_n)_*\xi) \subseteq \varphi_n(K)$. Since $\varphi_n(K)$ is compact, this support is compact. Hence $(\varphi_n)_*\xi \in C_c(\mathcal{G}_n, A)$ and $(\varphi_n)_*$ is a well-defined homomorphism.

- **Chain map property.** For $n \geq 1$ and $0 \leq i \leq n$, functoriality of the nerve gives $d_i^{\mathcal{G}} \circ \varphi_n = \varphi_{n-1} \circ d_i^{\mathcal{H}}$. By functoriality of pushforward for composition, Proposition 2.1.2 yields $(d_i^{\mathcal{G}})_* \circ (\varphi_n)_* = (\varphi_{n-1})_* \circ (d_i^{\mathcal{H}})_*$. Summing with signs gives

$$\partial_n^{\mathcal{G}} \circ (\varphi_n)_* = \sum_{i=0}^{n} (-1)^i (d_i^{\mathcal{G}})_* \circ (\varphi_n)_* = \sum_{i=0}^{n} (-1)^i (\varphi_{n-1})_* \circ (d_i^{\mathcal{H}})_* = (\varphi_{n-1})_* \circ \partial_n^{\mathcal{H}}.$$

- **Isomorphism on homology.** Choose an étale quasi-inverse $\psi : \mathcal{G} \to \mathcal{H}$ for $\varphi$, so that $\psi_n : \mathcal{G}_n \to \mathcal{H}_n$ is a local homeomorphism with finite fibres for all $n$. Define $(\psi_n)_*$ by the same fibrewise sum and obtain a chain map $(\psi_\bullet)_*$. Since $\varphi$ and $\psi$ are quasi-inverses, there are natural transformations $\psi \circ \varphi \Rightarrow \mathrm{id}_{\mathcal{H}}, \varphi \circ \psi \Rightarrow \mathrm{id}_{\mathcal{G}}$. Applying Proposition 2.5.8 to these similarities yields chain homotopies $(\psi_\bullet)_* \circ (\varphi_\bullet)_* \simeq \mathrm{id}_{C_\bullet(\mathcal{H};A)}, (\varphi_\bullet)_* \circ (\psi_\bullet)_* \simeq \mathrm{id}_{C_\bullet(\mathcal{G};A)}$. Thus $(\varphi_\bullet)_*$ and $(\psi_\bullet)_*$ are chain homotopy inverses, hence induce mutually inverse isomorphisms on homology.
- **Naturality.** Naturality in $A$ follows because any homomorphism $A \to B$ induces levelwise postcomposition maps $C_c(-, A) \to C_c(-, B)$, and these commute with all fibrewise-sum pushforwards.

□

The regular subgroupoid setting replaces the original inclusion $\mathcal{G}' \subseteq \mathcal{G}$ by the equivalent pair $\mathcal{H}' \subseteq \mathcal{H}$ from Setting 3.1.2, where the complement $\triangle$ becomes a clopen summand and the quotient complex is identified with compactly supported chains on the complement. Proposition 3.1.8 provides short exact sequences of Moore complexes for both pairs. To use these sequences for computations, one must know that the passage from $(\mathcal{G}, \mathcal{G}')$ to $(\mathcal{H}, \mathcal{H}')$ does not alter homology, and that the resulting long exact sequences match under the equivalence. The next theorem records this compatibility at the level of long exact sequences.

**Theorem 3.1.10.** Let $\mathcal{G}' \subseteq \mathcal{G}$ and $\mathcal{H}' \subseteq \mathcal{H}$ be as in Setting 3.1.2, so that $\mathcal{G}' \subseteq \mathcal{G}$ is regular and the pair $(\mathcal{G}, \mathcal{G}')$ is equivalent to $(\mathcal{H}, \mathcal{H}')$. Let $A$ be a discrete abelian group. Then for each $n \geq 0$ there are natural isomorphisms

$$H_n(\mathcal{G}'; A) \cong H_n(\mathcal{H}'; A), \qquad H_n(\mathcal{G}; A) \cong H_n(\mathcal{H}; A), \qquad H_n(\mathcal{G}/\mathcal{G}'; A) \cong H_n(\mathcal{H}/\mathcal{H}'; A),$$

such that the following diagram has exact horizontal lines and commutes:

$$\begin{array}{ccccccccc}
\cdots \longrightarrow & H_n(\mathcal{G}';A) & \xrightarrow{H_n((j^{\mathcal{G}}_\bullet)_*)} & H_n(\mathcal{G};A) & \xrightarrow{H_n((q^{\mathcal{G}}_\bullet)_*)} & H_n(\mathcal{G}/\mathcal{G}';A) & \xrightarrow{\partial^{\mathcal{G}}_n} & H_{n-1}(\mathcal{G}';A) & \longrightarrow \cdots \\
& {\scriptstyle H_n((\Phi'_\bullet)_*)}\downarrow & & {\scriptstyle H_n((\Phi_\bullet)_*)}\downarrow & & {\scriptstyle H_n((\overline{\Phi}_\bullet)_*)}\downarrow & & {\scriptstyle H_{n-1}((\Phi'_\bullet)_*)}\downarrow & \\
\cdots \longrightarrow & H_n(\mathcal{H}';A) & \xrightarrow{H_n((j^{\mathcal{H}}_\bullet)_*)} & H_n(\mathcal{H};A) & \xrightarrow{H_n((q^{\mathcal{H}}_\bullet)_*)} & H_n(\mathcal{H}/\mathcal{H}';A) & \xrightarrow{\partial^{\mathcal{H}}_n} & H_{n-1}(\mathcal{H}';A) & \longrightarrow \cdots
\end{array}$$

*Proof.* By Setting 3.1.2 there is an étale functor $\Phi : \mathcal{G} \to \mathcal{H}$ which restricts to $\Phi' : \mathcal{G}' \to \mathcal{H}'$ and induces a functor on complements

$$\overline{\Phi} : \mathcal{G} \setminus \mathcal{G}' \to \mathcal{H} \setminus \mathcal{H}'.$$

For each $n \geq 0$ this yields local homeomorphisms

$$\Phi_n : \mathcal{G}_n \to \mathcal{H}_n, \qquad \Phi'_n : \mathcal{G}'_n \to \mathcal{H}'_n, \qquad \overline{\Phi}_n : \mathcal{G}_n \setminus \mathcal{G}'_n \to \mathcal{H}_n \setminus \mathcal{H}'_n.$$

- **Chain-level short exact sequences.** Proposition 3.1.8 gives short exact sequences

$$0 \to C_c(\mathcal{G}'_\bullet, A) \xrightarrow{(j^{\mathcal{G}}_\bullet)_*} C_c(\mathcal{G}_\bullet, A) \xrightarrow{(q^{\mathcal{G}}_\bullet)_*} C_c(\mathcal{G}_\bullet, A)/C_c(\mathcal{G}'_\bullet, A) \to 0,$$
$$0 \to C_c(\mathcal{H}'_\bullet, A) \xrightarrow{(j^{\mathcal{H}}_\bullet)_*} C_c(\mathcal{H}_\bullet, A) \xrightarrow{(q^{\mathcal{H}}_\bullet)_*} C_c(\mathcal{H}_\bullet, A)/C_c(\mathcal{H}'_\bullet, A) \to 0.$$

- **Chain-level isomorphisms and compatibility.** Since $A$ is discrete, pushforward along $\Phi_n, \Phi'_n, \overline{\Phi}_n$ defines chain maps

$$\begin{aligned}
(\Phi_\bullet)_* &: C_c(\mathcal{G}_\bullet, A) \to C_c(\mathcal{H}_\bullet, A), \\
(\Phi'_\bullet)_* &: C_c(\mathcal{G}'_\bullet, A) \to C_c(\mathcal{H}'_\bullet, A), \\
(\overline{\Phi}_\bullet)_* &: C_c(\mathcal{G}_\bullet, A)/C_c(\mathcal{G}'_\bullet, A) \to C_c(\mathcal{H}_\bullet, A)/C_c(\mathcal{H}'_\bullet, A).
\end{aligned}$$

  By Proposition 3.1.9 these are chain isomorphisms. Moreover, $\Phi$ respects subgroupoids and complements, hence the induced chain maps satisfy

$$(\Phi_\bullet)_* \circ (j^{\mathcal{G}}_\bullet)_* = (j^{\mathcal{H}}_\bullet)_* \circ (\Phi'_\bullet)_*, \qquad (q^{\mathcal{H}}_\bullet)_* \circ (\Phi_\bullet)_* = (\overline{\Phi}_\bullet)_* \circ (q^{\mathcal{G}}_\bullet)_*,$$

  so we obtain a commutative diagram of short exact sequences of chain complexes

$$\begin{array}{ccccccccc}
0 & \longrightarrow & C_c(\mathcal{G}'_\bullet, A) & \xrightarrow{j^{\mathcal{G}}_*} & C_c(\mathcal{G}_\bullet, A) & \xrightarrow{q^{\mathcal{G}}_*} & C_c(\mathcal{G}_\bullet, A)/C_c(\mathcal{G}'_\bullet, A) & \longrightarrow & 0 \\
& & {\scriptstyle (\Phi'_\bullet)_*}\downarrow & & {\scriptstyle (\Phi_\bullet)_*}\downarrow & & {\scriptstyle (\overline{\Phi}_\bullet)_*}\downarrow & & \\
0 & \longrightarrow & C_c(\mathcal{H}'_\bullet, A) & \xrightarrow{j^{\mathcal{H}}_*} & C_c(\mathcal{H}_\bullet, A) & \xrightarrow{q^{\mathcal{H}}_*} & C_c(\mathcal{H}_\bullet, A)/C_c(\mathcal{H}'_\bullet, A) & \longrightarrow & 0.
\end{array}$$

- **Passing to homology.** Applying the homology functor to the two short exact sequences yields long exact sequences with connecting morphisms $\partial^{\mathcal{G}}_n$ and $\partial^{\mathcal{H}}_n$. Naturality of the connecting morphism for morphisms of short exact sequences of chain complexes implies that the above commutative diagram induces the commutative diagram of long exact sequences in the statement. Since $(\Phi_\bullet)_*$, $(\Phi'_\bullet)_*$, $(\overline{\Phi}_\bullet)_*$ are chain isomorphisms, the induced maps $H_n((\Phi_\bullet)_*), H_n((\Phi'_\bullet)_*), H_n((\overline{\Phi}_\bullet)_*)$ are isomorphisms for all $n \geq 0$.

□

### 3.1.2 Quotient Groupoid Case

In this section we turn to the second source of long exact sequences for Moore homology, namely the quotient situation. In contrast to the subgroupoid case, the relevant maps on nerves come from quotient-type surjections and are therefore typically not local homeomorphisms. Since our chain groups are $C_c(\mathcal{G}_n, A)$ and the Moore boundary is built from pushforward along local homeomorphisms, we need a substitute for pushforward along such quotient maps. Matui's approach is to single out a class of proper surjections between locally compact metric spaces for which fibres vary in a controlled way. This control is needed to define fibrewise summation on compactly supported functions, to prove continuity of the resulting function, and to verify compatibility with composition and with the simplicial face maps.

The metric input is phrased in terms of Hausdorff variation of compact fibres. For a metric space $(X, d)$ and a compact subset $K \subseteq X$ we write $\operatorname{diam}(K) := \sup\{d(x,y) \mid x, y \in K\}$. For compact $K, L \subseteq X$ we use the Hausdorff distance

$$d_H(K,L) := \max\Big\{\sup_{x\in K}\inf_{y\in L} d(x,y),\ \sup_{y\in L}\inf_{x\in K} d(x,y)\Big\},$$

and we record the estimate $|\operatorname{diam}(K) - \operatorname{diam}(L)| \leq 2\, d_H(K,L)$. Intuitively, regularity says that near a point of the target, either the fibre varies continuously in Hausdorff distance, or it becomes uniformly small. This is the precise hypothesis that makes fibrewise constructions compatible with compact supports.

After collecting the required consequences for such surjections, we apply them degreewise to simplicial spaces arising from an étale groupoid $\mathcal{G}$ together with a suitable normal wide subgroupoid of isotropy $\mathcal{N}$, so that the factor groupoid $\mathcal{G}/\mathcal{N}$ is defined and the quotient map behaves well on compactly supported chains. This yields a chain-level short exact sequence adapted to the quotient setting and, via Matui's explicit connecting morphism, a long exact sequence in Moore homology. Together with the subgroupoid long exact sequence, this provides the basic cut-and-paste technology used later for computations.

For a set $S$ we write $\#S$ for its cardinality.

**Definition 3.1.11** (Regular proper surjections [18, Def. 7.6, Def. 7.7])**.** Let $(X, d)$ and $(X', d')$ be locally compact metric spaces, and let $\pi : X \to X'$ be a continuous proper surjection.

1. We call $\pi$ regular if for every $x' \in X'$ and every $\varepsilon > 0$ there exists an open neighbourhood $U' \subseteq X'$ of $x'$ such that for every $y' \in U'$ one has

$$d_H(\pi^{-1}(x'), \pi^{-1}(y')) < \varepsilon \quad \text{or} \quad \operatorname{diam}(\pi^{-1}(y')) < \varepsilon.$$

2. Set $X'_\pi := \{x' \in X' \mid \#\pi^{-1}(x') > 1\}$ and $X_\pi := \pi^{-1}(X'_\pi)$. Define metrics

$$d'_\pi(x', y') := d_H(\pi^{-1}(x'), \pi^{-1}(y')) \qquad \text{for } x', y' \in X'_\pi,$$
$$d_\pi(x, y) := d(x, y) + d'_\pi(\pi(x), \pi(y)) \qquad \text{for } x, y \in X_\pi.$$

**Lemma 3.1.12** ([18, Proposition 7.8])**.** Let $(X, d)$ and $(X', d')$ be locally compact metric spaces and let $\pi : X \to X'$ be a continuous proper surjection. With $X'_\pi \subseteq X'$, $X_\pi := \pi^{-1}(X'_\pi)$, and the metrics $d'_\pi$ on $X'_\pi$ and $d_\pi$ on $X_\pi$ as in Definition 3.1.11, the inclusions $(X_\pi, d_\pi) \hookrightarrow (X, d)$, $(X'_\pi, d'_\pi) \hookrightarrow (X', d')$ are continuous.

*Proof.* Write $\iota : X_\pi \hookrightarrow X$ and $\iota' : X'_\pi \hookrightarrow X'$ for the inclusions.

- **Continuity of $\iota$.** For $x, y \in X_\pi$ one has $d_\pi(x, y) = d(x, y) + d'_\pi(\pi(x), \pi(y))$, hence $d(x, y) \leq d_\pi(x, y)$. Thus $\iota$ is 1-Lipschitz and therefore continuous.
- **Continuity of $\iota'$.** Fix $x' \in X'_\pi$ and let $U \subseteq X'$ be an open neighbourhood of $x'$. Pick $x \in \pi^{-1}(x')$. By continuity of $\pi$ there exists an open neighbourhood $V \subseteq X$ of $x$ with $\pi(V) \subseteq U$. Choose $\varepsilon > 0$ such that $B_d(x, \varepsilon) \subseteq V$.
  Let $y' \in X'_\pi$ with $d'_\pi(x', y') < \varepsilon$, i.e. $d_H(\pi^{-1}(x'), \pi^{-1}(y')) < \varepsilon$. By the definition of $d_H$, there exists $y \in \pi^{-1}(y')$ with $d(x, y) < \varepsilon$. Then $y \in B_d(x, \varepsilon) \subseteq V$, hence $y' = \pi(y) \in \pi(V) \subseteq U$. Since $U$ was arbitrary, $\iota'$ is continuous at $x'$, hence continuous.
- **Well-definedness.** Properness makes each fibre $\pi^{-1}(x')$ compact, so the Hausdorff distance $d_H(\pi^{-1}(x'), \pi^{-1}(y'))$ is well defined.

□

**Lemma 3.1.13** ([18, Lemma 4.4])**.** Let $(X, d)$ and $(X', d')$ be locally compact metric spaces and let $\pi : X \to X'$ be a continuous regular proper surjection. For every compact $K \subseteq X'$ and every $\delta > 0$, the set $K_\delta := K \cap \{x' \in X'_\pi \mid \operatorname{diam} \pi^{-1}(x') \geq \delta\}$ is compact in $(X'_\pi, d'_\pi)$.

*Proof.* Since $(X'_\pi, d'_\pi)$ is metric, it suffices to prove sequential compactness. Let $(x'_k)_{k \geq 1}$ be a sequence in $K_\delta$. Compactness of $K$ in $(X', d')$ yields a subsequence, still denoted $(x'_k)_k$, with $x'_k \to x'$ in $(X', d')$ for some $x' \in K$. Fix $\varepsilon > 0$ with $\varepsilon < \delta$. By regularity of $\pi$ at $x'$ there exists an open neighbourhood $U' \subseteq X'$ of $x'$ such that for every $y' \in U'$ one has

$$d_H(\pi^{-1}(x'), \pi^{-1}(y')) < \varepsilon \quad \text{or} \quad \operatorname{diam} \pi^{-1}(y') < \varepsilon.$$

For $k$ large we have $x'_k \in U'$. Since $x'_k \in K_\delta$, we have $\operatorname{diam} \pi^{-1}(x'_k) \geq \delta > \varepsilon$, so the second alternative is impossible. Hence $d_H(\pi^{-1}(x'), \pi^{-1}(x'_k)) < \varepsilon$ for all large $k$, and therefore $d_H(\pi^{-1}(x'), \pi^{-1}(x'_k)) \to 0$.

For compact subsets $A, B \subseteq X$ one has $|\operatorname{diam}(A) - \operatorname{diam}(B)| \leq 2d_H(A, B)$. Thus, for all large $k$,

$$\operatorname{diam} \pi^{-1}(x') \geq \operatorname{diam} \pi^{-1}(x'_k) - 2d_H(\pi^{-1}(x'), \pi^{-1}(x'_k)) \geq \delta - 2\varepsilon.$$

Since this holds for every $\varepsilon \in (0, \delta)$, it follows that $\operatorname{diam} \pi^{-1}(x') \geq \delta$. In particular $\pi^{-1}(x')$ contains at least two points, so $x' \in X'_\pi$, hence $x' \in K_\delta$.

Finally, for large $k$ both $x', x'_k$ lie in $X'_\pi$, and then $d'_\pi(x', x'_k) = d_H(\pi^{-1}(x'), \pi^{-1}(x'_k)) \to 0$. Hence $x'_k \to x'$ in $(X'_\pi, d'_\pi)$. $\square$

Lemma 3.1.13 isolates the part of $X'$ where the fibres have uniformly positive diameter. This is the key compactness input for proving local compactness of $X'_\pi$ in the Hausdorff fibre metric $d'_\pi$. Once local compactness is available, the metric definition of $d_\pi$ gives continuity and properness of $\pi : X_\pi \to X'_\pi$, and regularity provides openness.

**Proposition 3.1.14** ([18, Theorem 7.9])**.** Let $(X, d)$ and $(X', d')$ be locally compact metric spaces and let $\pi : X \to X'$ be a continuous proper surjection. Assume that $\pi$ is regular. Then $(X_\pi, d_\pi)$ and $(X'_\pi, d'_\pi)$ are locally compact metric spaces, and $\pi : (X_\pi, d_\pi) \to (X'_\pi, d'_\pi)$ is continuous, proper, and open.

*Proof.*

- **$X'_\pi$ is locally compact.** Fix $x' \in X'_\pi$. Choose a compact neighbourhood $K \subseteq X'$ of $x'$ with $x' \in \mathring{K}$. Set $\delta := \operatorname{diam} \pi^{-1}(x') > 0$ and

$$K_{\delta/2} := K \cap \{y' \in X'_\pi \mid \operatorname{diam} \pi^{-1}(y') \geq \delta/2\}.$$

  By Lemma 3.1.13, $K_{\delta/2}$ is compact in $(X'_\pi, d'_\pi)$.
  We show that $x'$ is an interior point of $K_{\delta/2}$ in $(X'_\pi, d'_\pi)$. Since the inclusion $(X'_\pi, d'_\pi) \hookrightarrow (X', d')$ is continuous by Lemma 3.1.12, there exists $\varepsilon_0 > 0$ such that $d'_\pi(x', y') < \varepsilon_0$ implies $y' \in \mathring{K}$. Let $\varepsilon := \min\{\varepsilon_0, \delta/4\}$ and assume $d'_\pi(x', y') < \varepsilon$. Then $y' \in K$. Moreover,

$$\operatorname{diam} \pi^{-1}(y') \geq \operatorname{diam} \pi^{-1}(x') - 2d_H(\pi^{-1}(x'), \pi^{-1}(y')) = \delta - 2d'_\pi(x', y') \geq \delta/2,$$

  so $y' \in K_{\delta/2}$. Thus $B_{d'_\pi}(x', \varepsilon) \subseteq K_{\delta/2}$, and $X'_\pi$ is locally compact.
- **Continuity of $\pi : X_\pi \to X'_\pi$.** For $x, y \in X_\pi$,

$$d'_\pi(\pi(x), \pi(y)) \leq d(x, y) + d'_\pi(\pi(x), \pi(y)) = d_\pi(x, y),$$

  so $\pi$ is 1-Lipschitz.
- **Properness of $\pi : X_\pi \to X'_\pi$.** Let $K \subseteq X'_\pi$ be compact. By Lemma 3.1.12, $K$ is compact in $X'$. Properness of $\pi : X \to X'$ gives compactness of $\pi^{-1}(K)$ in $X$. To show compactness in $(X_\pi, d_\pi)$, let $(x_k)_k$ be a sequence in $\pi^{-1}(K)$. After passing to a subsequence, $x_k \to x$ in $(X, d)$. Since $K$ is compact in $(X'_\pi, d'_\pi)$, after passing to a subsequence we have $\pi(x_k) \to z'$ in $(X'_\pi, d'_\pi)$. By Lemma 3.1.12 we also have $\pi(x_k) \to z'$ in $(X', d')$. Continuity of $\pi : (X, d) \to (X', d')$ yields $\pi(x_k) \to \pi(x)$ in $(X', d')$, hence $z' = \pi(x)$.
  Therefore $d'_\pi(\pi(x_k), \pi(x)) \to 0$, and

$$d_\pi(x_k, x) = d(x_k, x) + d'_\pi(\pi(x_k), \pi(x)) \to 0.$$

  Thus $\pi^{-1}(K)$ is compact in $X_\pi$, so $\pi$ is proper.

- **$X_\pi$ is locally compact.** Fix $x \in X_\pi$. Choose a compact neighbourhood $K \subseteq X'_\pi$ of $\pi(x)$. Then $\pi^{-1}(K)$ is compact in $X_\pi$ by properness, and it is a neighbourhood of $x$. Hence $X_\pi$ is locally compact.
- **Openness of $\pi : X_\pi \to X'_\pi$.** Let $V \subseteq X_\pi$ be open and let $x \in V$. Choose $\varepsilon > 0$ such that $B_{d_\pi}(x, \varepsilon) \subseteq V$. Let $z' \in X'_\pi$ with $d'_\pi(\pi(x), z') < \varepsilon/2$, that is $d_H(\pi^{-1}(\pi(x)), \pi^{-1}(z')) < \varepsilon/2$. Pick $y \in \pi^{-1}(z')$ with $d(x, y) < \varepsilon/2$. Then

$$d_\pi(x, y) = d(x, y) + d'_\pi(\pi(x), \pi(y)) < \varepsilon/2 + \varepsilon/2 = \varepsilon,$$

so $y \in V$ and $z' = \pi(y) \in \pi(V)$. Hence $B_{d'_\pi}(\pi(x), \varepsilon/2) \subseteq \pi(V)$, so $\pi(V)$ is open.

□

**Proposition 3.1.15.** For $i = 1, 2$, let $(X_i, d_i)$ and $(X'_i, d'_i)$ be locally compact metric spaces and let $\pi_i : X_i \to X'_i$ be a continuous proper surjection. Let $\varphi : X_1 \to X_2$ and $\varphi' : X'_1 \to X'_2$ be surjective local homeomorphisms such that $\pi_2 \circ \varphi = \varphi' \circ \pi_1$. Assume moreover that for every $x' \in X'_1$ the restriction $\varphi : \pi_1^{-1}(x') \to \pi_2^{-1}(\varphi'(x'))$ is bijective. Then $\pi_1$ is regular if and only if $\pi_2$ is regular.

If $\pi_1$ and $\pi_2$ are regular, then the restrictions $\varphi' : (X'_1)_{\pi_1} \to (X'_2)_{\pi_2}, \varphi : (X_1)_{\pi_1} \to (X_2)_{\pi_2}$ are local homeomorphisms, where $(X'_i)_{\pi_i} = \{x' \in X'_i \mid \#\pi_i^{-1}(x') > 1\}$ and $(X_i)_{\pi_i} = \pi_i^{-1}((X'_i)_{\pi_i})$, endowed with the metrics $d'_{\pi_i}$ and $d_{\pi_i}$ from Definition 3.1.11.

*Proof.*

- **A Hausdorff estimate under uniform continuity.** Let $(Y, \rho)$, $(Z, \sigma)$ be metric spaces, let $K \subseteq Y$ be compact, and let $F : K \to Z$ be continuous. Then $F$ is uniformly continuous, so for every $\varepsilon > 0$ there exists $\delta > 0$ such that $\rho(x, y) < \delta$ implies $\sigma(F(x), F(y)) < \varepsilon$ for all $x, y \in K$. Consequently, for nonempty compact $A, B \subseteq K$,

$$d_{H,\rho}(A, B) < \delta \;\Rightarrow\; d_{H,\sigma}(F(A), F(B)) < \varepsilon, \qquad \operatorname{diam}_\rho(A) < \delta \;\Rightarrow\; \operatorname{diam}_\sigma(F(A)) < \varepsilon.$$

Indeed, if $d_{H,\rho}(A, B) < \delta$, then for every $a \in A$ there exists $b \in B$ with $\rho(a, b) < \delta$, hence $\sigma(F(a), F(b)) < \varepsilon$, and taking suprema gives $\sup_{a \in A} \inf_{b \in B} \sigma(F(a), F(b)) \leq \varepsilon$; the other half is symmetric. The diameter implication is immediate.
- **$\pi_1$ regular implies $\pi_2$ regular.** Fix $x'_2 \in X'_2$ and $\varepsilon > 0$. Choose $x'_1 \in X'_1$ with $\varphi'(x'_1) = x'_2$. Let $K'_1 \subseteq X'_1$ be a compact neighbourhood of $x'_1$ and set $K_1 := \pi_1^{-1}(K'_1)$, compact by properness of $\pi_1$. Apply the estimate above to $F = \varphi|_{K_1}$ to obtain $\delta > 0$ such that for compact $A, B \subseteq K_1$,

$$d_{H,1}(A, B) < \delta \Rightarrow d_{H,2}(\varphi(A), \varphi(B)) < \varepsilon, \qquad \operatorname{diam}_1(A) < \delta \Rightarrow \operatorname{diam}_2(\varphi(A)) < \varepsilon.$$

By regularity of $\pi_1$ at $x'_1$ with parameter $\delta$, choose an open neighbourhood $U'_1 \subseteq K'_1$ of $x'_1$ such that for every $y'_1 \in U'_1$,

$$d_{H,1}(\pi_1^{-1}(x'_1), \pi_1^{-1}(y'_1)) < \delta \quad \text{or} \quad \operatorname{diam}_1(\pi_1^{-1}(y'_1)) < \delta.$$

Set $U'_2 := \varphi'(U'_1)$, open since $\varphi'$ is open. Fix $y'_2 \in U'_2$ and choose $y'_1 \in U'_1$ with $\varphi'(y'_1) = y'_2$.

The following square commutes:

$$\begin{array}{ccc} X_1 & \xrightarrow{\varphi} & X_2 \\ {\scriptstyle \pi_1}\downarrow & & \downarrow{\scriptstyle \pi_2} \\ X_1' & \xrightarrow{\varphi'} & X_2'. \end{array}$$

The commutative square and fibrewise bijectivity give

$$\varphi(\pi_1^{-1}(x_1')) = \pi_2^{-1}(x_2'), \qquad \varphi(\pi_1^{-1}(y_1')) = \pi_2^{-1}(y_2').$$

If the first alternative holds, then $d_{H,2}(\pi_2^{-1}(x_2'), \pi_2^{-1}(y_2')) < \varepsilon$. If the second alternative holds, then $\mathrm{diam}_2(\pi_2^{-1}(y_2')) < \varepsilon$. Thus $\pi_2$ is regular at $x_2'$.

- $\pi_2$ **regular implies** $\pi_1$ **regular.** Fix $x_1' \in X_1'$ and $\varepsilon > 0$, and set $x_2' := \varphi'(x_1')$. Choose an open neighbourhood $V_1' \subseteq X_1'$ of $x_1'$ such that $\varphi'|_{V_1'} : V_1' \to V_2'$ is a homeomorphism onto an open set $V_2' \subseteq X_2'$. Shrink $V_1'$ so that $\overline{V_1'}$ is compact, and set $K_1' := \overline{V_1'}$, $K_2' := \varphi'(K_1')$, $K_2 := \pi_2^{-1}(K_2')$, compact by properness of $\pi_2$.
We claim that $\varphi : \pi_1^{-1}(V_1') \to \pi_2^{-1}(V_2')$ is a homeomorphism. It is bijective by fibrewise bijectivity and the fact that $\varphi'|_{V_1'}$ is bijective: given $z \in \pi_2^{-1}(V_2')$, put $z_2' = \pi_2(z) \in V_2'$ and $z_1' = (\varphi'|_{V_1'})^{-1}(z_2') \in V_1'$, then there is a unique $x \in \pi_1^{-1}(z_1')$ with $\varphi(x) = z$. Since $\varphi$ is a local homeomorphism, a bijective restriction is a homeomorphism.
Apply the estimate above to $F = (\varphi|_{\pi_1^{-1}(K_1')})^{-1} : K_2 \to X_1$. Thus there exists $\delta > 0$ such that for compact $A, B \subseteq K_2$,

$$d_{H,2}(A,B) < \delta \Rightarrow d_{H,1}(F(A), F(B)) < \varepsilon, \qquad \mathrm{diam}_2(A) < \delta \Rightarrow \mathrm{diam}_1(F(A)) < \varepsilon.$$

By regularity of $\pi_2$ at $x_2'$ with parameter $\delta$, choose an open neighbourhood $U_2' \subseteq V_2'$ of $x_2'$ such that for every $y_2' \in U_2'$,

$$d_{H,2}(\pi_2^{-1}(x_2'), \pi_2^{-1}(y_2')) < \delta \quad \text{or} \quad \mathrm{diam}_2(\pi_2^{-1}(y_2')) < \delta.$$

Put $U_1' := (\varphi'|_{V_1'})^{-1}(U_2')$, an open neighbourhood of $x_1'$. Fix $y_1' \in U_1'$ and set $y_2' = \varphi'(y_1')$. Using

$$\pi_2^{-1}(x_2') = \varphi(\pi_1^{-1}(x_1')), \qquad \pi_2^{-1}(y_2') = \varphi(\pi_1^{-1}(y_1')),$$

the two alternatives imply, after applying $F$,

$$d_{H,1}(\pi_1^{-1}(x_1'), \pi_1^{-1}(y_1')) < \varepsilon \quad \text{or} \quad \mathrm{diam}_1(\pi_1^{-1}(y_1')) < \varepsilon.$$

Thus $\pi_1$ is regular at $x_1'$.

- **Local homeomorphisms on the nontrivial-fibre.** Fibrewise bijectivity implies $\#\pi_1^{-1}(x') = \#\pi_2^{-1}(\varphi'(x'))$ for all $x' \in X_1'$, hence

$$\varphi'((X_1')_{\pi_1}) = (X_2')_{\pi_2}, \qquad \varphi((X_1)_{\pi_1}) = (X_2)_{\pi_2}.$$

Fix $x_1' \in (X_1')_{\pi_1}$ and choose $V_1'$ as in the previous item, with $\overline{V_1'}$ compact and $\varphi'|_{V_1'}$ a homeomorphism onto $V_2'$. Then $\varphi : \pi_1^{-1}(V_1') \to \pi_2^{-1}(V_2')$ is a homeomorphism, so both $\varphi$ and $\varphi^{-1}$ are uniformly continuous on the compact sets $\pi_1^{-1}(\overline{V_1'})$ and $\pi_2^{-1}(\overline{V_2'})$. Applying the Hausdorff estimate to $\varphi$ and $\varphi^{-1}$ shows that the restriction

$$\varphi' : V_1' \cap (X_1')_{\pi_1} \to V_2' \cap (X_2')_{\pi_2}$$

is a homeomorphism for the metrics $d'_{\pi_1}$ and $d'_{\pi_2}$. Hence $\varphi' : (X_1')_{\pi_1} \to (X_2')_{\pi_2}$ is a local homeomorphism. Finally, since $d_{\pi_i}(x,y) = d_i(x,y) + d'_{\pi_i}(\pi_i(x), \pi_i(y))$ and $\pi_2 \circ \varphi = \varphi' \circ \pi_1$, the same neighbourhoods yield that $\varphi : (X_1)_{\pi_1} \to (X_2)_{\pi_2}$ is a local homeomorphism.

□

We now introduce the notion of reduction maps. Throughout this subsection we fix the following standing assumptions. Let $(X,d)$ and $(X',d')$ be locally compact metric spaces and let $\pi : X \to X'$ be a continuous proper surjection. Assume that $\pi$ is regular and fix a discrete abelian group $A$. By Proposition 3.1.14, the metric spaces $(X_\pi, d_\pi)$ and $(X'_\pi, d'_\pi)$ are locally compact, and the restricted map $\pi : X_\pi \to X'_\pi$ is continuous, proper, and open. If $X$ and $X'$ are totally disconnected, then $X_\pi \subseteq X$ and $X'_\pi \subseteq X'$ are totally disconnected as well.

**Lemma 3.1.16** ([14, Lemma 4.8])**.** Let $(X,d)$ and $(X',d')$ be locally compact metric spaces which are totally disconnected, let $A$ be a discrete abelian group, and let $\pi : X \to X'$ be a continuous proper surjection which is regular in the sense of Definition 3.1.11. Write $X_\pi \subseteq X$ and $X'_\pi \subseteq X'$ as in Definition 3.1.11, and regard $\pi$ as a map $X_\pi \to X'_\pi$. Then for every $f \in C_c(X,A)$ there exist a compact open subset $L' \subseteq X'_\pi$ and a continuous map $g : X'_\pi \to A$ such that $f(x) = g(\pi(x))$ for all $x \in X_\pi \setminus \pi^{-1}(L')$.

*Proof.* Let $K := \operatorname{supp}(f)$, which is compact. Since $A$ is discrete, $f$ is locally constant, hence $K = \{x \in X \mid f(x) \neq 0\}$ and $K$ is compact open in $X$. Choose $\delta_1 > 0$ such that $d(x,y) < \delta_1$ implies $f(x) = f(y)$ for all $x,y \in K$. Since $K$ is compact and $X \setminus K$ is closed and disjoint from $K$, the number $\delta_0 := \inf\{d(x,y) \mid x \in K,\ y \in X \setminus K\}$ is strictly positive. Set $\delta := \min\{\delta_0, \delta_1\}$. Then $d(x,y) < \delta$ implies $f(x) = f(y)$ for all $x,y \in X$.

Define $K' := \pi(K) \cap \{x' \in X'_\pi \mid \operatorname{diam}(\pi^{-1}(x')) \geq \delta\}$. Since $\pi$ is proper and $K$ is compact, $\pi(K)$ is compact in $X'$, and Lemma 3.1.13 shows that $K'$ is compact in $(X'_\pi, d'_\pi)$. Because $X'_\pi$ is locally compact and totally disconnected, it has a basis of compact open sets, so choose a compact open neighbourhood $L' \subseteq X'_\pi$ of $K'$. Since $X'_\pi$ is Hausdorff, $L'$ is clopen.

Put $B := X'_\pi \setminus L'$ and $E := \pi^{-1}(B) \subseteq X_\pi$. We claim that $f|_E$ is constant on each fibre of $\pi|_E : E \to B$. Fix $x' \in B$ and $x,y \in \pi^{-1}(x')$. If $x' \notin \pi(K)$, then $x,y \notin K$, hence $f(x) = f(y) = 0$. If $x' \in \pi(K)$, then $x' \notin K'$ because $K' \subseteq L'$, hence $\operatorname{diam}(\pi^{-1}(x')) < \delta$, so $d(x,y) < \delta$ and therefore $f(x) = f(y)$. By Proposition 3.1.14, the map $\pi : X_\pi \to X'_\pi$ is open. Hence $\pi|_E : E \to B$ is an open surjection and therefore a quotient map. Since $f|_E$ is constant on fibres, there exists a unique continuous map $g_B : B \to A$ such that $f|_E = g_B \circ (\pi|_E)$. Define $g : X'_\pi \to A$ by $g|_B := g_B$ and $g|_{L'} := 0$. Because $B$ and $L'$ are clopen and $A$ is discrete, $g$ is continuous. For $x \in X_\pi \setminus \pi^{-1}(L') = E$ we then have $f(x) = g(\pi(x))$, as required. □

The goal of this subsection is to isolate, from a compactly supported function $f \in C_c(X, A)$, the part which is genuinely transverse to a regular quotient map $\pi : X \to X'$. Regularity controls the variation of fibres in the Hausdorff metric and forces the region of large fibres to live in a compact open subset of $X'_\pi$. As a consequence, outside the preimage of a suitable compact open set in $X'_\pi$, the function $f$ becomes constant along fibres of $\pi$, hence factors through $\pi$. This is the key input for the reduction maps constructed below, and it is the point where regularity replaces excision in the factor groupoid situation, compare [14, Lemma 4.8, Definition 4.10].

**Lemma 3.1.17** ([14, Lemma 4.8])**.** Let $(X, d)$ and $(X', d')$ be locally compact metric spaces which are totally disconnected, let $A$ be a discrete abelian group, and let $\pi : X \to X'$ be a continuous proper surjection which is regular in the sense of Definition 3.1.11. Write $X_\pi \subseteq X$ and $X'_\pi \subseteq X'$ as in Definition 3.1.11, and regard $\pi$ as a map $X_\pi \to X'_\pi$. Then for every $f \in C_c(X, A)$ there exist a compact open subset $L' \subseteq X'_\pi$ and a continuous map $g : X'_\pi \to A$ such that $f(x) = g(\pi(x))$ for all $x \in X_\pi \setminus \pi^{-1}(L')$.

*Proof.* Let $K := \operatorname{supp}(f)$, which is compact. Since $A$ is discrete, $f$ is locally constant, hence $K = \{x \in X \mid f(x) \neq 0\}$ and $K$ is compact open in $X$. Choose $\delta_1 > 0$ such that $d(x, y) < \delta_1$ implies $f(x) = f(y)$ for all $x, y \in K$. Since $K$ is compact and $X \setminus K$ is closed and disjoint from $K$, the number

$$\delta_0 := \inf\{d(x, y) \mid x \in K,\ y \in X \setminus K\}$$

is strictly positive. Set $\delta := \min\{\delta_0, \delta_1\}$. Then $d(x, y) < \delta$ implies $f(x) = f(y)$ for all $x, y \in X$.

Define

$$K' := \pi(K) \cap \{x' \in X'_\pi \mid \operatorname{diam}(\pi^{-1}(x')) \geq \delta\}.$$

Since $\pi$ is proper and $K$ is compact, $\pi(K)$ is compact in $X'$, and Lemma 3.1.13 shows that $K'$ is compact in $(X'_\pi, d'_\pi)$. Because $X'_\pi$ is locally compact and totally disconnected, choose a compact open neighbourhood $L' \subseteq X'_\pi$ of $K'$. Since $X'_\pi$ is Hausdorff, $L'$ is clopen.

Put $B := X'_\pi \setminus L'$ and $E := \pi^{-1}(B) \subseteq X_\pi$. We claim that $f|_E$ is constant on each fibre of $\pi|_E : E \to B$. Fix $x' \in B$ and $x, y \in \pi^{-1}(x')$. If $x' \notin \pi(K)$, then $x, y \notin K$, hence $f(x) = f(y) = 0$. If $x' \in \pi(K)$, then $x' \notin K'$ because $K' \subseteq L'$, hence $\operatorname{diam}(\pi^{-1}(x')) < \delta$, so $d(x, y) < \delta$ and therefore $f(x) = f(y)$.

By Proposition 3.1.14, the map $\pi : X_\pi \to X'_\pi$ is open. Hence $\pi|_E : E \to B$ is an open surjection and therefore a quotient map. Since $f|_E$ is constant on fibres, there exists a unique continuous map $g_B : B \to A$ such that $f|_E = g_B \circ (\pi|_E)$. Define $g : X'_\pi \to A$ by $g|_B := g_B$ and $g|_{L'} := 0$. Because $B$ and $L'$ are clopen and $A$ is discrete, $g$ is continuous. For $x \in X_\pi \setminus \pi^{-1}(L') = E$ we then have $f(x) = g(\pi(x))$, as required. $\square$

For $f \in C_c(X, A)$ and a compact open subset $K' \subseteq X'_\pi$ define $f_{K'} \in C_c(X_\pi, A)$ by

$$f_{K'}(x) := \begin{cases} f(x), & x \in \pi^{-1}(K'), \\ 0, & x \notin \pi^{-1}(K'). \end{cases}$$

**Corollary 3.1.18** ([14, Corollary 4.9]). Let $f \in C_c(X, A)$, and let $L' \subseteq X'_\pi$ be as in Lemma 3.1.17. If $K'_1, K'_2 \subseteq X'_\pi$ are compact open subsets with $L' \subseteq K'_1$ and $L' \subseteq K'_2$, then $f_{K'_1} - f_{K'_2}$ lies in the subgroup $\pi^* C_c(X'_\pi, A) \subseteq C_c(X_\pi, A)$, where $\pi^*(h) := h \circ \pi$. In particular, the class of $f_{K'}$ in $C_c(X_\pi, A)/\pi^* C_c(X'_\pi, A)$ is independent of the choice of compact open $K' \supseteq L'$.

*Proof.* By Lemma 3.1.17 there exists a continuous $g : X'_\pi \to A$ such that $f = g \circ \pi$ on $X_\pi \setminus \pi^{-1}(L')$. Set $D := (K'_1 \setminus K'_2) \sqcup (K'_2 \setminus K'_1) \subseteq X'_\pi \setminus L'$. Define $h : X'_\pi \to A$ by

$$h(x') := \begin{cases} g(x'), & x' \in K'_1 \setminus K'_2, \\ -g(x'), & x' \in K'_2 \setminus K'_1, \\ 0, & x' \notin D. \end{cases}$$

The sets $K'_1 \setminus K'_2$ and $K'_2 \setminus K'_1$ are compact open, hence $h \in C_c(X'_\pi, A)$.

For $x \in X_\pi$, if $\pi(x) \in K'_1 \cap K'_2$ or $\pi(x) \notin K'_1 \cup K'_2$, then $(f_{K'_1} - f_{K'_2})(x) = 0 = h(\pi(x))$. If $\pi(x) \in K'_1 \setminus K'_2$, then $x \notin \pi^{-1}(L')$ and

$$(f_{K'_1} - f_{K'_2})(x) = f(x) = g(\pi(x)) = h(\pi(x)).$$

If $\pi(x) \in K'_2 \setminus K'_1$, then $x \notin \pi^{-1}(L')$ and

$$(f_{K'_1} - f_{K'_2})(x) = -f(x) = -g(\pi(x)) = h(\pi(x)).$$

Thus $f_{K'_1} - f_{K'_2} = h \circ \pi = \pi^*(h) \in \pi^* C_c(X'_\pi, A)$. □

Our goal is to construct, from a regular proper surjection $\pi : X \to X'$ between totally disconnected locally compact metric spaces, a canonical reduction of compactly supported $A$-valued functions on $X$ modulo pullbacks from $X'_\pi$. The key input is Lemma 3.1.17: outside the preimage of a suitable compact open subset $L' \subseteq X'_\pi$, every $f \in C_c(X, A)$ is constant along the fibres of $\pi$ and hence factors through $\pi$. Therefore one may truncate $f$ to $\pi^{-1}(K') \subseteq X_\pi$ for any compact open $K' \supseteq L'$, and the resulting class in $C_c(X_\pi, A)/\pi^* C_c(X'_\pi, A)$ is independent of the choice of $K'$. This is the reduction map $\Pi$ used later in the long exact sequence.

For $f \in C_c(X, A)$ and a compact open subset $K' \subseteq X'_\pi$, define $f_{K'} \in C_c(X_\pi, A)$ by

$$f_{K'}(x) := \begin{cases} f(x), & x \in \pi^{-1}(K'), \\ 0, & x \in X_\pi \setminus \pi^{-1}(K'). \end{cases}$$

**Lemma 3.1.19.** Let $A$ be a discrete abelian group. Fix $f \in C_c(X, A)$, and choose $L' \subseteq X'_\pi$ and $g : X'_\pi \to A$ as in Lemma 3.1.17, so that $L'$ is compact open and $f(x) = g(\pi(x))$ for all $x \in X_\pi \setminus \pi^{-1}(L')$. If $K'_1, K'_2 \subseteq X'_\pi$ are compact open with $L' \subseteq K'_1 \cap K'_2$, then

$$f_{K'_1} - f_{K'_2} \in \pi^* C_c(X'_\pi, A) \subseteq C_c(X_\pi, A),$$

where $\pi^*(h) := h \circ \pi$. In particular, the coset $f_{K'} + \pi^* C_c(X'_\pi, A)$ is independent of the choice of compact open $K' \supseteq L'$.

*Proof.* Set $D := (K'_1 \setminus K'_2) \sqcup (K'_2 \setminus K'_1)$, a compact open subset of $X'_\pi$. Define $h : X'_\pi \to A$ by

$$h(x') := \begin{cases} g(x'), & x' \in K'_1 \setminus K'_2, \\ -g(x'), & x' \in K'_2 \setminus K'_1, \\ 0, & x' \notin D. \end{cases}$$

Then $h \in C_c(X'_\pi, A)$. If $\pi(x) \in D$, then $\pi(x) \notin L'$ since $L' \subseteq K'_1 \cap K'_2$, hence $f(x) = g(\pi(x))$. If $\pi(x) \notin D$, then $\pi(x) \in K'_1 \cap K'_2$ or $\pi(x) \notin K'_1 \cup K'_2$, so $f_{K'_1}(x) - f_{K'_2}(x) = 0$. Therefore $f_{K'_1} - f_{K'_2} = h \circ \pi \in \pi^* C_c(X'_\pi, A)$. □

**Definition 3.1.20** (Reduction map)**.** Let $A$ be a discrete abelian group. For $f \in C_c(X, A)$ choose $L' \subseteq X'_\pi$ as in Lemma 3.1.17 and pick any compact open $K' \subseteq X'_\pi$ with $L' \subseteq K'$. Define

$$\Pi(f) := f_{K'} + \pi^* C_c(X'_\pi, A) \in C_c(X_\pi, A) / \pi^* C_c(X'_\pi, A),$$

where $\pi^*(h) := h \circ \pi$. By Lemma 3.1.19 this is independent of the choices and hence defines a well-defined homomorphism

$$\Pi : C_c(X, A) \to C_c(X_\pi, A) / \pi^* C_c(X'_\pi, A).$$

**Proposition 3.1.21** ([14, Proposition 4.11])**.** The reduction map $\Pi$ is a surjective homomorphism and $\ker(\Pi) = \pi^* C_c(X', A) \subseteq C_c(X, A)$.

*Proof.*

- **Homomorphism.** Fix $f_1, f_2 \in C_c(X, A)$. Choose compact open sets $L'_1, L'_2 \subseteq X'_\pi$ as in Lemma 3.1.17 for $f_1, f_2$, and choose a compact open $K' \subseteq X'_\pi$ with $L'_1 \cup L'_2 \subseteq K'$. Then $(f_1 + f_2)_{K'} = (f_1)_{K'} + (f_2)_{K'}$, hence $\Pi(f_1 + f_2) = \Pi(f_1) + \Pi(f_2)$.
- **Surjectivity.** Let $u \in C_c(X_\pi, A)$. Since $A$ is discrete and $X_\pi$ is totally disconnected, $u$ is locally constant with compact support, so $u = \sum_{j=1}^m a_j \chi_{C_j}$ for some $a_j \in A$ and pairwise disjoint compact open sets $C_j \subseteq X_\pi$. It suffices to lift $a\chi_C$ for a compact open $C \subseteq X_\pi$ and $a \in A$. Put $K' := \pi(C) \subseteq X'_\pi$. By Proposition 3.1.14, the map $\pi : X_\pi \to X'_\pi$ is open and proper, hence $K'$ is compact open. The sets $C$ and $\pi^{-1}(K') \setminus C$ are compact and disjoint in the totally disconnected locally compact space $X$, so there exists a compact open set $U \subseteq X$ with $C \subseteq U$ and $U \cap (\pi^{-1}(K') \setminus C) = \emptyset$. Let $f := a\chi_U \in C_c(X, A)$. Then on $X_\pi$ one has $f_{K'} = a\chi_{U \cap \pi^{-1}(K')} = a\chi_C$, so $\Pi(f) = a\chi_C + \pi^* C_c(X'_\pi, A)$. By additivity, $\Pi$ is surjective.
- **Kernel.** If $h \in C_c(X', A)$, then for any compact open $K' \subseteq X'_\pi$, $(h \circ \pi)_{K'} = (h|_{K'}) \circ \pi \in \pi^* C_c(X'_\pi, A)$, so $\Pi(h \circ \pi) = 0$ and $\pi^* C_c(X', A) \subseteq \ker(\Pi)$.
  Conversely, let $f \in \ker(\Pi)$. Choose $L' \subseteq X'_\pi$ and $g : X'_\pi \to A$ as in Lemma 3.1.17. Since $\Pi(f) = 0$, we may take $K' = L'$ and obtain $f_{L'} \in \pi^* C_c(X'_\pi, A)$. Hence $f$ is constant on each fibre of $\pi$ over $L'$. On $X_\pi \setminus \pi^{-1}(L')$ we have $f = g \circ \pi$, so $f$ is constant on fibres there as well. If $x' \in X' \setminus X'_\pi$, then $\#\pi^{-1}(x') = 1$ by definition of $X'_\pi$, so $f$ is trivially constant on that fibre. Therefore $f$ is constant on every fibre of $\pi$.

Define $\bar{f} : X' \to A$ by $\bar{f}(x') := f(x)$ for any $x \in \pi^{-1}(x')$. Then $f = \bar{f} \circ \pi$. Since $\pi$ is continuous, surjective, and open, it is a quotient map, hence $\bar{f}$ is continuous because $\bar{f} \circ \pi = f$ is continuous. Moreover, $\mathrm{supp}(\bar{f}) \subseteq \pi(\mathrm{supp}(f))$, which is compact since $\pi$ is proper. Thus $\bar{f} \in C_c(X', A)$ and $f \in \pi^* C_c(X', A)$.

□

The role of the spaces $X_\pi$ and $X'_\pi$ is to isolate the place where a proper surjection $\pi : X \to X'$ has nontrivial fibres. Under regularity, Proposition 3.1.14 shows that $(X_\pi, d_\pi)$ and $(X'_\pi, d'_\pi)$ are locally compact and that $\pi : X_\pi \to X'_\pi$ remains continuous, proper, and open. This is the input needed to form quotients of compactly supported functions in a controlled way. Lemma 3.1.17 is the key mechanism: for every $f \in C_c(X, A)$ the obstruction to factoring $f$ through $\pi$ is supported over a compact open region in $X'_\pi$. This motivates two steps.

- For a compact open $K' \subseteq X'_\pi$ one truncates $f$ to $f_{K'}$, keeping support over $\pi^{-1}(K') \subseteq X_\pi$.
- One then passes to the quotient by $\pi^* C_c(X'_\pi, A)$, which kills the fibrewise constant part. The reduction map $\Pi : C_c(X, A) \to C_c(X_\pi, A)/\pi^* C_c(X'_\pi, A)$ records the remaining fibrewise variation and is independent of auxiliary choices once $K'$ is large enough, see Definition 3.1.20 and Proposition 3.1.21.

In particular, Proposition 3.1.21 gives a canonical identification

$$C_c(X, A)/\pi^* C_c(X', A) \;\cong\; C_c(X_\pi, A)/\pi^* C_c(X'_\pi, A).$$

For later applications it is essential that this reduction procedure is natural with respect to local homeomorphisms. The next proposition records the compatibility needed when the construction is applied degreewise to nerves in the factor groupoid situation.

**Proposition 3.1.22** ([14, Proposition 4.12])**.** For $i = 1, 2$, let $(X_i, d_i)$ and $(X'_i, d'_i)$ be locally compact metric spaces which are totally disconnected, and let $\pi_i : X_i \to X'_i$ be continuous regular proper surjections. Let $\varphi : X_1 \to X_2$ and $\varphi' : X'_1 \to X'_2$ be surjective local homeomorphisms such that $\pi_2 \circ \varphi = \varphi' \circ \pi_1$, and assume that for every $x'_1 \in X'_1$ the restriction $\varphi : \pi_1^{-1}(x'_1) \to \pi_2^{-1}(\varphi'(x'_1))$ is bijective. Put $Y_i := (X_i)_{\pi_i}$ and $Y'_i := (X'_i)_{\pi_i}$. Let $\sigma$ be the map induced by $\varphi_*$,

$$\sigma : C_c(X_1, A)/\pi_1^* C_c(X'_1, A) \to C_c(X_2, A)/\pi_2^* C_c(X'_2, A),$$

and let $\tau$ be the map induced by $\pi_*|_{C_c(Y_1, A)}$,

$$\tau : C_c(Y_1, A)/\pi_1^* C_c(Y'_1, A) \to C_c(Y_2, A)/\pi_2^* C_c(Y'_2, A).$$

Then the reduction maps $\Pi_i$ satisfy $\tau \circ \Pi_1 = \Pi_2 \circ \sigma$ as maps $C_c(X_1, A)/\pi_1^* C_c(X'_1, A) \to C_c(Y_2, A)/\pi_2^* C_c(Y'_2, A)$.

*Proof.* We first record that $\sigma$ and $\tau$ are well defined. For $g \in C_c(X_1', A)$ and $y \in X_2$, using $\pi_2 \circ \varphi = \varphi' \circ \pi_1$ and the fibrewise bijectivity assumption, we have

$$\varphi_*(\pi_1^* g)(y) = \sum_{x \in \varphi^{-1}(y)} g(\pi_1(x)) = \sum_{x_1' \in (\varphi')^{-1}(\pi_2(y))} g(x_1') = (\varphi')_* g(\pi_2(y)) = \pi_2^*((\varphi')_* g)(y).$$

Hence $\varphi_*(\pi_1^* C_c(X_1', A)) \subseteq \pi_2^* C_c(X_2', A)$, and similarly $\varphi_*(\pi_1^* C_c(Y_1', A)) \subseteq \pi_2^* C_c(Y_2', A)$, so $\sigma$ and $\tau$ are induced by $\pi_*$.

Let $f \in C_c(X_1, A)$ and set $S := \{x \in X_1 \mid f(x) \neq 0\}$, so $S$ is compact. Choose a compact open set $D' \subseteq X_1'$ with $\pi_1(S) \subseteq D'$. By Proposition 3.1.15, the restrictions

$$\varphi|_{Y_1} : Y_1 \to Y_2, \qquad \varphi'|_{Y_1'} : Y_1' \to Y_2'$$

are local homeomorphisms and satisfy $\pi_2 \circ \varphi = \varphi' \circ \pi_1$ on $Y_1$. Choose compact open sets $L_1' \subseteq Y_1'$ and $L_2' \subseteq Y_2'$ as in Lemma 3.1.17 for $f$ and $\varphi_* f$, respectively. Choose a compact open set $C' \subseteq Y_2'$ with $L_2' \subseteq C'$, and define $E' := (\varphi')^{-1}(C') \cap D' \cap Y_1'$. Since $Y_1'$ is totally disconnected and locally compact, choose a compact open $K' \subseteq Y_1'$ with $L_1' \cup E' \subseteq K'$.

Enlarging cut-offs does not change reduction classes by Corollary 3.1.18, hence

$$\Pi_1(f + \pi_1^* C_c(X_1', A)) = f_{K'} + \pi_1^* C_c(Y_1', A),$$
$$\Pi_2(\pi_* f + \pi_2^* C_c(X_2', A)) = (\pi_* f)_{C'} + \pi_2^* C_c(Y_2', A).$$

It suffices to prove $\pi_*(f_{K'}) - (\pi_* f)_{C'} \in \pi_2^* C_c(Y_2', A)$. We prove the stronger identity $\pi_*(f_{E'}) = (\pi_* f)_{C'}$ in $C_c(Y_2, A)$. Since $f_{K'} - f_{E'} \in \pi_1^* C_c(Y_1', A)$, applying $\pi_*$ and using the well-definedness computation above yields $\pi_*(f_{K'}) - \pi_*(f_{E'}) \in \pi_2^* C_c(Y_2', A)$, so the claim follows.

- **Compactness device.** Let $C' \subseteq Y_2'$ be compact and let $D' \subseteq X_1'$ be compact. Then $E' := (\varphi')^{-1}(C') \cap D' \cap Y_1'$ is compact in $Y_1'$. Let $(y_k')_k$ be a sequence in $E'$. Since $D'$ is compact in $X_1'$, after passing to a subsequence we may assume $y_k' \to y'$ in $X_1'$. Since $\varphi'(y_k') \in C'$ and $C'$ is compact in $Y_2'$, after passing to a further subsequence we may assume $\varphi'(y_k') \to z'$ in $Y_2'$. By Lemma 3.1.12, convergence in $Y_2'$ implies convergence in $X_2'$, hence $\varphi'(y_k') \to z'$ in $X_2'$. Continuity of $\varphi'$ gives $\varphi'(y_k') \to \varphi'(y')$ in $X_2'$, hence $z' = \varphi'(y')$. Thus $\varphi'(y') \in C' \subseteq Y_2'$. By Proposition 3.1.15, this implies $y' \in Y_1'$. Moreover $y' \in D'$ since $D'$ is closed in $X_1'$, and $y' \in (\varphi')^{-1}(C')$. Hence $y' \in E'$.
  Since $\varphi'|_{Y_1'} : Y_1' \to Y_2'$ is a local homeomorphism, there exists an open neighbourhood $U' \subseteq Y_1'$ of $y'$ such that $\varphi'|_{U'}$ is a homeomorphism onto the open set $\varphi'(U') \subseteq Y_2'$. As $\varphi'(y_k') \to \varphi'(y')$ in $Y_2'$, we have $\varphi'(y_k') \in \varphi'(U')$ for $k \gg 0$, hence $y_k' \in U'$ for $k \gg 0$. Therefore $y_k' \to y'$ in $Y_1'$. Thus $E'$ is sequentially compact, hence compact.
- **Identify the representatives.** Fix $y \in Y_2$. If $\pi_2(y) \notin C'$, then for every $x \in \varphi^{-1}(y)$ one has

  $$\varphi'(\pi_1(x)) = \pi_2(\varphi(x)) = \pi_2(y) \notin C',$$

  so $\pi_1(x) \notin (\varphi')^{-1}(C')$ and hence $\pi_1(x) \notin E'$. Thus $f_{E'}(x) = 0$ for all $x \in \varphi^{-1}(y)$, so $\varphi_*(f_{E'})(y) = 0 = (\varphi_* f)_{C'}(y)$. If $\pi_2(y) \in C'$, then for every $x \in \varphi^{-1}(y)$ one has $\pi_1(x) \in$ $(\varphi')^{-1}(C')$. Moreover $f(x) = 0$ whenever $\pi_1(x) \notin D'$, since $\pi_1(S) \subseteq D'$. Hence the fibrewise sum defining $\varphi_* f(y)$ only receives contributions from $x \in \varphi^{-1}(y)$ with $\pi_1(x) \in D'$. For such $x$ we have $\pi_1(x) \in E'$ and $f_{E'}(x) = f(x)$. Therefore $\varphi_*(f_{E'})(y) = \varphi_* f(y) = (\varphi_* f)_{C'}(y)$. Thus $\pi_*(f_{E'}) = (\pi_* f)_{C'}$ in $C_c(Y_2, A)$.

Combining the reductions and the definitions of $\sigma$ and $\tau$, we obtain

$$\begin{aligned}\tau\Big(\Pi_1(f + \pi_1^* C_c(X_1', A))\Big) &= \pi_*(f_{K'}) + \pi_2^* C_c(Y_2', A)\\ &= (\pi_* f)_{C'} + \pi_2^* C_c(Y_2', A)\\ &= \Pi_2\Big(\sigma(f + \pi_1^* C_c(X_1', A))\Big).\end{aligned}$$

This proves $\tau \circ \Pi_1 = \Pi_2 \circ \sigma$. □

### 3.1.3 From Regularity on Units to Regularity on Nerves

We record the groupoid-level setting and its consequences following [14, Setting 4.13].

**Setting 3.1.23** ([14, Setting 4.13])**.** Let $\mathcal{G}$ and $\mathcal{G}'$ be étale groupoids and let $\pi : \mathcal{G} \to \mathcal{G}'$ be a continuous homomorphism such that:

1. $\pi$ is surjective,
2. $\pi$ is proper,
3. for every $x \in \mathcal{G}_0$, the restriction $\pi : r^{-1}(x) \to r^{-1}(\pi(x))$ is bijective,
4. $\mathcal{G}$ and $\mathcal{G}'$ are totally disconnected.

**Lemma 3.1.24** ([14, Lemma 4.15])**.** Let $\pi : \mathcal{G} \to \mathcal{G}'$ be as in Setting 3.1.23.

The following are equivalent:

1. $\pi : \mathcal{G}_0 \to \mathcal{G}_0'$ is regular,
2. $\pi : \mathcal{G} \to \mathcal{G}'$ is regular,
3. there exists $n \in \mathbb{N} \setminus \{1\}$ such that $\pi_n : \mathcal{G}_n \to (\mathcal{G}')_n$ is regular.

*Proof.*

- **1.⇔2.** Write $\pi_0 : \mathcal{G}_0 \to \mathcal{G}_0'$ for the map on units. Since $\mathcal{G}$ and $\mathcal{G}'$ are étale, the range maps $r : \mathcal{G} \to \mathcal{G}_0$ and $r : \mathcal{G}' \to \mathcal{G}_0'$ are local homeomorphisms. The commutative square

$$\begin{array}{ccc}\mathcal{G} & \xrightarrow{\pi} & \mathcal{G}' \\ {\scriptstyle r}\downarrow & & \downarrow{\scriptstyle r} \\ \mathcal{G}_0 & \xrightarrow{\pi_0} & \mathcal{G}_0'\end{array}$$

  satisfies the fibrewise bijectivity hypothesis along $r$ by assumption 3. Therefore Proposition 3.1.15 applied to this square yields 1.⇔2.
- **2.⇔3.** For $n \geq 2$, the projection $p_1 : \mathcal{G}_n \to \mathcal{G}$, $p_1(g_1, \dots, g_n) = g_1$, is a local homeomorphism, and similarly $p_1 : \mathcal{G}_n' \to \mathcal{G}'$. Moreover, $\pi_n$ is defined levelwise and satisfies $p_1 \circ \pi_n = \pi \circ p_1$, and the induced maps on fibres of $p_1$ are bijections. Applying Proposition 3.1.15 again to the following square gives 2.⇔3.:

$$\begin{array}{ccc} \mathcal{G}_n & \xrightarrow{\pi_n} & \mathcal{G}'_n \\ {\scriptstyle p_1}\downarrow & & \downarrow{\scriptstyle p_1} \\ \mathcal{G} & \xrightarrow{\pi} & \mathcal{G}'. \end{array}$$

□

**Proposition 3.1.25** ([14, Proposition 4.17])**.** Let $\pi : \mathcal{G} \to \mathcal{G}'$ be as in Setting 3.1.23 and assume $\pi$ is regular. Then the groupoids $\mathcal{G}_\pi$ and $\mathcal{G}'_\pi$ are étale for the topologies induced by the metrics $d_\pi$ and $d'_\pi$ from Definition 3.1.11.

*Proof.* Apply Proposition 3.1.14 to the regular proper surjection $\pi : \mathcal{G} \to \mathcal{G}'$. This yields that $\mathcal{G}_\pi$ and $\mathcal{G}'_\pi$ are locally compact Hausdorff and that the restricted map $\pi : \mathcal{G}_\pi \to \mathcal{G}'_\pi$ is continuous, proper, and open. By Lemma 3.1.24, regularity of $\pi$ implies regularity of $\pi_0 : \mathcal{G}_0 \to \mathcal{G}'_0$. Applying Proposition 3.1.14 again to $\pi_0$ shows that $(\mathcal{G}_0)_\pi$ and $(\mathcal{G}'_0)_\pi$ are locally compact and that $\pi : (\mathcal{G}_0)_\pi \to (\mathcal{G}'_0)_\pi$ is open. The range map $r : \mathcal{G} \to \mathcal{G}_0$ is a local homeomorphism. Restricting to $\mathcal{G}_\pi \subseteq \mathcal{G}$ and $(\mathcal{G}_0)_\pi \subseteq \mathcal{G}_0$ gives a continuous map $r : \mathcal{G}_\pi \to (\mathcal{G}_0)_\pi$. To see that $r$ is a local homeomorphism for the $d_\pi$-topology, fix $\gamma \in \mathcal{G}_\pi$ and choose an open bisection $U \subseteq \mathcal{G}$ with $\gamma \in U$. Then $r|_U : U \to r(U)$ is a homeomorphism. Since $\pi|_{r(U)}$ is a homeomorphism onto $\pi(r(U))$ by assumption 3. and étaleness, the set $U \cap \mathcal{G}_\pi$ is open in $\mathcal{G}_\pi$ and $r(U \cap \mathcal{G}_\pi) = r(U) \cap (\mathcal{G}_0)_\pi$ is open in $(\mathcal{G}_0)_\pi$. Hence $r : U \cap \mathcal{G}_\pi \to r(U) \cap (\mathcal{G}_0)_\pi$ is a homeomorphism. Thus $r : \mathcal{G}_\pi \to (\mathcal{G}_0)_\pi$ is a local homeomorphism. The same argument applies to $r : \mathcal{G}'_\pi \to (\mathcal{G}'_0)_\pi$. Therefore $\mathcal{G}_\pi$ and $\mathcal{G}'_\pi$ are étale. □

**Lemma 3.1.26** ([14, Lemma 4.18])**.** Let $\pi : \mathcal{G} \to \mathcal{G}'$ be as in Setting 3.1.23 and assume $\pi$ is regular. For any $n \in \mathbb{N} \setminus \{1\}$, the spaces $(\mathcal{G}'_\pi)_n$ and $(\mathcal{G}'_n)_{\pi_n}$ are canonically homeomorphic, and similarly $(\mathcal{G}_\pi)_n \cong (\mathcal{G}_n)_{\pi_n}$.

*Proof.* Fix $n \geq 2$. By definition,

$$\begin{aligned} (\mathcal{G}'_\pi)_n &= \{(g'_1, \dots, g'_n) \in (\mathcal{G}'_\pi)^n \mid s(g'_i) = r(g'_{i+1})\}, \\ (\mathcal{G}'_n)_{\pi_n} &= \{h' \in \mathcal{G}'_n \mid \#(\pi_n)^{-1}(h') > 1\}. \end{aligned}$$

We claim that these subsets of $(\mathcal{G}')^n$ coincide.

Let $(g'_1, \dots, g'_n) \in \mathcal{G}'_n$. Then $(g'_1, \dots, g'_n) \in (\mathcal{G}'_n)_{\pi_n}$ if and only if there exist two distinct elements of $\mathcal{G}_n$ mapping to $(g'_1, \dots, g'_n)$ under $\pi_n$. Since $\pi_n$ is defined componentwise,

$$\pi_n(g_1, \dots, g_n) = (\pi(g_1), \dots, \pi(g_n)),$$

and since $\pi$ is fibrewise bijective along the range map in Setting 3.1.23, the fibre $\pi_n^{-1}(g'_1,\dots,g'_n)$ is in bijection with $\pi^{-1}(g'_1)$ via projection to the first coordinate. In particular,

$$\#\pi_n^{-1}(g'_1,\dots,g'_n) > 1 \Leftrightarrow \#\pi^{-1}(g'_1) > 1.$$

On the other hand, $(g'_1,\dots,g'_n) \in (\mathcal{G}'_\pi)_n$ if and only if $g'_1 \in \mathcal{G}'_\pi$, that is $\#\pi^{-1}(g'_1) > 1$. Thus $(\mathcal{G}'_\pi)_n = (\mathcal{G}'_n)_{\pi_n}$ as subsets of $(\mathcal{G}')^n$. It remains to check that the induced topologies agree. Consider the map $\varphi : \mathcal{G}'_n \to \mathcal{G}'_0, (g'_1,\dots,g'_n) \mapsto r(g'_1)$. Since $\mathcal{G}'$ is étale, $\varphi$ is a local homeomorphism. Moreover, by Proposition 3.1.15 applied to $\varphi$ and its induced map on the base, the restriction $\varphi : (\mathcal{G}'_n)_{\pi_n} \to (\mathcal{G}'_0)_{\pi_0}$ is a local homeomorphism.

On the other hand, $(\mathcal{G}'_\pi)_n$ carries the relative product topology induced from $(\mathcal{G}'_\pi)^n$, and the same map $\varphi$ is a local homeomorphism $\varphi : (\mathcal{G}'_\pi)_n \to (\mathcal{G}'_0)_{\pi_0}$, since $r : \mathcal{G}'_\pi \to (\mathcal{G}'_0)_{\pi_0}$ is a local homeomorphism by Proposition 3.1.25. Hence both subspace topologies on the common subset $(\mathcal{G}')^n$ are characterised by the requirement that $\varphi$ be a local homeomorphism to $(\mathcal{G}'_0)_{\pi_0}$. Therefore they agree, and the identification is canonical.

The statement for $\mathcal{G}$ is identical. □

**Proposition 3.1.27** ([14, Proposition 4.19])**.** Let $\pi : \mathcal{G} \to \mathcal{G}'$ be as in Setting 3.1.23 and assume $\pi$ is regular. Let $A$ be a discrete abelian group. For each $n \geq 0$, let $\Pi_n$ be the reduction map associated to $\pi_n : \mathcal{G}_n \to \mathcal{G}'_n$ as in Definition 3.1.20. Then:

1. $\Pi_n$ induces an isomorphism

$$C_c(\mathcal{G}_n, A)/\pi_n^* C_c(\mathcal{G}'_n, A) \xrightarrow{\cong} C_c((\mathcal{G}_\pi)_n, A)/\pi_n^* C_c((\mathcal{G}'_\pi)_n, A),$$

2. the family $(\Pi_n)_{n\geq 0}$ defines an isomorphism of chain complexes

$$C_c(\mathcal{G}_\bullet, A)/\pi_\bullet^* C_c(\mathcal{G}'_\bullet, A) \xrightarrow{\cong} C_c((\mathcal{G}_\pi)_\bullet, A)/\pi_\bullet^* C_c((\mathcal{G}'_\pi)_\bullet, A),$$

where the differentials are the Moore differentials on the respective nerves.

*Proof.*

- **Degreewise isomorphism.** Fix $n \geq 0$. By Lemma 3.1.24, the map $\pi_n : \mathcal{G}_n \to \mathcal{G}'_n$ is a regular proper surjection between totally disconnected locally compact metric spaces. Apply Proposition 3.1.21 to $\pi_n$. Then $\Pi_n$ is surjective with kernel $\pi_n^* C_c(\mathcal{G}'_n, A)$, hence it induces an isomorphism

$$C_c(\mathcal{G}_n, A)/\pi_n^* C_c(\mathcal{G}'_n, A) \xrightarrow{\cong} C_c((\mathcal{G}_n)_{\pi_n}, A)/\pi_n^* C_c((\mathcal{G}'_n)_{\pi_n}, A).$$

By Lemma 3.1.26, there are canonical homeomorphisms $(\mathcal{G}_n)_{\pi_n} \cong (\mathcal{G}_\pi)_n$, $(\mathcal{G}'_n)_{\pi_n} \cong (\mathcal{G}'_\pi)_n$, so the claim follows.

- **Compatibility with Moore differentials.** Fix $n \geq 1$ and $0 \leq i \leq n$. The face maps $d_i : \mathcal{G}_n \to \mathcal{G}_{n-1}$ and $d_i : \mathcal{G}'_n \to \mathcal{G}'_{n-1}$ are local homeomorphisms since $\mathcal{G}$ and $\mathcal{G}'$ are étale, and the squares with $\pi_n$ commute. Apply Proposition 3.1.22 to the commutative square

$$\begin{array}{ccc} \mathcal{G}_n & \xrightarrow{d_i} & \mathcal{G}_{n-1} \\ {\scriptstyle \pi_n}\big\downarrow & & \big\downarrow{\scriptstyle \pi_{n-1}} \\ \mathcal{G}'_n & \xrightarrow[d_i]{} & \mathcal{G}'_{n-1}. \end{array}$$

This yields $(d_i)_* \circ \Pi_n = \Pi_{n-1} \circ (d_i)_*$. Taking the alternating sum over $i$ gives $\partial_n \circ \Pi_n = \Pi_{n-1} \circ \partial_n$ Therefore $(\Pi_n)_{n \geq 0}$ is a chain map. By the first item it is degreewise an isomorphism, hence a chain isomorphism.

□

**Theorem 3.1.28** ([14, Theorem 4.20])**.** Let $\pi : \mathcal{G} \to \mathcal{G}'$ be as in Setting 3.1.23 and assume $\pi$ is regular. Let $A$ be a discrete abelian group. Then there is a commutative diagram whose two rows are long exact sequences in Moore homology:

$$\begin{array}{ccccccccccc} \cdots & \longrightarrow & H_n(\mathcal{G}', A) & \xrightarrow{H_n(\pi^*)} & H_n(\mathcal{G}, A) & \xrightarrow{H_n(q)} & H_n(Q_\bullet) & \xrightarrow{\partial_n} & H_{n-1}(\mathcal{G}', A) & \longrightarrow & \cdots \\ & & {\scriptstyle H_n(\iota')}\big\downarrow & & {\scriptstyle H_n(\iota)}\big\downarrow & & {\scriptstyle H_n(\Pi_\bullet)}\big\downarrow & & {\scriptstyle H_{n-1}(\iota')}\big\downarrow & & \\ \cdots & \longrightarrow & H_n(\mathcal{G}'_\pi, A) & \xrightarrow{H_n(\pi^*)} & H_n(\mathcal{G}_\pi, A) & \xrightarrow{H_n(q_\pi)} & H_n(Q^\pi_\bullet) & \xrightarrow{\partial^\pi_n} & H_{n-1}(\mathcal{G}'_\pi, A) & \longrightarrow & \cdots \end{array}$$

Here

- $\pi^*$ denotes pullback along $\pi_n : \mathcal{G}_n \to \mathcal{G}'_n$ in each degree, hence a chain map $C_c(\mathcal{G}'_\bullet, A) \to C_c(\mathcal{G}_\bullet, A)$ and likewise $C_c((\mathcal{G}'_\pi)_\bullet, A) \to C_c((\mathcal{G}_\pi)_\bullet, A)$,
- $q$ and $q_\pi$ are the quotient maps of chain complexes,
- $\iota, \iota'$ are induced by inclusions $C_c((\mathcal{G}_\pi)_\bullet, A) \hookrightarrow C_c(\mathcal{G}_\bullet, A), C_c((\mathcal{G}'_\pi)_\bullet, A) \hookrightarrow C_c(\mathcal{G}'_\bullet, A)$,
- $H_n(\Pi_\bullet)$ is the isomorphism induced by Proposition 3.1.27,
- $\partial_n$ and $\partial^\pi_n$ are the connecting morphisms,
- $Q_\bullet := C_c(\mathcal{G}_\bullet, A)/\pi^* C_c(\mathcal{G}'_\bullet, A)$ and $Q^\pi_\bullet := C_c((\mathcal{G}_\pi)_\bullet, A)/\pi^* C_c((\mathcal{G}'_\pi)_\bullet, A)$.

In particular, the vertical map on the quotient term is an isomorphism.

*Proof.*

- **Short exact sequence and connecting morphisms for $\mathcal{G}$.** For each $n \geq 0$, pullback along $\pi_n$ defines a homomorphism

$$\pi_n^* : C_c(\mathcal{G}'_n, A) \to C_c(\mathcal{G}_n, A), \qquad (\pi_n^* f)(g) := f(\pi_n(g)).$$

Since $\pi_n$ is surjective, $\pi_n^*$ is injective. Since $\pi_n \circ d_i = d_i \circ \pi_{n+1}$ for all face maps $d_i$, these maps assemble to a chain map $\pi^* : C_c(\mathcal{G}'_\bullet, A) \to C_c(\mathcal{G}_\bullet, A)$. By definition of $Q_\bullet$, for each $n \geq 0$ there is a short exact sequence

$$0 \to C_c(\mathcal{G}'_n, A) \xrightarrow{\pi_n^*} C_c(\mathcal{G}_n, A) \xrightarrow{q_n} Q_n \to 0,$$

hence a short exact sequence of chain complexes

$$0 \to C_c(\mathcal{G}'_\bullet, A) \xrightarrow{\pi^*} C_c(\mathcal{G}_\bullet, A) \xrightarrow{q} Q_\bullet \to 0.$$

This yields a long exact sequence in homology with connecting morphisms $\partial_n$. Concretely, writing the Moore differential on all complexes as $d_n$, the connecting morphism

$$\partial_n : H_n(Q_\bullet) \to H_{n-1}(\mathcal{G}', A)$$

is defined by the following lift and boundary procedure. Choose a cycle $c \in Q_n$ representing a class in $H_n(Q_\bullet)$. Choose $b \in C_c(\mathcal{G}_n, A)$ with $q_n(b) = c$. Then $0 = d_{n-1}(c) = q_{n-1}(d_{n-1}(b))$, so $d_{n-1}(b) \in \pi_{n-1}^* C_c(\mathcal{G}'_{n-1}, A)$. Since $\pi_{n-1}^*$ is injective, there is a unique $a \in C_c(\mathcal{G}'_{n-1}, A)$ such that $\pi_{n-1}^*(a) = d_{n-1}(b)$. Set $\partial_n([c]) := [a]$.
This is well defined and gives exactness.

- **Short exact sequence and connecting morphisms for $\mathcal{G}_\pi$.** Assume $\pi$ is regular. Applying the same construction degreewise to $\pi_n : (\mathcal{G}_\pi)_n \to (\mathcal{G}'_\pi)_n$ gives a short exact sequence

$$0 \to C_c((\mathcal{G}'_\pi)_\bullet, A) \xrightarrow{\pi^*} C_c((\mathcal{G}_\pi)_\bullet, A) \xrightarrow{q_\pi} Q_\bullet^\pi \to 0,$$

hence a long exact sequence with connecting morphisms $\partial_n^\pi$, by the same lift.

- **Vertical maps and commutativity.** Since the topology on $(\mathcal{G}_\pi)_n$ and $(\mathcal{G}'_\pi)_n$ is finer than the subspace topology from $\mathcal{G}_n$ and $\mathcal{G}'_n$, every element of $C_c((\mathcal{G}_\pi)_n, A)$ and $C_c((\mathcal{G}'_\pi)_n, A)$ is also an element of $C_c(\mathcal{G}_n, A)$ and $C_c(\mathcal{G}'_n, A)$. Thus the inclusions of chain complexes

$$C_c((\mathcal{G}_\pi)_\bullet, A) \hookrightarrow C_c(\mathcal{G}_\bullet, A), \qquad C_c((\mathcal{G}'_\pi)_\bullet, A) \hookrightarrow C_c(\mathcal{G}'_\bullet, A)$$

induce homology maps $H_n(\iota)$ and $H_n(\iota')$. By Proposition 3.1.27, the reduction maps $\Pi_n$ assemble to a chain isomorphism $\Pi_\bullet : Q_\bullet \xrightarrow{\cong} Q_\bullet^\pi$, hence induce isomorphisms $H_n(\Pi_\bullet)$ on homology. The long exact sequence construction is natural with respect to morphisms of short exact sequences of chain complexes. Finally, apply naturality and homology to

$$\begin{array}{ccccccccc}
0 & \longrightarrow & C_c(\mathcal{G}'_\bullet, A) & \xrightarrow{\pi^*} & C_c(\mathcal{G}_\bullet, A) & \xrightarrow{q} & Q_\bullet & \longrightarrow & 0 \\
 & & \big\downarrow{\scriptstyle \iota'} & & \big\downarrow{\scriptstyle \iota} & & \big\downarrow{\scriptstyle \Pi_\bullet} & & \\
0 & \longrightarrow & C_c((\mathcal{G}'_\pi)_\bullet, A) & \xrightarrow{\pi^*} & C_c((\mathcal{G}_\pi)_\bullet, A) & \xrightarrow{q_\pi} & Q_\bullet^\pi & \longrightarrow & 0.
\end{array}$$

□

## 3.2 The Universal Coefficient Theorem

The UCT is purely algebraic and applies to chain complexes of free abelian groups with coefficients introduced via tensor products. In order to bring our groupoid homology with coefficients in a topological group $A$ into this framework, we now identify it with the homology of $C_c(\mathcal{G}_\bullet, \mathbb{Z}) \otimes_{\mathbb{Z}} A$ by constructing an explicit chain isomorphism $C_c(\mathcal{G}_n, \mathbb{Z}) \otimes_{\mathbb{Z}} A \cong C_c(\mathcal{G}_n, A)$ in each degree. We will restrict the hypothesis to discrete topological groups. As it turns out, homology using Moore complexes does not yield a UCT for general topological abelian groups.

**Proposition 3.2.1.** Let $\mathcal{G}$ be an ample groupoid and $A$ a discrete abelian group. For $n \geq 0$ write

$$C_n(\mathcal{G};\mathbb{Z}) := C_c(\mathcal{G}_n, \mathbb{Z}), \quad C_n(\mathcal{G};A) := C_c(\mathcal{G}_n, A),$$
$$\partial_n := \sum_{i=0}^{n} (-1)^i (d_i)_* : C_n(\mathcal{G};\mathbb{Z}) \to C_{n-1}(\mathcal{G};\mathbb{Z}),$$
$$\partial_n^A := \sum_{i=0}^{n} (-1)^i (d_i)_* : C_n(\mathcal{G};A) \to C_{n-1}(\mathcal{G};A)$$

be the Moore differentials. Then for each $n \geq 0$ there is a natural isomorphism in **Ab**

$$\Phi_n : C_c(\mathcal{G}_n, \mathbb{Z}) \otimes_{\mathbb{Z}} A \xrightarrow{\cong} C_c(\mathcal{G}_n, A), \qquad \Phi_n(f \otimes a)(x) := f(x) \cdot a,$$

and the family $\Phi_\bullet = (\Phi_n)_{n\geq 0}$ is a chain isomorphism

$$\Phi_\bullet : (C_c(\mathcal{G}_\bullet, \mathbb{Z}) \otimes_{\mathbb{Z}} A,\ \partial_\bullet \otimes \mathrm{id}_A) \xrightarrow{\cong} (C_c(\mathcal{G}_\bullet, A),\ \partial_\bullet^A).$$

In particular, $H_n(\mathcal{G};A) \ \cong\ H_n(C_c(\mathcal{G}_\bullet, \mathbb{Z}) \otimes_{\mathbb{Z}} A)$ for all $n \geq 0$.

*Proof.* Fix $n \geq 0$. Since $\mathcal{G}$ is ample, $\mathcal{G}_n$ is locally compact Hausdorff totally disconnected, hence has a basis of compact open sets. Because $A$ is discrete, every $\xi \in C_c(\mathcal{G}_n, A)$ is locally constant and $\{x \mid \xi(x) \neq 0\}$ is compact open in $\mathcal{G}_n$.

- **Definition of $\Phi_n$.** Define $\beta_n : C_c(\mathcal{G}_n, \mathbb{Z}) \times A \to C_c(\mathcal{G}_n, A)$ by $\beta_n(f, a)(x) := f(x) \cdot a$. This is biadditive and $\mathrm{supp}(\beta_n(f,a)) \subseteq \mathrm{supp}(f)$, hence $\beta_n(f, a) \in C_c(\mathcal{G}_n, A)$. By the universal property of $\otimes_{\mathbb{Z}}$ there is a unique homomorphism $\Phi_n : C_c(\mathcal{G}_n, \mathbb{Z}) \otimes_{\mathbb{Z}} A \to C_c(\mathcal{G}_n, A)$ with $\Phi_n(f \otimes a) = \beta_n(f, a)$.
- **Explicit inverse.** Let $\xi \in C_c(\mathcal{G}_n, A)$. Since $\xi$ is locally constant with compact open support, we can write

  $$\xi = \sum_{j=1}^{m} a_j \chi_{U_j},$$

  with $a_j \in A$ and pairwise disjoint compact open $U_j \subseteq \mathcal{G}_n$. Define

  $$\Psi_n(\xi) := \sum_{j=1}^{m} \chi_{U_j} \otimes a_j \in C_c(\mathcal{G}_n, \mathbb{Z}) \otimes_{\mathbb{Z}} A.$$

If $\xi = \sum_j a_j \chi_{U_j} = \sum_k b_k \chi_{V_k}$ are two such decompositions, refine by $W_{jk} := U_j \cap V_k$, still pairwise disjoint compact open. On $W_{jk} \neq \emptyset$ one has $a_j = b_k$, and bilinearity gives

$$\sum_j \chi_{U_j} \otimes a_j = \sum_{j,k} \chi_{W_{jk}} \otimes a_j = \sum_{j,k} \chi_{W_{jk}} \otimes b_k = \sum_k \chi_{V_k} \otimes b_k,$$

so $\Psi_n$ is well defined. Then $\Phi_n(\Psi_n(\xi)) = \xi$ is immediate. Conversely, for a pure tensor $f \otimes a$ write $f = \sum_{j=1}^{m} m_j \chi_{U_j}$ with pairwise disjoint compact open $U_j$, and compute

$$\Psi_n(\Phi_n(f \otimes a)) = \Psi_n\Big(\sum_j (m_j \cdot a)\, \chi_{U_j}\Big) = \sum_j \chi_{U_j} \otimes (m_j \cdot a) = \sum_j (m_j \chi_{U_j}) \otimes a = f \otimes a.$$

Hence $\Psi_n \circ \Phi_n = \mathrm{id}_{C_c(\mathcal{G}_n, \mathbb{Z}) \otimes_{\mathbb{Z}} A}$ and $\Phi_n$ is an isomorphism with inverse $\Psi_n$.

- **Compatibility with pushforward.** Let $p : Y \to Z$ be a local homeomorphism between locally compact Hausdorff spaces. For $f \in C_c(Y, \mathbb{Z})$ and $z \in Z$, the set $p^{-1}(z) \cap \mathrm{supp}(f)$ is finite because it is compact and discrete. Thus $p_*$ is defined by finite sums on both $\mathbb{Z}$- and $A$-valued chains. For $a \in A$ and $z \in Z$,

$$(\Phi_Z(p_* f \otimes a))(z) = (p_* f)(z) \cdot a = \Big( \sum_{p(y)=z} f(y) \Big) \cdot a = \sum_{p(y)=z} f(y) \cdot a = (p_* \Phi_Y(f \otimes a))(z).$$

Hence $\Phi_Z \circ (p_* \otimes \mathrm{id}_A) = p_* \circ \Phi_Y$. Applying this with $p = d_i : \mathcal{G}_n \to \mathcal{G}_{n-1}$ and summing with signs yields $\Phi_{n-1} \circ (\partial_n \otimes \mathrm{id}_A) = \partial_n^A \circ \Phi_n$ for all $n \geq 1$.

Therefore $\Phi_\bullet$ is a chain isomorphism, and the claimed identification on homology follows. □

We now define the cohomology theory that pairs naturally with the integral Moore complex and hence fits the classical cohomological universal coefficient formalism. Throughout, $\partial_\bullet$ denotes the Moore differential on the integral chain complex $(C_c(\mathcal{G}_\bullet, \mathbb{Z}), \partial_\bullet)$, constructed earlier from pushforwards along the face maps of the nerve.

**Proposition 3.2.2.** Let $A$ be an abelian group. For $n \geq 0$ set $C^n(\mathcal{G}; A) := \mathrm{Hom}_{\mathbb{Z}}(C_c(\mathcal{G}_n, \mathbb{Z}), A)$. Define $\delta^n : C^n(\mathcal{G}; A) \to C^{n+1}(\mathcal{G}; A)$ by $\delta^n(\varphi) := \varphi \circ \partial_{n+1}$, for $\varphi \in C^n(\mathcal{G}; A)$. Then $\delta^{n+1} \circ \delta^n = 0$ for all $n \geq 0$, so $(C^\bullet(\mathcal{G}; A), \delta^\bullet)$ is a cochain complex. Its cohomology groups are $H^n(\mathcal{G}; A) := \ker(\delta^n) / \operatorname{im}(\delta^{n-1})$, for $n \geq 0$, where $\operatorname{im}(\delta^{-1}) := 0$. Thus $H^n(\mathcal{G}; A)$ is the cohomology of the cochain complex $\mathrm{Hom}_{\mathbb{Z}}(C_c(\mathcal{G}_\bullet, \mathbb{Z}), A)$.

*Proof.* For $n \geq 0$ the map $\delta^n$ is well-defined since it is the composition of group homomorphisms. For $\varphi \in C^n(\mathcal{G}; A)$ one has

$$\delta^{n+1}(\delta^n(\varphi)) = (\varphi \circ \partial_{n+1}) \circ \partial_{n+2} = \varphi \circ (\partial_{n+1} \circ \partial_{n+2}) = 0,$$

because $\partial_{n+1} \circ \partial_{n+2} = 0$ in the chain complex $(C_c(\mathcal{G}_\bullet, \mathbb{Z}), \partial_\bullet)$. □

### 3.2.1 UCT for Moore Homology

We have homology with coefficients as the homology of $C_c(\mathcal{G}_\bullet, \mathbb{Z}) \otimes_{\mathbb{Z}} A$ and cohomology with coefficients as the cohomology of the dual complex $\mathrm{Hom}_{\mathbb{Z}}(C_c(\mathcal{G}_\bullet, \mathbb{Z}), A)$. Since each $C_c(\mathcal{G}_n, \mathbb{Z})$ is a free abelian group, the chain complex $C_c(\mathcal{G}_\bullet, \mathbb{Z})$ lies exactly in the algebraic setting where the classical UCT applies. In the next step we recall these UCT statements for an arbitrary chain complex of free abelian groups. We will then modify to $C_c(\mathcal{G}_\bullet, \mathbb{Z})$ to obtain the UCT-sequences for $H_n(\mathcal{G}; A)$ and $H^n(\mathcal{G}; A)$.

**Theorem 3.2.3.** Let $\mathcal{G}$ be an ample étale groupoid and let $A$ be a discrete abelian group. Write $H_n(\mathcal{G}) := H_n(\mathcal{G}; \mathbb{Z}) = H_n(C_c(\mathcal{G}_\bullet, \mathbb{Z}), \partial_\bullet)$. Then for each $n \geq 0$ there is a natural short exact sequence of abelian groups

$$0 \to H_n(\mathcal{G}) \otimes_{\mathbb{Z}} A \xrightarrow{\iota_n^{\mathcal{G}}} H_n(\mathcal{G}; A) \xrightarrow{\kappa_n^{\mathcal{G}}} \mathrm{Tor}_1^{\mathbb{Z}}(H_{n-1}(\mathcal{G}), A) \to 0.$$

*Proof.* For each $n \geq 0$, ampleness of $\mathcal{G}$ implies that $\mathcal{G}_n$ is locally compact, Hausdorff, totally disconnected, and has a basis of compact open subsets. By Lemma 2.4.9, every $f \in C_c(\mathcal{G}_n, \mathbb{Z})$ is a finite $\mathbb{Z}$-linear combination of characteristic functions $\chi_U$ of compact open subsets $U \subseteq \mathcal{G}_n$. In particular, $C_c(\mathcal{G}_n, \mathbb{Z})$ is a free abelian group, hence $C_c(\mathcal{G}_\bullet, \mathbb{Z})$ is a chain complex of free $\mathbb{Z}$-modules. Since $A$ is discrete, Proposition 3.2.1 provides a natural chain isomorphism

$$\Phi_\bullet \colon (C_c(\mathcal{G}_\bullet, \mathbb{Z}) \otimes_{\mathbb{Z}} A,\ \partial_\bullet \otimes \mathrm{id}_A) \xrightarrow{\cong} (C_c(\mathcal{G}_\bullet, A),\ \partial_\bullet^A), \qquad \Phi_n(f \otimes a) = [g \mapsto f(g) \cdot a].$$

Hence $\alpha_n^{\mathcal{G}} := H_n(\Phi_\bullet)$ is a natural isomorphism $\alpha_n^{\mathcal{G}} \colon H_n(C_c(\mathcal{G}_\bullet, \mathbb{Z}) \otimes_{\mathbb{Z}} A) \xrightarrow{\cong} H_n(\mathcal{G}; A)$.

Now apply the classical homological universal coefficient theorem for chain complexes of free abelian groups to the chain complex $C_c(\mathcal{G}_\bullet, \mathbb{Z})$ and the $\mathbb{Z}$-module $A$: for each $n \geq 0$ there is a natural short exact sequence

$$0 \to H_n(\mathcal{G}) \otimes_{\mathbb{Z}} A \xrightarrow{\iota_n} H_n(C_c(\mathcal{G}_\bullet, \mathbb{Z}) \otimes_{\mathbb{Z}} A) \xrightarrow{\kappa_n} \mathrm{Tor}_1^{\mathbb{Z}}(H_{n-1}(\mathcal{G}), A) \to 0,$$

see [22, Theorem 3.6.2]. Transporting this sequence along the isomorphism $\alpha_n^{\mathcal{G}}$ yields the claimed short exact sequence, with $\iota_n^{\mathcal{G}} := \alpha_n^{\mathcal{G}} \circ \iota_n, \kappa_n^{\mathcal{G}} := \kappa_n \circ (\alpha_n^{\mathcal{G}})^{-1}$. Naturality in $A$ follows from naturality of the UCT sequence and functoriality of $\Phi_\bullet$ with respect to homomorphisms $A \to B$ via postcomposition on $C_c(\mathcal{G}_n, A)$. $\square$

The UCT for ample étale groupoids established above is, in general, a discrete-coefficients result: its proof relies on the identification

$$C_c(\mathcal{G}_n, A) \cong C_c(\mathcal{G}_n, \mathbb{Z}) \otimes_{\mathbb{Z}} A,$$

which holds whenever every continuous compactly supported $A$-valued function on $\mathcal{G}_n$ is locally constant. This is automatic if $A$ is discrete, but it may also hold for certain non-discrete $A$ when $\mathcal{G}_n$ is discrete. In contrast, for non-discrete $\mathcal{G}_n$ and non-discrete $A$ one should not expect such an identification, and therefore no general UCT for the homology defined via $C_c(\mathcal{G}_\bullet, A)$.

**Corollary 3.2.4.** Let $X$ be a locally compact, totally disconnected Hausdorff space with a basis of compact open sets, and let $A$ be a topological abelian group. Consider the homomorphism

$$\Phi: \ C_c(X, \mathbb{Z}) \otimes_{\mathbb{Z}} A \to C_c(X, A), \qquad f \otimes a \mapsto (x \mapsto f(x) \cdot a),$$

where $f(x) \in \mathbb{Z}$ acts on $a \in A$ by repeated addition. In particular, for compact open $U \subseteq X$, $\Phi(\chi_U \otimes a) = a \cdot \chi_U$. Then the following are equivalent:

1. every $\xi \in C_c(X, A)$ is locally constant,
2. $\Phi$ is surjective,
3. $\Phi$ is an isomorphism of abelian groups.

In particular, if $A$ is discrete, then $\Phi$ is an isomorphism. The converse need not hold.

*Proof.* Define the $\mathbb{Z}$-bilinear map

$$\beta: \ C_c(X, \mathbb{Z}) \times A \to C_c(X, A), \qquad \beta(f, a)(x) := f(x) \cdot a.$$

Since $f$ is $\mathbb{Z}$-valued and locally constant, $\beta(f, a)$ is locally constant for every $a \in A$, hence continuous for any topology on $A$. Moreover $\operatorname{supp}(\beta(f, a)) \subseteq \operatorname{supp}(f)$, so $\beta(f, a)$ is compactly supported. Thus $\beta$ induces a unique group homomorphism

$$\Phi: \ C_c(X, \mathbb{Z}) \otimes_{\mathbb{Z}} A \to C_c(X, A), \qquad \Phi(f \otimes a) = \beta(f, a).$$

- **Every element of** $\operatorname{im}(\Phi)$ **is locally constant.** Write $h = \Phi(\sum_{j=1}^m f_j \otimes a_j) = \sum_{j=1}^m \beta(f_j, a_j)$. Each $f_j : X \to \mathbb{Z}$ is locally constant because $\mathbb{Z}$ is discrete, hence

  $$F: \ X \to \mathbb{Z}^m, \qquad x \mapsto (f_1(x), \dots, f_m(x)),$$

  is locally constant. Define

  $$\lambda: \ \mathbb{Z}^m \to A, \qquad (n_1, \dots, n_m) \mapsto \sum_{j=1}^m n_j \cdot a_j.$$

  Then $h = \lambda \circ F$, so $h$ is locally constant.
- **1.⇒2.** Let $\xi \in C_c(X, A)$. Assume $\xi$ is locally constant. Its support $\operatorname{supp}(\xi) = \{x \in X \mid \xi(x) \neq 0\}$ is compact. For each $x \in \operatorname{supp}(\xi)$ choose an open neighbourhood $O_x$ on which $\xi$ is constant. Using the basis of compact open sets, choose compact open $U_x \subseteq O_x$ with $x \in U_x$. By compactness, pick $x_1, \dots, x_N \in \operatorname{supp}(\xi)$ such that $\operatorname{supp}(\xi) \subseteq U_{x_1} \cup \dots \cup U_{x_N}$. Set $K := \bigcup_{i=1}^N U_{x_i}$, which is compact open.
  Since $X$ is Hausdorff, every compact open set is clopen. Hence we can refine the finite cover $(U_{x_i})_{i=1}^N$ to a finite family of pairwise disjoint compact open sets $(V_\ell)_{\ell=1}^L$ with $\bigcup_{\ell=1}^L V_\ell = K$.

Because $\xi$ is constant on each $U_{x_i}$, it is constant on each $V_\ell$; write $\xi|_{V_\ell} \equiv a_\ell \in A$. Then $\xi = 0$ on $X \setminus K$, and therefore

$$\xi = \sum_{\ell=1}^{L} a_\ell \chi_{V_\ell} = \Phi\Big(\sum_{\ell=1}^{L} \chi_{V_\ell} \otimes a_\ell\Big).$$

Thus $\Phi$ is surjective.

- **2.⇒1.** If $\Phi$ is surjective, then every $\xi \in C_c(X, A)$ lies in $\mathrm{im}(\Phi)$, hence is locally constant.
- **2.⇒3.** Assume $\Phi$ is surjective. By 2.⇒1., every $\xi \in C_c(X, A)$ is locally constant. As in 1.⇒2., every $\xi \in C_c(X, A)$ admits a decomposition

$$\xi = \sum_{i=1}^{n} a_i \chi_{U_i},$$

where $U_1, \dots, U_n \subseteq X$ are pairwise disjoint compact open sets and $a_i \in A$. Define

$$\Psi : C_c(X, A) \to C_c(X, \mathbb{Z}) \otimes_{\mathbb{Z}} A, \qquad \Psi\Big(\sum_{i=1}^{n} a_i \chi_{U_i}\Big) := \sum_{i=1}^{n} \chi_{U_i} \otimes a_i.$$

  - **Well-definedness of $\Psi$.** Suppose also $\xi = \sum_{j=1}^{m} b_j \chi_{V_j}$ with pairwise disjoint compact open $V_j$. Put $W_{ij} := U_i \cap V_j$, which are pairwise disjoint compact open sets, and note that in $C_c(X, \mathbb{Z})$,

$$\chi_{U_i} = \sum_{j=1}^{m} \chi_{W_{ij}}, \qquad \chi_{V_j} = \sum_{i=1}^{n} \chi_{W_{ij}}.$$

    On each $W_{ij} \neq \emptyset$, the function $\xi$ is constant with value $a_i = b_j$, hence

$$\sum_{i=1}^{n} \chi_{U_i} \otimes a_i = \sum_{i,j} \chi_{W_{ij}} \otimes a_i = \sum_{i,j} \chi_{W_{ij}} \otimes b_j = \sum_{j=1}^{m} \chi_{V_j} \otimes b_j.$$

    Thus $\Psi$ is well-defined.
  - **Inverse property.** For $\xi = \sum_i a_i \chi_{U_i}$,

$$\Phi(\Psi(\xi)) = \sum_i \Phi(\chi_{U_i} \otimes a_i) = \sum_i a_i \chi_{U_i} = \xi,$$

    so $\Phi \circ \Psi = \mathrm{id}_{C_c(X,A)}$.

Conversely, let $f \in C_c(X, \mathbb{Z})$. Since $X$ is totally disconnected and $f$ is locally constant with compact support, there exists a finite family of pairwise disjoint compact open subsets $U_1, \dots, U_m \subseteq X$ and integers $b_1, \dots, b_m \in \mathbb{Z}$ such that

$$f = \sum_{j=1}^{m} b_j \chi_{U_j}.$$

Indeed, for each $x \in \mathrm{supp}(f)$ choose a compact open neighbourhood $V_x$ on which $f$ is constant. Compactness of $\mathrm{supp}(f)$ yields finitely many such sets $V_{x_1}, \dots, V_{x_m}$ covering $\mathrm{supp}(f)$. Replace this cover by the disjoint refinement given by successive differences, and read off the constants on each piece. Hence every element of $C_c(X, \mathbb{Z}) \otimes_{\mathbb{Z}} A$ is a finite sum of pure tensors $\chi_U \otimes a$ with $U$ compact open. Hence elements of $C_c(X, \mathbb{Z}) \otimes_{\mathbb{Z}} A$ are finite sums of tensors $\chi_U \otimes a$ with $U$ compact open. For such a tensor,

$$\Psi(\Phi(\chi_U \otimes a)) = \Psi(a\chi_U) = \chi_U \otimes a.$$

By additivity, $\Psi \circ \Phi = \mathrm{id}_{C_c(X,\mathbb{Z})\otimes_{\mathbb{Z}} A}$. Thus $\Phi$ is an isomorphism with inverse $\Psi$.

- **(3)⇒(2).** Trivial.

Finally, if $A$ is discrete then every continuous $X \to A$ is locally constant, so 1. holds. □

**Example 3.2.5.**

- **Non-discrete $A$, still $\Phi$ an isomorphism because $X$ is discrete.** Let $X := \mathbb{N}$ with the discrete topology and let $A := (\mathbb{R}, +)$ with its usual topology. Then every $\xi \in C_c(X, A)$ has finite support because compact subsets of $X$ are finite. Hence $\xi$ is locally constant. By Corollary 3.2.4, $\Phi$ is an isomorphism although $A$ is not discrete.
- **Non-discrete $X$ and non-discrete $A$, $\Phi$ not surjective.** Let $X := \{0, 1\}^{\mathbb{N}}$ be the Cantor space with the product topology, and let $A := (\mathbb{R}, +)$ with its usual topology. For $n \in \mathbb{N}$ and $\varepsilon = (\varepsilon_1, \dots, \varepsilon_n) \in \{0, 1\}^n$ set

$$Z(\varepsilon) := \{x = (x_k)_{k\geq 1} \in X \mid x_1 = \varepsilon_1, \dots, x_n = \varepsilon_n\} = \prod_{k=1}^{n} \{\varepsilon_k\} \times \prod_{k=n+1}^{\infty} \{0, 1\}.$$

Then $\mathcal{B} := \{Z(\varepsilon) \mid n \in \mathbb{N},\ \varepsilon \in \{0, 1\}^n\} \cup \{\emptyset\}$ is a basis of compact open subsets of $X$. Define

$$\xi : X \to \mathbb{R}, \qquad \xi(x) := \sum_{n=1}^{\infty} 2^{-n} x_n.$$

Then $\xi \in C_c(X, \mathbb{R})$ but $\xi$ is not locally constant. Consequently, condition 1. of Corollary 3.2.4 fails for $(X, A)$, and hence the canonical map

$$\Phi : C_c(X, \mathbb{Z}) \otimes_{\mathbb{Z}} \mathbb{R} \to C_c(X, \mathbb{R}), \qquad \chi_U \otimes a \mapsto a\chi_U,$$

is not surjective. Since $X$ is compact, every continuous map $X \to \mathbb{R}$ has compact support equal to $X$, so it suffices to show that $\xi$ is continuous and not locally constant. For $N \in \mathbb{N}$ define the partial sums

$$\xi_N(x) := \sum_{n=1}^{N} 2^{-n} x_n.$$

Each $\xi_N$ depends only on the first $N$ coordinates, hence is continuous in the product topology. Moreover, for all $x \in X$,

$$|\xi(x) - \xi_N(x)| = \sum_{n>N} 2^{-n} x_n \leq \sum_{n>N} 2^{-n} = 2^{-N}.$$

Thus $\xi_N \to \xi$ uniformly on $X$, and $\xi$ is continuous. Hence $\xi \in C(X, \mathbb{R}) = C_c(X, \mathbb{R})$.

Fix $x = (x_n)_{n\geq 1} \in X$ and let $U$ be any neighbourhood of $x$. Since $\mathcal{B}$ is a basis, there exist $N \in \mathbb{N}$ and $\varepsilon = (x_1, \dots, x_N) \in \{0,1\}^N$ such that $x \in Z(\varepsilon) \subseteq U$. Define $y, z \in X$ by

$$y_n := x_n \text{ for } 1 \leq n \leq N, \quad y_{N+1} := 0, \quad y_n := 0 \text{ for } n \geq N+2,$$
$$z_n := x_n \text{ for } 1 \leq n \leq N, \quad z_{N+1} := 1, \quad z_n := 0 \text{ for } n \geq N+2.$$

Then $y, z \in Z(\varepsilon) \subseteq U$, and

$$\xi(z) - \xi(y) = 2^{-(N+1)}(z_{N+1} - y_{N+1}) = 2^{-(N+1)} \neq 0,$$

so $\xi(z) \neq \xi(y)$. Hence $\xi$ is not constant on $U$. Since $U$ was arbitrary, $\xi$ is not locally constant. Every element of $\operatorname{im}(\Phi)$ is locally constant by Corollary 3.2.4. Since $\xi \in C_c(X, \mathbb{R})$ is not locally constant, $\xi \notin \operatorname{im}(\Phi)$. Thus $\Phi$ is not surjective.

In our coefficient discussion the homomorphism

$$\Phi_X : C_c(X, \mathbb{Z}) \otimes_{\mathbb{Z}} A \to C_c(X, A), \qquad \chi_U \otimes a \mapsto a\chi_U,$$

is the conceptual bridge that would allow us to transport the classical algebraic universal coefficient theorem to the chain complexes arising from locally compact, totally disconnected spaces. Indeed, if $\Phi_X$ were an isomorphism, then homology with coefficients in $A$ could be computed purely algebraically from the free $\mathbb{Z}$-complex $C_c(X, \mathbb{Z})$ by tensoring with $A$.

The point of the next results is that this bridge typically collapses as soon as $A$ carries a non-discrete topology, and the obstruction is already visible on the level of values. Every element of $C_c(X, \mathbb{Z}) \otimes_{\mathbb{Z}} A$ is a finite sum of pure tensors. Expanding each $f \in C_c(X, \mathbb{Z})$ as a finite $\mathbb{Z}$-linear combination of characteristic functions of compact open sets, one may assume that an arbitrary tensor has the form

$$z = \sum_{j=1}^{m} \chi_{K_j} \otimes b_j, \qquad K_j \subseteq X \text{ compact open, } b_j \in A.$$

Then

$$\Phi_X(z) = \sum_{j=1}^{m} b_j \chi_{K_j}, \qquad \Phi_X(z)(x) = \sum_{\{j \,:\, x \in K_j\}} b_j,$$

so $\Phi_X(z)(x)$ is a sum over a subset of the finite set $\{b_1, \dots, b_m\}$. Consequently, every function in $\operatorname{im}(\Phi_X)$ has finite image, contained in

$$\Big\{ \sum_{j\in J} b_j \;\Big|\; J \subseteq \{1, \dots, m\} \Big\} \subseteq A.$$

This finiteness constraint is purely algebraic. It does not depend on finer topological properties of $A$ such as total disconnectedness or first countability.

In contrast, once $A$ admits a nontrivial sequence $a_n \to 0$, any non-discrete zero-dimensional space $X$ supports continuous compactly supported $A$-valued functions with infinitely many distinct values. One concentrates the $a_n$ on a clopen partition shrinking to a non-isolated point and uses $a_n \to 0$ to enforce continuity at the limit point. Such a function cannot lie in $\mathrm{im}(\Phi_X)$ by the finite-image obstruction, and therefore $\Phi_X$ fails to be surjective. The Cantor-space example makes this failure completely explicit and shows that even for the most regular $X$, non-discrete coefficients already break the tensor-product model. This is exactly why the UCT for ample étale groupoids obtained below is intrinsically a discrete-coefficients statement when homology is defined using compactly supported continuous $A$-valued chains.

**Lemma 3.2.6.** Let $X$ be a locally compact Hausdorff space and let $A$ be an abelian group. Consider the homomorphism

$$\Phi_X : C_c(X, \mathbb{Z}) \otimes_{\mathbb{Z}} A \to C_c(X, A), \qquad f \otimes a \mapsto (x \mapsto f(x) \cdot a),$$

where $f(x) \in \mathbb{Z}$ acts on $a \in A$ by repeated addition. Then $\xi \in \mathrm{im}(\Phi_X)$ has finite image $\xi(X)$.

*Proof.* Let $z = \sum_{j=1}^{m} f_j \otimes a_j \in C_c(X, \mathbb{Z}) \otimes_{\mathbb{Z}} A$. For each $x \in X$ we have

$$\Phi_X(z)(x) = \sum_{j=1}^{m} f_j(x) \cdot a_j.$$

Fix $j$. Since $f_j \in C_c(X, \mathbb{Z})$, the support $\mathrm{supp}(f_j)$ is compact and $f_j$ is continuous into the discrete space $\mathbb{Z}$. Hence $f_j(\mathrm{supp}(f_j))$ is compact in $\mathbb{Z}$, therefore finite. Because $f_j$ vanishes on $X \setminus \mathrm{supp}(f_j)$, it follows that $f_j(X)$ is finite. Consequently,

$$\Phi_X(z)(X) \subseteq \Big\{ \sum_{j=1}^{m} n_j \cdot a_j \;\Big|\; n_j \in f_j(X) \Big\},$$

and the right-hand side is finite since each $f_j(X)$ is finite. Thus $\Phi_X(z)$ has finite image. $\square$

**Corollary 3.2.7.** Let $X$ be a locally compact, totally disconnected Hausdorff space with a basis of compact open sets. Assume that $X$ is non-discrete. Assume that there exists a non-isolated point $x \in X$ and a decreasing sequence $(U_n)_{n \geq 1}$ of compact open neighbourhoods of $x$ such that

$$U_{n+1} \subseteq U_n \quad \text{for all } n \geq 1, \qquad \bigcap_{n \geq 1} U_n = \{x\}.$$

Let $A$ be a topological abelian group. Assume that there exists a sequence $(a_n)_{n \geq 1}$ in $A \setminus \{0\}$ with $a_n \to 0$. Then the homomorphism

$$\Phi_X : C_c(X, \mathbb{Z}) \otimes_{\mathbb{Z}} A \to C_c(X, A), \qquad \chi_U \otimes a \mapsto a \cdot \chi_U,$$

is not surjective.

*Proof.* For $n \geq 1$ set $V_n := U_n \setminus U_{n+1}$. Since each $U_n$ is compact open, it is clopen, hence each $V_n$ is compact open. The sets $(V_n)_{n \geq 1}$ are pairwise disjoint and satisfy $U_1 \setminus \{x\} = \bigsqcup_{n \geq 1} V_n$, $X \setminus U_1$ is open.

Define $\xi : X \to A$ by

$$\xi(y) := \begin{cases} a_n, & y \in V_n \text{ for some } n \geq 1, \\ 0, & y \notin \bigsqcup_{n\geq 1} V_n. \end{cases}$$

Then $\xi$ is continuous on $X \setminus \{x\}$ since it is constant on each clopen set $V_n$ and on $X \setminus U_1$. It remains to check continuity at $x$. Let $W \subseteq A$ be an open neighbourhood of 0. Choose $n_0 \geq 1$ such that $a_n \in W$ for all $n \geq n_0$. We claim that $\xi(U_{n_0}) \subseteq W$. Indeed, if $y \in U_{n_0}$ and $y \neq x$, then $y \in U_1 \setminus \{x\} = \bigsqcup_{n\geq 1} V_n$, hence $y \in V_n$ for a unique $n \geq 1$. Since $y \in U_{n_0}$ and $V_n \subseteq U_n$, this forces $n \geq n_0$, and therefore $\xi(y) = a_n \in W$. Also $\xi(x) = 0 \in W$. Thus $U_{n_0} \subseteq \xi^{-1}(W)$ and $\xi$ is continuous at $x$. Moreover, $\xi$ has compact support since $\xi = 0$ on $X \setminus U_1$, hence $\operatorname{supp}(\xi) \subseteq U_1$ and $U_1$ is compact. Finally, $\xi(X)$ contains the infinite set $\{a_n \mid n \geq 1\}$, hence $\xi$ has infinite image. By Lemma 3.2.6, every element in $\operatorname{im}(\Phi_X)$ has finite image. Therefore $\xi \notin \operatorname{im}(\Phi_X)$. □

*Remark* 3.2.8.

1. **On the space** $X$**.** Corollary 3.2.7 assumes the existence of a non-isolated point $x \in X$ admitting a countable neighbourhood basis consisting of compact open sets. This holds, for example, if $X$ is non-discrete, locally compact, Hausdorff, totally disconnected, and second countable. Indeed, in a second countable locally compact Hausdorff space every point has a countable neighbourhood basis. If, moreover, $X$ is totally disconnected, then every neighbourhood of $x$ contains a compact open neighbourhood of $x$, so the neighbourhood basis can be chosen to consist of compact open sets.
2. **On the coefficient group** $A$**.** The hypothesis on $A$ used in Corollary 3.2.7 is the existence of a sequence $(a_n)_{n\in\mathbb{N}} \subset A \setminus \{0\}$ with $a_n \to 0$. This is strictly weaker than first countability at 0. It is also stronger than the purely closure-theoretic statement $0 \in \overline{A \setminus \{0\}}$, and stronger than the existence of a countable subset $S \subseteq A \setminus \{0\}$ with $0 \in \overline{S}$. In general, even if such a countable $S$ exists, it need not contain a sequence converging to 0.
   A sufficient topological condition ensuring equivalence is the Fréchet–Urysohn property at 0: if $A$ is Fréchet–Urysohn at 0 then $0 \in \overline{A \setminus \{0\}}$ implies the existence of a sequence in $A \setminus \{0\}$ converging to 0. In particular, if $A$ is first countable at 0, then 0 is non-isolated if and only if there exists a sequence in $A \setminus \{0\}$ converging to 0.

**Corollary 3.2.9.** The short exact sequence of Theorem 3.2.3

$$0 \to H_n(\mathcal{G}) \otimes_{\mathbb{Z}} A \xrightarrow{\iota_n^{\mathcal{G}}} H_n(\mathcal{G}; A) \xrightarrow{\kappa_n^{\mathcal{G}}} \operatorname{Tor}_1^{\mathbb{Z}}(H_{n-1}(\mathcal{G}), A) \to 0$$

splits in **Ab**. In general, such a splitting is not natural and hence not canonical.

*Proof.* By Proposition 3.2.1 there is a chain isomorphism $\Phi_\bullet : C_c(\mathcal{G}_\bullet, \mathbb{Z}) \otimes_{\mathbb{Z}} A \to C_c(\mathcal{G}_\bullet, A)$, hence an isomorphism $H_n(\Phi_\bullet) : H_n(C_c(\mathcal{G}_\bullet, \mathbb{Z}) \otimes_{\mathbb{Z}} A) \xrightarrow{\cong} H_n(\mathcal{G}; A)$. The algebraic UCT for the chain complex of free abelian groups $C_c(\mathcal{G}_\bullet, \mathbb{Z})$ yields a split short exact sequence

$$0 \to H_n(\mathcal{G}, \mathbb{Z}) \otimes_{\mathbb{Z}} A \xrightarrow{\iota_n} H_n(C_c(\mathcal{G}_\bullet, \mathbb{Z}) \otimes_{\mathbb{Z}} A) \xrightarrow{\kappa_n} \operatorname{Tor}_1^{\mathbb{Z}}(H_{n-1}(\mathcal{G}, \mathbb{Z}), A) \to 0$$

whose splitting is in general non-natural. Transporting this splitting along $H_n(\Phi_\bullet)$ yields a splitting of the displayed sequence. $\square$

Corollary 3.2.4 shows that the tensor-product model $C_c(X,\mathbb{Z})\otimes_{\mathbb{Z}} A$ captures $C_c(X,A)$ automatically for discrete coefficient groups $A$, but fails for many natural non-discrete coefficients as soon as $X$ supports continuous compactly supported functions that are not locally constant. Thus, for ample étale groupoids $\mathcal{G}$, the universal coefficient theorem in the form of Theorem 3.2.3 is intrinsically a discrete-coefficients statement when homology is defined via the Moore complexes $C_c(\mathcal{G}_\bullet,A)$. In particular, for discrete $A$, the abelian groups $H_n(\mathcal{G};A)$ are controlled by the integral homology $H_\bullet(\mathcal{G})$ through the functors $-\otimes_{\mathbb{Z}} A$ and $\mathrm{Tor}_1^{\mathbb{Z}}(-,A)$. This rigidity is the mechanism by which the UCT short exact sequences can be compared to Matui's long exact sequences arising from short exact sequences of Moore complexes.

**Theorem 3.2.10.** Let $\mathcal{G},\mathcal{G}',\mathcal{G}''$ be ample étale groupoids. Assume that there is a short exact sequence of chain complexes of free abelian groups

$$0\to C_c(\mathcal{G}'_\bullet,\mathbb{Z})\xrightarrow{i_\bullet} C_c(\mathcal{G}_\bullet,\mathbb{Z})\xrightarrow{p_\bullet} C_c(\mathcal{G}''_\bullet,\mathbb{Z})\to 0.$$

Let $A$ be a discrete abelian group. Then for each $n\geq 0$ the universal coefficient sequences of Theorem 3.2.3 fit into a commutative diagram with exact rows

$$
\begin{array}{ccccccccc}
0 & \longrightarrow & H_n(\mathcal{G}')\otimes_{\mathbb{Z}} A & \xrightarrow{\iota_n^{\mathcal{G}'}} & H_n(\mathcal{G}';A) & \xrightarrow{\kappa_n^{\mathcal{G}'}} & \mathrm{Tor}_1^{\mathbb{Z}}(H_{n-1}(\mathcal{G}'),A) & \longrightarrow & 0\\
 & & {\scriptstyle H_n(i)\otimes\mathrm{id}_A}\downarrow & & {\scriptstyle H_n(i)}\downarrow & & {\scriptstyle \mathrm{Tor}_1^{\mathbb{Z}}(H_{n-1}(i))}\downarrow & & \\
0 & \longrightarrow & H_n(\mathcal{G})\otimes_{\mathbb{Z}} A & \xrightarrow{\iota_n^{\mathcal{G}}} & H_n(\mathcal{G};A) & \xrightarrow{\kappa_n^{\mathcal{G}}} & \mathrm{Tor}_1^{\mathbb{Z}}(H_{n-1}(\mathcal{G}),A) & \longrightarrow & 0\\
 & & {\scriptstyle H_n(p)\otimes\mathrm{id}_A}\downarrow & & {\scriptstyle H_n(p)}\downarrow & & {\scriptstyle \mathrm{Tor}_1^{\mathbb{Z}}(H_{n-1}(p))}\downarrow & & \\
0 & \longrightarrow & H_n(\mathcal{G}'')\otimes_{\mathbb{Z}} A & \xrightarrow{\iota_n^{\mathcal{G}''}} & H_n(\mathcal{G}'';A) & \xrightarrow{\kappa_n^{\mathcal{G}''}} & \mathrm{Tor}_1^{\mathbb{Z}}(H_{n-1}(\mathcal{G}''),A) & \longrightarrow & 0.
\end{array}
$$

Here $H_n(i)$ and $H_n(p)$ are induced by $i_\bullet$ and $p_\bullet$ on integral Moore homology. The maps $H_n(i):=H_n(i;A)$ and $H_n(p):=H_n(p;A)$ are induced on homology by the chain maps

$$i_\bullet^A: C_c(\mathcal{G}'_\bullet,A)\to C_c(\mathcal{G}_\bullet,A),\qquad p_\bullet^A: C_c(\mathcal{G}_\bullet,A)\to C_c(\mathcal{G}''_\bullet,A),$$

obtained by tensoring $i_\bullet,p_\bullet$ with $\mathrm{id}_A$ and composing with the chain isomorphisms $\Phi_\bullet$ from Proposition 3.2.1. The right vertical maps are the functorial maps induced on $\mathrm{Tor}_1^{\mathbb{Z}}(-,A)$.

*Proof.* Since $\mathcal{G},\mathcal{G}',\mathcal{G}''$ are ample, Lemma 2.4.9 implies that each group $C_c(\mathcal{G}_n,\mathbb{Z})$, $C_c(\mathcal{G}'_n,\mathbb{Z})$, $C_c(\mathcal{G}''_n,\mathbb{Z})$ is free abelian, hence the algebraic UCT applies to these chain complexes. Applying the algebraic UCT to the chain maps $i_\bullet$ and $p_\bullet$ yields a commutative diagram of short exact sequences for the integral complexes.

Next, Proposition 3.2.1 provides chain isomorphisms

$$\Phi_\bullet^{\mathcal{G}'} : C_c(\mathcal{G}'_\bullet, \mathbb{Z}) \otimes_{\mathbb{Z}} A \xrightarrow{\cong} C_c(\mathcal{G}'_\bullet, A),$$
$$\Phi_\bullet^{\mathcal{G}} : C_c(\mathcal{G}_\bullet, \mathbb{Z}) \otimes_{\mathbb{Z}} A \xrightarrow{\cong} C_c(\mathcal{G}_\bullet, A),$$
$$\Phi_\bullet^{\mathcal{G}''} : C_c(\mathcal{G}''_\bullet, \mathbb{Z}) \otimes_{\mathbb{Z}} A \xrightarrow{\cong} C_c(\mathcal{G}''_\bullet, A).$$

By construction these satisfy the chain-level commutation relations

$$\Phi_\bullet^{\mathcal{G}} \circ (i_\bullet \otimes \mathrm{id}_A) = i_\bullet^A \circ \Phi_\bullet^{\mathcal{G}'}, \qquad \Phi_\bullet^{\mathcal{G}''} \circ (p_\bullet \otimes \mathrm{id}_A) = p_\bullet^A \circ \Phi_\bullet^{\mathcal{G}}.$$

Passing to homology identifies the middle terms in the algebraic UCT diagram with $H_n(\mathcal{G}'; A), H_n(\mathcal{G}; A), H_n(\mathcal{G}''; A)$ and identifies the induced maps with $H_n(i; A)$ and $H_n(p; A)$. The remaining vertical maps are those induced functorially on tensor products and on $\mathrm{Tor}_1^{\mathbb{Z}}(-, A)$. This yields the displayed commutative $3 \times 3$ diagram with exact rows. □

We have thus obtained a homological universal coefficient theorem for ample étale groupoids with discrete coefficients. In particular, for a discrete abelian group $A$ the short exact sequence

$$0 \to H_n(\mathcal{G}) \otimes_{\mathbb{Z}} A \xrightarrow{\iota_n^{\mathcal{G}}} H_n(\mathcal{G}; A) \xrightarrow{\kappa_n^{\mathcal{G}}} \mathrm{Tor}_1^{\mathbb{Z}}(H_{n-1}(\mathcal{G}), A) \to 0$$

shows that $H_n(\mathcal{G}; A)$ is determined, up to the Tor-term, by the integral homology groups $H_\bullet(\mathcal{G})$.

### 3.2.2 UCT for Cohomology

On the cohomological side, coefficients enter via the contravariant functor $\mathrm{Hom}_{\mathbb{Z}}(-, A)$ applied to the integral Moore chain complex $C_c(\mathcal{G}_\bullet, \mathbb{Z})$. This yields the classical cohomological universal coefficient short exact sequence with $\mathrm{Ext}_{\mathbb{Z}}^1$ and $\mathrm{Hom}_{\mathbb{Z}}$ in place of $\mathrm{Tor}_1^{\mathbb{Z}}$ and $\otimes_{\mathbb{Z}}$.

**Theorem 3.2.11** (Cohomological UCT)**.** Assume that each $C_c(\mathcal{G}_n, \mathbb{Z})$ is a free abelian group. Then for every $n \geq 0$ there is a natural short exact sequence of abelian groups

$$0 \to \mathrm{Ext}_{\mathbb{Z}}^1(H_{n-1}(\mathcal{G}), A) \xrightarrow{\kappa_{\mathcal{G}}^n} H^n(\mathcal{G}; A) \xrightarrow{\rho_{\mathcal{G}}^n} \mathrm{Hom}_{\mathbb{Z}}(H_n(\mathcal{G}), A) \to 0,$$

where $H_{-1}(\mathcal{G}) = 0$ and $H^{-1}(\mathcal{G}; A) = 0$ by convention. The sequence splits in general, but the splitting is not canonical.

*Proof.* Write $C_n := C_c(\mathcal{G}_n, \mathbb{Z})$ with differential $\partial_n$. Set $Z_n := \ker(\partial_n)$ and $B_n := \mathrm{im}(\partial_{n+1})$, so $H_n(\mathcal{G}) = Z_n / B_n$. Let $C^n(\mathcal{G}; A) := \mathrm{Hom}_{\mathbb{Z}}(C_n, A)$ and define the coboundary by $\delta^n(\xi) := \xi \circ \partial_{n+1}$. Then $\delta^{n+1} \circ \delta^n = 0$ since $\partial_{n+1} \circ \partial_{n+2} = 0$, and $H^n(\mathcal{G}; A) = H^n(C^\bullet(\mathcal{G}; A), \delta^\bullet)$.

- **Definition of $\rho^n_{\mathcal{G}}$.** Let $[\xi] \in H^n(\mathcal{G}; A)$ be represented by a cocycle $\xi \in \mathrm{Hom}_{\mathbb{Z}}(C_n, A)$, so $\xi \circ \partial_{n+1} = 0$. Then $\xi$ vanishes on $B_n = \mathrm{im}(\partial_{n+1})$, hence factors uniquely through $Z_n/B_n = H_n(\mathcal{G})$. Define $\rho^n_{\mathcal{G}}([\xi]) \in \mathrm{Hom}_{\mathbb{Z}}(H_n(\mathcal{G}), A)$ by

$$\rho^n_{\mathcal{G}}([\xi])\colon z + B_n \mapsto \xi(z) \qquad \text{for } z \in Z_n.$$

This is well defined and depends only on the cohomology class $[\xi]$.

- **Surjectivity of $\rho^n_{\mathcal{G}}$.** Since $C_{n-1}$ is free, the subgroup $B_{n-1} \subseteq C_{n-1}$ is free. Hence $B_{n-1}$ is projective and the short exact sequence

$$0 \to Z_n \to C_n \xrightarrow{\partial_n} B_{n-1} \to 0$$

splits. Fix a splitting $C_n = Z_n \oplus S_n$ such that $\partial_n|_{S_n}\colon S_n \xrightarrow{\cong} B_{n-1}$. Let $\phi \in \mathrm{Hom}_{\mathbb{Z}}(H_n(\mathcal{G}), A)$. Let $\widetilde{\phi}\colon Z_n \to A$ be the composite $Z_n \to Z_n/B_n = H_n(\mathcal{G}) \xrightarrow{\phi} A$, and extend $\widetilde{\phi}$ to $\xi \in \mathrm{Hom}_{\mathbb{Z}}(C_n, A)$ by $\xi|_{S_n} = 0$. Then $\delta^n(\xi) = \xi \circ \partial_{n+1} = 0$ because $\partial_{n+1}(C_{n+1}) = B_n \subseteq Z_n$ and $\widetilde{\phi}(B_n) = 0$. Thus $[\xi] \in H^n(\mathcal{G}; A)$ and $\rho^n_{\mathcal{G}}([\xi]) = \phi$, proving surjectivity.

- **Identification of $\ker(\rho^n_{\mathcal{G}})$.** A class $[\xi] \in H^n(\mathcal{G}; A)$ lies in $\ker(\rho^n_{\mathcal{G}})$ if and only if it has a cocycle representative $\xi$ with $\xi|_{Z_n} = 0$. With the fixed splitting $C_n = Z_n \oplus S_n$, such a cocycle is uniquely determined by its restriction $\xi|_{S_n} \in \mathrm{Hom}_{\mathbb{Z}}(S_n, A)$. Using $\partial_n|_{S_n}\colon S_n \xrightarrow{\cong} B_{n-1}$, we identify $\mathrm{Hom}_{\mathbb{Z}}(S_n, A) \cong \mathrm{Hom}_{\mathbb{Z}}(B_{n-1}, A)$. If we change $\xi$ by a coboundary $\delta^{n-1}(\eta) = \eta \circ \partial_n$, then under this identification the restriction $\xi|_{S_n}$ changes by $\eta|_{B_{n-1}}$ viewed as an element of $\mathrm{Hom}_{\mathbb{Z}}(B_{n-1}, A)$. Therefore

$$\ker(\rho^n_{\mathcal{G}}) \cong \mathrm{Hom}_{\mathbb{Z}}(B_{n-1}, A) \big/ \mathrm{im}\big(\mathrm{Hom}_{\mathbb{Z}}(C_{n-1}, A) \to \mathrm{Hom}_{\mathbb{Z}}(B_{n-1}, A)\big).$$

Apply $\mathrm{Hom}_{\mathbb{Z}}(-, A)$ to the short exact sequence $0 \to B_{n-1} \to C_{n-1} \to C_{n-1}/B_{n-1} \to 0$. Since $C_{n-1}$ is free, $\mathrm{Ext}^1_{\mathbb{Z}}(C_{n-1}, A) = 0$, and the long exact Hom–Ext sequence identifies the above quotient with $\mathrm{Ext}^1_{\mathbb{Z}}(C_{n-1}/B_{n-1}, A)$. Next, since $C_{n-2}$ is free, the subgroup $B_{n-2} \subseteq C_{n-2}$ is free. Hence the short exact sequence $0 \to Z_{n-1} \to C_{n-1} \to B_{n-2} \to 0$ splits. Quotienting by $B_{n-1} \subseteq Z_{n-1}$ yields a split short exact sequence

$$0 \to H_{n-1}(\mathcal{G}) = Z_{n-1}/B_{n-1} \to C_{n-1}/B_{n-1} \to B_{n-2} \to 0,$$

hence $C_{n-1}/B_{n-1} \cong H_{n-1}(\mathcal{G}) \oplus B_{n-2}$. Since $B_{n-2}$ is free, $\mathrm{Ext}^1_{\mathbb{Z}}(B_{n-2}, A) = 0$, so

$$\mathrm{Ext}^1_{\mathbb{Z}}(C_{n-1}/B_{n-1}, A) \cong \mathrm{Ext}^1_{\mathbb{Z}}(H_{n-1}(\mathcal{G}), A).$$

Composing the identifications yields a natural injection $\kappa^n_{\mathcal{G}}\colon \mathrm{Ext}^1_{\mathbb{Z}}(H_{n-1}(\mathcal{G}), A) \hookrightarrow H^n(\mathcal{G}; A)$ with image $\ker(\rho^n_{\mathcal{G}})$, proving exactness.

- **Naturality.** Let $f_\bullet\colon C_\bullet \to D_\bullet$ be a chain map between chain complexes of free abelian groups. Then precomposition induces a cochain map $f^\bullet\colon \mathrm{Hom}_{\mathbb{Z}}(D_\bullet, A) \to \mathrm{Hom}_{\mathbb{Z}}(C_\bullet, A)$, and the constructions of $\rho$ and $\kappa$ are functorial with respect to $f_\bullet$. Hence the short exact sequence is natural in $\mathcal{G}$ and in $A$.

- **Splitting.** Choosing splittings $C_n = Z_n \oplus S_n$ for all $n$ yields a splitting of the short exact sequence. These choices are not canonical, hence neither is the resulting splitting.

□

**Corollary 3.2.12** (UCT criteria in the Moore framework)**.** Let $\mathcal{G}$ be an ample étale groupoid and write $C_c(\mathcal{G}_\bullet, \mathbb{Z}) = (C_c(\mathcal{G}_n, \mathbb{Z}), \partial_n)_{n\geq 0}$ for the integral Moore chain complex.

1. **Discreteness obstruction for tensor comparison.** Let $A$ be a topological abelian group. Assume that there exists a sequence $(a_n)_{n\geq 1}$ in $A \setminus \{0\}$ with $a_n \to 0$. Then for the Cantor space $X = \{0,1\}^{\mathbb{N}}$ the canonical map

$$\Phi_X : C_c(X, \mathbb{Z}) \otimes_{\mathbb{Z}} A \to C_c(X, A), \qquad \chi_U \otimes a \mapsto a \cdot \chi_U,$$

   is not surjective. In particular, there is no functorial identification $C_c(\mathcal{G}_n, \mathbb{Z}) \otimes_{\mathbb{Z}} A \cong C_c(\mathcal{G}_n, A)$ for all ample $\mathcal{G}$, and hence no functorial $\otimes$–Tor UCT of the classical form for Moore homology with such non-discrete coefficients.
2. **Algebraic input for UCT.** If $A$ is a discrete abelian group, then the chain-level identification $C_c(\mathcal{G}_n, \mathbb{Z}) \otimes_{\mathbb{Z}} A \cong C_c(\mathcal{G}_n, A)$ holds for all $n$, so the homological UCT of Theorem 3.2.3 applies. Moreover, if each $C_c(\mathcal{G}_n, \mathbb{Z})$ is free, then the cohomological UCT of Theorem 3.2.11 applies as well. In both cases the UCT sequences split in general non-canonically.

*Proof.*

1. This is the failure mechanism established in Corollary 3.2.4 and made explicit in the Cantor-space Example 3.2.5: one constructs $\xi \in C_c(X, A)$ with infinite image using a clopen partition shrinking to a non-isolated point and the convergence $a_n \to 0$. Every element in $\mathrm{im}(\Phi_X)$ has finite image, hence $\xi \notin \mathrm{im}(\Phi_X)$.
2. Discreteness of $A$ is exactly what makes the tensor maps $C_c(\mathcal{G}_n, \mathbb{Z}) \otimes_{\mathbb{Z}} A \to C_c(\mathcal{G}_n, A)$ isomorphisms in each degree, so Theorem 3.2.3 applies. If each $C_c(\mathcal{G}_n, \mathbb{Z})$ is free, then Theorem 3.2.11 applies to the dual complex $\mathrm{Hom}_{\mathbb{Z}}(C_c(\mathcal{G}_\bullet, \mathbb{Z}), A)$.

□

In the Moore-complex approach, the classical $\otimes$–Tor and Ext–Hom universal coefficient sequences rest on two independent inputs.

1. **Topological.** Coefficients must be discrete in order that compactly supported $A$-valued chains are finite sums of characteristic functions with coefficients in $A$, so that the canonical tensor comparison map is available as in Corollary 3.2.4.
2. **Algebraic.** Degreewise freeness of the integral Moore chain complex $C_c(\mathcal{G}_\bullet, \mathbb{Z})$ is the hypothesis that yields the short exact UCT sequences of Theorems 3.2.3 and 3.2.11.

## 3.3 Moore–Mayer–Vietoris Sequence

The Moore–Mayer–Vietoris sequence for groupoid homology is the homological analogue of gluing along a cover. In the ample setting, the gluing data live on the unit space. An admissible cover $U_1, U_2 \subseteq \mathcal{G}_0$ determines the reductions $\mathcal{G}|_{U_1}$, $\mathcal{G}|_{U_2}$, and $\mathcal{G}|_{U_1 \cap U_2}$. The Moore–Mayer–Vietoris long exact sequence expresses $H_\bullet(\mathcal{G}; A)$ in terms of the homology of these three reduced ample groupoids.

**Definition 3.3.1** (Admissible Mayer–Vietoris cover)**.** Let $\mathcal{G}$ be an ample groupoid. An admissible Mayer–Vietoris cover of $\mathcal{G}$ is a pair of clopen subsets $U_1, U_2 \subseteq \mathcal{G}_0$ such that

1. $U_1 \cup U_2 = \mathcal{G}_0$,
2. each $U_i$ is saturated for $i = 1, 2$: for $x \in U_i$ and $y \in \mathcal{G}_0$, if $x \sim_{\mathcal{G}} y$ then $y \in U_i$.

Here $x \sim_{\mathcal{G}} y$ means that there exists $\gamma \in \mathcal{G}$ with $s(\gamma) = y$ and $r(\gamma) = x$. For such $U_1, U_2$ we consider the reductions

$$\begin{aligned} \mathcal{G}|_{U_1} &\coloneqq \{\gamma \in \mathcal{G} \mid r(\gamma) \in U_1,\ s(\gamma) \in U_1\}, \\ \mathcal{G}|_{U_2} &\coloneqq \{\gamma \in \mathcal{G} \mid r(\gamma) \in U_2,\ s(\gamma) \in U_2\}, \\ \mathcal{G}|_{U_1 \cap U_2} &\coloneqq \{\gamma \in \mathcal{G} \mid r(\gamma) \in U_1 \cap U_2,\ s(\gamma) \in U_1 \cap U_2\}, \end{aligned}$$

each endowed with the groupoid structure obtained by restricting the range, source, unit, inverse, and multiplication maps of $\mathcal{G}$.

**Lemma 3.3.2.** Let $\mathcal{G}$ be an ample groupoid and let $U \subseteq \mathcal{G}_0$ be open. Then $\mathcal{G}|_U$ is an ample groupoid with unit space $U$, and the inclusion $\mathcal{G}|_U \hookrightarrow \mathcal{G}$ is an open embedding of topological groupoids. If moreover $U$ is clopen, then $\mathcal{G}|_U$ is clopen in $\mathcal{G}$. In particular, if $\mathcal{G}$ is locally compact and Hausdorff, then $\mathcal{G}|_U$ is locally compact and Hausdorff, and if $\mathcal{G}$ is totally disconnected then $\mathcal{G}|_U$ is totally disconnected.

*Proof.* Since $\mathcal{G}$ is ample, it is étale, locally compact, Hausdorff, and totally disconnected. By definition $\mathcal{G}|_U = r^{-1}(U) \cap s^{-1}(U)$, hence $\mathcal{G}|_U$ is open in $\mathcal{G}$. All structure maps of $\mathcal{G}|_U$ are restrictions of those of $\mathcal{G}$, hence are continuous and satisfy the groupoid axioms. Since $r, s$ are local homeomorphisms and $U$ is open in $\mathcal{G}_0$, the restrictions $r|_{\mathcal{G}|_U}, s|_{\mathcal{G}|_U} \colon \mathcal{G}|_U \to U$ are again local homeomorphisms. Thus $\mathcal{G}|_U$ is étale. As an open subspace of a locally compact Hausdorff totally disconnected space, $\mathcal{G}|_U$ is locally compact, Hausdorff, and totally disconnected. Therefore $\mathcal{G}|_U$ is ample. The inclusion $\mathcal{G}|_U \hookrightarrow \mathcal{G}$ is the inclusion of an open subspace, hence an open embedding of topological groupoids. If $U$ is clopen, then $r^{-1}(U)$ and $s^{-1}(U)$ are clopen in $\mathcal{G}$, so $\mathcal{G}|_U = r^{-1}(U) \cap s^{-1}(U)$ is clopen. $\square$

*Remark* 3.3.3. Let $\mathcal{G}$ be ample and let $U \subseteq \mathcal{G}_0$ be clopen. Then a compact open bisection of $\mathcal{G}|_U$ is precisely a set of the form $B \cap (\mathcal{G}|_U)$ with $B \subseteq \mathcal{G}$ a compact open bisection.

**Definition 3.3.4.** Let $U_1, U_2 \subseteq \mathcal{G}_0$ be an admissible Mayer–Vietoris cover. For $n \geq 0$ write

$$\mathcal{G}_n^{(1)} \coloneqq (\mathcal{G}|_{U_1})_n, \qquad \mathcal{G}_n^{(2)} \coloneqq (\mathcal{G}|_{U_2})_n, \qquad \mathcal{G}_n^{(12)} \coloneqq (\mathcal{G}|_{U_1 \cap U_2})_n$$

for the $n$-simplices in the nerves of $\mathcal{G}|_{U_1}$, $\mathcal{G}|_{U_2}$, and $\mathcal{G}|_{U_1 \cap U_2}$.

**Lemma 3.3.5.** Let $U_1, U_2 \subseteq \mathcal{G}_0$ be an admissible Mayer–Vietoris cover.

Then for every $n \geq 0$ the following hold.

1. The subsets $\mathcal{G}_n^{(1)}$ and $\mathcal{G}_n^{(2)}$ are clopen in $\mathcal{G}_n$, and $\mathcal{G}_n^{(12)} = \mathcal{G}_n^{(1)} \cap \mathcal{G}_n^{(2)}$.
2. One has $\mathcal{G}_n^{(1)} \cup \mathcal{G}_n^{(2)} = \mathcal{G}_n$.

*Proof.* For $n = 0$ the claims are immediate since $\mathcal{G}_0^{(1)} = U_1$, $\mathcal{G}_0^{(2)} = U_2$, and $\mathcal{G}_0^{(12)} = U_1 \cap U_2$.

Assume $n \geq 1$. Define continuous maps

$$r_n \colon \mathcal{G}_n \to \mathcal{G}_0, \quad r_n(g_1, \dots, g_n) := r(g_1),$$
$$s_n \colon \mathcal{G}_n \to \mathcal{G}_0, \quad s_n(g_1, \dots, g_n) := s(g_n).$$

Let $U \subseteq \mathcal{G}_0$ be saturated. We claim that $(\mathcal{G}|_U)_n = r_n^{-1}(U) = s_n^{-1}(U)$ as subsets of $\mathcal{G}_n$. If $(g_1, \dots, g_n) \in (\mathcal{G}|_U)_n$, then $r(g_1) \in U$, hence $(g_1, \dots, g_n) \in r_n^{-1}(U)$. Conversely, if $(g_1, \dots, g_n) \in r_n^{-1}(U)$, then $r(g_1) \in U$. For each $k$, the units $r(\gamma_k)$ and $s(\gamma_k)$ lie in the $\mathcal{G}$-orbit of $r(g_1)$, hence belong to $U$ by saturation. Thus $\gamma_k \in \mathcal{G}|_U$ for all $k$, so $(g_1, \dots, g_n) \in (\mathcal{G}|_U)_n$. The identity with $s_n^{-1}(U)$ is proved similarly. Applying this to $U_1$ and $U_2$ yields

$$\mathcal{G}_n^{(i)} = (\mathcal{G}|_{U_i})_n = r_n^{-1}(U_i) = s_n^{-1}(U_i).$$

Since each $U_i$ is clopen, each $\mathcal{G}_n^{(i)}$ is clopen in $\mathcal{G}_n$. Moreover,

$$\mathcal{G}_n^{(12)} = (\mathcal{G}|_{U_1 \cap U_2})_n = r_n^{-1}(U_1 \cap U_2) = r_n^{-1}(U_1) \cap r_n^{-1}(U_2) = \mathcal{G}_n^{(1)} \cap \mathcal{G}_n^{(2)}.$$

Finally, $\mathcal{G}_n = r_n^{-1}(\mathcal{G}_0) = r_n^{-1}(U_1 \cup U_2) = r_n^{-1}(U_1) \cup r_n^{-1}(U_2) = \mathcal{G}_n^{(1)} \cup \mathcal{G}_n^{(2)}$. $\square$

**Definition 3.3.6.** Let $U_1, U_2 \subseteq \mathcal{G}_0$ be an admissible Mayer–Vietoris cover and let $A$ be a discrete abelian group. For $n \geq 0$ define

$$C_n(\mathcal{G}|_{U_1}; A) := C_c(\mathcal{G}_n^{(1)}, A), \quad C_n(\mathcal{G}|_{U_2}; A) := C_c(\mathcal{G}_n^{(2)}, A), \quad C_n(\mathcal{G}|_{U_1 \cap U_2}; A) := C_c(\mathcal{G}_n^{(12)}, A),$$

where compact support is taken with respect to the subspace topology on $\mathcal{G}_n^{(1)}, \mathcal{G}_n^{(2)}, \mathcal{G}_n^{(12)} \subseteq \mathcal{G}_n$.

For $n \geq 1$ set

$$\partial_n^{A, U_1} := \sum_{j=0}^{n} (-1)^j (d_j)_* \colon C_c(\mathcal{G}_n^{(1)}, A) \to C_c(\mathcal{G}_{n-1}^{(1)}, A),$$
$$\partial_n^{A, U_2} := \sum_{j=0}^{n} (-1)^j (d_j)_* \colon C_c(\mathcal{G}_n^{(2)}, A) \to C_c(\mathcal{G}_{n-1}^{(2)}, A),$$
$$\partial_n^{A, U_1 \cap U_2} := \sum_{j=0}^{n} (-1)^j (d_j)_* \colon C_c(\mathcal{G}_n^{(12)}, A) \to C_c(\mathcal{G}_{n-1}^{(12)}, A),$$

where $d_j \colon \mathcal{G}_n \to \mathcal{G}_{n-1}$ are the face maps of the nerve of $\mathcal{G}$ and we use the restricted maps

$$d_j \colon \mathcal{G}_n^{(1)} \to \mathcal{G}_{n-1}^{(1)}, \qquad d_j \colon \mathcal{G}_n^{(2)} \to \mathcal{G}_{n-1}^{(2)}, \qquad d_j \colon \mathcal{G}_n^{(12)} \to \mathcal{G}_{n-1}^{(12)}.$$

For $n = 0$ set $\partial_0^{A, U_1} = \partial_0^{A, U_2} = \partial_0^{A, U_1 \cap U_2} := 0$. The resulting homology groups are denoted by

$$H_n(\mathcal{G}|_{U_1}; A), \qquad H_n(\mathcal{G}|_{U_2}; A), \qquad H_n(\mathcal{G}|_{U_1 \cap U_2}; A).$$

**Lemma 3.3.7.** With notation as above, the pairs

$$(C_\bullet(\mathcal{G}|_{U_1};A),\partial_\bullet^{A,U_1}),\qquad (C_\bullet(\mathcal{G}|_{U_2};A),\partial_\bullet^{A,U_2}),\qquad (C_\bullet(\mathcal{G}|_{U_1\cap U_2};A),\partial_\bullet^{A,U_1\cap U_2})$$

are chain complexes. Concretely, for all $n\geq 1$ one has

$$\partial_{n-1}^{A,U_1}\circ\partial_n^{A,U_1}=0,\qquad \partial_{n-1}^{A,U_2}\circ\partial_n^{A,U_2}=0,\qquad \partial_{n-1}^{A,U_1\cap U_2}\circ\partial_n^{A,U_1\cap U_2}=0.$$

*Proof.* Since $\mathcal{G}$ is ample, it is étale. By Lemma 3.3.2, each reduction $\mathcal{G}|_{U_1}$, $\mathcal{G}|_{U_2}$, and $\mathcal{G}|_{U_1\cap U_2}$ is ample, hence étale. Its nerve is obtained by applying the nerve construction to that reduced groupoid, and each face map is the restriction of the corresponding face map of $\mathcal{G}$. The simplicial identities among the face maps therefore restrict to the reduced nerves. Since pushforward along local homeomorphisms is functorial and respects composition, the Moore differential computation shows that the alternating sums of the restricted pushforwards square to zero. □

### 3.3.1 Moore–Mayer–Vietoris at Chain Level

The Moore–Mayer–Vietoris construction starts from a purely chain-level observation. A clopen saturated cover $U_1,U_2\subseteq\mathcal{G}_0$ cuts every nerve space $\mathcal{G}_n$ into two clopen pieces

$$\mathcal{G}_n^{(1)}=(\mathcal{G}|_{U_1})_n,\qquad \mathcal{G}_n^{(2)}=(\mathcal{G}|_{U_2})_n,\qquad \mathcal{G}_n=\mathcal{G}_n^{(1)}\cup\mathcal{G}_n^{(2)},\qquad \mathcal{G}_n^{(12)}=\mathcal{G}_n^{(1)}\cap\mathcal{G}_n^{(2)}.$$

Because the pieces are clopen, compactly supported $A$-valued chains on each reduction extend by zero to compactly supported chains on $\mathcal{G}_n$. This makes it possible to regard chains on $\mathcal{G}|_{U_1}$ and $\mathcal{G}|_{U_2}$ as global chains that vanish off the corresponding clopen region. The overlap $\mathcal{G}|_{U_1\cap U_2}$ measures the compatibility of these extensions. A chain on the overlap can be inserted into the two pieces with opposite signs, and hence disappears after gluing. Conversely, if two chains on the pieces glue to zero, then each must vanish off the overlap and the remaining values on the overlap must cancel. This yields, in every degree, the short exact sequence

$$0\to C_n(\mathcal{G}|_{U_1\cap U_2};A)\xrightarrow{\alpha_n}C_n(\mathcal{G}|_{U_1};A)\oplus C_n(\mathcal{G}|_{U_2};A)\xrightarrow{\beta_n}C_n(\mathcal{G};A)\to 0,$$

which is the only algebraic input needed to extract the Moore–Mayer–Vietoris long exact sequence in homology. The rest of the argument is then formal homological algebra, applied to this degreewise short exact sequence of chain complexes.

**Lemma 3.3.8** (Clopen Mayer–Vietoris for compactly supported chains)**.** Let $X$ be a locally compact Hausdorff space and let $A$ be a topological abelian group. Let $X_1,X_2\subseteq X$ be clopen subsets with $X=X_1\cup X_2$ and set $X_{12}:=X_1\cap X_2$. For $i\in\{1,2,12\}$ write $C_c(X_i,A)$ for compactly supported continuous maps on $X_i$. Define group homomorphisms

$$\begin{aligned}\alpha\colon C_c(X_{12},A)\to C_c(X_1,A)\oplus C_c(X_2,A),&\qquad \xi\mapsto(\xi^{(1)},-\xi^{(2)}),\\ \beta\colon C_c(X_1,A)\oplus C_c(X_2,A)\to C_c(X,A),&\qquad (\xi_1,\xi_2)\mapsto\widetilde{\xi_1}+\widetilde{\xi_2},\end{aligned}$$

where $\xi^{(i)} \in C_c(X_i, A)$ denotes extension by zero of $\xi$ along the clopen inclusion $X_{12} \hookrightarrow X_i$, and $\widetilde{\xi_i} \in C_c(X, A)$ denotes extension by zero of $\xi_i$ along the clopen inclusion $X_i \hookrightarrow X$.

Then we have an exact sequence

$$0 \to C_c(X_{12}, A) \xrightarrow{\alpha} C_c(X_1, A) \oplus C_c(X_2, A) \xrightarrow{\beta} C_c(X, A) \to 0.$$

*Proof.* All extension-by-zero maps are well defined because the relevant inclusions are clopen.

- **Injectivity of $\alpha$.** If $\alpha(\xi) = 0$, then $\xi^{(1)} = 0$ in $C_c(X_1, A)$. Evaluating on $x \in X_{12} \subseteq X_1$ gives $\xi(x) = \xi^{(1)}(x) = 0$, hence $\xi = 0$.
- **Surjectivity of $\beta$.** Let $\eta \in C_c(X, A)$. Set $\xi_1 := \eta|_{X_1} \in C_c(X_1, A)$ and $\xi_2 := \eta|_{X_2 \setminus X_1} \in C_c(X_2 \setminus X_1, A)$. Since $X_2 \setminus X_1$ is clopen in $X_2$, extending $\xi_2$ by zero to $X_2$ yields an element, still denoted $\xi_2 \in C_c(X_2, A)$. Then pointwise on $X$ one has $\eta = \widetilde{\xi_1} + \widetilde{\xi_2}$, because on $X_1 \setminus X_2$ only $\widetilde{\xi_1}$ contributes, on $X_2 \setminus X_1$ only $\widetilde{\xi_2}$ contributes, and on $X_{12}$ one has $\xi_2 = 0$ by construction. Thus $\beta(\xi_1, \xi_2) = \eta$.
- **Exactness at the middle term.** First, $\beta \circ \alpha = 0$ holds pointwise because both components of $\alpha(\xi)$ extend $\xi$ by zero to $X$ with opposite signs. Conversely, suppose $(\xi_1, \xi_2) \in C_c(X_1, A) \oplus C_c(X_2, A)$ satisfies $\beta(\xi_1, \xi_2) = 0$. Evaluating on $x \in X_1 \setminus X_2$ yields $0 = \widetilde{\xi_1}(x) = \xi_1(x)$, so $\xi_1$ vanishes on $X_1 \setminus X_{12}$. Similarly, evaluating on $x \in X_2 \setminus X_1$ yields $\xi_2(x) = 0$, so $\xi_2$ vanishes on $X_2 \setminus X_{12}$. Thus the restrictions $\xi_i|_{X_{12}} \in C_c(X_{12}, A)$ are defined and satisfy

$$0 = \beta(\xi_1, \xi_2)(x) = \xi_1(x) + \xi_2(x) \quad \text{for all } x \in X_{12},$$

  hence $\xi_2|_{X_{12}} = -\xi_1|_{X_{12}}$. Set $\xi := \xi_1|_{X_{12}} \in C_c(X_{12}, A)$. Then by the vanishing outside $X_{12}$ one has $\xi^{(1)} = \xi_1$ and $\xi^{(2)} = \xi_2$. Therefore $\alpha(\xi) = (\xi_1, \xi_2)$, so $\ker(\beta) \subseteq \operatorname{im}(\alpha)$.

□

**Corollary 3.3.9.** Let $\mathcal{G}$ be an ample groupoid and let $U_1, U_2 \subseteq \mathcal{G}_0$ be a clopen saturated cover in the sense of Definition 3.3.1. Let $A$ be a topological abelian group. For every $n \geq 0$ the sequence

$$0 \to C_n(\mathcal{G}|_{U_1 \cap U_2}; A) \xrightarrow{\alpha_n} C_n(\mathcal{G}|_{U_1}; A) \oplus C_n(\mathcal{G}|_{U_2}; A) \xrightarrow{\beta_n} C_n(\mathcal{G}; A) \to 0$$

is exact, where $\alpha_n, \beta_n$ are the Mayer–Vietoris chain maps from Lemma 3.3.8.

*Proof.* Fix $n \geq 0$ and set

$$X := \mathcal{G}_n, \qquad X_1 := \mathcal{G}_n^{(1)} = (\mathcal{G}|_{U_1})_n, \qquad X_2 := \mathcal{G}_n^{(2)} = (\mathcal{G}|_{U_2})_n, \qquad X_{12} := X_1 \cap X_2 = \mathcal{G}_n^{(12)}.$$

Since $U_1, U_2$ are clopen and saturated, Lemma 3.3.5 implies that $X_1, X_2$ are clopen in $X$, that $X = X_1 \cup X_2$, and that $X_{12} = X_1 \cap X_2$. With the identifications

$$C_n(\mathcal{G}; A) = C_c(X, A), \qquad C_n(\mathcal{G}|_{U_i}; A) = C_c(X_i, A), \qquad C_n(\mathcal{G}|_{U_1 \cap U_2}; A) = C_c(X_{12}, A),$$

the maps $\alpha_n, \beta_n$ coincide with the canonical clopen Mayer–Vietoris maps

$$\begin{aligned} \alpha\colon C_c(X_{12}, A) &\to C_c(X_1, A) \oplus C_c(X_2, A), & \xi &\mapsto (\xi, -\xi), \\ \beta\colon C_c(X_1, A) \oplus C_c(X_2, A) &\to C_c(X, A), & (\xi_1, \xi_2) &\mapsto \widetilde{\xi_1} + \widetilde{\xi_2}, \end{aligned}$$

where $\widetilde{\xi_i}$ denotes extension by zero from $X_i$ to $X$ along the clopen inclusion. Hence the displayed sequence is identified with the short exact sequence of Lemma 3.3.8. $\square$

### 3.3.2 Moore–Mayer–Vietoris Long Exact Homology Sequence

The Mayer–Vietoris principle reconstructs homology from two subobjects and their overlap. In the ample groupoid setting, the gluing data live on the unit space. A clopen saturated cover $\mathcal{G}_0 = U_1 \cup U_2$ yields three reduced groupoids $\mathcal{G}|_{U_1}$, $\mathcal{G}|_{U_2}$, and $\mathcal{G}|_{U_1 \cap U_2}$ and three Moore complexes of compactly supported chains. The basic input is a short exact sequence of Moore complexes, obtained by support decompositions along clopen subsets. Once this short exact sequence is established, the long exact sequence is obtained from the connecting homomorphism construction for short exact sequences of chain complexes of abelian groups.

Hausdorffness of $A$ is used only to ensure that $C_c(-, A)$ is well defined on locally compact Hausdorff spaces. After the Moore complexes are defined, the argument below takes place in the category of abelian groups.

**Theorem 3.3.10** (Moore–Mayer–Vietoris LES)**.** Let $\mathcal{G}$ be an ample groupoid, let $A$ be a topological Hausdorff abelian group, and let $U_1, U_2 \subseteq \mathcal{G}_0$ be clopen saturated subsets with $U_1 \cup U_2 = \mathcal{G}_0$. Then the homology groups $H_n(\mathcal{G}; A) := H_n(C_c(\mathcal{G}_\bullet, A))$ fit into a natural long exact sequence

$$\begin{array}{ccccc} \cdots \longleftarrow & H_{n-1}(\mathcal{G}|_{U_1}; A) \oplus H_{n-1}(\mathcal{G}|_{U_2}; A) & \xleftarrow{H_{n-1}(\alpha_\bullet)} & H_{n-1}(\mathcal{G}|_{U_1 \cap U_2}; A) & \\ & & \partial_n & & \\ H_n(\mathcal{G}; A) \xleftarrow{H_n(\beta_\bullet)} & H_n(\mathcal{G}|_{U_1}; A) \oplus H_n(\mathcal{G}|_{U_2}; A) & \xleftarrow{H_n(\alpha_\bullet)} & H_n(\mathcal{G}|_{U_1 \cap U_2}; A) & \\ & & \partial_{n+1} & & \\ H_{n+1}(\mathcal{G}; A) \xleftarrow{H_{n+1}(\beta_\bullet)} & H_{n+1}(\mathcal{G}|_{U_1}; A) \oplus H_{n+1}(\mathcal{G}|_{U_2}; A) & \longleftarrow \cdots & & \end{array}$$

where $H_n(\alpha_\bullet)$ and $H_n(\beta_\bullet)$ are induced by the chain maps

$$\begin{aligned} \alpha_\bullet &\colon C_\bullet(\mathcal{G}|_{U_1 \cap U_2}; A) \to C_\bullet(\mathcal{G}|_{U_1}; A) \oplus C_\bullet(\mathcal{G}|_{U_2}; A), \\ \beta_\bullet &\colon C_\bullet(\mathcal{G}|_{U_1}; A) \oplus C_\bullet(\mathcal{G}|_{U_2}; A) \to C_c(\mathcal{G}_\bullet, A) \end{aligned}$$

from Definition 3.3.6. The connecting homomorphisms $\partial_n\colon H_n(\mathcal{G}; A) \to H_{n-1}(\mathcal{G}|_{U_1 \cap U_2}; A)$ are defined explicitly in the proof.

*Proof.* By Corollary 3.3.9 there is a short exact sequence of Moore chain complexes

$$0 \to C_\bullet(\mathcal{G}|_{U_1 \cap U_2}; A) \xrightarrow{\alpha_\bullet} C_\bullet(\mathcal{G}|_{U_1}; A) \oplus C_\bullet(\mathcal{G}|_{U_2}; A) \xrightarrow{\beta_\bullet} C_c(\mathcal{G}_\bullet, A) \to 0. \tag{3.3.1}$$

Write

$$C^{12}_\bullet := C_\bullet(\mathcal{G}|_{U_1\cap U_2}; A), \quad C^{1,2}_\bullet := C_\bullet(\mathcal{G}|_{U_1}; A) \oplus C_\bullet(\mathcal{G}|_{U_2}; A), \quad C_\bullet := C_c(\mathcal{G}_\bullet, A),$$
$$\partial^{12}_n \colon C^{12}_n \to C^{12}_{n-1}, \qquad \partial^{1,2}_n := \partial^1_n \oplus \partial^2_n \colon C^{1,2}_n \to C^{1,2}_{n-1}, \qquad \partial^{\mathcal{G}}_n \colon C_n \to C_{n-1}.$$

The chain map identities are $\alpha_{n-1} \circ \partial^{12}_n = \partial^{1,2}_n \circ \alpha_n, \beta_{n-1} \circ \partial^{1,2}_n = \partial^{\mathcal{G}}_n \circ \beta_n$. Exactness of (3.3.1) in each degree yields that $\alpha_n$ is injective and $\beta_n$ is surjective for every $n$.

- **Definition of the connecting homomorphism $\partial_n$.** Fix $n \geq 0$ and let $[c] \in H_n(\mathcal{G}; A)$ be represented by a cycle $c \in C_n$ with $\partial^{\mathcal{G}}_n(c) = 0$. Choose $b \in C^{1,2}_n$ with $\beta_n(b) = c$, which exists since $\beta_n$ is surjective. Then $\beta_{n-1}(\partial^{1,2}_n(b)) = \partial^{\mathcal{G}}_n(\beta_n(b)) = \partial^{\mathcal{G}}_n(c) = 0$, so $\partial^{1,2}_n(b) \in \ker(\beta_{n-1})$. Exactness in degree $n-1$ gives $\ker(\beta_{n-1}) = \operatorname{im}(\alpha_{n-1})$, so there exists a unique $a \in C^{12}_{n-1}$ such that
$$\alpha_{n-1}(a) = \partial^{1,2}_n(b). \tag{3.3.2}$$
We show that $a$ is a cycle. Using that $\alpha_\bullet$ is a chain map and that $\partial^{1,2}_{n-1} \circ \partial^{1,2}_n = 0$, $\alpha_{n-2}(\partial^{12}_{n-1}(a)) = \partial^{1,2}_{n-1}(\alpha_{n-1}(a)) = \partial^{1,2}_{n-1}(\partial^{1,2}_n(b)) = 0$. Injectivity of $\alpha_{n-2}$ implies $\partial^{12}_{n-1}(a) = 0$. Define $\partial_n([c]) := [a] \in H_{n-1}(\mathcal{G}|_{U_1\cap U_2}; A)$.
- **Independence of the choice of $b$.** Let $b, b' \in C^{1,2}_n$ satisfy $\beta_n(b) = \beta_n(b') = c$. Then $b' - b \in \ker(\beta_n) = \operatorname{im}(\alpha_n)$, so choose $u \in C^{12}_n$ with $b' = b + \alpha_n(u)$. Let $a, a' \in C^{12}_{n-1}$ be defined by (3.3.2) for $b$ and $b'$. Then
$$\alpha_{n-1}(a') = \partial^{1,2}_n(b') = \partial^{1,2}_n(b) + \partial^{1,2}_n(\alpha_n(u)) = \alpha_{n-1}(a) + \alpha_{n-1}(\partial^{12}_n(u)),$$
since $\alpha_\bullet$ is a chain map. Injectivity of $\alpha_{n-1}$ yields $a' = a + \partial^{12}_n(u)$, hence $[a'] = [a]$. Thus $\partial_n([c])$ is independent of the lift $b$.
- **Independence of the representative of $[c]$.** Let $c, c' \in C_n$ be cycles with $[c] = [c']$ in $H_n(\mathcal{G}; A)$, so $c' - c = \partial^{\mathcal{G}}_{n+1}(d)$ for some $d \in C_{n+1}$. Choose $b, b' \in C^{1,2}_n$ with $\beta_n(b) = c$ and $\beta_n(b') = c'$. Choose $e \in C^{1,2}_{n+1}$ with $\beta_{n+1}(e) = d$, using surjectivity of $\beta_{n+1}$. Then
$$\beta_n(b' - b - \partial^{1,2}_{n+1}(e)) = c' - c - \partial^{\mathcal{G}}_{n+1}(d) = 0,$$
so $b' - b - \partial^{1,2}_{n+1}(e) \in \ker(\beta_n) = \operatorname{im}(\alpha_n)$. Choose $u \in C^{12}_n$ with $b' = b + \partial^{1,2}_{n+1}(e) + \alpha_n(u)$. Let $a, a' \in C^{12}_{n-1}$ be the elements defined from $b, b'$ by (3.3.2). Then
$$\alpha_{n-1}(a') = \partial^{1,2}_n(b') = \partial^{1,2}_n(b) + \partial^{1,2}_n(\partial^{1,2}_{n+1}(e)) + \partial^{1,2}_n(\alpha_n(u)) = \alpha_{n-1}(a) + \alpha_{n-1}(\partial^{12}_n(u)),$$
so $a' = a + \partial^{12}_n(u)$ and $[a'] = [a]$. Hence $\partial_n$ depends only on $[c]$.
- **Exactness.** The equalities needed for exactness follow from the construction.
  1. $\beta_\bullet \circ \alpha_\bullet = 0$ implies $H_n(\beta_\bullet) \circ H_n(\alpha_\bullet) = 0$ for all $n$.
  2. $\operatorname{im}(\mathbf{H_n}(\boldsymbol{\alpha_\bullet})) = \ker(\mathbf{H_n}(\boldsymbol{\beta_\bullet}))$. Let $[b] \in H_n(C^{1,2}_\bullet)$ with $H_n(\beta_\bullet)([b]) = 0$. Choose a cycle $b \in C^{1,2}_n$ representing $[b]$. Then $\beta_n(b)$ is a boundary in $C_\bullet$, so choose $d \in C_{n+1}$ with $\beta_n(b) = \partial^{\mathcal{G}}_{n+1}(d)$. Choose $e \in C^{1,2}_{n+1}$ with $\beta_{n+1}(e) = d$. Then $b - \partial^{1,2}_{n+1}(e) \in \ker(\beta_n) =$ $\mathrm{im}(\alpha_n)$, so choose $u \in C_n^{12}$ with $b - \partial_{n+1}^{1,2}(e) = \alpha_n(u)$. A direct computation gives $\partial_n^{12}(u) = 0$ since

$$\alpha_{n-1}(\partial_n^{12}(u)) = \partial_n^{1,2}(\alpha_n(u)) = \partial_n^{1,2}(b) - \partial_n^{1,2}(\partial_{n+1}^{1,2}(e)) = 0,$$

and $\alpha_{n-1}$ is injective. Hence $[u] \in H_n(C_\bullet^{12})$ and $H_n(\alpha_\bullet)([u]) = [b]$. The reverse inclusion follows from $H_n(\beta_\bullet) \circ H_n(\alpha_\bullet) = 0$.

3. $\mathrm{im}(\mathbf{H_n}(\boldsymbol{\beta_\bullet})) = \ker(\boldsymbol{\partial_n})$. If $[b] \in H_n(C_\bullet^{1,2})$ is represented by a cycle $b$, then $\partial_n([\beta_n(b)]) = 0$ because one may take this $b$ in the definition and obtain $a = 0$ from (3.3.2). Conversely, let $[c] \in H_n(C_\bullet)$ with $\partial_n([c]) = 0$. Choose $c$ and a lift $b$ with $\beta_n(b) = c$. Let $a$ be defined by (3.3.2). The condition $\partial_n([c]) = 0$ means $[a] = 0$, so choose $u \in C_n^{12}$ with $a = \partial_n^{12}(u)$. Then

$$\partial_n^{1,2}(b - \alpha_n(u)) = \partial_n^{1,2}(b) - \alpha_{n-1}(\partial_n^{12}(u)) = \alpha_{n-1}(a) - \alpha_{n-1}(a) = 0,$$

so $b - \alpha_n(u)$ is a cycle in $C_n^{1,2}$ and $\beta_n(b - \alpha_n(u)) = c$. Thus $[c] \in \mathrm{im}(H_n(\beta_\bullet))$.

4. $\mathrm{im}(\boldsymbol{\partial_n}) = \ker(\mathbf{H_{n-1}}(\boldsymbol{\alpha_\bullet}))$. If $\partial_n([c]) = [a]$, then $\alpha_{n-1}(a) = \partial_n^{1,2}(b)$ is a boundary, hence $H_{n-1}(\alpha_\bullet)([a]) = 0$. Conversely, let $[a] \in H_{n-1}(C_\bullet^{12})$ with $H_{n-1}(\alpha_\bullet)([a]) = 0$. Choose a cycle $a \in C_{n-1}^{12}$ representing $[a]$. Then $\alpha_{n-1}(a)$ is a boundary in $C_\bullet^{1,2}$, so choose $b \in C_n^{1,2}$ with $\partial_n^{1,2}(b) = \alpha_{n-1}(a)$. Set $c := \beta_n(b)$. Then $\partial_n^{\mathcal{G}}(c) = 0$ since

$$\partial_n^{\mathcal{G}}(\beta_n(b)) = \beta_{n-1}(\partial_n^{1,2}(b)) = \beta_{n-1}(\alpha_{n-1}(a)) = 0.$$

Using this $b$ in the definition of $\partial_n([c])$ yields the unique element $a$ again by injectivity of $\alpha_{n-1}$. Thus $[a] \in \mathrm{im}(\partial_n)$.

The displayed equalities give exactness at every term. Naturality follows because the connecting map construction is functorial for morphisms of short exact sequences of chain complexes of abelian groups.

□

*Remark* 3.3.11. The connecting homomorphism $\partial_n$ is explicit. Let $c \in C_c(\mathcal{G}_n, A)$ be a cycle. Choose $b \in C_n(\mathcal{G}|_{U_1}; A) \oplus C_n(\mathcal{G}|_{U_2}; A)$ with $\beta_n(b) = c$. Compute $\partial_n^{1,2}(b)$. Exactness in degree $n-1$ implies $\partial_n^{1,2}(b) \in \mathrm{im}(\alpha_{n-1})$. Let $a \in C_{n-1}(\mathcal{G}|_{U_1 \cap U_2}; A)$ be the unique element with $\alpha_{n-1}(a) = \partial_n^{1,2}(b)$. Then $\partial_n([c]) = [a]$. The class $\partial_n([c])$ vanishes precisely when $c$ admits a lift $b$ that is a cycle in the middle complex.

In the Moore–Mayer–Vietoris construction above, the admissibility hypothesis requires the clopen pieces to be saturated. This hypothesis is used only to ensure that for every simplicial degree $n$ one has a global cover $\mathcal{G}_n = (\mathcal{G}|_{U_1})_n \cup (\mathcal{G}|_{U_2})_n$, so that for each simplicial degree $n$ the maps induced by the inclusions

$$0 \to C_c((\mathcal{G}|_{U_1 \cap U_2})_n, A) \xrightarrow{\iota} C_c((\mathcal{G}|_{U_1})_n, A) \oplus C_c((\mathcal{G}|_{U_2})_n, A) \xrightarrow{\rho} C_c(\mathcal{G}_n, A) \to 0$$

given by $\iota(c) = (c, -c)$ and $\rho(c_1, c_2) = c_1 + c_2$ form a short exact sequence. Injectivity of $\iota$ and $\rho \circ \iota = 0$ are immediate. Exactness in the middle follows since a compactly supported function on $\mathcal{G}_n$ vanishing on $(\mathcal{G}|_{U_1})_n$ and $(\mathcal{G}|_{U_2})_n$ vanishes everywhere. Surjectivity of $\rho$ holds because any $c \in C_c(\mathcal{G}_n, A)$ decomposes as $c = c\chi_{(\mathcal{G}|_{U_1})_n} + c\chi_{(\mathcal{G}|_{U_2})_n}$, and the intersection term adjusts the double count on $(\mathcal{G}|_{U_1 \cap U_2})_n$. Therefore there is a short exact sequence of Moore chain complexes and hence a long exact sequence in homology. In many applications, however, the relevant chains are compactly supported. Thus one does not need a global cover of all of $\mathcal{G}_n$. Instead, it suffices that a given cycle representative $c \in C_n(\mathcal{G}; A)$ is supported in the region where the cover behaves well, namely $\operatorname{supp}(c) \subseteq (\mathcal{G}|_{U_1})_n \cup (\mathcal{G}|_{U_2})_n$. Under this support condition the Mayer–Vietoris splitting and extension by zero arguments apply verbatim on the clopen subspace $(\mathcal{G}|_{U_1})_n \cup (\mathcal{G}|_{U_2})_n \subseteq \mathcal{G}_n$, even when $U_1$ and $U_2$ are not saturated. Corollary 3.3.12 records this support local replacement of saturation and yields a canonical long exact sequence controlling the homology classes represented by such compactly supported cycles.

**Corollary 3.3.12.** Let $\mathcal{G}$ be an ample groupoid and let $U_1, U_2 \subseteq \mathcal{G}_0$ be clopen with $U_1 \cup U_2 = \mathcal{G}_0$. Let $A$ be a discrete abelian group. For every $n \geq 0$ set $(\mathcal{G}|_{U_1})_n \cup (\mathcal{G}|_{U_2})_n \subseteq \mathcal{G}_n$ and define

$$C_n^{U_1,U_2}(\mathcal{G}; A) := C_c((\mathcal{G}|_{U_1})_n \cup (\mathcal{G}|_{U_2})_n, A), \qquad C_\bullet^{U_1,U_2}(\mathcal{G}; A) := (C_n^{U_1,U_2}(\mathcal{G}; A), \partial_n)_{n \geq 0},$$

where $\partial_\bullet$ is the Moore boundary of $C_\bullet(\mathcal{G}; A)$ restricted to $C_\bullet^{U_1,U_2}(\mathcal{G}; A)$.

Then the following holds:

1. For every $n \geq 0$ the sequence

$$0 \to C_n(\mathcal{G}|_{U_1 \cap U_2}; A) \xrightarrow{\alpha_n} C_n(\mathcal{G}|_{U_1}; A) \oplus C_n(\mathcal{G}|_{U_2}; A) \xrightarrow{\beta_n} C_n^{U_1,U_2}(\mathcal{G}; A) \to 0$$

is exact, where $\alpha_n(\xi) := (\xi, -\xi)$ and $\beta_n(\xi_1, \xi_2) := \tilde{\xi}_1 + \tilde{\xi}_2$, with $\tilde{\xi}_i$ denoting extension by zero along the clopen inclusion $(\mathcal{G}|_{U_i})_n \hookrightarrow (\mathcal{G}|_{U_1})_n \cup (\mathcal{G}|_{U_2})_n$.

2. The maps $\alpha_\bullet, \beta_\bullet$ form a short exact sequence of Moore chain complexes

$$0 \to C_\bullet(\mathcal{G}|_{U_1 \cap U_2}; A) \xrightarrow{\alpha_\bullet} C_\bullet(\mathcal{G}|_{U_1}; A) \oplus C_\bullet(\mathcal{G}|_{U_2}; A) \xrightarrow{\beta_\bullet} C_\bullet^{U_1,U_2}(\mathcal{G}; A) \to 0.$$

3. Define $C_\bullet^{U_1,U_2} := C_\bullet^{U_1,U_2}(\mathcal{G}; A)$. There is an induced long exact homology sequence

$$\cdots \longleftarrow H_{n-1}(\mathcal{G}|_{U_1}; A) \oplus H_{n-1}(\mathcal{G}|_{U_2}; A) \xleftarrow{H_{n-1}(\alpha_\bullet)} H_{n-1}(\mathcal{G}|_{U_1 \cap U_2}; A) \xleftarrow{\partial_n}$$

$$H_n(C_\bullet^{U_1,U_2}) \xleftarrow{H_n(\beta_\bullet)} H_n(\mathcal{G}|_{U_1}; A) \oplus H_n(\mathcal{G}|_{U_2}; A) \xleftarrow{H_n(\alpha_\bullet)} H_n(\mathcal{G}|_{U_1 \cap U_2}; A) \xleftarrow{\partial_{n+1}}$$

$$H_{n+1}(C_\bullet^{U_1,U_2}) \xleftarrow{H_{n+1}(\beta_\bullet)} H_{n+1}(\mathcal{G}|_{U_1}; A) \oplus H_{n+1}(\mathcal{G}|_{U_2}; A) \longleftarrow \cdots$$

4. Let $c \in C_n(\mathcal{G};A)$ be an $n$-cycle with $\operatorname{supp}(c) \subseteq (\mathcal{G}|_{U_1})_n \cup (\mathcal{G}|_{U_2})_n$. Then $c \in C_n^{U_1,U_2}(\mathcal{G};A)$ and the class $[c] \in H_n(\mathcal{G};A)$ lies in the image of the map induced by the inclusion

$$\iota_\bullet \colon C_\bullet^{U_1,U_2}(\mathcal{G};A) \hookrightarrow C_\bullet(\mathcal{G};A).$$

*Proof.* Fix $n \geq 0$. $U_1, U_2$ are clopen in $\mathcal{G}_0$, $\mathcal{G}$ is ample, the subsets $(\mathcal{G}|_{U_1})_n$, $(\mathcal{G}|_{U_2})_n$, $(\mathcal{G}|_{U_1\cap U_2})_n$ are clopen in $\mathcal{G}_n$. In particular, $(\mathcal{G}|_{U_1})_n$ and $(\mathcal{G}|_{U_2})_n$ are clopen in $(\mathcal{G}|_{U_1})_n \cup (\mathcal{G}|_{U_2})_n$, and $(\mathcal{G}|_{U_1\cap U_2})_n = (\mathcal{G}|_{U_1})_n \cap (\mathcal{G}|_{U_2})_n$.

1. **Injectivity of $\alpha_n$.** Injectivity of $\alpha_n$ is immediate and $\beta_n \circ \alpha_n = 0$.
2. **Surjectivity of $\beta_n$.** To prove surjectivity of $\beta_n$, let $\eta \in C_c(Y_n, A)$. Consider the clopen partition

$$Y_n = \Big((\mathcal{G}|_{U_1})_n \setminus (\mathcal{G}|_{U_1\cap U_2})_n\Big) \sqcup \Big((\mathcal{G}|_{U_2})_n \setminus (\mathcal{G}|_{U_1\cap U_2})_n\Big) \sqcup (\mathcal{G}|_{U_1\cap U_2})_n.$$

   Define $\eta_1 \in C_c((\mathcal{G}|_{U_1})_n, A)$ and $\eta_2 \in C_c((\mathcal{G}|_{U_2})_n, A)$ by

$$\begin{aligned} \eta_1|_{(\mathcal{G}|_{U_1})_n\setminus(\mathcal{G}|_{U_1\cap U_2})_n} &:= \eta|_{(\mathcal{G}|_{U_1})_n\setminus(\mathcal{G}|_{U_1\cap U_2})_n}, \quad \eta_1|_{(\mathcal{G}|_{U_1\cap U_2})_n} := 0,\\ \eta_2|_{(\mathcal{G}|_{U_2})_n\setminus(\mathcal{G}|_{U_1\cap U_2})_n} &:= \eta|_{(\mathcal{G}|_{U_2})_n\setminus(\mathcal{G}|_{U_1\cap U_2})_n}, \quad \eta_2|_{(\mathcal{G}|_{U_1\cap U_2})_n} := \eta|_{(\mathcal{G}|_{U_1\cap U_2})_n}. \end{aligned}$$

   Then $\eta_1, \eta_2$ are compactly supported and locally constant since all pieces are clopen. By construction, $\tilde{\eta}_1 + \tilde{\eta}_2 = \eta$ on each of the three clopen pieces of $Y_n$, hence on all of $Y_n$. Thus $\beta_n(\eta_1, \eta_2) = \eta$. To identify the kernel, let

$$(\xi_1, \xi_2) \in C_c((\mathcal{G}|_{U_1})_n, A) \oplus C_c((\mathcal{G}|_{U_2})_n, A)$$

   satisfy $\tilde{\xi}_1 + \tilde{\xi}_2 = 0$ in $C_c(Y_n, A)$. Restricting to $(\mathcal{G}|_{U_1})_n \setminus (\mathcal{G}|_{U_1\cap U_2})_n$ gives $\xi_1 = 0$ there. Restricting to $(\mathcal{G}|_{U_2})_n \setminus (\mathcal{G}|_{U_1\cap U_2})_n$ gives $\xi_2 = 0$ there. Restricting to $(\mathcal{G}|_{U_1\cap U_2})_n$ gives $\xi_1|_{(\mathcal{G}|_{U_1\cap U_2})_n} + \xi_2|_{(\mathcal{G}|_{U_1\cap U_2})_n} = 0$. Set $\zeta := \xi_1|_{(\mathcal{G}|_{U_1\cap U_2})_n} \in C_c((\mathcal{G}|_{U_1\cap U_2})_n, A)$. Then $(\xi_1, \xi_2) = (\zeta, -\zeta) = \alpha_n(\zeta)$, hence $\ker(\beta_n) = \operatorname{im}(\alpha_n)$. This proves exactness.
3. $C_\bullet^{U_1,U_2}(\mathcal{G};A)$ **is a subcomplex.** It remains to check that $C_\bullet^{U_1,U_2}(\mathcal{G};A)$ is a subcomplex of $C_\bullet(\mathcal{G};A)$. For every $i \in \{1,2\}$ and every face map $d_j \colon \mathcal{G}_n \to \mathcal{G}_{n-1}$ one has

$$d_j((\mathcal{G}|_{U_i})_n) \subseteq (\mathcal{G}|_{U_i})_{n-1}, \qquad d_j((\mathcal{G}|_{U_1\cap U_2})_n) \subseteq (\mathcal{G}|_{U_1\cap U_2})_{n-1},$$

   since deleting one arrow from an $n$-simplex in the nerve does not introduce units outside the same restriction. Therefore each pushforward $(d_j)_*$ sends

$$C_c((\mathcal{G}|_{U_i})_n, A) \to C_c((\mathcal{G}|_{U_i})_{n-1}, A), \qquad C_c((\mathcal{G}|_{U_1\cap U_2})_n, A) \to C_c((\mathcal{G}|_{U_1\cap U_2})_{n-1}, A).$$

   Moreover, by extension by zero along the clopen inclusions

$$(\mathcal{G}|_{U_1})_n \hookrightarrow Y_n, \qquad (\mathcal{G}|_{U_2})_n \hookrightarrow Y_n, \qquad (\mathcal{G}|_{U_1\cap U_2})_n \hookrightarrow Y_n,$$

the same holds with $C_c(Y_n, A)$ in place of $C_c((\mathcal{G}|_{U_i})_n, A)$. Hence the Moore boundary $\partial_n = \sum_{j=0}^{n}(-1)^j (d_j)_*$ restricts to a homomorphism $\partial_n \colon C_c(Y_n, A) \to C_c(Y_{n-1}, A)$, so the degreewise short exact sequence from 1. is a short exact sequence of chain complexes.

4. **Support reduction.** The condition $\operatorname{supp}(c) \subseteq Y_n$ means that $c$ vanishes on $\mathcal{G}_n \setminus Y_n$. Hence $c$ is a compactly supported locally constant function on $Y_n$, so $c \in C_c(Y_n, A) = C_n^{U_1,U_2}(\mathcal{G}; A)$. Since $\iota_\bullet$ is a chain map, it induces a homomorphism on homology and carries the class of $c$ in $H_n(C_\bullet^{U_1,U_2}(\mathcal{G}; A))$ to the class of $c$ in $H_n(\mathcal{G}; A)$.

□

## 3.4 Computing Homology of a SFT Groupoid

We demonstrate, that the Moore–Mayer–Vietoris long exact sequence from Theorem 3.3.10 is useful for explicit computations. One covers the unit space by saturated clopen subsets, computes the homology of the corresponding reductions, and then recovers the homology of the whole groupoid via exactness. In this example we combine Moore–Mayer–Vietoris with the UCT from Theorem 3.2.3 to exhibit how torsion in degree 0 affects homology with finite field coefficients.

Let $A \in \operatorname{Mat}(N \times N, \mathbb{N}_0)$ be a square matrix with no zero row and no zero column. Fix a finite directed graph $E_A$ with vertex set $\{1, \dots, N\}$ whose adjacency matrix is $A$ – allowing multiple edges. Let $E_A^\infty = \left\{(e_n)_{n\geq 0} \in E_A^{\mathbb{N}_0} \,\middle|\, r(e_n) = s(e_{n+1}) \text{ for all } n \geq 0\right\}$ be the space of infinite directed paths, endowed with the product topology, and let $\sigma \colon E_A^\infty \to E_A^\infty, (e_0, e_1, e_2, \dots) \mapsto (e_1, e_2, e_3, \dots)$ be the left shift. The range and source maps are $r(x, n, y) = x$ and $s(x, n, y) = y$, with units $\chi_x = (x, 0, x)$, inverses $(x, n, y)^{-1} = (y, -n, x)$, and multiplication $(x, n, y) \cdot (y, m, z) = (x, n+m, z)$ whenever $s(x, n, y) = r(y, m, z)$. Then $\sigma$ is a local homeomorphism, and the associated Deaconu–Renault groupoid $\mathcal{G}_A$ has unit space $(\mathcal{G}_A)_0 = E_A^\infty$ and arrow space

$$(\mathcal{G}_A)_1 = \{(x, n, y) \in E_A^\infty \times \mathbb{Z} \times E_A^\infty \mid \exists\, k, \ell \in \mathbb{N}_0 :\ n = k - \ell, \sigma^k(x) = \sigma^\ell(y)\}.$$

Equip $(\mathcal{G}_A)_1$ with the étale topology generated by the compact open bisections

$$Z(\alpha, \beta) := \{(\alpha z, |\alpha| - |\beta|, \beta z) \mid z \in E_A^\infty\},$$

where $\alpha, \beta$ range over finite paths with common range and $|\alpha|$ denotes length. Then $r$ and $s$ restrict to homeomorphisms on each $Z(\alpha, \beta)$ and these sets form a basis. $\mathcal{G}_A$ is second countable, locally compact, Hausdorff, totally disconnected, and étale; in particular it is ample. This is the same construction of the Deaconu–Renault groupoid of a local homeomorphism as in [1, §2.5].

For SFT groupoids, the integral homology of $\mathcal{G}_A$ is given in terms of $\mathbb{1} - A^\mathsf{T}$. We have

$$\begin{aligned} H_0(\mathcal{G}_A) &\cong \operatorname{coker}(\mathbb{1} - A^\mathsf{T}), \\ H_1(\mathcal{G}_A) &\cong \ker(\mathbb{1} - A^\mathsf{T}), \\ H_n(\mathcal{G}_A) &= 0 \ \text{ for } n \geq 2, \end{aligned}$$

where $\mathbb{1} - A^{\mathsf{T}}$ acts on $\mathbb{Z}^N$ [13, Theorem 4.14]. We consider now the matrices

$$A = \begin{pmatrix} 2 & 1 \\ 1 & 0 \end{pmatrix}, \qquad B = \begin{pmatrix} 2 & 1 \\ 1 & 2 \end{pmatrix}, \qquad C = (3),$$

and compute the integral homology of $\mathcal{G}_A$, $\mathcal{G}_B$, and $\mathcal{G}_C$. For $A$ we have

$$\mathbb{1} - A^{\mathsf{T}} = \begin{pmatrix} -1 & -1 \\ -1 & 1 \end{pmatrix}, \qquad \det(\mathbb{1} - A^{\mathsf{T}}) = (-1) \cdot 1 - (-1) \cdot (-1) = -2.$$

Hence $\mathbb{1} - A^{\mathsf{T}}$ has full rank over $\mathbb{Z}$ and $\ker(\mathbb{1} - A^{\mathsf{T}}) = 0$. Moreover, the Smith normal form is

$$\begin{pmatrix} -1 & -1 \\ -1 & 1 \end{pmatrix} \overset{R_1 \leftarrow (-1) \cdot R_1}{\rightsquigarrow} \begin{pmatrix} 1 & 1 \\ -1 & 1 \end{pmatrix} \overset{R_2 \leftarrow R_2 + R_1}{\rightsquigarrow} \begin{pmatrix} 1 & 1 \\ 0 & 2 \end{pmatrix} \overset{C_2 \leftarrow C_2 + (-1) \cdot C_1}{\rightsquigarrow} \begin{pmatrix} 1 & 0 \\ 0 & 2 \end{pmatrix} = \operatorname{diag}(1, 2),$$

so $\operatorname{coker}(\mathbb{1} - A^{\mathsf{T}}) \cong \mathbb{Z}/2\mathbb{Z}$. Therefore

$$\begin{aligned} H_0(\mathcal{G}_A) &\cong \mathbb{Z}/2\mathbb{Z}, \\ H_1(\mathcal{G}_A) &= 0, \\ H_n(\mathcal{G}_A) &= 0 \text{ for } n \geq 2. \end{aligned}$$

For $B$ we have

$$\mathbb{1} - B^{\mathsf{T}} = \begin{pmatrix} -1 & -1 \\ -1 & -1 \end{pmatrix}.$$

The condition $(\mathbb{1} - B^{\mathsf{T}})(x, y)^{\mathsf{T}} = 0$ is $-x - y = 0$, hence $\ker(\mathbb{1} - B^{\mathsf{T}}) \cong \mathbb{Z}$ generated by $(1, -1)$. The image is generated by $(1, 1)$, which is primitive in $\mathbb{Z}^2$, so $\operatorname{coker}(\mathbb{1} - B^{\mathsf{T}}) \cong \mathbb{Z}^2/\langle (1, 1) \rangle_{\mathbb{Z}} \cong \mathbb{Z}$. Thus, we have for homology

$$\begin{aligned} H_0(\mathcal{G}_B) &\cong \mathbb{Z}, \\ H_1(\mathcal{G}_B) &\cong \mathbb{Z}, \\ H_n(\mathcal{G}_B) &= 0 \text{ for } n \geq 2. \end{aligned}$$

For $C$ we have $\mathbb{1} - C^{\mathsf{T}} = -2$, so $\ker(\mathbb{1} - C^{\mathsf{T}}) = 0$ and $\operatorname{coker}(\mathbb{1} - C^{\mathsf{T}}) \cong \mathbb{Z}/2\mathbb{Z}$. Hence

$$\begin{aligned} H_0(\mathcal{G}_C) &\cong \mathbb{Z}/2\mathbb{Z}, \\ H_1(\mathcal{G}_C) &= 0, \\ H_n(\mathcal{G}_C) &= 0 \text{ for } n \geq 2. \end{aligned}$$

Next set $\mathcal{G} := \mathcal{G}_A \sqcup \mathcal{G}_B \sqcup \mathcal{G}_C$, the disjoint union groupoid. Then $\mathcal{G}$ is ample, and levelwise its nerve decomposes as $\mathcal{G}_n = (\mathcal{G}_A)_n \sqcup (\mathcal{G}_B)_n \sqcup (\mathcal{G}_C)_n$. Consequently the associated compactly supported Moore chain complex splits as a direct sum, and therefore

$$H_n(\mathcal{G}) \cong H_n(\mathcal{G}_A) \oplus H_n(\mathcal{G}_B) \oplus H_n(\mathcal{G}_C) \text{ for } n \geq 0.$$

In particular,

$$
\begin{aligned}
H_0(\mathcal{G}) &\cong \mathbb{Z} \oplus (\mathbb{Z}/2\mathbb{Z})^2, \\
H_1(\mathcal{G}) &\cong \mathbb{Z}, \\
H_n(\mathcal{G}) &= 0 \text{ for } n \geq 2.
\end{aligned}
$$

We now illustrate Moore–Mayer–Vietoris on a saturated clopen cover of $\mathcal{G}_0$. Define

$$
\begin{aligned}
U_1 &:= (\mathcal{G}_A)_0 \sqcup (\mathcal{G}_B)_0, \\
U_2 &:= (\mathcal{G}_B)_0 \sqcup (\mathcal{G}_C)_0.
\end{aligned}
$$

These subsets are clopen because they are unions of clopen components of the disjoint union space $\mathcal{G}_0$ and saturated as there are no arrows between distinct components, so any union of components is a union of orbits. One has $U_1 \cup U_2 = \mathcal{G}_0$ by construction and $U_1 \cap U_2 = (\mathcal{G}_B)_0$ since $(\mathcal{G}_B)_0$ is the unique component contained in both unions. The reductions are

$$
\begin{aligned}
\mathcal{G}|_{U_1} &= \mathcal{G}_A \sqcup \mathcal{G}_B, \\
\mathcal{G}|_{U_2} &= \mathcal{G}_B \sqcup \mathcal{G}_C, \\
\mathcal{G}|_{U_1 \cap U_2} &= \mathcal{G}_B.
\end{aligned}
$$

Applying Theorem 3.3.10 yields the long exact sequence

$$
\cdots \to H_n(\mathcal{G}_B) \xrightarrow{\alpha_n} H_n(\mathcal{G}_A \sqcup \mathcal{G}_B) \oplus H_n(\mathcal{G}_B \sqcup \mathcal{G}_C) \xrightarrow{\beta_n} H_n(\mathcal{G}) \xrightarrow{\partial_n} H_{n-1}(\mathcal{G}_B) \to \cdots.
$$

Using the canonical identifications

$$
\begin{aligned}
H_n(\mathcal{G}_A \sqcup \mathcal{G}_B) &\cong H_n(\mathcal{G}_A) \oplus H_n(\mathcal{G}_B), \\
H_n(\mathcal{G}_B \sqcup \mathcal{G}_C) &\cong H_n(\mathcal{G}_B) \oplus H_n(\mathcal{G}_C), \\
H_n(\mathcal{G}) &\cong H_n(\mathcal{G}_A) \oplus H_n(\mathcal{G}_B) \oplus H_n(\mathcal{G}_C),
\end{aligned}
$$

the map $\alpha_n$ is the difference of the two inclusions of the $\mathcal{G}_B$-summand and is given by $\alpha_n([b]) = ([0], [b], [-b], [0])$. The map $\beta_n$ is induced by the two inclusions $\mathcal{G}|_{U_1} \hookrightarrow \mathcal{G}$ and $\mathcal{G}|_{U_2} \hookrightarrow \mathcal{G}$; under the above identifications it is $\beta_n([a], [b_1], [b_2], [c]) = ([a], [b_1 + b_2], [c])$. In particular, $\beta_n$ is surjective and $\ker(\beta_n) = \{([0], [b], [-b], [0]) \mid b \in H_n(\mathcal{G}_B)\} = \operatorname{im}(\alpha_n)$. Exactness therefore forces $\partial_n = [0]$ for all $n$, and the Moore–Mayer–Vietoris sequence recovers the direct-sum decomposition of $H_n(\mathcal{G})$.

Finally we compute homology with finite coefficients via UCT. Fix a prime $p$. Since $H_2(\mathcal{G}) = 0$ and $H_1(\mathcal{G}) \cong \mathbb{Z}$ is torsion-free, the UCT implies $H_n(\mathcal{G}; \mathbb{Z}/p\mathbb{Z}) = 0$ for all $n \geq 2$. For $n = 0$:

$$
H_0(\mathcal{G}; \mathbb{Z}/p\mathbb{Z}) \cong H_0(\mathcal{G}) \otimes_{\mathbb{Z}} \mathbb{Z}/p\mathbb{Z} \cong \begin{cases} \mathbb{Z}/p\mathbb{Z}, & \text{for } p \text{ odd}, \\ (\mathbb{Z}/2\mathbb{Z})^3, & \text{for } p = 2. \end{cases}
$$

For $n = 1$ the UCT from Theorem 3.2.3 yields a short exact sequence

$$0 \to H_1(\mathcal{G}) \otimes_{\mathbb{Z}} \mathbb{Z}/p\mathbb{Z} \xrightarrow{\iota_1} H_1(\mathcal{G}; \mathbb{Z}/p\mathbb{Z}) \xrightarrow{\beta_1} \mathrm{Tor}_1^{\mathbb{Z}}(H_0(\mathcal{G}), \mathbb{Z}/p\mathbb{Z}) \to 0.$$

Here $\iota_1$ is the change-of-coefficients map induced by $C_\bullet(\mathcal{G}; \mathbb{Z}) \otimes \mathbb{Z}/p\mathbb{Z} \to C_\bullet(\mathcal{G}; \mathbb{Z}/p\mathbb{Z})$, and $\beta_1$ is the Bockstein connecting morphism associated to $0 \to \mathbb{Z} \xrightarrow{\cdot p} \mathbb{Z} \xrightarrow{\pi} \mathbb{Z}/p\mathbb{Z} \to 0$. In our example $H_1(\mathcal{G}) \cong \mathbb{Z}$ and $H_0(\mathcal{G}) \cong \mathbb{Z} \oplus (\mathbb{Z}/2\mathbb{Z})^2$. If $p$ is odd, then $\mathrm{Tor}_1^{\mathbb{Z}}(H_0(\mathcal{G}), \mathbb{Z}/p\mathbb{Z}) = 0$, so $\beta_1 = 0$ and $\iota_1$ is an isomorphism. If $p = 2$, then $H_1(\mathcal{G}; \mathbb{Z}/2\mathbb{Z}) \cong (\mathbb{Z}/2\mathbb{Z})^3$ and $\mathrm{Tor}_1^{\mathbb{Z}}(H_0(\mathcal{G}), \mathbb{Z}/2\mathbb{Z}) \cong (\mathbb{Z}/2\mathbb{Z})^2$. With coordinates corresponding to the decomposition $\mathcal{G} = \mathcal{G}_A \sqcup \mathcal{G}_B \sqcup \mathcal{G}_C$, the maps are $\iota_1([1]) = ([0], [1], [0])$ and $\beta_1([a], [b], [c]) = ([a], [c])$. Here $H_1(\mathcal{G}) \cong \mathbb{Z}$ gives $H_1(\mathcal{G}) \otimes_{\mathbb{Z}} \mathbb{Z}/p\mathbb{Z} \cong \mathbb{Z}/p\mathbb{Z}$, and

$$\mathrm{Tor}_1^{\mathbb{Z}}(\mathbb{Z}, \mathbb{Z}/p\mathbb{Z}) = 0,$$

$$\mathrm{Tor}_1^{\mathbb{Z}}(\mathbb{Z}/2\mathbb{Z}, \mathbb{Z}/p\mathbb{Z}) \cong \begin{cases} 0, & \text{for } p \text{ odd,} \\ \mathbb{Z}/2\mathbb{Z}, & \text{for } p = 2. \end{cases}$$

Since $H_0(\mathcal{G}) \cong \mathbb{Z} \oplus (\mathbb{Z}/2\mathbb{Z})^2$, this gives

$$\mathrm{Tor}_1^{\mathbb{Z}}(H_0(\mathcal{G}), \mathbb{Z}/p\mathbb{Z}) \cong \begin{cases} 0, & \text{for } p \text{ odd,} \\ (\mathbb{Z}/2\mathbb{Z})^2, & \text{for } p = 2. \end{cases}$$

Thus $H_1(\mathcal{G}; \mathbb{Z}/p\mathbb{Z}) \cong \mathbb{Z}/p\mathbb{Z}$ for odd $p$. For $p = 2$, the chain complex defining $H_\bullet(\mathcal{G}; \mathbb{Z}/2\mathbb{Z})$ is a complex of $\mathbb{Z}/2\mathbb{Z}$-vector spaces, hence $H_1(\mathcal{G}; \mathbb{Z}/2\mathbb{Z})$ is itself a $\mathbb{Z}/2\mathbb{Z}$-vector space; the above short exact sequence therefore splits though non-canonically in $\mathbf{Vect}_{\mathbb{Z}/2\mathbb{Z}}$. Consequently,

$$H_1(\mathcal{G}; \mathbb{Z}/p\mathbb{Z}) \cong \begin{cases} \mathbb{Z}/p\mathbb{Z}, & \text{for } p \text{ odd,} \\ (\mathbb{Z}/2\mathbb{Z})^3, & \text{for } p = 2. \end{cases}$$